\documentclass[10pt,reqno,a4paper]{amsart} 

\usepackage[a4paper,
left=1in,right=1in,top=1.2in,bottom=1.2in,%
footskip=.25in]{geometry}

\usepackage{bm}
\usepackage{latexsym, amsmath, amstext, amssymb, amsfonts, amscd, bm, array, multirow, amsbsy, mathrsfs, stmaryrd}
\usepackage{amsthm}
\usepackage{t1enc}
\usepackage[mathscr]{eucal}
\usepackage{indentfirst}
\usepackage{pb-diagram}
\usepackage{graphicx}
\usepackage{fancyhdr}
\usepackage{fancybox}
\usepackage{color}
\usepackage{tikz-cd}
\usetikzlibrary{arrows}
\usepackage[all]{xy}
\usepackage{tikz}
\usepackage{xparse}
\usetikzlibrary{matrix}
\usepackage{upgreek}
\usepackage[utf8]{inputenc}
\usepackage[T1]{fontenc}
\usepackage{bbm}
\usepackage{enumitem}
\setlist[enumerate]{leftmargin=*,noitemsep}
\setlist[itemize]{leftmargin=*,noitemsep}
\usepackage{float}
\usepackage{soul} 
\usepackage{hyperref}
\hypersetup{colorlinks=false,pdfborderstyle={/S/U/W 0}}

\newcommand{\doi}[1]{\url{https://doi.org/#1}}
\usepackage[noadjust]{cite}

\makeatletter
\def\@defaultbiblabelstyle#1{[#1]}
\makeatother

\usepackage{url}

\theoremstyle{plain}
\newtheorem{thm}{Theorem}[section]

\newtheorem{prop}[thm]{Proposition}
\newtheorem{lemma}[thm]{Lemma}
\newtheorem{cor}[thm]{Corollary}

\theoremstyle{definition}
\newtheorem{definition}[thm]{Definition}

\theoremstyle{remark}
\newtheorem{remark}[thm]{Remark}
\newtheorem{example}[thm]{Example}

\usepackage[math]{anttor}
\usepackage{orcidlink}

\allowdisplaybreaks
 
\theoremstyle{remark}
\newtheorem{ep}[thm]{Example}

\usetikzlibrary{arrows.meta}
 
\hypersetup{
  pdftitle={The Cauchy problem for skew-torsion-parallel spinors with applications to supergravity},
  pdfauthor={C. S. Shahbazi},
  pdfkeywords={Spin geometry, metric connections with torsion, parallel spinors, Lorentzian four-manifolds, Cauchy problem, supergravity, bundle gerbes}
}

\def\mA{\mathrm{A}}
\def\mT{\mathrm{T}}

\newcommand{\End}{\mathrm{End}}

\newcommand{\frt}{\mathfrak{t}}

\newcommand{\frG}{\mathfrak{G}}

\usepackage{tikz}
\usetikzlibrary{arrows.meta}

\DeclareFontFamily{U}{rsf}{}
\DeclareFontShape{U}{rsf}{m}{n}{<5> <6> rsfs5 <7> <8> <9> rsfs7 <10-> rsfs10}{}
\DeclareMathAlphabet\Scr{U}{rsf}{m}{n}

\def\dd{\mathrm{d}}

\def\l\Xi{\overrightarrow{\Xi}}
\def\r\Xi{\overleftarrow{\Xi}}

\def\l{\partial^l}

\def\Conf{\mathrm{Conf}}

\def\Sol{\mathrm{Sol}}

\newcommand{\Tr}{\mathrm{Tr}}

\def\cR{{\mathcal R}}

\def\cB{\Scr B}

\def\cH{\mathcal{H}}

\def\Cl{\mathrm{Cl}}

\def\U{\mathrm{U}}

\def\cD{\mathcal{D}}

\def\cA{\mathcal{A}}
\def\cE{\mathcal{E}}
\def\cI{\mathcal{I}}
\def\cP{\mathcal{P}}

\def\cT{\mathcal{T}}
\def\cF{\mathcal{F}}

\def\G_2{\mathrm{G_2}}

\def\cL{\mathcal{L}}

\def\cX{\mathcal{X}}

\def\frc{\mathfrak{c}}
\def\frw{\mathfrak{w}}

\def\G{\mathrm{G}}

\def\cM{\mathcal{M}}
\def\cW{\mathcal{W}}
\def\cL{\mathcal{L}}
\def\cH{\mathcal{H}}

\def\fra{\mathfrak{a}}

\def\cU{\mathcal{U}}

\def\fre{\mathfrak{e}}

\newcommand{\Ker}{\operatorname{Ker}}

\newcommand{\Ric}{\mathrm{Ric}}

\newcommand{\surj}{\to\kern-1.8ex\to}

\newcommand{\frf}{\mathfrak{f}}

\newcommand{\arxiv}[1]{{\tt
		\href{http://www.arXiv.org/abs/#1}{arXiv:#1}}}

\title[Cauchy problem for skew-torsion-parallel spinors]{The Cauchy problem for skew-torsion-parallel spinors with applications to supergravity}
\author{C. S. Shahbazi}
\date{}

\begin{document}

\begin{abstract}
We solve the Cauchy problem for an irreducible real spinor parallel with respect to a metric connection with skew-symmetric torsion on a globally hyperbolic Lorentzian four-manifold, proving that it is well-posed modulo an obstruction that we identify explicitly, and obtaining an initial-data characterization of such spinors. Our proof relies on the theory of spinorial forms, which we use to explicitly integrate the coupled flow for the spinor and underlying metric via a rotated pull-back action of the flow generated by a vector field on the Cauchy hypersurface $\Sigma$. We characterize the corresponding constraint equations as an explicit exterior differential system and describe all its solutions in terms of contact structures and flat foliations, obtaining a complete classification in the compact case. In particular, we conclude that every oriented compact three-manifold arises as the Cauchy hypersurface of a globally hyperbolic Lorentzian four-manifold equipped with a skew-torsion parallel spinor. Building on these results, we study the interaction between the four-dimensional Neveu-Schwarz supergravity evolution flow on an abelian bundle gerbe and the spinorial evolution flow determined by its supersymmetric conditions, which we prove to be well-posed and compatible. This establishes the well-posedness of the supersymmetric Neveu-Schwarz system in four dimensions and provides an initial-data characterization of supersymmetry. Furthermore, we obtain a strong rigidity result stating that every globally hyperbolic supersymmetric Neveu-Schwarz solution with a compact Cauchy hypersurface is flat and has vanishing bundle-gerbe curvature.
\end{abstract}

\maketitle

\setcounter{tocdepth}{1}
\tableofcontents


\section{Introduction}
\label{chapter:introduction} 



\subsection*{Main results and goals} 


This article is devoted to the study of irreducible real spinors parallel with respect to a metric connection with totally skew torsion on a Lorentzian four-manifold, which we call four-dimensional skew-torsion parallel spinors, focusing on the case in which the underlying Lorentzian metric is globally hyperbolic. Whereas the mathematical theory of such spinors is well developed in Riemannian signature \cite{AgricolaFriedrich,AgricolaSrni,FriedrichIvanov,IvanovPapadopoulos}, its Lorentzian counterpart remains in its infancy, with some notable exceptions \cite{ErnstGalaev,ErnstGalaevTwistorFree}. This is despite the extensive occurrence of Lorentzian spinors parallel with respect to metric connections with torsion in the theory of supersymmetric solutions in supergravity, see \cite{Gran:2018ijr,Ortin,Tomasiello} and their references and citations for more details. \medskip

Our main goal is to establish the well-posedness of the Cauchy problem for a skew-torsion parallel spinor in four dimensions, characterize the corresponding initial data, and apply the results to study the Cauchy problem for the supersymmetry conditions of Neveu-Schwarz supergravity \cite{Ortin} and its interaction with the Cauchy problem for its bosonic equations. This requires, as a preliminary step, reducing the skew-torsion parallelicity condition into a differential system of evolution equations and time-dependent constraints for a family of Riemannian metrics and the appropriate geometric structure induced by the parallel spinor on the given Cauchy hypersurface. To efficiently perform this reduction and bypass the subtleties associated with the dependence of the spinor bundle on the underlying metric, we adopt the theory of spinorial polyforms developed in \cite{Cortes:2019xmk}, which, in the particular case of four Lorentzian dimensions, allows us to equivalently study any parallelicity condition on an irreducible real spinor in terms of a differential system for a naturally associated \emph{isotropic parallelism} \cite{ShahbaziTesis}, defined as a particular equivalence class of global coframes on the underlying four-manifold $M$. A representative in an isotropic parallelism is a global null tetrad in the sense of Newman and Penrose \cite{NewmanPenrose}, and the aforementioned equivalence relation is the \emph{gauge freedom} in the spin-coefficient formalism of Geroch, Held and Penrose \cite{GerochHeldPenrose}. In this sense, an isotropic parallelism encodes the square \cite{Atiyah} of a nowhere vanishing irreducible spinor in four Lorentzian dimensions \cite{Cortes:2019xmk}, see \cite{GilGarciaShahbazi} for the general complex case. The novelty of the isotropic parallelisms lies in the global information they contain. A spinor requires an a priori choice of metric and spin structure, whereas an isotropic parallelism determines by itself the Lorentzian metric, the time orientation and the isomorphism class of the spin structure for which $(M,g)$ admits a nowhere-vanishing irreducible real spinor. The sign ambiguity of the spinor is absorbed in the equivalence relation, and the dependence of the spinor bundle on the metric, which is the main technical obstacle in the formulation of spinorial evolution and moduli problems, disappears altogether. First-order spinorial equations become explicit exterior differential systems for isotropic parallelisms for which the curvature is read off from a pair of derived one-forms, and global invariants can be extracted via the theory of exterior differential systems \cite{BryantBook}. The fact that isotropic parallelisms encode the underlying metric as a variable is fundamental to studying evolution and moduli problems, and this becomes apparent in the supersymmetric conditions in supergravity, which always involve the metric as a variable.\medskip

Using the formalism of isotropic parallelisms introduced in Section \ref{sec:torsionparallel}, in Section \ref{sec:CauchyKillingSkew} we proceed to study the globally hyperbolic case, reformulating the problem as a differential system of evolution equations and time-dependent constraints for a family of global coframes, a triplet of functions, and a family of one-forms from which the full four-dimensional Lorentzian manifold and corresponding skew-torsion parallel spinor is reconstructed. In Subsection \ref{sec:ReductionGloballyHyperbolicparallelisms}, we begin this reduction procedure by reducing an isotropic parallelism whose associated Lorentzian metric is globally hyperbolic. This is equivalent to reducing a pair $(g,\varepsilon)$, where $g$ is a globally hyperbolic Lorentzian metric and $\varepsilon$ is a nowhere vanishing irreducible real spinor on $(M,g)$. In Subsection \ref{sec:evolutionskewtorsion} we apply this reduction to the specific case of a skew-torsion parallel spinor, and we identify the system's first-order evolution and constraint equations as an explicit exterior differential system, which we present in Proposition \ref{prop:Thetaevolutionsystem}. In Subsection \ref{sec:skewtorsionwellposedness} we build on this constrained evolution system to prove, in Theorem \ref{thm:wellposedFinalSkew}, well-posedness of the corresponding initial value problem, including the existence and uniqueness of the four-dimensional germ. Our method uses an educated choice of the lapse function, with respect to which we can integrate the evolution equations explicitly in terms of a deformation of the flow of a vector field on the Cauchy hypersurface, given by the spatial projection of the Dirac current of the skew-torsion parallel spinor. Once existence and uniqueness are proven in this natural gauge, Lemma \ref{lemma:uniqueevol} guarantees uniqueness for any other choice of lapse function. Furthermore, Corollary \ref{cor:timeindependentdevelopment} removes the \emph{time dependence} of the metric and torsion of the globally hyperbolic development constructed in Theorem \ref{thm:wellposedFinalSkew}, identifying every globally hyperbolic development of a skew-torsion parallel spinor in four dimensions, near its Cauchy hypersurface, with an open subset of an algebraically special wave spacetime admitting a null Killing vector field with constant optical twist. This identifies, in full generality, such developments and leads to an initial-data characterization of skew-torsion parallel spinors, which we present via an exterior differential system in Corollary \ref{cor:initialdataskew}. This exterior differential system is simple enough that it can be presented here in an abbreviated form without further explanation.

\begin{thm}
A connected and oriented three-manifold $\Sigma$ occurs as the Cauchy hypersurface of a Lorentzian four-manifold equipped with a skew-torsion parallel spinor if and only if it admits a global coframe $(e_u,e_l,e_n)$ satisfying the following differential system:
\begin{equation*} 
\dd  e_l\wedge e_u=0,  \qquad \dd_\Sigma e_n\wedge e_u=0,  
\end{equation*}
together with the condition that $\dd e_u \wedge e_u \in \Omega^3(\Sigma)$ be either identically zero or nowhere vanishing on $\Sigma$.
\end{thm}  

If $\Sigma$ admits a solution $(e_u,e_l,e_n)$ to this differential system, then $(e_u,e_l,e_n)$, together with a choice of nowhere vanishing function in the case $\dd_{\Sigma} e_u \wedge e_u \equiv 0$, determines uniquely the corresponding skew torsion and germ of globally hyperbolic extension, as explained in Theorem \ref{thm:wellposedFinalSkew}. Conversely, every globally hyperbolic Lorentzian four-manifold equipped with a skew-torsion parallel spinor determines, on each smooth Cauchy hypersurface, a natural solution to the previous system. Furthermore, we prove that if a certain conformal transformation of the Riemannian metric on the Cauchy hypersurface is complete, then the globally hyperbolic development constructed in Theorem \ref{thm:wellposedFinalSkew} is maximal. After proving well-posedness and identifying the time-independent model, and still within Section \ref{sec:CauchyKillingSkew}, we study the three-dimensional Cauchy hypersurfaces that admit admissible initial data in terms of the aforementioned exterior differential system, yielding an explicit characterization in terms of an integration constant $r$ that depends on the given skew-torsion and determines the geometric properties of the initial data. If $r\neq 0$ we encounter a particular class of contact structures that we describe in detail, whereas if $r=0$ we find certain foliations with flat leaves. This is the content of Proposition \ref{prop:geometryconstraintbranches}. In particular, we prove that every compact oriented three-manifold without boundary admits admissible initial data, and we relate it to a particular class of contact structures. \medskip

In Section \ref{sec:Globallyhyperbolicsusy}, we proceed to apply the results of Section \ref{sec:CauchyKillingSkew} to the supersymmetric Neveu-Schwarz system. This system is a coupled gauge system defined on an abelian bundle gerbe by the bosonic equations of motion of Neveu-Schwarz supergravity together with the corresponding supersymmetry conditions evaluated on a bosonic configuration. Crucially, these supersymmetry conditions involve an irreducible real spinor parallel with respect to a metric connection with completely skew torsion given by the curvature of a curving on an abelian bundle gerbe, which is the system's gauge variable, together with a circle-valued function, the so-called \emph{dilaton}. Thus, the supersymmetric Neveu-Schwarz system combines second-order Einstein- and Maxwell-type equations with a pair of first-order spinorial equations. In Subsection \ref{sec:NSNSsystem}, we introduce the Neveu-Schwarz system and its supersymmetry conditions, providing a global model based on an abelian bundle gerbe, and describing the \emph{dilaton} of the theory in terms of a certain pullback of a circle-valued function. In Subsection \ref{sec:reductionbundlegerbe}, we reduce the bundle gerbe to a Cauchy hypersurface as a preliminary step toward reducing the full Neveu-Schwarz system, which we perform in detail in \ref{sec:NSNSsystemreduction}. This yields a system of second-order evolution equations and time-dependent constraints that generalize the vacuum Einstein flow and must be supplemented by the evolution and time-dependent constraints imposed by the supersymmetry conditions, including the skew-torsion parallelity condition for a spinor. Proposition \ref{prop:normalhyperbolicflow} presents the complete second-order system. The evolution equations appear in \eqref{eq:normalevolution}, while the time-dependent constraints are given in \eqref{eq:normalconstraint}. In Subsection \ref{sec:susyflowCauchy} we obtain the evolution flow, which we call the \emph{supersymmetric Neveu-Schwarz evolution flow}, determined by the supersymmetry conditions of Neveu-Schwarz supergravity, using for this the machinery developed in Section \ref{sec:torsionparallel}. Every supersymmetric Neveu-Schwarz evolution flow determines a natural candidate flow for the second-order equations of motion. In Theorem \ref{thm:compatibilityflows} we obtain the compatibility condition between the supersymmetric Neveu-Schwarz evolution flow and the constraint equations of the second-order evolution problem. In Subsection \ref{sec:susyinitialdevelopment} we address the Cauchy problem posed by the supersymmetry conditions of four-dimensional Neveu-Schwarz supergravity by applying the formalism developed in Section \ref{sec:CauchyKillingSkew}, proving in Theorem \ref{thm:wellposedFinalSusy} that it is well-posed. The evolution in a natural gauge for the lapse function simplifies and is fully determined by the flow of the spatial projection of the Dirac current of the supersymmetry parameter. Furthermore, we show that the admissible initial data corresponds to the $r=0$ case of the general skew-torsion parallel case, and hence, if compact, can only exist on torus bundles over $S^1$. Altogether, this establishes the well-posedness of the Cauchy problem determined by the supersymmetry conditions of Neveu-Schwarz supergravity. In Subsection \ref{sec:compatibleinitialdata}, we observe that every initial datum for the supersymmetric Neveu-Schwarz evolution flow defines a natural candidate initial data for the second-order Cauchy problem determined by the corresponding equations of motion. Therefore, it makes sense to consider initial data that satisfies the constraint equations for \emph{both} evolution problems. In Theorem \ref{thm:susyflowcompatibility}, we prove that if the supersymmetric Neveu-Schwarz evolution flow satisfies a particular subset of the second-order constraint equations at zero time, given by the Maxwell constraint and a component of the momentum constraint, then it satisfies all evolution and constraint equations of Neveu-Schwarz supergravity, and therefore defines a supersymmetric solution of the theory. This implies that the supersymmetry conditions provide partial \emph{first-order integrability} of the second-order flow when starting on initial data admissible to both evolution problems. Furthermore, since the Cauchy problem determined by the second-order equations of motion is well-posed by the results of \cite{Schiller:2026vls}, this leads to an initial data characterization of supersymmetry that we present in Corollary \ref{cor:NSNSsusyinitialdataII}, and which, to the best of our knowledge, is the first of its kind in the literature. Finally, and rather surprisingly, in Theorem \ref{thm:compactrigidity} we prove a strong rigidity result when the underlying Cauchy hypersurface is compact, which states that every globally hyperbolic supersymmetric Neveu-Schwarz solution admitting a compact Cauchy hypersurface is flat and flux-less, its supersymmetry parameter is Levi-Civita parallel, and its dilaton curvature is a constant multiple of its Dirac current. \medskip

The proof of Theorem \ref{thm:compactrigidity} combines the results of Section \ref{sec:CauchyKillingSkew} with a leafwise Bochner--Hodge argument on the compact Cauchy hypersurface. By the universal property established in the proof of Corollary \ref{cor:timeindependentdevelopment}, every globally hyperbolic development of the induced data embeds into its time-independent development, which is maximal since the Cauchy hypersurface is compact. Integrating a certain three-form constructed from the Maxwell constraint and the restriction of the Dirac current to the Cauchy hypersurface implies that the restriction of the dilaton curvature is proportional to the Dirac current. This restriction is then closed and nowhere vanishing, and its kernel defines a foliation by flat surfaces equipped with a parallel orthonormal frame. The remaining scalar constraint is quadratic in the leafwise part of the second fundamental form. Lemma \ref{lemma:leafwiseaveraging} shows that the weighted integral pairing of two leafwise closed one-forms depends only on their leafwise parallel parts, obtained by $L^2$-orthogonal projection. The decomposition follows from the $L^2$-orthogonal projection onto the closure of the space of leafwise exact forms. The flat Weitzenböck identity, applied after smoothing along the leaves by the induced $\mathbb{R}^2$-action, then shows that the remainder has vanishing tangential derivatives. Since the argument takes place in $L^2(\Sigma,\omega)$ and does not integrate over an individual leaf, it applies equally when the leaves are dense. Applied to the rows of the leafwise directions of the second fundamental form, which are leafwise closed by the Codazzi equations, Lemma \ref{lemma:leafwiseaveraging} shows that the leafwise average of its determinant equals the determinant of its leafwise average. Integrated over the Cauchy hypersurface against a suitable positive weight, the scalar constraint becomes a sum of squares, which forces both the flux and the averaged tensor to vanish. Finally, the coefficients encoding the curvature form a symmetric, trace-free matrix with leafwise closed rows and vanishing leafwise average. A second application of Lemma \ref{lemma:leafwiseaveraging} shows that the integral of its determinant vanishes, which for such a matrix forces it to vanish, implying that the development is flat.

\begin{table}[htbp]
\centering
\small
\renewcommand{\arraystretch}{1.2}
\begin{tabular}{|>{\raggedright\arraybackslash}p{0.23\textwidth}|>{\raggedright\arraybackslash}p{\dimexpr0.77\textwidth-4\tabcolsep-3\arrayrulewidth\relax}|}
\hline
\textbf{Result} & \textbf{Description} \\
\hline
Theorem \ref{thm:wellposedFinalSkew}
&
Establishes existence and uniqueness of the germ of a globally hyperbolic development carrying a skew-torsion parallel spinor, for initial data satisfying an explicit exterior differential system.
\\
\hline
Corollary \ref{cor:timeindependentdevelopment}
&
Provides an explicit development with time-independent metric and torsion. Completeness yields a maximal globally hyperbolic development on $\mathbb{R}\times\Sigma$.
\\
\hline
Corollary \ref{cor:initialdataskew}
&
Expresses the initial constraints in terms of a coframe and a function, with contact $(r\neq 0)$ and integrable $(r=0)$ branches, and reconstructs the torsion.
\\
\hline
Corollary \ref{cor:compactconstraintbranches}
&
Every connected, closed and oriented three-manifold admits contact-branch initial data with nonzero torsion. Integrable-branch initial data exist precisely on oriented two-torus bundles over the circle.
\\
\hline
Theorem \ref{thm:compatibilityflows}
&
Establishes an identity among the Hamiltonian, momentum, and Maxwell constraints evaluated on the Neveu--Schwarz supersymmetry flow. This identity is crucial to clarify the compatibility between the supersymmetric Neveu-Schwarz flow and its second-order counterpart.
\\
\hline
Theorem \ref{thm:wellposedFinalSusy}
&
Solves the supersymmetric Cauchy problem for Neveu--Schwarz supergravity on an abelian bundle gerbe, with uniqueness of the metric, dilaton, and flux, up to automorphisms fixing the initial data and Cauchy hypersurface.
\\
\hline
Corollary \ref{cor:NSNSsusyinitialdataII}
&
Initial data characterization of supersymmetric Neveu--Schwarz solutions by an exterior differential system combining the supersymmetry constraints with certain components of the Maxwell and momentum constraints.
\\
\hline
Theorem \ref{thm:compactrigidity}
&
Every globally hyperbolic supersymmetric solution with compact Cauchy hypersurface is flat and flux-less, with parallel dilaton curvature.
\\
\hline
\end{tabular} 
\end{table}
 
\begin{figure}[htbp]
\centering
\begin{tikzpicture}[
    x=\textwidth,
    y=1cm,
    result/.style={
        draw,
        rounded corners=2pt,
        align=center,
        text width=0.26\textwidth,
        minimum height=1.45cm,
        inner sep=5pt,
        font=\small
    },
    flow/.style={
        -{Stealth[length=2mm]},
        semithick,
        rounded corners=3pt
    }
]

\node[result] (skew) at (-0.34,0) {
    \textbf{Cauchy problem for skew-torsion parallel spinors}\\[4pt]
    {\footnotesize Theorem~\ref{thm:wellposedFinalSkew}}
};

\node[result] (supergravity) at (0.34,0) {
    \textbf{Neveu--Schwarz supergravity}
};

\node[result] (explicit) at (-0.34,-2.8) {
    \textbf{Explicit developments}\\[4pt]
    {\footnotesize Corollary~\ref{cor:timeindependentdevelopment}}
};

\node[result] (susy) at (0,-2.8) {
    \textbf{Cauchy problem for supersymmetric configurations}\\[4pt]
    {\footnotesize Theorem~\ref{thm:wellposedFinalSusy}}
};

\node[result] (constraints) at (0.34,-2.8) {
    \textbf{Constraint identities and propagation}\\[4pt]
    {\footnotesize
    Theorems~\ref{thm:compatibilityflows}
    and~\ref{thm:susyflowcompatibility}}
};

\node[result] (topology) at (-0.34,-5.6) {
    \textbf{Initial data and compact topology}\\[4pt]
    {\footnotesize
    Corollaries~\ref{cor:initialdataskew}
    and~\ref{cor:compactconstraintbranches}}
};

\node[result] (solutions) at (0,-5.6) {
    \textbf{Initial data for supersymmetric Neveu--Schwarz solutions}\\[4pt]
    {\footnotesize Corollary~\ref{cor:NSNSsusyinitialdataII}}
};

\node[result] (rigidity) at (0,-8.1) {
    \textbf{Compact rigidity}\\[4pt]
    {\footnotesize Theorem~\ref{thm:compactrigidity}}
};

\draw[flow] (skew.south) -- (explicit.north);
\draw[flow] (skew.east) -| (susy.north);
\draw[flow] (skew.east) -- (-0.17,0) |- (topology.east);

\draw[flow] (supergravity.south) -- (constraints.north);

\draw[flow] (susy.south) -- (solutions.north);
\draw[flow] (constraints.south) |- (solutions.east);

\draw[flow] (solutions.south) -- (rigidity.north);

\end{tikzpicture}
\caption{Conceptual structure of the paper and principal dependencies.}
\label{fig:paperstructure}
\end{figure}
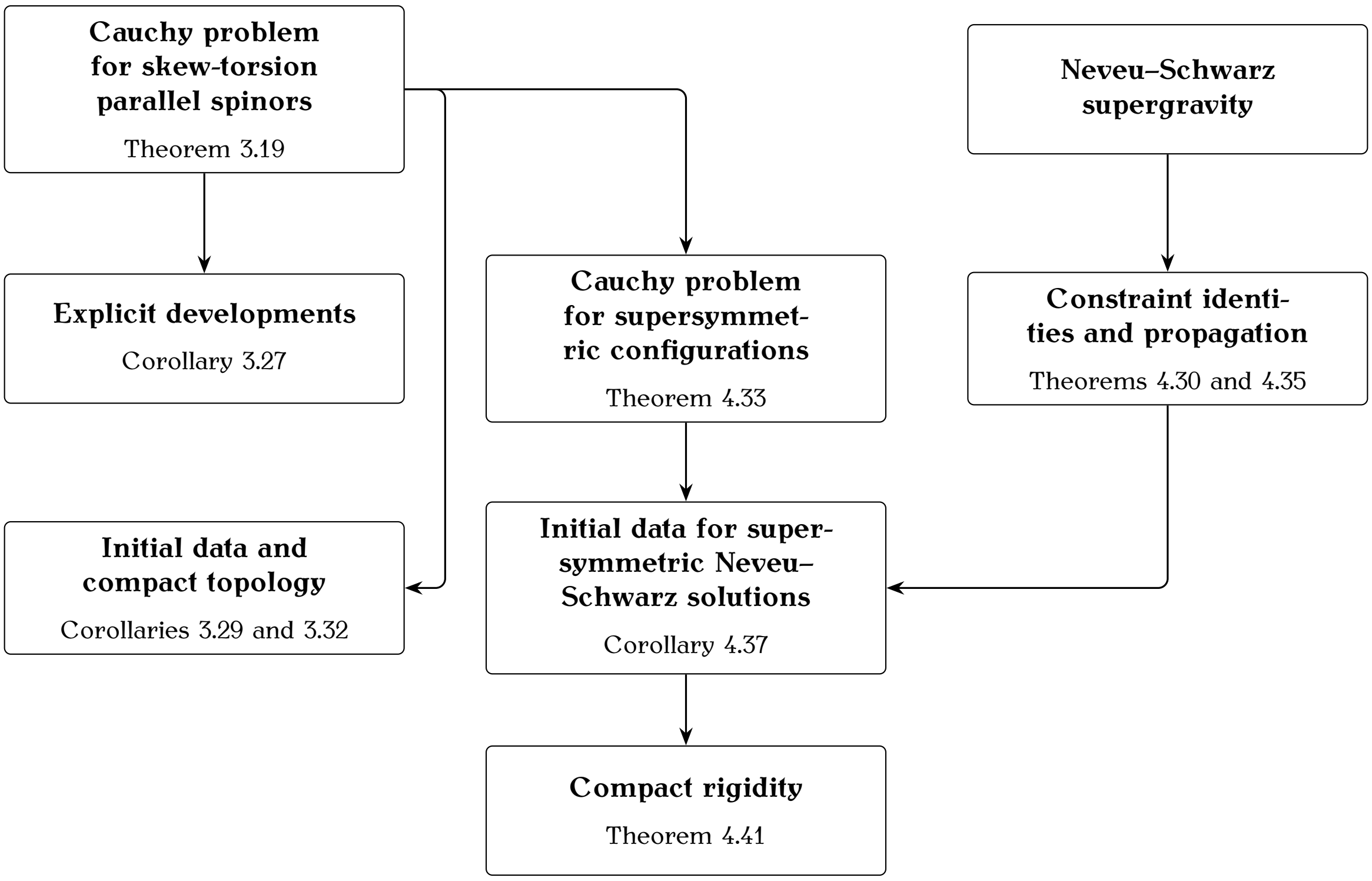
 
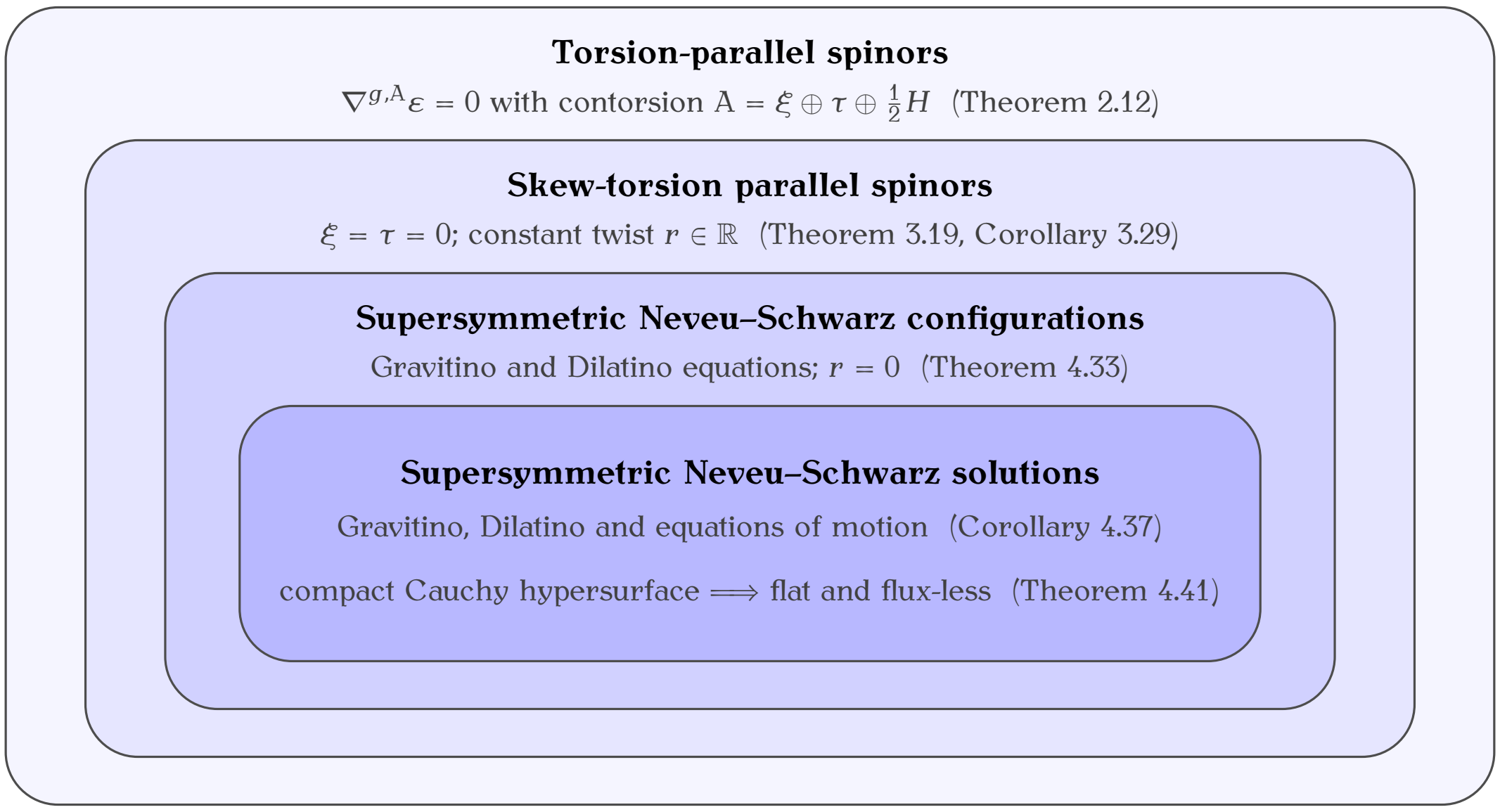
\begin{figure}[htbp]
\centering
\begin{tikzpicture}[
    x=1cm, y=1cm,
    ring/.style={draw=black!70, semithick, rounded corners=14pt},
    lab/.style={align=center, font=\small, inner sep=1pt},
    cond/.style={align=center, font=\footnotesize, inner sep=1pt, text=black!75}
]
\filldraw[ring, fill=blue!4]  (-7.0,-3.75) rectangle (7.0,3.75);
\filldraw[ring, fill=blue!10] (-6.25,-3.3)  rectangle (6.25,2.5);
\filldraw[ring, fill=blue!18] (-5.5,-2.85) rectangle (5.5,1.25);
\filldraw[ring, fill=blue!28] (-4.8,-2.4)  rectangle (4.8,0.0);

\node[lab]  at (0,3.3)  {\textbf{Torsion-parallel spinors}};
\node[cond] at (0,2.85) {$\nabla^{g,\mA}\varepsilon = 0$ with contorsion $\mA = \xi\oplus\tau\oplus\tfrac{1}{2}H$ \ (Theorem~\ref{thm:existencenullcoframeII})};

\node[lab]  at (0,2.05) {\textbf{Skew-torsion parallel spinors}};
\node[cond] at (0,1.6)  {$\xi = \tau = 0$; constant twist $r\in\mathbb{R}$ \ (Theorem~\ref{thm:wellposedFinalSkew}, Corollary~\ref{cor:initialdataskew})};

\node[lab]  at (0,0.8)  {\textbf{Supersymmetric Neveu--Schwarz configurations}};
\node[cond] at (0,0.35) {Gravitino and Dilatino equations; $r = 0$ \ (Theorem~\ref{thm:wellposedFinalSusy})};

\node[lab]  at (0,-0.65) {\textbf{Supersymmetric Neveu--Schwarz solutions}};
\node[cond] at (0,-1.15) {Gravitino, Dilatino and equations of motion \ (Corollary~\ref{cor:NSNSsusyinitialdataII})};
\node[cond] at (0,-1.75) {compact Cauchy hypersurface $\Longrightarrow$ flat and flux-less \ (Theorem~\ref{thm:compactrigidity})};
\end{tikzpicture}
\caption{The four classes of four-dimensional parallel spinors considered in this article, ordered by inclusion. Each class is cut out of the previous one by a condition on the torsion, indicated together with the results of the article that concern it.}
\label{fig:nestedclasses}
\end{figure}


\subsection*{Relation to prior work} 


The results of Section \ref{sec:CauchyKillingSkew} include as a particular case the resolution of the Cauchy problem for a Levi-Civita parallel spinor of real type in four dimensions. In arbitrary dimension, Baum, Leistner, and Lischewski first solved this problem in the real-analytic category \cite{BaumLeistnerLischewski} using the Cauchy--Kovalevskaya theorem. In the smooth category, Lischewski \cite{Lischewski} proved that every Riemannian manifold carrying an imaginary generalized Killing spinor whose Dirac current satisfies an additional algebraic constraint embeds as a Cauchy hypersurface into a Lorentzian manifold on which the spinor extends to a parallel null spinor by solving a first-order symmetric hyperbolic system relative to a background metric. Leistner and Lischewski \cite{LeistnerLischewski} subsequently proved the well-posedness of the Cauchy problem for a parallel null vector field, through hyperbolic evolution equations in which the Ricci tensor appears as an additional unknown. Together, these works establish the well-posedness of the Cauchy problem for parallel null spinors in arbitrary dimension. In contrast to the methods used in these works, our integration procedure is closely related to recent constructions in the parallel and Killing cases. In the parallel case, Ammann, Kr\"oncke, and M\"uller \cite{AmmannKronckeMuller} showed how to solve the constraint equations of \cite{BaumLeistnerLischewski} in general dimension. Gl\"ockle \cite{Glockle} proved that, under a local rigidity assumption on the initial data, this development is locally unique among spacetimes satisfying the dominant energy condition, and Ammann, Gl\"ockle, and Kr\"oncke \cite{AmmannGlockleKroncke} used it to parametrize pp-waves with compact codimension-two hypersurfaces, extending the spinor by parallel transport.   \medskip

In the Killing case, Hirsch and Zhang \cite{HirschZhang} have recently characterized the asymptotically anti-de Sitter spin initial data sets that saturate a \emph{BPS bound}. A key step in their analysis \cite[Theorem 6.2]{HirschZhang} shows that an initial data set carrying a spinor that satisfies the restriction of the imaginary Killing spinor equation to a spacelike hypersurface embeds into its Killing development, whose lapse and shift are given by the Dirac current of the spinor, and that the spinor extends to an imaginary Killing spinor on it, see also \cite{HirschZhangpp} for the asymptotically flat case, in which initial data sets with vanishing mass are contained in pp-wave spacetimes. There, the spinor is extended off the initial hypersurface by an ordinary differential equation generated by a hidden symmetry of the imaginary Killing spinor equation \cite[Theorem 6.1]{HirschZhang}, which becomes trivial when the spinor is null, a case that leads in their setting to Siklos waves \cite[Theorem 8.11]{HirschZhang}. Our procedure exhibits a similar structure in the presence of skew-torsion, which introduces several novelties, including an obstruction to the existence of a globally hyperbolic development.   \medskip

The second-order Cauchy problem posed by the equations of motion of Neveu-Schwarz supergravity was solved only recently by Schiller \cite{Schiller:2026vls} in every dimension, see \cite{CB85} for a pioneering study of the Cauchy problem in eleven-dimensional supergravity. The starting point in \cite{Schiller:2026vls} is the generalized Einstein equations for a generalized metric on an exact Courant algebroid \cite{GarciaFernandezStreets}, in the case of closed divergence. This generalized Einstein condition leads to the Neveu-Schwarz equations for a Lorentzian metric, a closed three-form in a fixed cohomology class, and a closed one-form that plays the role of the dilaton curvature. The proof of well-posedness combines a DeTurck gauge for the metric with a generalized Lorenz gauge for the $b$-field, which turns the system into a system of wave equations, and follows the approach of Ringstr\"om \cite{Ringstrom} to the Einstein--scalar field system, which builds on the classical results of Choquet-Bruhat and Geroch \cite{ChoquetBruhatGeroch}. Since the equations involve the $b$-field and the dilaton only through their curvatures, and the divergence is only required to be closed, Schiller's theorem applies, at the level of the metric, the $b$-field curvature, and the dilaton curvature, to the Neveu-Schwarz system on an abelian bundle gerbe with a circle-valued dilaton considered in this article. Hence, his result and ours are complementary. His theorem concerns the second-order equations of motion and is proven by hyperbolic methods, whereas ours concerns the first-order supersymmetry conditions, whose evolution we integrate explicitly. Theorem \ref{thm:susyflowcompatibility} connects the two: on initial data admissible to both problems, the development of Theorem \ref{thm:wellposedFinalSusy} solves the equations of motion, so that, by the uniqueness part of Schiller's theorem, it coincides near the Cauchy hypersurface with the development of the second-order Cauchy problem. This allows us to characterize supersymmetric solutions purely in terms of initial data in Corollary \ref{cor:NSNSsusyinitialdataII}.\medskip

The present work continues the study of spinorial evolution flows on globally hyperbolic Lorentzian four-manifolds initiated in \cite{Murcia:2020zig,Murcia:2021dur,MurciaKundt}. The four-dimensional case considered in this article constitutes a preliminary step toward studying the general supersymmetric Cauchy problem in Heterotic and Type-II supergravity in ten dimensions, which is the natural dimension for these theories.\medskip

\noindent
Throughout this article, Lorentzian metrics have signature $(-+++)$, and we identify vectors and one-forms with the metric whenever no confusion arises. The Clifford algebra is defined with the plus-sign convention. 


\subsection*{Acknowledgments}


This work is based on my doctoral thesis at the University of Hamburg, supervised by Professor Vicente Cortés and defended in April 2025. I was supported by the Leonardo grant LEO22-2-2155 of the BBVA Foundation and by the research grant PID2023-152822NBI00 of the Ministry of Science of the Government of Spain.


\section{Torsion-parallel spinors and isotropic parallelisms}
\label{sec:torsionparallel}


In this section we develop the general theory of irreducible real spinors parallel with respect to a metric connection with torsion, in terms of the equivalent differential system satisfied by their associated \emph{square} differential forms. In this dimension and signature, this leads to an explicit exterior differential system for \emph{isotropic parallelisms}, which we introduce below.  We recall only the necessary background for this purpose, referring to \cite{Cortes:2019xmk,ShahbaziTesis} for the general theory of spinorial forms.  


\subsection{Torsion-parallel spinors}
\label{sec:metric-torsion-spinors}


Let $(M,g)$ be an oriented and time-oriented Lorentzian four-manifold and let $\nabla^g$ denote its Levi--Civita connection.  Every metric connection on $TM$ can be written as follows:
\begin{equation*} 
\nabla^{g,\mA}_X Y = \nabla^g_X Y + \mA_X Y
\end{equation*}

\noindent
for a unique tensor $\mA\in\Omega^1(M,\End_g(TM))$ called the \emph{contorsion tensor} of $\nabla^{g,\mA}$ on $(M,g)$. Here $\End_g(TM)$ denotes the skew-symmetric endomorphisms of the Euclidean vector bundle $(TM,g)$. Thus, for every vector field $X\in \mathfrak{X}(M)$, the endomorphism $\mA_X$ is skew-adjoint with respect to $g$, that is:
\begin{equation*}
g(\mA_XY,Z)+g(Y,\mA_XZ)=0\, , \qquad Y , Z \in \mathfrak{X}(M)
\end{equation*}

\noindent
Using $g$, we will identify at our convenience $\End_g(TM)$ with $\wedge^2T^*M$, and hence we will occasionally regard $\mA$ as an element in $\Omega^1(M,\wedge^2T^*M) = \Gamma(\wedge^1 M \otimes \wedge^2 M)$.  The torsion of $\nabla^{g,\mA}$ is then given by the familiar expression:
\begin{equation*} 
\mT(X,Y)= \nabla^{g,\mA}_X Y - \nabla^{g,\mA}_Y X - [X,Y] = \mA_X Y-\mA_Y X.
\end{equation*}

\noindent
For a connection compatible with a given metric, passage from contorsion to torsion is invertible. More explicitly, if $\mT(X,Y,Z) : = g(\mT(X,Y),Z)$, then:
\begin{equation*} 
2g(\mA_XY,Z)=\mT(X,Y,Z)-\mT(Y,Z,X)+\mT(Z,X,Y).
\end{equation*}

\noindent
In particular, a metric connection is uniquely determined by its torsion. The $\mathrm{O}(3,1)$-module $\wedge^1 M \otimes \wedge^2 M$ has the following standard irreducible decomposition:
\begin{equation*} 
\wedge^1 M \otimes \wedge^2 M = \wedge^1 M \oplus\cT M\oplus\wedge^3 M,
\end{equation*}

\noindent
where $\cT M$ is called the \emph{twistorial component} of the contorsion, which is defined as follows. A tensor $\tau \in \wedge^1 M \otimes \wedge^2 M $ belongs to $\cT M$ if and only if it satisfies the following algebraic conditions for every $X, Y, Z \in TM$: 
\begin{equation*} 
\tau(X, Y, Z) = -\tau(X, Z, Y)\, , \quad \tau(X, Y, Z) + \tau(Y, Z, X) + \tau(Z, X, Y) = 0 \, , \quad \Tr_g(\tau(- , - , Z)) = 0
\end{equation*}
    
\noindent
Given a contorsion tensor $\mathrm{A} \in \Gamma(\wedge^1 M \otimes \wedge^2 M)$, there are unique tensors:
\begin{equation*}
\xi\in\Omega^1(M),\qquad \tau\in\Gamma(\cT M),\qquad H\in\Omega^3(M)
\end{equation*}

\noindent
such that:
\begin{equation}
\label{eq:contorsioncomponents}
\mA_XY = g(X,Y)\xi^{\sharp_g}-\xi(Y)X+\tau_X Y + \frac12 H(X,Y)^{\sharp_g}.
\end{equation}

\noindent
For ease of notation, we will write this decomposition succinctly as:
\begin{equation*}
\mA = \xi \oplus \tau \oplus \frac{1}{2} H    
\end{equation*}

\noindent
With this convention, the totally skew-symmetric component of the torsion is precisely $H$. For later use, if $\beta\in\Omega^1(M)$ we set:
\begin{equation}
\label{eq:Htaubeta} 
H_\beta:=\iota_{\beta^{\sharp_g}}H, \qquad \tau_\beta(X,Y) := g (\tau_Y X -\tau_X Y,\beta^{\sharp_g} ).
\end{equation}  

\noindent
For further reference, we note that our curvature convention is:
\begin{equation*} 
\cR^{g,\mA}_{X,Y}Z : = \nabla^{g,\mA}_X \nabla^{g,\mA}_Y Z - \nabla^{g,\mA}_Y \nabla^{g,\mA}_X Z - \nabla^{g,\mA}_{[X,Y]} Z.
\end{equation*}

\noindent
In terms of the contorsion tensor $\mA \in  \Gamma(T^*M\otimes\wedge^2 M)$, this equation exapnds as follows:
\begin{equation} 
\cR^{g,\mA}_{X,Y} = \cR^o_{x,y} + (\nabla^g_X\mA)_Y - (\nabla^g_Y\mA)_X + [\mA_X,\mA_Y].
\end{equation}

\noindent
where $\cR^g$ is the Riemann curvature tensor of $g$. \medskip

\noindent
Let $(M,g)$ be a Lorentzian four-manifold. A bundle of irreducible spinors on $(M,g)$ will be denoted by $(S ,\Gamma)$, where $S$ is a rank-four real vector bundle on $(M,g)$ and $\Gamma \colon \Cl(M,g) \to \End(S)$ denotes Clifford multiplication, which we present as a morphism of bundles of unital and associative algebras. We equip such a spinor bundle $(S,\Gamma)$ with a skew-symmetric admissible pairing $\cB$ of negative adjoint type, as explained in \cite{Cortes:2019xmk}. By the results of \cite{Lazaroiu:2016vov}, it follows that such an irreducible real spinor bundle $(S,\Gamma)$ exists if and only if it is associated to a spin structure on $(M,g)$. Furthermore, $\cB$ is guaranteed to exist if such a spin structure reduces to its identity component. We will refer to triples of the form $(S,\Gamma,\cB)$ as \emph{paired spinor bundles}. Since $(S,\Gamma)$ is necessarily associated to a spin structure, the metric connection $\nabla^{g,\mA}$ has a canonical lift to $S$ that we denote by the same symbol for ease of notation. We can write:
\begin{equation*} 
\nabla^{g,\mA}_X\varepsilon = \nabla^g_X\varepsilon-\frac12\mA_X\cdot\varepsilon := \nabla^g_X\varepsilon-\frac12 \Gamma(\mA_X) \varepsilon 
\end{equation*}
where $\mA_X$ is regarded as a two-form for every $X\in \mathfrak{X}(M)$. 

\begin{definition} 
\label{def:torsionparallel}
A \emph{torsion parallel spinor} on $M$ with contorsion $\mA \in \Gamma(\wedge^1 M \otimes \wedge^2 M)$ is a pair $(g , \varepsilon)$ consisting of a Lorentzian metric $g$ on $M$ and a non-zero section $\varepsilon\in \Gamma(S)$ of an irreducible paired spinor bundle $(S,\Gamma,\cB)$ on $(M,g)$ satisfying $\nabla^{g,\mA}_X\varepsilon = 0$.
\end{definition}

\noindent
In this definition, neither the Lorentzian metric nor spin structure is fixed in advance, a fact that is crucial for supergravity applications.  


\subsection{Nowhere vanishing spinors and isotropic parallelisms}

 
In this subsection, we introduce \emph{isotropic parallelism} as an equivalent formulation of the square of a real, irreducible, nowhere-vanishing spinor, particularly well adapted to studying first-order parallelicity conditions. We will closely follow \cite{ShahbaziTesis}, in which the theory is developed in detail. The square of a nowhere-vanishing irreducible real spinor $\varepsilon\in\Gamma(S)$ is given by:
\begin{equation*}
\alpha_\varepsilon=u+u\wedge l \in \Gamma(\wedge M),
\end{equation*}

\noindent
where $u$ is a nowhere-vanishing null one-form and $l$ is a unit spacelike one-form orthogonal to $u$.  This is the four-dimensional specialization of \cite[Theorem~3.34]{Cortes:2019xmk}, where the result is explicitly proven.  Note that replacing $l$ by $l+fu$, for any $f\in C^\infty(M)$, leaves $\alpha_\varepsilon$ unchanged, whence $l$ is not uniquely determined by $\varepsilon$. This leads us to introduce the following definition.

\begin{definition} 
A \emph{parabolic pair} on $(M,g)$ is a pair $(u,[l]_u)$, where $u$ is a nowhere-vanishing null one-form and $[l]_u$ is the equivalence class determined by a unit-norm one-form $l\in \Omega^1(M)$ orthogonal to $u$, under the equivalence relation $l\sim l+fu$.
\end{definition}

\noindent
Note that every representative $l^{\prime}\in [l]_u$ satisfies:
\begin{equation*} 
g(u,l^{\prime})=0,\qquad g(l^{\prime},l^{\prime})=1.
\end{equation*}

\noindent
as a direct consequence of $u$ being istropic and $l$ being of unit-norm and orthogonal to $u$. The nowhere vanishing null one-form $u$ defines the rank-two screen bundle in the usual way:
\begin{equation*} 
\frG_u:= \frac{\mathrm{Ker}(u)}{\langle u^{\sharp_g}\rangle},
\end{equation*}
where $\mathrm{Ker}(u)\subset T^*M$ denotes the distribution defined as the kernel of $u$.  The Lorentzian metric descends to a Riemannian metric $q_u$ on $\frG_u$.  Orientation and time orientation induce an orientation on $\frG_u$, whence Hodge duality is well defined for sections of $\frG_u$. In particular:
\begin{equation*}
[n]_u:=\ast_{q_u}[l]_u
\end{equation*}
is a second globally defined unit section. Hence the screen bundle associated with a parabolic pair is canonically trivialized by $([l]_u,[n]_u)$.\medskip

\noindent
A one-form $v$ is said to be \emph{conjugate} to $u$ if it is null and satisfies $g(u,v)=1$.  It can be seen that conjugate one-forms exist globally.  For example, if $\frt \in \Omega^1(M)$ is a nowhere vanishing timelike one-form, then the one-form:
\begin{equation*}
v=\frac{\frt}{g(u,\frt)}-\frac{g(\frt,\frt)}{2g(u,\frt)^2}u
\end{equation*}
is conjugate to $u$.  Once $v$ is chosen, there are unique representatives $l\in[l]_u$ and $n\in[n]_u$ orthogonal to $v$ for which $(u,v,l,n)$ is positively oriented.  They form a global isotropic coframe in terms of which the metric reads:
\begin{equation}
\label{eq:metricfromisotropiccoframe}
g=u\odot v+l\otimes l+n\otimes n.
\end{equation}

\noindent
This discussion leads to the following correspondence.

\begin{prop}
Let $(u,[l]_u)$ be a parabolic pair on $(M,g)$. Every choice of isotropic one-form $v\in \Omega^1(M)$ conjugate to $u$ canonically determines a positively oriented coframe $(u,v,l,n)$ on $M$, where $l\in [l]_u$ and $n\in [n]_u$ are orthogonal unit vectors uniquely determined by the orthogonality condition $l,n\in ( \langle \mathbb{R}\, u\rangle  \oplus \langle \mathbb{R}\, v \rangle)^{\perp_g}$ and the positive orientation.
\end{prop}

\noindent
We will refer to the global coframes $(u,v,l,n)$ occurring in the previous proposition as \emph{isotropic coframes} associated to $(u,[l]_u)$ or its corresponding spinor. Given a nowhere vanishing isotropic one-form $u\in \Omega^1(M)$, conjugate one-forms are abundant. Let $v, v^{\prime} \in \Omega^1(M)$ be both conjugate to $u$. A calculation gives:
\begin{equation*}
v^{\prime} = v - \frac{1}{2} \vert \frw\vert_g^2 u + \frw
\end{equation*} 

\noindent
for a unique one-form $\frw\in\Omega^1(M)$ orthogonal to both $u$ and $v$. Elaborating on this, we arrive at the following result, proven in detail in \cite{ShahbaziTesis}.

\begin{prop}
\label{prop:primerparallelizationglobal}
Let $(u,[l]_u)$ be a parabolic pair on $(M,g)$. Given isotropic coframes $(u,v,l,n)$ and $(u,v^{\prime},l^{\prime},n^{\prime})$ associated to $(u,[l]_u)$, there exists a unique one-form $\frw\in \Gamma( \langle \mathbb{R}\, u\rangle  \oplus \langle \mathbb{R}\, v \rangle)^{\perp_h}$ such that:
\begin{equation}
\label{eq:equivalencerelation}
(u,v^{\prime},l^{\prime},n^{\prime}) = (u,v-\frac{1}{2} \vert\frw\vert_g^2 u + \frw,l - \frw(l) u,n - \frw(n) u)
\end{equation}
	
\noindent
In particular, the set of isotropic coframes associated to $(u,[l]_u)$ is a torsor over the infinite-dimensional vector space $\Gamma( \langle \mathbb{R}\, u\rangle  \oplus \langle \mathbb{R}\, v \rangle)^{\perp_g}$.
\end{prop}   

\noindent
A choice of one-form $v\in\Omega^1(M)$ conjugate to a given parabolic pair $(u,[l]_u)$ defines a \emph{realization} of the abstract screen bundle $\mathfrak{G}_u$ as the codimension-two distribution in $T^{\ast}M$ spanned by $l$ and $n$. That is:
\begin{eqnarray*}
\mathfrak{G}_u \xrightarrow{\simeq} ( \langle \mathbb{R}\, u\rangle  \oplus \langle \mathbb{R}\, v \rangle)^{\perp_g} 
\end{eqnarray*}

\noindent
The following definition identifies the intrinsic characterization of a pair $(g,\varepsilon)$ consisting of a Lorentzian metric $g$ and a nowhere-vanishing irreducible real spinor $\varepsilon$, which we will use throughout this article.

\begin{definition}
An \emph{isotropic parallelism} is an equivalence class $[u,v,l,n]$ of isotropic coframes on $M$ under the equivalence relation defined by Equation \eqref{eq:equivalencerelation}.
\end{definition}  
 
\noindent
Denote by $\mathfrak{F}(M)$ the category of isotropic parallelisms on $M$, where arrows are diffeomorphisms of $M$ mapping isotropic parallelisms to each other. We have a natural functor:
\begin{equation*}
\mathfrak{F}(M) \to \mathrm{Lor}(M)\, , \qquad [u,v,l,n] \mapsto g = u\odot v + l\otimes l + n\otimes n
\end{equation*}

\noindent
where $\mathrm{Lor}(M)$ is considered as the category of Lorentzian metrics and isometries on $M$. Such $g$ is the Lorentzian metric canonically associated to $[u,v,l,n]$, and, as explained in \cite{ShahbaziTesis}, leads to the following result.

\begin{prop}
\label{prop:global-spinor-parabolic}
Let $M$ be an oriented four-manifold.  Every isotropic parallelism $[u,v,l,n]$ determines an isomorphism class of spin structures on $(M,g)$ and, on any representative of that class, a nowhere-vanishing irreducible real spinor, unique up to sign, whose signed square is:
\begin{equation*}
\alpha = u + u\wedge l    
\end{equation*}

\noindent
Conversely, a nowhere-vanishing irreducible real spinor on a time-oriented Lorentzian four-manifold $(M,g)$ determines an isotropic parallelism.   
\end{prop}

\begin{proof}
Let $[u,v,l,n]$ be an isotropic parallelism on $M$. By \cite[Theorem~3.34]{Cortes:2019xmk}, the polyform: 
\begin{equation*}
\alpha=u+u\wedge l \in \Gamma (\wedge^1 M \oplus \wedge^2 M)    
\end{equation*}

\noindent
pointwise satisfies the algebraic conditions for a nonzero signed spinor square with respect to the Lorentzian metric \eqref{eq:metricfromisotropiccoframe}. Since the parabolic pair gives a global coframe on $M$, it follows that $(M,g)$ is spin. Starting from any compatible spin structure, \cite[Proposition~4.22]{Cortes:2019xmk} selects the unique isomorphism class of spin structures for which $\alpha$ is the square of a global spinor.  The global reconstruction and the two-to-one character of the signed-square map give uniqueness up to sign \cite[Section~4.4 and Proposition~4.18]{Cortes:2019xmk}. The converse follows from the fiberwise square characterization, combined with the construction of the associated isotropic parallelism, as explained in \cite{ShahbaziTesis}.
\end{proof}

\noindent
As a consequence, we have the following correspondence, which highlights the convenience of considering the underlying metric as a variable when studying parallelicity equations for spinors. This is particularly fundamental for supergravity applications, and it is why we introduced torsion parallel spinors in Definition \ref{def:torsionparallel} as pairs $(g,\varepsilon)$ instead of considering $\varepsilon$ on a fixed Lorentzian manifold $(M,g)$.

\begin{prop}
There is a two-to-one correspondence between pairs $(g,\varepsilon)$ consisting of a time-oriented Lorentzian metric $g$ on $M$ and a nowhere vanishing irreducible real spinor $\varepsilon$ on $(M,g)$, and the set of isotropic parallelisms $[u,v,l,n]$ on $M$.
\end{prop}

\noindent
The \emph{two-to-one} qualifier in the previous proposition simply reflects that both $(g,\varepsilon)$ and $(g,-\varepsilon)$ map to the same isotropic parallelism. We recall the following formulas for the Hodge operator of \eqref{eq:metricfromisotropiccoframe}: 
\begin{eqnarray*}
& \ast_h u = -u\wedge l\wedge n \, , \qquad \ast_h v = v\wedge l\wedge n \, , \qquad \ast_h l = u\wedge v\wedge n \, , \qquad \ast_h n = - u\wedge v\wedge l\\
& \ast_h(u\wedge v) = - l\wedge n\, , \quad \ast_h(u\wedge l) = u\wedge n \, , \quad \ast_h(v\wedge l) = - v\wedge n\, , \\
& \ast_h \nu_g = -1\, , \quad \ast_h^2\vert_{\wedge^1 V^{\ast}} = 1\, , \quad \ast_h^2\vert_{\wedge^2 V^{\ast}} = - 1\, , \quad\iota_w \ast_h \beta = \ast_h(\beta\wedge w)  
\end{eqnarray*}

\noindent
for every $w\in \Omega^1(M)$ and every $\beta \in \Gamma(\wedge T^{\ast}M)$. These formulas will be extensively used in later computations.


\subsection{Torsion-parallel spinors as an exterior differential system}
\label{sec:torsion-spinor-characterization}


In this subsection, we translate the torsion-parallelity condition for a pair $(g,\varepsilon)$ into an equivalent exterior differential system for the associated isotropic parallelism. Our starting point is the following result.   

\begin{prop}
\label{prop:existencenullcoframegeneral}
An oriented four-manifold $M$ admits a torsion parallel spinor with contorsion $\mA \in \Gamma(\wedge^1 M \otimes \wedge^2 M)$ if and only if it admits an isotropic parallelism $[u,v,l,n]$ satisfying the following differential system:
\begin{eqnarray}
\label{eq:existencenullcoframegeneral}
\nabla^{g,\mA} u =0, \qquad  \nabla^{g,\mA} v =-\kappa\otimes l-\rho\otimes n,\qquad \nabla^{g,\mA} l =\kappa\otimes u, \qquad  \nabla^{g,\mA} n =\rho\otimes u ,
\end{eqnarray}

\noindent
for certain one-forms $\kappa,\rho\in\Omega^1(M)$, where $(u,v,l,n)\in [u,v,l,n]$ is a representative and $g$ is given in terms of $[u,v,l,n]$ as in Equation \eqref{eq:metricfromisotropiccoframe}. 
\end{prop}
 
\begin{proof}
Follows from \cite[Lemma~4.24, Theorem~4.26 and Remark~4.27]{Cortes:2019xmk}.
\end{proof}

\noindent
If a representative $(u,v,l,n) \in [u,v,l,n]$ satisfies the differential system \eqref{eq:existencenullcoframegeneral} for certain one-forms $\kappa,\rho \in \Omega^1(M)$, then, any other representative $(u,v^{\prime},l^{\prime},n^{\prime})\in [u,v,l,n]$ also satisfies the differential system \eqref{eq:existencenullcoframegeneral} for generally different one-forms $\kappa^{\prime},\rho^{\prime} \in \Omega^1(M)$. These can be characterized as follows.
 
\begin{lemma} 
If $(u,v,l,n)$ and $(u,v^{\prime},l^{\prime},n^{\prime})$ both satisfy Equation \eqref{eq:existencenullcoframegeneral} and are related as in Equation\eqref{eq:equivalencerelation}, then:
\begin{equation}
\label{eq:kapparhogauge}
\kappa'=\kappa-\dd (\frw(l)) ,\qquad \rho'=\rho-\dd (\frw(n)).
\end{equation}
where we have expanded $\frw = \frw(l^{\sharp})\, l + \frw(n^{\sharp})\, n = \frw(l)\, l + \frw(n)\, n $.
\end{lemma}

\noindent
Consequently, equations \eqref{eq:existencenullcoframegeneral} are intrinsic to the isotropic-parallelism class, even though the one-forms $\kappa$ and $\rho$ may depend on the representative. 

\begin{definition}
A \emph{torsion isotropic parallelism} with contorsion $\mA \in \Gamma(\wedge^1 M \otimes \wedge^2 M)$ is an isotropic parallelism satisfying equations \eqref{prop:existencenullcoframegeneral}.
\end{definition} 

\noindent
The following lemma will be used to translate the previous covariant system into an exterior differential system.

\begin{lemma}
\label{lemma:skewsymconnection}
Let $\mA = \xi \oplus \tau \oplus \frac12H$ as in Equation \eqref{eq:contorsioncomponents}.  For every $\beta\in\Omega^1(M)$, we have:
\begin{equation*} 
(\nabla^{g,\mA}_X\beta)(Y) - (\nabla^{g,\mA}_Y\beta)(X) =\dd\beta(X,Y)+(\beta\wedge\xi)(X,Y) + \tau_\beta(X,Y)-H_\beta(X,Y).
\end{equation*}
\end{lemma}

\begin{proof}
Follows from inserting \eqref{eq:contorsioncomponents} into $(\nabla^{g,\mA}_X\beta)(Y)=(\nabla^g_X\beta)(Y)-\beta(\mA_XY)$, alternating in $X,Y$, and using definitions \eqref{eq:Htaubeta}.
\end{proof}

\noindent
We arrive now at our final characterization of torsion parallel spinors in terms of an equivalent exterior differential system for their associated isotropic parallelisms.  

\begin{thm}
\label{thm:existencenullcoframeII}
An oriented four-manifold $M$ admits a torsion parallel spinor with contorsion $\mA = \xi \oplus \tau \oplus \frac{1}{2} H$ if and only if it admits an isotropic parallelism $[u,v,l,n]$ satisfying the following exterior differential system:
\begin{eqnarray}
& \dd u = H_{u} + \xi\wedge u - \tau_u  \, , \qquad \dd v = H_{v} + \xi\wedge v - \tau_v - \kappa\wedge l - \rho\wedge n   \label{eq:duvtorsiongeneral}\\
& \dd l = H_{l} + \xi\wedge l - \tau_l + \kappa\wedge u \, , \qquad \dd n = H_{n} + \xi\wedge n - \tau_n + \rho\wedge u   \label{eq:dlntorsiongeneral}
\end{eqnarray}

\noindent
for certain one-forms $\kappa,\rho\in\Omega^1(M)$, where $(u,v,l,n)\in [u,v,l,n]$ is a representative and the musical isomorphism is defined with respect to the Lorentzian metric $g$ given in Equation \eqref{eq:metricfromisotropiccoframe}.  
\end{thm}

\begin{proof}
Suppose first that $(g,\varepsilon)$ is a torsion-parallel spinor on $M$, and let $[u,v,l,n]$ be its associated torsion isotropic parallelism.  Choose a representative $(u,v,l,n) \in [u,v,l,n]$. Proposition \ref{prop:existencenullcoframegeneral} implies that $(u,v,l,n)$ satisfies all equations in \eqref{eq:existencenullcoframegeneral}. Alternating those four equations and using Lemma~\ref{lemma:skewsymconnection} gives all equations in \eqref{eq:duvtorsiongeneral} and \eqref{eq:dlntorsiongeneral}. \medskip

Conversely, suppose that the exterior system holds for a certain isotropic coframe $(u,v,l,n)$ and pair of one-forms $\kappa,\rho \in \Omega^1(M)$. We need to prove that all equations in \eqref{eq:existencenullcoframegeneral} are satisfied. This occurs if and only if their symmetric and skew-symmetric projections are satisfied. The skew-symmetric projection of the equations in \eqref{eq:existencenullcoframegeneral} yields \eqref{eq:duvtorsiongeneral} and \eqref{eq:dlntorsiongeneral}. On the other hand, the symmetrization of \eqref{eq:existencenullcoframegeneral} gives the following system of equations:
\begin{eqnarray*}
& \cL_{u} g = 2 \xi(u) g	- u\odot \xi  - \tau^s_u\, , \qquad \cL_v g = 2 \xi(v) g	- v\odot \xi - \tau^s_v - \kappa\odot l - \rho\odot n \\
& \cL_l g = 2 \xi(l) g	- l\odot \xi - \tau^s_l + \kappa\odot u\, , \qquad \cL_n g = 2 \xi(n) g	- n\odot \xi - \tau^s_n  + \rho\odot u 
\end{eqnarray*}

\noindent
where $\tau^s_{\beta}(w_1 , w_2) = \tau(w_1,\beta,w_2) + \tau(w_2,\beta,w_1)$, $w_1 , w_2 \in \mathfrak{X}(M)$. To obtain this symmetrization, we have used the following identity:
\begin{eqnarray*}
& (\nabla^{g,\mA}_{w_1} \beta) (w_2) + (\nabla^{g,\mA}_{w_2} \beta) (w_1) = \cL_{\beta} g (w_1,w_2) + (\beta \odot \xi)(w_1 , w_2) \\
& - 2 \beta(\xi) g(w_1,w_2)  + \tau (w_1,\beta,w_2) + \tau(w_2,\beta,w_1) 
\end{eqnarray*}

\noindent
applied to each of the elements in $(u,v,l,n)$. Explicitly computing using equations \eqref{eq:duvtorsiongeneral} and \eqref{eq:dlntorsiongeneral}, we obtain:
\begin{eqnarray*}
& \cL_{u} g = \dd u (u)\odot \dd v  + \dd u \odot \dd v(u) + l \odot \dd l (u) + n \odot \dd n (u)\\
& = (\xi(u) u - \tau(u,u)) \odot v + u\odot (H(v,u) + \xi(u) v - \xi -\tau(v,u) - \kappa(u) l - \rho(u) n) \\
& + l \odot (H(l,u) + \xi(u) l - \tau(l,u) + \kappa(u) u) + n \odot (H(n,u) + \xi(u) n - \tau(n,u) + \rho(u) u)\\
& =  2 \xi(u) g	- u\odot \xi  - \tau^s_u\, .
\end{eqnarray*}
	
\noindent
This yields the symmetrization of the first equation in \eqref{eq:existencenullcoframegeneral}. Computing similarly for $\cL_v g$, $\cL_l g$ and $\cL_n g$, we obtain the symmetrization of the remaining equations in \eqref{eq:existencenullcoframegeneral} and thus we conclude. 
\end{proof}

 \noindent
 Hence, torsion parallel spinors can be equivalently studied in terms of isotropic parallelisms satisfying equations \eqref{eq:duvtorsiongeneral} and \eqref{eq:dlntorsiongeneral}.


\subsection{Curvature and cohomological invariants}
\label{sec:generalcurvature}


Let $[u,v,l,n]$ be a torsion isotropic parallelism. The one-forms $\kappa,\rho\in \Omega^1(M)$ for which a given representative $(u,v,l,n)$ satisfies the system \eqref{eq:existencenullcoframegeneral} are uniquely determined by the representative. They depend, however, on the representative $(u,v,l,n)\in [u,v,l,n]$ as prescribed in Equation \eqref{eq:kapparhogauge}. Hence, every torsion isotropic parallelism, equivalently, every torsion parallel spinor, defines the following pair of intrinsically defined quantities:
\begin{eqnarray}
\label{eq:class}
\mathfrak{F}(M) \to \frac{\Omega^1(M)}{\dd C^{\infty}(M)}\times \frac{\Omega^1(M)}{\dd C^{\infty}(M)}\, , \qquad [u,v,l,n] \mapsto ([\kappa],[\rho]) 
\end{eqnarray}

\noindent
where $\mathfrak{F}(M)$ denotes the set of torsion isotropic parallelisms on $M$. 

\begin{definition}
Let $[u,v,l,n]$ be a torsion isotropic parallelism. Using the notation introduced above, the pair $[\mathfrak{o}] := [\kappa] \oplus [\rho]$ is the \emph{rank-one invariant} of $[u,v,l,n]$. 
\end{definition}

 \noindent
 Representatives $\mathfrak{o} \in [\mathfrak{o}]$ will be considered as $\mathbb{R}^2$-valued one-forms, that is $\mathfrak{o}\in \Omega^1(M,\mathbb{R}^2)$. Every class $[\mathfrak{o}]$ is natural under pullback, but for a general, nonclosed one-form it is not invariant under diffeomorphisms isotopic to the identity.  Indeed, if $\varphi_s$ is a family of isotopies generated by the family of vector fields $X_{s}\in \mathfrak{X}(M)$ for some real parameter $s$, then:
\begin{equation*}
\partial_s \varphi_s^{\ast}\mathfrak{o} = \varphi_s^{\ast} (\dd\iota_{X_s}\kappa+\iota_{X_s}\dd\kappa ),
\end{equation*}

\noindent
and consequently the second summand may not be exact. Hence, \eqref{eq:class} is \emph{gauge-invariant} data intrinsically associated with the given torsion isotropic parallelism and is equivariant under diffeomorphisms. In the flat case considered below, it becomes a genuine invariant on the quotient by diffeomorphisms isotopic to the identity.

\begin{prop}
\label{prop:curvaturetorsion}
Let $[u,v,l,n]$ be a torsion isotropic parallelism with contorsion $\mA\in \Gamma(\wedge^1 M \otimes \wedge^2 M)$ and let $(u,v,l,n)\in [u,v,l,n]$ be any representative.  Then:
\begin{equation}
\label{eq:torsioncurvaturegeneral}
\cR^{g,\mA}_{X,Y} =-\dd\kappa(X,Y)\,u\wedge l - \dd\rho(X,Y)\,u\wedge n
\end{equation}
for all $X,Y\in\mathfrak X(M)$.
\end{prop}

\begin{proof}
We apply the operator $(\nabla^{g,\mA}_X \nabla^{g,\mA}_Y - \nabla^{g,\mA}_Y \nabla^{g,\mA}_X - \nabla^{g,\mA}_{[X,Y]})$ to the four equations in \eqref{eq:existencenullcoframegeneral}.  The quadratic terms cancel, yielding:
\begin{eqnarray*}
\cR^D_{X,Y}u =0\, , \quad \cR^D_{X,Y}v =-\dd\kappa(X,Y)l-\dd\rho(X,Y)n\, , \quad \cR^D_{X,Y}l = \dd\kappa(X,Y)u\, , \quad \cR^D_{X,Y}n = \dd\rho(X,Y)u\, 
\end{eqnarray*}
These four expressions are equivalent to \eqref{eq:torsioncurvaturegeneral}.  Gauge invariance is also explicit:
$\dd\kappa$ and $\dd\rho$ are unchanged by \eqref{eq:kapparhogauge}, while $u\wedge l$ and $u\wedge n$ are unchanged
by \eqref{eq:equivalencerelation}.
\end{proof} 

\begin{definition}
\label{def:torsion-flat}
A torsion isotropic parallelism is \emph{torsion-flat} if its preserving metric connection $\nabla^{g,\mA}$ is flat, namely if $\cR^{g,\mA} = 0$.
\end{definition}

\noindent
Note that by \emph{flatness} we refer to flatness of the connection with torsion $\nabla^{g ,\mA}$, not to vanishing of the torsion tensor or the Riemann curvature tensor of $g$.

\begin{cor}
\label{cor:torsionflat}
For a torsion-parallel spinor $(g,\varepsilon)$, equivalently, its associated torsion isotropic parallelism $[u,v,l,n]$, the following conditions are equivalent:
\begin{enumerate}
\item $\nabla^{g,\mA}$ is flat;
\item $\dd\kappa=0$ and $\dd\rho=0$;
\item Its rank-one invariants \eqref{eq:class} descend to de Rham cohomology $[\mathfrak{o}] = [\kappa] \oplus [\rho] \in H^1 (M,\mathbb R^2)$, which define a cohomological invariant modulo isotopies.
\end{enumerate}
\end{cor}

\begin{proof}
The first two conditions are equivalent by Proposition~\ref{prop:curvaturetorsion}, because $u\wedge l$ and $u\wedge n$ are everywhere linearly independent.  The equivalence with the third follows directly from Equation \eqref{eq:kapparhogauge}.
\end{proof}  

\noindent
Let $[u,v,l,n]$ be a flat torsion isotropic parallelism. The cohomology class $[\mathfrak{o}]  \in H^1 (M,\mathbb R^2)$ records the holonomy of the flat connection $\nabla^{g,\mA}$, in the sense that if $[\mathfrak{o}] = 0$ is trivial, then $\nabla^{g,\mA}$ preserves a conformal transformation of every representative $(u,v,l,n) \in [u,v,l,n]$ and therefore has trivial holonomy. More generally, for every smooth map $f \colon S^1 \to M$ we have a \emph{Torelli period} type of map:
\begin{equation*}
f\longmapsto (\int_f\kappa,\int_f\rho)\in\mathbb R^2.
\end{equation*}
which vanishes precisely when $[\mathfrak{o}] = 0$.  If $\mathfrak F_o(M)$ denotes the set of flat torsion isotropic parallelisms, we obtain a map:
\begin{equation*} 
\mathfrak F_o(M)\longrightarrow H^1(M,\mathbb R^2),
\qquad [u,v,l,n]\longmapsto\mathfrak \mathfrak{o},
\end{equation*}
which is invariant under diffeomorphisms connected to the identity.\medskip

\noindent
Tracing Equation \eqref{eq:torsioncurvaturegeneral} gives the Ricci and scalar curvatures of $\nabla^{g,\mA}$ 
\begin{eqnarray*}
& \Ric^{g,\mA} (Y) = (\dd\kappa(Y,l)+\dd\rho(Y,n) )u -\dd\kappa(Y,u)l-\dd\rho(Y,u)n \\
&s^{g,\mA} = 2 (\dd\kappa(u,l)+\dd\rho(u,n) ) 
\end{eqnarray*} 

\noindent
We will use these equations in the following.


\section{The Cauchy problem for parallel skew-torsion spinors}
\label{sec:CauchyKillingSkew}


In this section, we pose and solve the Cauchy problem for skew-torsion parallel spinors on a globally hyperbolic Lorentzian four-manifold, considered as a constrained initial evolution problem for both the spinor and the underlying metric. To do so, we use the equivalent formulation in terms of isotropic parallelisms obtained in the previous section, which encodes both the underlying Lorentzian metric and the parallel spinor. We first establish a canonical form of a globally hyperbolic isotropic parallelism. We then decompose the skew-torsion equations into evolution equations, spatial constraints, and a cohomological condition. We then integrate the complete system explicitly via a deformation of the flow of the spatial projection of the Dirac current, obtained in an appropriately chosen gauge for the lapse function.  
 


\subsection{Globally hyperbolic isotropic parallelisms}
\label{sec:ReductionGloballyHyperbolicparallelisms}


Let $M$ be an oriented four-manifold and let $\mathfrak{F}(M)$ denote the category of isotropic parallelisms on $M$. As explained in the previous section, there is a natural functor:
\begin{equation*}
\mathfrak{F}(M) \to \mathrm{Lor}(M)\, , \qquad [u,v,l,n] \mapsto g = u\odot v + l\otimes l + n\otimes n
\end{equation*}

\noindent
where $\mathrm{Lor}(M)$ is considered as the category of Lorentzian metrics and isometries on $M$.  

\begin{definition}
\label{def:globallyhyperbolic}
An isotropic parallelism $[u,v,l,n]\in \mathfrak{F}(M)$ is \emph{globally hyperbolic} if its associated metric:
\begin{eqnarray*}
g = u\odot v + l\otimes l + n\otimes n
\end{eqnarray*}

\noindent
is a globally hyperbolic metric on $M$.
\end{definition}  

\noindent
We denote by $\mathfrak{F}_o(M)$ the full subcategory of globally hyperbolic isotropic parallelisms $\mathfrak{F}(M)$. Let $[u,v,l,n]$ be a globally hyperbolic isotropic parallelism on $M$ with associated globally hyperbolic metric $g$, and let $\Sigma \subset (M,g)$ be a fixed Cauchy hypersurface. A celebrated theorem of Bernal and S\'anchez \cite{Bernal:2003jb,Bernal:2004gm} states that in this case $(M,g)$ has the following isometry type:
\begin{equation}
\label{eq:globahyp}
(M,g) = (\cI\times \Sigma, -\lambda^2_t \dd t\otimes \dd t + h_t) 
\end{equation}

\noindent
where $t$ is the canonical coordinate of the interval $\cI \subset \mathbb{R}$ containing the origin, $\lambda_t$ is a smooth family of nowhere vanishing functions on $\Sigma$ parametrized by $t\in\cI$ and $h_t$ is a family of Riemannian metrics on $\Sigma$, again parametrized by $t\in\cI$. By the previous discussion, it is enough to fix the topological type of $M = \cI\times \Sigma$ for an oriented three-dimensional manifold $\Sigma$. Hence, from now on, we consider the identification \eqref{eq:globahyp} to be fixed, and we denote by $\mathfrak{F}_o(\Sigma)$ the set of globally hyperbolic isotropic parallelisms on $M = \cI\times \Sigma$. We set:
\begin{equation*}
\Sigma_t := \left\{ t\right\}\times \Sigma \hookrightarrow M\, , \qquad \Sigma := \left\{ 0\right\}\times \Sigma \hookrightarrow M 
\end{equation*}

\noindent
and define:
\begin{equation*}
\mathfrak{t}_t = \lambda_t\, \dd t 
\end{equation*}

\noindent
to be the past-directed unit time-like one-form orthogonal to $T^{\ast}\Sigma_t$ for every $t\in \cI$. We will consider $\Sigma\hookrightarrow M$, endowed with the induced Riemannian metric:
\begin{equation*}
h := h_0\vert_{T\Sigma\times T\Sigma} 
\end{equation*}

\noindent
to be the Cauchy hypersurface of $(M,g)$. The \emph{shape operator} or scalar second fundamental form $\Theta_t \in \Gamma (T^{\ast}\Sigma_t \odot T^{\ast}\Sigma_t)$ of the embedded manifold $\Sigma_t\hookrightarrow M$ is defined in the usual way as follows:
\begin{equation}
\label{eq:secondfundamentalform}
\Theta_t  := \nabla^g \mathfrak{t}_t\vert_{T\Sigma_t\times T\Sigma_t} 
\end{equation}  

\noindent
Using standard theory of immersed Riemannian manifolds, this definition can be shown to be equivalent to:
\begin{equation*}
\Theta_t = - \frac{1}{2\lambda_t} \partial_t h_t \in \Gamma(T^{\ast}\Sigma_t\odot T^{\ast}\Sigma_t)
\end{equation*}

\noindent
Moreover, it can be seen that:
\begin{equation*}
\nabla^g \alpha \vert_{T\Sigma_t\times TM} = \nabla^{h_t} \alpha + \Theta_t(\alpha^{\sharp})\otimes \mathfrak{t}_t\, , \qquad \forall\,\,\alpha\in \Omega^1(\Sigma_t) 
\end{equation*}

\noindent
where $\nabla^{h_t}$ denotes the Levi-Civita connection on $(\Sigma_t,h_t)$ and $\Theta_t(\alpha) := \Theta_t(\alpha^{\sharp})$ is by definition the evaluation of $\Theta_t$ on the metric dual of $\alpha$. For every representative  $(u,v,l,n)\in [u,v,l,n]\in \mathfrak{F}_o(M)$ we write:
\begin{equation}
\label{eq:hyperbolicsplittingframe}
u = u^o_t\, \mathfrak{t}_t + u^{\perp}_t\, , \qquad v = v^o_t\, \mathfrak{t}_t + v^{\perp}_t  \, , \qquad l = l^o_t\, \mathfrak{t}_t + l^{\perp}_t \, , \qquad n = n^o_t\, \mathfrak{t}_t + n^{\perp}_t 
\end{equation}

\noindent
where the superscript $\perp$ denotes orthogonal projection to $T^{\ast}\Sigma_t$ with respect to $g$ and where we have defined:
\begin{equation*}
u^o_t = - g(u,\mathfrak{t}_t)\, , \qquad v^o_t = - g(v,\mathfrak{t}_t)\, , \qquad l^o_t = - g(l,\mathfrak{t}_t)\, ,\qquad n^o_t = - g(n,\mathfrak{t}_t) 
\end{equation*}

\noindent
This splitting is not preserved within the equivalence class $[u,v,l,n]$. Instead, we have:
\begin{eqnarray*}
& \bar{v}^o_t = v^0_t -\frac{1}{2} \vert \frw\vert^2_g u^o_t + \frw^o_t\, , \quad \bar{v}^{\perp}_t = v^{\perp}_t - \frac{1}{2} \vert \frw\vert^2_g u^{\perp}_t + \frw^{\perp}_t\\
& \bar{l}_t^{o} = l^o_t - \frw(l) u^o_t\, , \quad  \bar{l}_t^{\perp} = l^{\perp}_t - \frw (l) u^{\perp}_t\, , \quad \bar{n}_t^{o} = n^o_t - \frw(n) u^o_t\, , \quad  \bar{n}_t^{\perp} = n^{\perp}_t - \frw (n) u^{\perp}_t
\end{eqnarray*}

\noindent
where we have set:
\begin{equation}
\label{eq:parametrizedw}
(u,\bar{v},\bar{l},\bar{n}) = \mathfrak{w}\cdot (u,v,l,n)  = (u,v-\frac{1}{2} \vert \frw\vert^2_g u + \frw,l - \frw (l) u,n - \frw (n) u)\, , \quad \frw\in \Gamma (\langle \mathbb{R} u\rangle \oplus \langle \mathbb{R} v\rangle)^{\perp_g}
\end{equation}
 
\noindent
and where we have split $(u,\bar{v},\bar{l},\bar{n})$ analogously to the splitting of $(u,v,l,n)$ in Equation \eqref{eq:hyperbolicsplittingframe}. Hence, every isotropic parallelism $[u,v,l,n]$ defines an equivalence class of functions and frames on $\Sigma$ through its globally hyperbolic reduction. Rather than considering this equivalence class, we focus on the privileged representative most convenient for studying torsion parallel spinors on globally hyperbolic Lorentzian four-manifolds. 

\begin{lemma}
\label{lemma:ghcoframe}
Let $[u,v,l,n]$ be a globally hyperbolic isotropic parallelism on $M$. Then, there exists a unique representative $(u,v,l,n)\in[u,v,l,n]$ such that:
\begin{equation*} 
u = u^o_t\, \mathfrak{t}_t + u^{\perp}_t\, , \qquad v = \frac{1}{2(u^o_t)^2}\,( -u^o_t \mathfrak{t}_t + u^{\perp}_t)  \, , \qquad l =   l^{\perp}_t \, , \qquad n = n^{\perp}_t 
\end{equation*}

\noindent
in the splitting \eqref{eq:hyperbolicsplittingframe}.
\end{lemma}

\begin{proof}
Let $[u,v,l,n]$ be a globally hyperbolic istropic parallelims, and let $(u,v,l,n)$ be a given representative. Any other representative $(u,\bar{v},\bar{l},\bar{n}) \in [u,v,l,n]$ can be written as in Equation \eqref{eq:parametrizedw} for a unique $\frw\in\Gamma (\langle \mathbb{R} u\rangle \oplus \langle \mathbb{R} v\rangle)^{\perp_g}$. We choose $\frw$ such that:
\begin{equation*}
0 = \bar{l}_t^{o} = l^o_t - \frw(l) u^o_t\, , \qquad 0 = \bar{n}_t^{o} = n^o_t - \frw(n) u^o_t
\end{equation*}

\noindent
This already determines $\frw$ uniquely in terms of $(u,v,l,n)$  as follows:
\begin{equation*}
\frw = \frac{1}{u^o_t} (l^o_t l + n^o_t n)
\end{equation*}

\noindent
Hence $(u^o_t)^2 \vert \frw \vert_g^2 =  (l^o_t)^2   + (n^o_t)^2$. Plugging these expressions for $\frw$ and its norm into the equation for $\bar{v}$ given in \eqref{eq:parametrizedw}, we obtain:
\begin{eqnarray}
\label{eq:splittingchange}
\bar{v}^o_t = v^o_t + \frac{1}{2  u^o_t} ((l^o_t)^2   + (n^o_t)^2)  \, , \qquad \bar{v}^{\perp}_t = v^{\perp}_t - \frac{1}{2 (u^o_t)^2}((l^o_t)^2   + (n^o_t)^2) u^{\perp}_t + \frac{1}{u^o_t} (l^o_t l^{\perp}_t + n^o_t n^{\perp}_t)
\end{eqnarray}

\noindent
On the other hand, by assumption we must have:
\begin{eqnarray*}
& g = u\odot \bar{v}  + \bar{l}\otimes \bar{l} + \bar{n}\otimes \bar{n}  = (2 u^o_t v^o_t + (l^o_t)^2 + (n^o_t)^2) \frt_t\otimes \frt_t + \frt_t \odot (v^o_t u^{\perp}_t + u^o_t v^{\perp}_t + l^o_t l^{\perp}_t + n^o_t n^{\perp}_t)\\
& + u^{\perp}_t \odot v^{\perp}_t + (u^o_t)^{-2} (1 + 2 u^o_t v^o_t) u^{\perp}_t\otimes u^{\perp}_t + (u^o_t)^{-1} u^{\perp}_t\odot (l^o_t l^{\perp}_t + n^o_t n^{\perp}_t) + l^{\perp}_t \otimes l^{\perp}_t + n^{\perp}_t \otimes n^{\perp}_t = - \frt_t\otimes \frt_t + h_t
\end{eqnarray*}

\noindent
and thus:
\begin{equation*}
2 u^o_t v^o_t + (l^o_t)^2 + (n^o_t)^2 = -1\, , \qquad v^o_t u^{\perp}_t + u^o_t v^{\perp}_t + l^o_t l^{\perp}_t + n^o_t n^{\perp}_t = 0
\end{equation*}

\noindent
Plugging this equation back into \eqref{eq:splittingchange}, we obtain:
\begin{equation*}
\bar{v}^o_t = -  \frac{1}{2  u^o_t} \, , \qquad \bar{v}^{\perp}_t =  \frac{1}{2 (u^o_t)^2}  u^{\perp}_t  
\end{equation*}

\noindent
and thus we conclude after relabeling the variables to drop the \emph{bar}. 
\end{proof}

\noindent
We will refer to the representative $(u,v,l,n)\in [u,v,l,n]$ of the previous lemma as the \emph{adapted} representative in $[u,v,l,n]$. Hence, every globally hyperbolic isotropic parallelism $[u,v,l,n]$ naturally defines a family of tuples:
\begin{equation*}
(\lambda_t , u^o_t , u^{\perp}_t , l^{\perp}_t , n^{\perp}_t ) 
\end{equation*}

\noindent
on $\Sigma$, as explained in the previous discussion. Note that metric coefficient $\lambda_t$, or alternatively the unit time-like one-form $\frt_t$, needs to be considered as part of the globally hyperbolic reduction of $[u,v,l,n]$ in order to be able to recover the globally hyperbolic metric associated to $[u,v,l,n]$ from its globally hyperbolic reduction. By direct substitution, we obtain:
\begin{eqnarray*}
g = u\odot v + l \otimes l + n\otimes n = -\frt_t\otimes \frt_t + \frac{1}{(u^o_t)^2} u^{\perp}_t \otimes  u^{\perp}_t +  l^{\perp}_t\otimes  l^{\perp}_t +  n^{\perp}_t\otimes  n^{\perp}_t
\end{eqnarray*}

\noindent
Note that $u^{\perp}_{t}$ is not necessarily of unit norm, whereas $(u^o_t)^{-1} u^{\perp}_{t}$ is because of $u$ being isotropic. Since $l^{\perp}_t$ and $n^{\perp}_t$ are already unit-norm one-forms on $(\Sigma,h_t)$, it is convenient to introduce the following notation:
\begin{equation*}
e_{tu} := \frac{1}{u^o_t} u^{\perp}_{t} \, , \qquad e_{tl}:= l^{\perp}_{t} \, , \qquad e_{tn} := n^{\perp}_{t} 
\end{equation*}

\noindent
in terms of which the metric $g$ associated to $[u,v,l,n]$ adopts the standard form:
\begin{equation*}
g =  -\frt_t\otimes \frt_t + e_{tu} \otimes  e_{tu} + e_{tl} \otimes  e_{tl} + e_{tn} \otimes  e_{tn}
\end{equation*}

\noindent
In particular, $(e_{tu},e_{tl},e_{tn})$ is a global orthonormal coframe on $\Sigma$. For ease of notation, in the following we will set $\fre_t := (e_{tu},e_{tl},e_{tn})$. Furthermore, it will be convenient for future computations to introduce a family of functions $\fra_t$ on $\Sigma$ as follows:
\begin{equation*}
u_t^o = e^{\fra_t}\, , \qquad \forall\,\, t\in \cI
\end{equation*}

\noindent
Note that $u_t^o$ is nowhere vanishing for every $t\in\cI$ and furthermore, without loss of generality, we can take it to be strictly positive, perhaps after changing the admissible bilinear pairing from $\cB$ to $-\cB$.  

\begin{definition}
The tuple $(\lambda_t,\fra_t,\fre_t)$ is the \emph{globally hyperbolic reduction} of the isotropic parallelism $[u,v,l,n]$ on the three-manifold $\Sigma$.
\end{definition}

\noindent
This reduction of the globally hyperbolic isotropic parallelism $[u,v,l,n]$ will play a fundamental role in our study of the Cauchy problem for skew-torsion parallel spinors, which we consider next.


\subsection{The evolution problem of skew-torsion parallel spinors}
\label{sec:evolutionskewtorsion}


In this section, we characterize the evolution problem for a globally hyperbolic skew-torsion parallel spinor in terms of the equivalent system for the associated globally hyperbolic isotropic skew-torsion parallelism. This evolution flow provides a general framework for studying the evolution problem associated with globally hyperbolic supersymmetric solutions of the Neveu-Schwarz system, which, as explained later, relies on a particular class of skew-torsion parallel spinors. Our starting point is Theorem \ref{thm:existencenullcoframeII}, which, when adapted to skew-torsion isotropic parallelisms, states that a four-manifold $M$ admits a skew-torsion isotropic parallelism with torsion $H = \ast_g \alpha \in \Omega^3(M)$ if and only if it admits an isotropic parallelism satisfying the following differential system:
\begin{eqnarray}
& \dd u = \ast_g(\alpha\wedge u)  \, , \qquad \dd v = \ast_g(\alpha\wedge v)  - \kappa\wedge l - \rho\wedge n   \label{eq:duvtorsiongeneralII}\\
& \dd l = \ast_g(\alpha\wedge l) + \kappa\wedge u \, , \qquad \dd n = \ast_g(\alpha\wedge n) + \rho\wedge u   \label{eq:dlntorsiongeneralII}
\end{eqnarray}

\noindent
where $\alpha = \ast_g H$. 

\begin{remark}
\label{remark:torsiondetermined}
The system \eqref{eq:duvtorsiongeneralII} and \eqref{eq:dlntorsiongeneralII} determines $\alpha$, $\kappa$ and $\rho$ uniquely in terms of the representative $(u,v,l,n)$. Indeed, since $(v^{\sharp_g},u^{\sharp_g},l^{\sharp_g},n^{\sharp_g})$ is the frame dual to $(u,v,l,n)$, the first equation in \eqref{eq:duvtorsiongeneralII} reads:
\begin{equation*}
\dd u = \alpha(n^{\sharp_g})\, u\wedge l - \alpha(l^{\sharp_g})\, u\wedge n + \alpha(u^{\sharp_g})\, l\wedge n\, ,
\end{equation*}

\noindent
and determines $\alpha$ modulo $u$, whereas the remaining component is recovered from the other three equations as:
\begin{equation*}
\alpha(v^{\sharp_g}) = \dd l(v^{\sharp_g},n^{\sharp_g}) - \dd n(v^{\sharp_g},l^{\sharp_g}) + \dd v(l^{\sharp_g},n^{\sharp_g})\, ,
\end{equation*}

\noindent
after which $\kappa$ and $\rho$ are read off from $\dd l$, $\dd n$ and $\dd v$. In particular, the skew-torsion of a skew-torsion parallel spinor is determined by its isotropic parallelism, and:
\begin{equation}
\label{eq:screenrotationinvariant}
u\wedge\dd u = -\alpha(u^{\sharp_g})\, \ast_g u\, .
\end{equation}

\noindent
Hence, the distribution $\Ker(u)\subset TM$ is integrable if and only if $\alpha(u^{\sharp_g})= 0$.
\end{remark}

\noindent
Our goal now is to \emph{reduce} the differential system \eqref{eq:duvtorsiongeneralII} and \eqref{eq:dlntorsiongeneralII} when evaluated on the adapted representative $(u,v,l,n)\in [u,v,l,n]$ of a globally hyperbolic isotropic parallelism $[u,v,l,n]$. Recall that by Lemma \ref{lemma:ghcoframe} we can write:
\begin{equation*}
u = e^{\fra_t} (\mathfrak{t}_t + e_{tu})\, , \qquad v = \frac{e^{-\fra_t}}{2}\,( -  \mathfrak{t}_t + e_{tu})  \, , \qquad l =   e_{tl} \, , \qquad n = e_{tn} 
\end{equation*}

\noindent
Using these equations, we compute the exterior derivative of the adapted representative $(u,v,l,n) \in [u,v,l,n]$ of a skew-torsion globally hyperbolic isotropic parallelism $[u,v,l,n]$ in terms of its globally hyperbolic reduction, namely in terms of objects and differentials defined on the three-manifold $\Sigma$.

\begin{lemma}
\label{lemma:derivativeadaptedcoframe}
Let $(u,v,l,n)\in [u,v,l,n]$ be the adapted representative of a globally hyperbolic parallelism. The following formulas hold:
\begin{eqnarray*}
& \dd u = e^{\fra_t} ( \frt_t \wedge (\lambda_t^{-1}\partial_t \fra_t e_{tu} - \dd_{\Sigma} \fra_t - \lambda_t^{-1} \dd_{\Sigma} \lambda_t + \lambda_t^{-1} \partial_t e_{tu}) + \dd_{\Sigma} e_{tu} + \dd_{\Sigma} \fra_t \wedge e_{tu})\\
& \dd v = \frac{1}{2} e^{-\fra_t} ( \frt_t \wedge (\lambda_t^{-1} \partial_t e_{tu} - \lambda_t^{-1}\partial_t \fra_t e_{tu} - \dd_{\Sigma} \fra_t + \lambda_t^{-1} \dd_{\Sigma} \lambda_t) + \dd_{\Sigma} e_{tu} - \dd_{\Sigma} \fra_t \wedge e_{tu})\\
& \dd e_{tl} = \frt_t \wedge \lambda_t^{-1} \partial_t e_{tl} + \dd_{\Sigma} e_{tl} \, , \qquad \dd e_{tn} = \frt_t \wedge \lambda_t^{-1} \partial_t e_{tn} + \dd_{\Sigma} e_{tn} 
\end{eqnarray*}
	
\noindent
where $(\lambda_t,\fra_t,\fre_t)$ is the globally hyperbolic reduction of $[u,v,l,n]$.
\end{lemma}

\begin{proof}
We compute:
\begin{eqnarray*}
& \dd u = e^{\fra_t} \dd \fra_t \wedge (\mathfrak{t}_t + e_{tu}) + e^{\fra_t} (\dd\mathfrak{t}_t + \dd e_{tu}) = e^{\fra_t} (\dd_{\Sigma} \fra_t \wedge \mathfrak{t}_t + \partial_t \fra_t \dd t\wedge e_{tu} + \dd_{\Sigma} \fra_t \wedge e_{tu} )\\
& + e^{\fra_t} (\dd_{\Sigma} \lambda_t \wedge \dd t + \dd_{\Sigma} e_{tu} + \dd t \wedge \partial_t e_{tu})
\end{eqnarray*}

\noindent
After isolating for type of tensor we obtain the first equation in the lemma. Similarly, we have:
\begin{eqnarray*}
& \dd v = \frac{1}{2}e^{-\fra_t}  \dd \fra_t \wedge (\mathfrak{t}_t - e_{tu}) +  \frac{e^{-\fra_t}}{2} (- \dd\mathfrak{t}_t + \dd e_{tu}) =  \frac{e^{-\fra_t}}{2} ( \dd_{\Sigma} \fra_t \wedge \mathfrak{t}_t - \dd_{\Sigma} \fra_t\wedge e_{tu} - \partial_t \fra_t \dd t \wedge e_{tu})\\
& +  \frac{e^{-\fra_t}}{2} (\dd_{\Sigma} e_{tu} + \dd t \wedge \partial_t e_{tu} - \lambda_t^{-1}\dd_{\Sigma}\lambda_t \wedge \mathfrak{t}_t)
\end{eqnarray*}

\noindent
which gives the second equation in the lemma. The third and fourth equations are proven similarly and thus we conclude. 
\end{proof}

\noindent
Let $(\lambda_t,\fra_t,\fre_t)$ be the globally hyperbolic reduction of a skew-torsion globally hyperbolic parallelism $[u,v,l,n]$ with torsion $H = \ast_g \alpha$. We split $\alpha$ as expected:
\begin{equation}
\label{eq:splittingxialpha}
\alpha = \alpha_t^o \frt_t +  \alpha^{\perp}_t
\end{equation}

\noindent
in terms of uniquely determined families of functions $\alpha_t^o$ and one-forms $\alpha^{\perp}_t$ on $\Sigma$. For further reference in upcoming calculations, we recall the following identities:
\begin{equation}
\label{eq:hodgehyperbolic}
\ast_g(\frt_t \wedge \beta) = -\ast_{\fre_t}\beta\, , \qquad \ast_{g}\beta = (-1)^{\vert\beta \vert} \frt_t \wedge \ast_{\fre_t} \beta\, , \qquad \beta \in \Gamma(\wedge M) 
\end{equation}

\noindent
relating the Hodge dual $\ast_g$ associated with $g$ to the Hodge dual $\ast_{\fre_t}$ associated with $h_t$. Using these relations, the following lemma follows by a direct computation.

\begin{lemma}
\label{lemma:dualphasplitting}
Let $(u,v,l,n)$ be the adapted representative of a globally hyperbolic
isotropic parallelism $[u,v,l,n]$ on $\cI\times\Sigma$. Then:
\begin{eqnarray*}
& \ast_g(\alpha\wedge u) = e^{\fra_t}  (\frt_t \wedge \ast_{\fre_t}(\alpha^{\perp}_t \wedge e_{tu})  - \alpha^o_t \ast_{\fre_t} e_{tu} + \ast_{\fre_t}\alpha^{\perp}_t) \\
& \ast_g(\alpha\wedge v) = \frac{1}{2} e^{-\fra_t} (\frt_t \wedge \ast_{\fre_t}(\alpha^{\perp}_t \wedge e_{tu})  - \alpha^o_t \ast_{\fre_t} e_{tu} - \ast_{\fre_t}\alpha^{\perp}_t) \\
& \ast_g(\alpha\wedge e_{tl}) =  \frt_t \wedge \ast_{\fre_t}(\alpha^{\perp}_t \wedge e_{tl})  - \alpha^o_t \ast_{\fre_t} e_{tl} \\
& \ast_g(\alpha\wedge e_{tn}) =  \frt_t \wedge \ast_{\fre_t}(\alpha^{\perp}_t \wedge e_{tn})  - \alpha^o_t \ast_{\fre_t} e_{tn} 
\end{eqnarray*}

\noindent
where $\alpha = \alpha^o_t \frt_t + \alpha^{\perp}_t \in \Omega^1(\cI\times \Sigma)$.
\end{lemma}

\noindent
The following lemma is the first step in our characterization of skew torsion globally hyperbolic parallelisms.

\begin{lemma}
\label{lemma:skewghevolution1}
A globally hyperbolic parallelism $[u,v,l,n]$ on $\cI\times \Sigma$ is skew-torsion with torsion $H =  \ast_g\alpha$ if and only if there exists a family of functions $(\kappa^o_t,\rho^o_t)$ and a family of one-forms $(\kappa^{\perp}_t , \rho^{\perp}_t)$ on $\Sigma$ such that the globally hyperbolic reduction $(\lambda_t,\fra_t,\fre_t)$ of $[u,v,l,n]$ satisfies the following differential equations on $\Sigma$:
\begin{eqnarray*}
& \lambda_t^{-1} \partial_t e_{tu} = \ast_{\fre_t} (\alpha^{\perp}_t \wedge e_{tu}) - \lambda_t^{-1}\partial_t \fra_t e_{tu} + \dd_{\Sigma} \fra_t + \lambda_t^{-1} \dd_{\Sigma} \lambda_t\\
& \lambda_t^{-1} \partial_t e_{tl} = \ast_{\fre_t} (\alpha^{\perp}_t \wedge e_{tl}) +  e^{\fra_t} (\kappa^o_t e_{tu} - \kappa^{\perp}_t) \, , \quad  \lambda_t^{-1} \partial_t e_{tn} = \ast_{\fre_t} (\alpha^{\perp}_t\wedge e_{tn})  + e^{\fra_t} (\rho^o_t e_{tu} - \rho^{\perp}_t) \\
& \dd_{\Sigma}  e_{tu}  = e_{tu}\wedge \dd_{\Sigma} \fra_t + \ast_{\fre_t}(\alpha^{\perp}_t  - \alpha^o_t e_{tu}) \\
& \dd_{\Sigma} e_{tl} = e^{\fra_t} \kappa^{\perp}_t\wedge e_{tu} - \alpha^o_t \ast_{\fre_t} e_{tl}\, , \quad  \dd_{\Sigma} e_{tn} = e^{\fra_t} \rho^{\perp}_t\wedge e_{tu}  - \alpha^o_t \ast_{\fre_t} e_{tn} 
\end{eqnarray*}
	
\noindent
together with the equations:
\begin{eqnarray*}
 e^{\fra_t}(\kappa^o_t e_{tl} +   \rho^o_t e_{tn} )= \lambda_t^{-1}\partial_t \fra_t e_{tu}  - \lambda_t^{-1} \dd_{\Sigma} \lambda_t   \, , \quad    \dd_{\Sigma} \fra_t \wedge e_{tu} =  \ast_{\fre_t}\alpha^{\perp}_t     +    e^{\fra_t} (\kappa^{\perp}_t \wedge e_{tl} + \rho^{\perp}_t \wedge e_{tn}  )
\end{eqnarray*}

\noindent
further restricting the evolution of $(\lambda_t,\fra_t,\fre_t)$ and $(\kappa^o_t,\rho^o_t, \kappa^{\perp}_t , \rho^{\perp}_t)$.
\end{lemma}
 
\begin{proof}
Let $[u,v,l,n]$ be a globally hyperbolic parallelism on $\cI\times \Sigma$ relative to $\kappa , \rho \in \Omega^1(\cI\times\Sigma)$ and with torsion $H = \ast_{g} \alpha$. We split $\alpha$ as prescribed in Equation \eqref{eq:splittingxialpha}. We split $\kappa$ and $\rho$ similarly to $\alpha$, namely, we set:
\begin{equation*}
\kappa = \kappa^o_t \frt_t + \kappa^{\perp}_t\, ,\qquad \rho = \rho^o_t \frt_t + \rho^{\perp}_t
\end{equation*}
	
\noindent
where $(\kappa^o_t , \rho^o_t)$ are families of functions on $\Sigma$ and $(\kappa^{\perp}_t , \rho^{\perp}_t)$ are families of one-forms on $\Sigma$. By Theorem \ref{thm:existencenullcoframeII} it follows that $[u,v,l,n]$ is a globally hyperbolic skew-torsion parallelism if and only if its adapted representative $(u,v,l,n)\in [u,v,l,n]$ on $\cI\times \Sigma$ satisfies \eqref{eq:duvtorsiongeneralII} and \eqref{eq:dlntorsiongeneralII}. Using Lemmas \ref{lemma:derivativeadaptedcoframe} and \ref{lemma:dualphasplitting} we obtain that the first equation in \eqref{eq:duvtorsiongeneralII} is equivalent to:
\begin{equation*}
\frt_t \wedge (\lambda_t^{-1}\partial_t \fra_t e_{tu} - \dd_{\Sigma} \fra_t - \lambda_t^{-1} \dd_{\Sigma} \lambda_t + \lambda_t^{-1} \partial_t e_{tu}) + \dd_{\Sigma} e_{tu} + \dd_{\Sigma} \fra_t \wedge e_{tu} =  \frt_t  \wedge  \ast_{\fre_t} (\alpha^{\perp}_t \wedge e_{tu})  +  \ast_{\fre_t}(\alpha^{\perp}_t   - \alpha^o_t e_{tu}) 
\end{equation*}
	
\noindent
Isolating by tensor type, this equation becomes equivalent to:
\begin{eqnarray}
\label{eq:partialeuproof}
\lambda_t^{-1} \partial_t e_{tu} = \ast_{\fre_t} (\alpha^{\perp}_t \wedge e_{tu}) - \lambda_t^{-1}\partial_t \fra_t e_{tu} + \dd_{\Sigma} \fra_t + \lambda_t^{-1} \dd_{\Sigma} \lambda_t\, , \,\,  e^{-\fra_t} \dd_{\Sigma} (e^{\fra_t} e_{tu} ) = \ast_{\fre_t}(\alpha^{\perp}_t  - \alpha^o_t e_{tu}) 
\end{eqnarray}
	
\noindent
Proceeding similarly for both equations in \eqref{eq:dlntorsiongeneralII}, we obtain that they are equivalent to the following equations:
\begin{eqnarray*}
& \frt_t \wedge \lambda_t^{-1} \partial_t e_{tl} + \dd_{\Sigma} e_{tl} =  \frt_t  \wedge (\ast_{\fre_t} (\alpha^{\perp}_t  \wedge e_{tl}) + e^{\fra_t} (\kappa^o_t e_{tu} - \kappa^{\perp}_t)) + e^{\fra_t} \kappa^{\perp}_t\wedge e_{tu} - \alpha^o_t \ast_{\fre_t} e_{tl} \\
& \frt_t \wedge \lambda_t^{-1} \partial_t e_{tn} + \dd_{\Sigma} e_{tn} = \frt_t  \wedge (\ast_{\fre_t} (\alpha^{\perp}_t \wedge e_{tn}) +   e^{\fra_t} (\rho^o_t e_{tu} - \rho^{\perp}_t)) +  e^{\fra_t} \rho^{\perp}_t\wedge e_{tu}  - \alpha^o_t \ast_{\fre_t} e_{tn}  
\end{eqnarray*}
	
\noindent
Isolating by tensor type, we obtain the following equivalent system:
\begin{eqnarray*}
& \lambda_t^{-1} \partial_t e_{tl} = \ast_{\fre_t} (\alpha^{\perp}_t \wedge e_{tl}) +  e^{\fra_t} (\kappa^o_t e_{tu} - \kappa^{\perp}_t) \, , \qquad  \lambda_t^{-1} \partial_t e_{tn} = \ast_{\fre_t} (\alpha^{\perp}_t\wedge e_{tn})  + e^{\fra_t} (\rho^o_t e_{tu} - \rho^{\perp}_t) \\
& \dd_{\Sigma} e_{tl} = e^{\fra_t} \kappa^{\perp}_t\wedge e_{tu} - \alpha^o_t \ast_{\fre_t} e_{tl}\, , \qquad \dd_{\Sigma} e_{tn} = e^{\fra_t} \rho^{\perp}_t\wedge e_{tu}  - \alpha^o_t \ast_{\fre_t} e_{tn}   
\end{eqnarray*}
	
\noindent
We consider now the second equation in  \eqref{eq:duvtorsiongeneralII}, which is the last that remains to be considered. Using Lemmas \ref{lemma:derivativeadaptedcoframe} and \ref{lemma:dualphasplitting} it follows that it is equivalent to:
\begin{eqnarray*}
& \frac{1}{2} e^{-\fra_t}  \frt_t \wedge (\lambda_t^{-1} \partial_t e_{tu} - \lambda_t^{-1}\partial_t \fra_t e_{tu} - \dd_{\Sigma} \fra_t + \lambda_t^{-1} \dd_{\Sigma} \lambda_t) + \frac{1}{2} e^{-\fra_t}  \dd_{\Sigma} e_{tu} - \frac{1}{2} e^{-\fra_t}  \dd_{\Sigma} \fra_t \wedge e_{tu} \\
& =  \frt_t  \wedge ( \frac{1}{2} e^{-\fra_t} \ast_{\fre_t} (\alpha^{\perp}_t \wedge e_{tu}) - \kappa^o_t e_{tl} - \rho^o_t e_{tn}) -   \frac{1}{2} e^{-\fra_t} \ast_{\fre_t}(\alpha^{\perp}_t + \alpha^o_t e_{tu})  - \kappa^{\perp}_t \wedge e_{tl} - \rho^{\perp}_t \wedge e_{tn}
\end{eqnarray*}
	
\noindent
Isolating by tensor type, we obtain:
\begin{eqnarray}
&  \lambda_t^{-1} \partial_t e_{tu}    = \lambda_t^{-1}\partial_t \fra_t e_{tu} + \dd_{\Sigma} \fra_t - \lambda_t^{-1} \dd_{\Sigma} \lambda_t +   \ast_{\fre_t} (\alpha^{\perp}_t \wedge e_{tu}) -2 e^{\fra_t}\kappa^o_t e_{tl} - 2 e^{\fra_t} \rho^o_t e_{tn}  \\
&  \dd_{\Sigma} e_{tu}  = \dd_{\Sigma} \fra_t \wedge e_{tu} - \ast_{\fre_t}(\alpha^{\perp}_t +  \alpha^o_t e_{tu})   +  2 e^{\fra_t}( e_{tl} \wedge \kappa^{\perp}_t +   e_{tn} \wedge \rho^{\perp}_t)   \label{eq:auxiliaryeuII}
\end{eqnarray}
	
\noindent
The first equation prescribes the \emph{time-derivative} of $e_{tu}$, and therefore needs to be compared with the first equation in \eqref{eq:partialeuproof}. Similarly, the second equation prescribes the exterior derivative of $e_{tu}$ and therefore needs to be compared with the second equation in \eqref{eq:partialeuproof}. This results, together with the previous equations, in the \emph{constraints} appearing in the statement of the lemma. The converse follows by tracing back the previous steps and hence we conclude. 
\end{proof}

\noindent
The previous lemma can be conveniently refined as follows.
\begin{prop}
\label{prop:evolutionsystem}
A globally hyperbolic parallelism $[u,v,l,n]$ on $\cI\times \Sigma$ is skew-torsion with torsion $H =  \ast_g\alpha$ if and only if there exists a family of one-forms $(\kappa^{\perp}_t , \rho^{\perp}_t)$ on $\Sigma$ such that the globally hyperbolic reduction $(\lambda_t,\fra_t,\fre_t)$ of $[u,v,l,n]$ satisfies the following system of evolution equations on $\Sigma$:
\begin{eqnarray}
& \partial_t \fra_t =  \dd_{\Sigma} \lambda_t( e_{tu}) \nonumber \\
&  \lambda_t^{-1} \partial_t e_{tu} = \ast_{\fre_t} (\alpha^{\perp}_t \wedge e_{tu}) + \dd_{\Sigma} \fra_t + \lambda_t^{-1} \dd_{\Sigma} \lambda_t (e_{tl}) e_{tl} + \lambda_t^{-1} \dd_{\Sigma} \lambda_t (e_{tn}) e_{tn} \nonumber \\
& \lambda_t^{-1} \partial_t e_{tl} = \ast_{\fre_t} (\alpha^{\perp}_t \wedge e_{tl}) - (\lambda_t^{-1} \dd_{\Sigma} \lambda_t (e_{tl}) e_{tu} + e^{\fra_t}\kappa^{\perp}_t) \label{eq:evolutionequationsI} \\
& \lambda_t^{-1} \partial_t e_{tn} = \ast_{\fre_t} (\alpha^{\perp}_t\wedge e_{tn})  - (\lambda_t^{-1} \dd_{\Sigma} \lambda_t (e_{tn}) e_{tu} + e^{\fra_t}\rho^{\perp}_t) \nonumber
\end{eqnarray}

\noindent
together with the following time-dependent constraint equations:
\begin{eqnarray}
&  \dd_{\Sigma} \fra_t \wedge e_{tu} =  \ast_{\fre_t}\alpha^{\perp}_t     +    e^{\fra_t} \kappa^{\perp}_t \wedge e_{tl} +   e^{\fra_t} \rho^{\perp}_t \wedge e_{tn}   \nonumber \\
& \dd_{\Sigma}  e_{tu}  =  - e^{\fra_t} \kappa^{\perp}_t \wedge e_{tl} - e^{\fra_t} \rho^{\perp}_t \wedge e_{tn} -   \alpha^o_t \ast_{\fre_t} e_{tu} \nonumber \\
& \dd_{\Sigma} e_{tl} = e^{\fra_t} \kappa^{\perp}_t\wedge e_{tu} - \alpha^o_t \ast_{\fre_t} e_{tl} \label{eq:constraintequationsI}\\
& \dd_{\Sigma} e_{tn} = e^{\fra_t} \rho^{\perp}_t\wedge e_{tu}  - \alpha^o_t \ast_{\fre_t} e_{tn} \nonumber
\end{eqnarray}

\noindent
where $\alpha = \alpha^o_t \frt_t + \alpha^{\perp}_t$. 
\end{prop}  

\begin{remark}
From the differential system \eqref{eq:evolutionequationsI} and \eqref{eq:constraintequationsI} it can be easily shown that the exterior derivative of $e_{tu}$ can be equivalently written as follows:
\begin{equation*}
\dd_{\Sigma}  e_{tu}  = e_{tu}\wedge \dd_{\Sigma} \fra_t + \ast_{\fre_t}(\alpha^{\perp}_t  - \alpha^o_t e_{tu})
\end{equation*}

\noindent
as it follows from a simple combination with the equation for the exterior derivative of $\fra_t$. 
\end{remark}

\begin{proof}
We consider the following equation:
\begin{equation*}
 e^{\fra_t}\kappa^o_t e_{tl} + e^{\fra_t} \rho^o_t e_{tn} = \lambda_t^{-1}\partial_t \fra_t e_{tu}  - \lambda_t^{-1} \dd_{\Sigma} \lambda_t   
\end{equation*}

\noindent
in Lemma \ref{lemma:skewghevolution1}.  Evaluating it on $e^{\sharp}_{tu}$ gives $\partial_t\fra_t=\dd_\Sigma\lambda_t(e^{\sharp}_{tu})$, which in turn implies:
\begin{equation}
\label{eq:kapparhonormalcomponents}
\kappa_t^o=-e^{-\fra_t}\lambda_t^{-1} \dd_\Sigma\lambda_t(e^{\sharp}_{tl}), \qquad \rho_t^o=-e^{-\fra_t}\lambda_t^{-1} \dd_\Sigma\lambda_t(e^{\sharp}_{tn}).
\end{equation}

\noindent
Substitution of these identities into the three evolution equations of Lemma \ref{lemma:skewghevolution1} gives \eqref{eq:evolutionequationsI}. Note that the spatial equations are unchanged and give \eqref{eq:constraintequationsI}.
\end{proof}

\noindent
Hence, the differential system characterizing skew-torsion parallel spinors splits, when transformed into the equivalent system for skew-torsion globally hyperbolic parallelisms, into two natural blocks, namely block \eqref{eq:evolutionequationsI} consisting of the evolution equations for $(\fra_t,\fre_t)$ together with block \eqref{eq:constraintequationsI} which contains the natural \emph{time-dependent constraint equations} for the aforementioned evolution equations. This is therefore an evolution problem for $(\fra_t,\fre_t)$ for which we can consider $\lambda_t$ as given data. We elaborate on the previous lemma in the following in order to obtain a geometric interpretation of the globally hyperbolic skew-torsion isotropic parallelisms in terms of the second fundamental form of the underlying Cauchy hypersurface.

\begin{lemma}
\label{lemma:ThetatI}
Let $(\lambda_t,\fra_t,\fre_t)$ be the globally hyperbolic reduction of a skew-torsion globally hyperbolic parallelism $[u,v,l,n]$ relative to $(\kappa^{\perp}_t , \rho^{\perp}_t)$. Then, the second fundamental form $\Theta_t \in \Gamma(T^{\ast}\Sigma \odot T^{\ast}\Sigma)$ of $(\Sigma,h_t)$ is given by:
\begin{eqnarray}
\label{eq:ThetatI}
-2 \Theta_t = \dd_{\Sigma} \fra_t \odot e_{tu} - e^{\fra_t} (\kappa^{\perp}_t \odot e_{tl} + \rho^{\perp}_t \odot e_{tn})
\end{eqnarray}

\noindent
for every $t\in \cI$.
\end{lemma}

\begin{proof}
Using equation \eqref{eq:evolutionequationsI} we compute:
\begin{eqnarray*}
& -2\Theta_t  = \frac{1}{\lambda_t} (e_{tu} \odot \partial_t e_{tu}  + e_{tl}\odot \partial_t e_{tl} + e_{tn} \odot\partial_t e_{tn}) \\
& = (\ast_{\fre_t} (\alpha^{\perp}_t \wedge e_{tu})  + \dd_{\Sigma} \fra_t  + \lambda_t^{-1} \dd_{\Sigma} \lambda_t (e_{tl}) e_{tl}  + \lambda_t^{-1} \dd_{\Sigma} \lambda_t (e_{tn}) e_{tn}) \odot e_{tu} \\
& \quad + (\ast_{\fre_t} (\alpha^{\perp}_t \wedge e_{tl}) - (\lambda_t^{-1} \dd_{\Sigma} \lambda_t (e_{tl}) e_{tu} + e^{\fra_t}\kappa^{\perp}_t))\odot e_{tl} \\
& + (\ast_{\fre_t} (\alpha^{\perp}_t\wedge e_{tn})  - (\lambda_t^{-1} \dd_{\Sigma} \lambda_t (e_{tn}) e_{tu} + e^{\fra_t}\rho^{\perp}_t)) \odot e_{tn}\\
& = \dd_{\Sigma} \fra_t \odot e_{tu} - e^{\fra_t} (\kappa^{\perp}_t \odot e_{tl} + \rho^{\perp}_t \odot e_{tn})
\end{eqnarray*}

\noindent
and thus we conclude. 
\end{proof}

\begin{prop}
\label{prop:kapparhoisolation}
Let $(\lambda_t,\fra_t,\fre_t)$ be the globally hyperbolic reduction of a skew-torsion globally hyperbolic parallelism $[u,v,l,n]$ relative to $(\kappa^{\perp}_t , \rho^{\perp}_t)$. Then:
\begin{eqnarray}
& \Theta_t (e_{tu}) = -\dd_{\Sigma} \fra_t - \frac{1}{2}\ast_{\fre_t} (\alpha^{\perp}_t \wedge e_{tu})\nonumber\\
& \Theta_t (e_{tl}) = e^{\fra_t} \kappa^{\perp}_t - \frac{1}{2} \ast_{\fre_t} (\alpha^{\perp}_t \wedge e_{tl}) \label{eq:ThetaexpressionII}\\
& \Theta_t (e_{tn}) = e^{\fra_t} \rho^{\perp}_t - \frac{1}{2} \ast_{\fre_t} (\alpha^{\perp}_t \wedge e_{tn})\nonumber
\end{eqnarray}

\noindent
for every $t\in \cI$, where $\Theta_t$ is the second fundamental form of $(\Sigma,h_t)$.
\end{prop}

\begin{proof}
Evaluating Equation \eqref{eq:ThetatI} on $e_{tu}^{\sharp}$ we obtain:
\begin{eqnarray*}
-2 \Theta_t(e_{tu}) = \dd_{\Sigma} \fra_t   + \dd_{\Sigma} \fra_t (e_{tu}) e_{tu} - e^{\fra_t} (\kappa^{\perp}_t (e_{tu}) e_{tl} + \rho^{\perp}_t(e_{tu}) e_{tn})
\end{eqnarray*}

\noindent
Combining this equation with the evaluation of the first equation in \eqref{eq:constraintequationsI} on $e_{tu}^{\sharp}$, we obtain the first equation in \eqref{eq:ThetaexpressionII}. Evaluating Equation \eqref{eq:ThetatI} on $e_{tl}^{\sharp}$ we obtain:
\begin{equation*}
-2 \Theta_t (e_{tl})= \dd_{\Sigma} \fra_t(e_{tl})  e_{tu} - e^{\fra_t} (\kappa^{\perp}_t (e_{tl})   e_{tl} + \rho^{\perp}_t(e_{tl})  e_{tn} + \kappa^{\perp}_t)
\end{equation*}

\noindent
Combining this equation with the evaluation of the first equation in \eqref{eq:constraintequationsI} on $e_{tl}^{\sharp}$, we obtain the second equation in \eqref{eq:ThetaexpressionII}. The third equation in \eqref{eq:ThetaexpressionII} is obtained similarly by evaluating Equation \eqref{eq:ThetatI} on $e_{tn}^{\sharp}$ and then combining the result with the evaluation of the first equation in \eqref{eq:constraintequationsI} on $e_{tn}^{\sharp}$.
\end{proof}

\noindent
The previous proposition allows to obtain an equivalent expression for the evolution problem defined by the differential system \eqref{eq:evolutionequationsI} and \eqref{eq:constraintequationsI} which is obtained by isolating $\kappa^{\perp}_t$ and $\rho^{\perp}_t$ in \eqref{eq:ThetaexpressionII} and plugging the result back into \eqref{eq:evolutionequationsI} and \eqref{eq:constraintequationsI}.

\begin{cor}
\label{cor:Thetaevolutionsystem}
A globally hyperbolic parallelism $[u,v,l,n]$ on $\cI\times \Sigma$ is skew-torsion with torsion $H =  \ast_g\alpha$ if and only if its globally hyperbolic reduction $(\lambda_t,\fra_t,\fre_t)$ satisfies the following system of evolution equations on $\Sigma$:
\begin{eqnarray}
& \partial_t \fra_t =  \dd_{\Sigma} \lambda_t( e_{tu}) \nonumber \\
&  \lambda_t^{-1} \partial_t e_{tu} + \Theta_t(e_{tu}) = \frac{1}{2}\ast_{\fre_t} (\alpha^{\perp}_t \wedge e_{tu}) +   \lambda_t^{-1} \dd_{\Sigma} \lambda_t (e_{tl}) e_{tl} + \lambda_t^{-1} \dd_{\Sigma} \lambda_t (e_{tn}) e_{tn} \nonumber \\
& \lambda_t^{-1} \partial_t e_{tl} + \Theta_t(e_{tl})  = \frac{1}{2} \ast_{\fre_t} (\alpha^{\perp}_t \wedge e_{tl}) -  \lambda_t^{-1} \dd_{\Sigma} \lambda_t (e_{tl}) e_{tu}   \label{eq:evolutionequationsII} \\
& \lambda_t^{-1} \partial_t e_{tn} + \Theta_t(e_{tn}) = \frac{1}{2} \ast_{\fre_t} (\alpha^{\perp}_t\wedge e_{tn})  -  \lambda_t^{-1} \dd_{\Sigma} \lambda_t (e_{tn}) e_{tu}   \nonumber
\end{eqnarray}

\noindent
together with the following time-dependent constraint equations:
\begin{eqnarray}
&  \dd_{\Sigma} \fra_t   =  - \Theta_t (e_{tu}) - \frac{1}{2}  \ast_{\fre_t} (\alpha^{\perp}_t \wedge e_{tu})  \nonumber \\
& \dd_{\Sigma} e_{tu} = \Theta_t(e_{tu})\wedge e_{tu} + \ast_{\fre_t} (\alpha^{\perp}_t   - \alpha^o_t e_{tu}) - \frac{1}{2} e_{tu} \wedge \ast_{\fre_t} (\alpha^{\perp}_t \wedge e_{tu}) \nonumber \\
& \dd_{\Sigma} e_{tl} = \Theta_t(e_{tl})\wedge e_{tu} +  (\frac{1}{2}\alpha^{\perp}_t (e_{tu}) - \alpha^o_t) \ast_{\fre_t} e_{tl} \label{eq:constraintequationsII} \\
& \dd_{\Sigma} e_{tn} = \Theta_t(e_{tn})\wedge e_{tu} +  (\frac{1}{2} \alpha^{\perp}_t (e_{tu}) - \alpha^o_t) \ast_{\fre_t} e_{tn} \nonumber
\end{eqnarray}

\noindent
where $\Theta_t$ is the second fundamental form of $(\Sigma,h_t)$ and $\alpha = \alpha^o_t \frt_t + \alpha^{\perp}_t$. 
\end{cor}

\begin{remark}
The first and second equations in \eqref{eq:constraintequationsII} give the following equation:
\begin{equation*}
e^{-\fra_t}\dd_{\Sigma} (e^{\fra_t} e_{tu}) =  \ast_{\fre_t} (\alpha^{\perp}_t   - \alpha^o_t e_{tu})
\end{equation*}

\noindent
which will be useful in the following. 
\end{remark}

\noindent
The previous reformulation of the evolution problem posed by a globally hyperbolic skew-torsion spinor is particularly well-suited for defining the associated initial data and constraint equations. Still, it can be refined further by eliminating the family $\fra_t$ in favor of a condition in de Rham cohomology.

\begin{prop}
\label{prop:Thetaevolutionsystem}
Let $(\lambda_t,\fre_t)$ determine a globally hyperbolic metric $g=-\lambda_t^2\dd t\otimes\dd t+h_{\fre_t}$ on $\cI\times\Sigma$. Then $(\lambda_t,\fre_t)$ is the reduction of a skew-torsion parallelism with torsion $H = \ast_g\alpha$, for a suitable family of functions $\fra_t$, if and only if it satisfies the following system of evolution equations on $\Sigma$:
\begin{eqnarray}
	&  \lambda_t^{-1} \partial_t e_{tu} + \Theta_t(e_{tu}) = \frac{1}{2}\ast_{\fre_t} (\alpha^{\perp}_t \wedge e_{tu}) +   \lambda_t^{-1} \dd_{\Sigma} \lambda_t (e_{tl}) e_{tl} + \lambda_t^{-1} \dd_{\Sigma} \lambda_t (e_{tn}) e_{tn} \nonumber \\
	& \lambda_t^{-1} \partial_t e_{tl} + \Theta_t(e_{tl})  = \frac{1}{2} \ast_{\fre_t} (\alpha^{\perp}_t \wedge e_{tl}) -  \lambda_t^{-1} \dd_{\Sigma} \lambda_t (e_{tl}) e_{tu}   \label{eq:evolutionequationsIII} \\
	& \lambda_t^{-1} \partial_t e_{tn} + \Theta_t(e_{tn}) = \frac{1}{2} \ast_{\fre_t} (\alpha^{\perp}_t\wedge e_{tn})  -  \lambda_t^{-1} \dd_{\Sigma} \lambda_t (e_{tn}) e_{tu}   \nonumber
\end{eqnarray}

\noindent
together with the following time-dependent constraint equations:
\begin{eqnarray}
& \dd_{\Sigma} e_{tu} = \Theta_t(e_{tu})\wedge e_{tu} + \ast_{\fre_t} (\alpha^{\perp}_t   - \alpha^o_t e_{tu}) - \frac{1}{2} e_{tu} \wedge \ast_{\fre_t} (\alpha^{\perp}_t \wedge e_{tu}) \nonumber \\
& \dd_{\Sigma} e_{tl} = \Theta_t(e_{tl})\wedge e_{tu} +  (\frac{1}{2} \alpha^{\perp}_t (e_{tu}) - \alpha^o_t) \ast_{\fre_t} e_{tl} \label{eq:constraintequationsIII} \\
& \dd_{\Sigma} e_{tn} = \Theta_t(e_{tn})\wedge e_{tu} +  (\frac{1}{2} \alpha^{\perp}_t (e_{tu}) - \alpha^o_t) \ast_{\fre_t} e_{tn} \nonumber 
\end{eqnarray}

\noindent
and the following conditions:
\begin{equation}
\label{eq:cohomologicalccondition}
0 = [\Theta_t (e_{tu}) + \frac{1}{2}\ast_{\fre_t} (\alpha^{\perp}_t \wedge e_{tu})]\in H^1(\Sigma,\mathbb{R})\, , \quad \partial_t (\Theta_t (e_{tu}) + \frac{1}{2}\ast_{\fre_t} (\alpha^{\perp}_t \wedge e_{tu})) = -\dd_{\Sigma}(\dd_{\Sigma}\lambda_t (e_{tu}))
\end{equation}

\noindent
where $\Theta_t$ is the second fundamental form of $(\Sigma,h_t)$ and $\alpha = \alpha^o_t \frt_t + \alpha^{\perp}_t$. 
\end{prop}

\begin{remark}
The first equation in \eqref{eq:constraintequationsIII} can be equivalently written as follows:
\begin{equation*}
\dd_{\Sigma} e_{tu} = \Theta_t(e_{tu})\wedge e_{tu} +  (\alpha^{\perp}_t (e_{tu})   - \alpha^o_t  ) \ast_{\fre_t} e_{tu}  +\frac{1}{2} (\alpha^{\perp}_t(e_{tl}) \ast_{\fre_t} e_{tl} + \alpha^{\perp}_t(e_{tn}) \ast_{\fre_t} e_{tn}) 
\end{equation*}

\noindent
which is sometimes useful in computations. 
\end{remark}

\begin{proof}
Let $(\lambda_t,\fra_t,\fre_t)$ be a solution to the differential system \eqref{eq:evolutionequationsII} and \eqref{eq:constraintequationsII}. Then, $(\lambda_t , \fre_t)$ automatically satisfies equations \eqref{eq:evolutionequationsIII} and \eqref{eq:constraintequationsIII} and thus we only need to deal with the conditions contained in Equation \eqref{eq:cohomologicalccondition}. By the first equation in \eqref{eq:ThetaexpressionII} the first equation in \eqref{eq:cohomologicalccondition} follows. Taking the time derivative of the first equation in \eqref{eq:ThetaexpressionII} and combining the result with the first equation in \eqref{eq:evolutionequationsII}  we obtain the second equation in \eqref{eq:cohomologicalccondition}. For the converse, suppose that $(\lambda_t,\fre_t)$ satisfies equations \eqref{eq:evolutionequationsIII}, \eqref{eq:constraintequationsIII} and \eqref{eq:cohomologicalccondition}. We only need to prove that the first equation \eqref{eq:evolutionequationsII} and the first equation in \eqref{eq:constraintequationsII} both hold. By the first equation in \eqref{eq:cohomologicalccondition}, there exist a family of functions $\bar{a}_t$ such that:
\begin{equation*}
\dd_{\Sigma} \bar{a}_t + \Theta_t (e_{tu}) + \frac{1}{2}\ast_{\fre_t} (\alpha^{\perp}_t \wedge e_{tu}) = 0
\end{equation*}

\noindent
This is the first equation in \eqref{eq:constraintequationsII} with $\fra_t$ replaced by $\bar{a}_t$. On the other hand, taking its time derivative and using the second equation in \eqref{eq:cohomologicalccondition} we obtain:
\begin{equation*}
\dd_{\Sigma} \partial_t \bar{a}_t = \dd_{\Sigma}(\dd_{\Sigma}\lambda_t (e_{tu}))
\end{equation*}

\noindent
and thus:
\begin{equation*}
\partial_t \bar{a}_t = \dd_{\Sigma}\lambda_t (e_{tu}) + c_t
\end{equation*}

\noindent
where $c_t$ is a family of constants on $\Sigma$. Defining $\fra_t = \bar{a}_t - \int c_t \,\dd t$ we recover the first equations in \eqref{eq:evolutionequationsII} and \eqref{eq:constraintequationsII} and thus we conclude. 
\end{proof}


\subsection{Well-posedness of the Cauchy problem}
\label{sec:skewtorsionwellposedness}


We now formulate and solve the Cauchy problem for skew-torsion parallel spinors in terms of their associated skew-torsion isotropic parallelisms, identifying a minimal set of initial data that we will construct in the process. 

\begin{lemma}
\label{lemma:refor1}
 A pair $(\lambda_t,\fre_t)$ is a solution of the system \eqref{eq:evolutionequationsIII}, \eqref{eq:constraintequationsIII}, and \eqref{eq:cohomologicalccondition} if and only if:
 \begin{eqnarray}
 & \partial_t \frf_t + \dd\lambda_t(e_{tu}^{\sharp}) = 0\, , \nonumber \\  
 & \partial_t \fre_t= \mathcal{L}_{\lambda_t e^{\sharp}_{tu}} \fre_t - (\dd\lambda_t(\fre_t^{\sharp}) + \lambda_t \dd\frf_t(\fre_t^{\sharp}) ) e_{tu} - \lambda_t (\alpha^{\perp}_t (e^{\sharp}_{tu}) - \alpha^o_t) ( \ast_{\fre_t} \fre_t)(e^{\sharp}_{tu} ) \nonumber  \\   
& \dd_{\Sigma} e_{tu} = \dd\frf_t \wedge e_{tu} + \ast_{\fre_t} (\alpha^{\perp}_t   - \alpha^o_t e_{tu})   \label{eq:evolutionconstraintskew} \\
& \dd_{\Sigma} e_{tl} = \Theta_t(e_{tl})\wedge e_{tu} +  (\frac{1}{2} \alpha^{\perp}_t (e_{tu}) - \alpha^o_t) \ast_{\fre_t} e_{tl} \nonumber \\
& \dd_{\Sigma} e_{tn} = \Theta_t(e_{tn})\wedge e_{tu} +  (\frac{1}{2} \alpha^{\perp}_t (e_{tu}) - \alpha^o_t) \ast_{\fre_t} e_{tn} \nonumber  
 \end{eqnarray}
 for a family of smooth functions $\frf_t\in C^{\infty}(\Sigma)$.
\end{lemma}
    
\begin{proof}
Let $(\lambda_t,\fre_t)$ be a solution of equations \eqref{eq:evolutionequationsIII}, \eqref{eq:constraintequationsIII}, and \eqref{eq:cohomologicalccondition}. Evaluating the constraint equations \eqref{eq:constraintequationsIII} on $e_{tu}^{\sharp}$, we obtain:
\begin{align*}
& (\dd e_{tu})(e_{tu}^{\sharp}) = \Theta_t (e_{tu}^{\sharp} , e_{tu}^{\sharp}) e_{tu} -\Theta_t(e_{tu}^{\sharp}) + \frac{1}{2} \ast_{\fre_t} (\alpha^{\perp}_t \wedge e_{tu})\\
& (\dd e_{tl})(e_{tu}^{\sharp}) = \Theta_t (e_{tu}^{\sharp} , e_{tl}^{\sharp}) e_{tu} - \Theta_t (e_{tl}^{\sharp}) - (\frac{1}{2}\alpha^{\perp}_t (e_{tu})  - \alpha^o_t) \ast_{\fre_t}(e_{tu}\wedge e_{tl}) \\
& (\dd e_{tn})(e_{tu}^{\sharp}) = \Theta_t (e_{tu}^{\sharp} , e_{tn}^{\sharp}) e_{tu} - \Theta_t (e_{tn}^{\sharp})- (\frac{1}{2}\alpha^{\perp}_t (e_{tu}) - \alpha^o_t) \ast_{\fre_t}(e_{tu}\wedge e_{tn}) 
\end{align*}
On the other hand, by the cohomological condition \eqref{eq:cohomologicalccondition}, there exists a family of smooth functions $\frf_t$ on $\Sigma$ such that:
\begin{equation}
\label{eq:frfdef}
\dd \frf_t = \Theta_t (e_{tu}) + \frac{1}{2}\ast_{\fre_t} (\alpha^{\perp}_t \wedge e_{tu})
\end{equation}

\noindent
which, plugged into the previous equations, yields:
\begin{align*}
& \Theta_t(e_{tu}^{\sharp})  = \dd\frf_t(e_{tu}^{\sharp}) e_{tu}  - (\dd e_{tu})(e_{tu}^{\sharp}) + \frac{1}{2}\ast_{\fre_t} (\alpha^{\perp}_t \wedge e_{tu})  \\
& \Theta_t (e_{tl}^{\sharp})  =( \dd\frf_t(e_{tl}^{\sharp}) - \frac{1}{2} \alpha^{\perp}_t(e_{tn})) e_{tu}   - (\dd e_{tl})(e_{tu}^{\sharp}) - (\frac{1}{2} \alpha^{\perp}_t (e_{tu}) - \alpha^o_t) e_{tn} \\
& \Theta_t (e_{tn}^{\sharp})   = ( \dd\frf_t(e_{tn}^{\sharp}) + \frac{1}{2} \alpha^{\perp}_t(e_{tl})) e_{tu}  -  (\dd e_{tn})(e_{tu}^{\sharp}) + (\frac{1}{2} \alpha^{\perp}_t (e_{tu}) - \alpha^o_t) e_{tl}
\end{align*}

\noindent
Plugging the first equation into the first equation of \eqref{eq:evolutionequationsIII}, the second equation into the second equation of  \eqref{eq:evolutionequationsIII}, and the third equation into the third equation of  \eqref{eq:evolutionequationsIII}, we conclude:
\begin{eqnarray*}
	&   \partial_t e_{tu} + \lambda_t (\dd\frf_t(e_{tu}^{\sharp}) e_{tu}  - (\dd e_{tu})(e_{tu}^{\sharp})) =    \dd_{\Sigma} \lambda_t (e_{tl}) e_{tl} +   \dd_{\Sigma} \lambda_t (e_{tn}) e_{tn} \nonumber \\
	&   \partial_t e_{tl}  + \lambda_t \dd\frf_t(e_{tl}^{\sharp})  e_{tu}   - \lambda_t (\dd e_{tl})(e_{tu}^{\sharp}) + \lambda_t\alpha^o_t e_{tn}  = \lambda_t \alpha^{\perp}_t(e_{tu}) e_{tn}   -   \dd_{\Sigma} \lambda_t (e_{tl}) e_{tu}     \\
	&  \partial_t e_{tn} + \lambda_t  \dd\frf_t(e_{tn}^{\sharp})  e_{tu}  -  \lambda_t (\dd e_{tn})(e_{tu}^{\sharp})   - \lambda_t \alpha^o_t e_{tl} =   - \lambda_t \alpha^{\perp}_t (e_{tu}) e_{tl} -   \dd_{\Sigma} \lambda_t (e_{tn}) e_{tu}   
\end{eqnarray*}
Using Cartan's magic formula to identify the Lie derivative, these equations can be uniformly written as in the second equation of \eqref{eq:evolutionconstraintskew}. The first equation in \eqref{eq:evolutionconstraintskew} on the other hand, follows directly from the second equation in \eqref{eq:cohomologicalccondition} after redefining $\frf_t$, which by construction is only well-defined modulo the addition of a smooth family of time-dependent constants. The third line in Equation \eqref{eq:evolutionconstraintskew} is the first equation in \eqref{eq:constraintequationsIII} after substituting Equation \eqref{eq:frfdef}. The fourth and fifth lines in \eqref{eq:evolutionconstraintskew} are the second and third equations in \eqref{eq:constraintequationsIII}.   

 \noindent
For the converse, let $(\frf_t,\lambda_t,\fre_t)$ be a solution of equations \eqref{eq:evolutionconstraintskew}. We need to prove that $(\frf_t,\lambda_t,\fre_t)$ satisfies equations  \eqref{eq:evolutionequationsIII}, \eqref{eq:constraintequationsIII} and \eqref{eq:cohomologicalccondition} in the original system. By tracing back the previous steps, we only need to prove that the equations in \eqref{eq:cohomologicalccondition} hold. We compute:
\begin{eqnarray*}
& -2\lambda_t \Theta_t (\fre_{tu}^{\sharp}) =  (\partial_t h_{\fre_t})(\fre_{tu}^{\sharp}) =  \partial_t e_{tu} + (\partial_t e_{tu})(\fre_{tu}^{\sharp}) e_{tu} + (\partial_t e_{tl})(\fre_{tu}^{\sharp}) e_{tl} + (\partial_t e_{tn})(\fre_{tu}^{\sharp}) e_{tn}\\
& =  \cL_{\lambda_t e_{tu}} e_{tu} + (\cL_{\lambda_t e_{tu}} e_{tu}) (e_{tu}) e_{tu} - \dd_{\Sigma}\lambda_t - \lambda_t \dd_{\Sigma} \frf_t - (\dd_{\Sigma}\lambda_t(e_{tu}) + \lambda_t \dd_{\Sigma} \frf_t (e_{tu}))e_{tu} \\
&  = \lambda_t (\dd_{\Sigma} e_{tu} )(e_{tu}) - \lambda_t \dd_{\Sigma} \frf_t (e_{tu}) e_{tu} - \lambda_t \dd_{\Sigma} \frf_t = -2\lambda_t \dd_{\Sigma} \frf_t + \lambda_t\ast_{\fre_t} (\alpha^{\perp}_t \wedge e_{tu})
\end{eqnarray*}
where we have used the third equation in \eqref{eq:evolutionconstraintskew}. Hence we recover Equation \eqref{eq:frfdef}. Combining its time derivative with the first equation in \eqref{eq:evolutionconstraintskew} we on the second equation in \eqref{eq:cohomologicalccondition}, and hence we conclude.
\end{proof}

\begin{prop}
\label{prop:Thetaevolution}
Let $(\lambda_t,\fre_t,\frf_t)$ be a solution of all equations in \eqref{eq:evolutionconstraintskew}. Then:
\begin{eqnarray}
\Theta_t = -\frac{1}{2\lambda_t} \partial_t h_{\fre_t} =  -\frac{1}{2\lambda_t} \cL_{\lambda_t e_{tu}} h_{\fre_t}  + \frac{1}{2} e_{tu}\odot (\dd_{\Sigma} \log(\lambda_t) + \dd_{\Sigma}\frf_t)
\end{eqnarray}

\noindent
where $h_{\fre_t} = e_{tu}\otimes e_{tu} + e_{tl}\otimes e_{tl} + e_{tn}\otimes e_{tn}$.
\end{prop}

\begin{proof}
We compute:
\begin{eqnarray*}
& \partial_t h_{\fre_t} = e_{tu} \odot \partial_t e_{tu} + e_{tl} \odot \partial_t e_{tl} + e_{tn} \odot \partial_t e_{tn} \\
& = e_{tu} \odot \cL_{\lambda_t e_{tu}} e_{tu} + e_{tl} \odot \cL_{\lambda_t e_{tu}} e_{tl} + e_{tn} \odot \cL_{\lambda_t e_{tu}} e_{tn} - e_{tu}\odot ( \dd_{\Sigma}\lambda_t + \lambda_t \dd_{\Sigma} \frf_t) \\
& = \cL_{\lambda_t e_{tu}} h_{\fre_t} - e_{tu}\odot (\dd_{\Sigma}\lambda_t + \lambda_t \dd_{\Sigma}\frf_t)
\end{eqnarray*}
where we have used the evolution equations in \eqref{eq:evolutionconstraintskew}.
\end{proof} 

\begin{lemma}
\label{lemma:refor2}
 A triple $(\lambda_t,\fre_t,\frf_t)$ is a solution of all equations in \eqref{eq:evolutionconstraintskew} if and only if:
 \begin{eqnarray}
 & \partial_t \frf_t + \dd\lambda_t(e_{tu}^{\sharp}) = 0\, , \quad \partial_t \fre_t= \mathcal{L}_{\lambda_t e^{\sharp}_{tu}} \fre_t - (\dd\lambda_t(\fre_t^{\sharp}) + \lambda_t \dd\frf_t(\fre_t^{\sharp}) ) e_{tu} - \lambda_t (\alpha^{\perp}_t (e^{\sharp}_{tu}) - \alpha^o_t) ( \ast_{\fre_t} \fre_t)(e^{\sharp}_{tu} ) \nonumber  \\   
& \dd_{\Sigma} e_{tu} = \dd_{\Sigma}\frf_t \wedge e_{tu} + \ast_{\fre_t} (\alpha^{\perp}_t   - \alpha^o_t e_{tu})  \, , \qquad \dd_{\Sigma} e_{tl} \wedge e_{tu} = 0 \, , \qquad \dd_{\Sigma} e_{tn} \wedge e_{tu} = 0 \label{eq:evolutionconstraintII} \\
& (\cL_{e_{tu}^{\sharp}}e_{tn})(e_{tl}^{\sharp}) - (\cL_{e_{tu}^{\sharp}}e_{tl})(e_{tn}^{\sharp}) = \alpha^{\perp}_t(e_{tu}) - 2\alpha^o_t\nonumber
\end{eqnarray}
 where $\frf_t\in C^{\infty}(\Sigma)$ is a family of smooth functions.
\end{lemma}

\begin{proof}
Assume that $(\lambda_t,\fre_t,\frf_t)$ is a solution of all equations in \eqref{eq:evolutionconstraintskew}. The first line in \eqref{eq:evolutionconstraintII} corresponds to the first and second lines in \eqref{eq:evolutionconstraintskew}. Wedging lines four and five in \eqref{eq:evolutionconstraintskew} with $e_{tu}$ we obtain, together with the third line of   \eqref{eq:evolutionconstraintskew}, all the equations in the second line of \eqref{eq:evolutionconstraintII}. To recover the third line in \eqref{eq:evolutionconstraintII}, we compute:
\begin{eqnarray*}
(\cL_{e_{tu}^{\sharp}}e_{tn})(e_{tl}^{\sharp}) - (\cL_{e_{tu}^{\sharp}}e_{tl})(e_{tn}^{\sharp}) = (\dd_{\Sigma} e_{tn}) (e_{tu} , e_{tl}) - (\dd_{\Sigma} e_{tl})(e_{tu},e_{tn}) = \alpha^{\perp}_t (e_{tu}) - 2\alpha^o_t
\end{eqnarray*}
where we have used the fourth and fifth lines in \eqref{eq:evolutionconstraintskew}.  

\noindent
For the converse, assume that $(\lambda_t,\fre_t,\frf_t)$ is a solution of all equations in \eqref{eq:evolutionconstraintII}. The evolution equations correspond to those of the system \eqref{eq:evolutionconstraintskew}, whereas the first equation in the second line of \eqref{eq:evolutionconstraintII} corresponds to the third line in \eqref{eq:evolutionconstraintskew}. Evaluating the second and third equations in the second line of \eqref{eq:evolutionconstraintII} in $e^{\sharp}_{tu}$, we obtain:
\begin{equation}
\label{eq:auxI}
\dd_{\Sigma} e_{tl} + (\dd_{\Sigma} e_{tl})(e_{tu})\wedge e_{tu} = 0\, , \qquad \dd_{\Sigma} e_{tn} + (\dd_{\Sigma} e_{tn})(e_{tu})\wedge e_{tu} = 0
\end{equation}

\noindent
On the other hand, we compute:
\begin{eqnarray*}
& -2\lambda_t \Theta_t (e_{tl}) \wedge e_{tu} = (\partial_t h_{\fre_t})(e_{tl})\wedge e_{tu} = (\partial_t e_{tl} + (\partial_t e_{tl} )(e_{tl}) e_{tl} + (\partial_t e_{tn} )(e_{tl}) e_{tn} ) \wedge e_{tu} \\ 
&  = \lambda_t ((\dd_{\Sigma} e_{tl})(e_{tu})  + (\dd_{\Sigma} e_{tl})(e_{tu},e_{tl}) e_{tl}  + (\dd_{\Sigma} e_{tn})(e_{tu},e_{tl}) e_{tn} )  \wedge e_{tu}\\
& =  2 \lambda_t (\dd_{\Sigma} e_{tl})(e_{tu}) \wedge e_{tu} - \lambda_t (\alpha^{\perp}_t (e_{tu}) - 2 \alpha^o_t) e_{tu}\wedge e_{tn} = - 2 \lambda_t  \dd_{\Sigma} e_{tl} - \lambda_t (\alpha^{\perp}_t (e_{tu}) - 2 \alpha^o_t) e_{tu}\wedge e_{tn}
\end{eqnarray*}
which therefore gives the fourth line in \eqref{eq:evolutionconstraintskew}. Here, we have used the evolution equations in \eqref{eq:evolutionconstraintII}, together with the third line of \eqref{eq:evolutionconstraintII} and the first equation in \eqref{eq:auxI}. The fifth line of \eqref{eq:evolutionconstraintskew} follows similarly, using the second equation in \eqref{eq:auxI} instead.
\end{proof} 

\noindent
We are ready to establish the well-posedness of the Cauchy problem determined by a skew-torsion parallel spinor on a Lorentzian four-manifold. Direct inspection of Lemma \ref{lemma:refor1} shows that the \emph{initial data} that needs to be specified to solve the Cauchy problem of the evolution and time-dependent constraint equations a priori consists of tuples of the form $(\fre,\frf,\alpha^o,\alpha^{\perp})$, where $\fre$ is a global coframe, $\frf , \alpha^o \in C^{\infty}(\Sigma)$ are functions, and $\alpha^{\perp} \in \Omega^1(\Sigma)$ is a one-form.

\begin{thm}
\label{thm:wellposedFinalSkew}
Let $\Sigma$ be a connected and oriented three-manifold admitting a tuple $(\fre,\frf,\alpha^o,\alpha^{\perp})$ satisfying the following differential system: 
 \begin{eqnarray}
& \dd_{\Sigma} ( e^{-\frf} e_{u}) = e^{-\frf} \ast_{\fre} (\alpha^{\perp}   - \alpha^o e_{u})  \, , \qquad \dd_{\Sigma} e_{l} \wedge e_{u} = 0 \, , \qquad \dd_{\Sigma} e_{n} \wedge e_{u} = 0 \label{eq:constraintIIthm1} \\
& \dd_{\Sigma} (e^{-\frf}(\alpha^{\perp}(e_u^{\sharp}) - \alpha^o)) = 0\, , \qquad   (\cL_{e_{u}^{\sharp}}e_{n})(e_{l}^{\sharp}) - (\cL_{e_{u}^{\sharp}}e_{l})(e_{n}^{\sharp}) = \alpha^{\perp}(e_{u}^{\sharp}) - 2\alpha^o  \label{eq:constraintIIthm2}
\end{eqnarray}
where $\fre^{\sharp} = (e^{\sharp}_u,e^{\sharp}_l,e^{\sharp}_n)$ is the frame dual to $\fre$. Then, there exists a unique germ of a globally hyperbolic Lorentzian four-manifold $(M,g)$ equipped with a skew-torsion parallel spinor with skew-torsion $\alpha \in \Omega^1(M)$ and an isometric embedding $\iota\colon (\Sigma,h_{\fre})\hookrightarrow (M,g)$ such that $(\Sigma,h_{\fre})$ is a Cauchy hypersurface in $(M,g)$ and:
\begin{equation*}
\iota^{\ast}\alpha = \alpha^{\perp}\, , \qquad \alpha^o = - \alpha(N)
\end{equation*}

\noindent
for the past-directed unit normal vector $N$, and the globally hyperbolic reduction of the associated isotropic parallelism has initial values $(\fra_0,\fre_0)=(-\frf,\fre)$. Conversely, every smooth spacelike Cauchy hypersurface $\Sigma$ of a globally hyperbolic Lorentzian four-manifold $(M,g)$ equipped with a skew-torsion parallel spinor of real type carries a solution $(\fre,\frf,\alpha^o,\alpha^{\perp})$ of the above system \eqref{eq:constraintIIthm1} and \eqref{eq:constraintIIthm2}. 
\end{thm}   

\begin{remark}
 Uniqueness is understood up to an orientation- and time-orientation-preserving isometry fixing the embedded initial hypersurface and preserving $\alpha$ and the isotropic parallelism, hence the spinor up to an overall sign.
\end{remark}

\begin{proof}
Let $(\fre,\frf,\alpha^o,\alpha^{\perp})$ be a solution of all equations in \eqref{eq:constraintIIthm1} and \eqref{eq:constraintIIthm2}. First, we need to prove that there exists a smooth family $(\fre_t,\frf_t,\alpha^o_t,\alpha^{\perp}_t)$ such that:
\begin{equation*}
(\fre_t,\frf_t,\alpha^o_t,\alpha^{\perp}_t) \vert_{t=0} = (\fre,\frf,\alpha^o,\alpha^{\perp})
\end{equation*}

\noindent
and furthermore there exists a lapse function $\lambda_t$ for which the tuple $(\fre_t,\frf_t,\alpha^o_t,\alpha^{\perp}_t, \lambda_t)$ satisfies all equations in \eqref{eq:evolutionconstraintII}. Using the freedom in the choice of the lapse function to integrate the flow, we make the educated choice $\lambda_t = e^{-\frf_t}$. With this choice, the system of evolution and constraint equations \eqref{eq:evolutionconstraintII} that needs to be solved reduces to: 
\begin{eqnarray}
 & \partial_t \frf_t =  \cL_{\cX_t}\frf_t \, , \qquad \partial_t e_{tu} = \mathcal{L}_{\cX_t} e_{tu}\, , \qquad  \partial_t \zeta_{t} = \mathcal{L}_{\cX_t} \zeta_t   - i (\alpha^{\perp}_t (\cX_t) - e^{-\frf_t}\alpha^o_t)\, \zeta_t   \nonumber \\ 
& \dd_{\Sigma} ( e^{-\frf_t} e_{tu}) = e^{-\frf_t} \ast_{\fre_t} (\alpha^{\perp}_t   - \alpha^o_t e_{tu})  \, , \qquad \dd_{\Sigma} e_{tl} \wedge e_{tu} = 0 \, , \qquad \dd_{\Sigma} e_{tn} \wedge e_{tu} = 0  \label{eq:auxiliarysystem} \\
& (\cL_{\cX_t}e_{tn})(e_{tl}^{\sharp}) - (\cL_{\cX_t}e_{tl})(e_{tn}^{\sharp}) = \alpha^{\perp}_t(\cX_t) - 2e^{-\frf_t} \alpha^o_t\nonumber
\end{eqnarray}

\noindent
where $\cL$ denotes the Lie derivative symbol, and we have introduced the vector field $\cX_t = e^{-\frf_t} e^{\sharp_{\fre_t}}_{tu}$ and the complexified one-form $\zeta_t = e_{tl} + \mathrm{i} e_{tn}$. To construct the desired solution starting in the initial data $(\fre,\frf,\alpha^o,\alpha^{\perp})$, we will reinterpret the system in terms of the flow of $\cX_t$, which will determine how to evolve any given initial data and will reveal a hidden obstruction in the system. Suppose that $(\fre_t,\frf_t,\alpha^o_t,\alpha^{\perp}_t)$ is a solution of all equations in \eqref{eq:auxiliarysystem}. A direct calculation via substitution gives:
\begin{equation*}
\partial_t h_{\fre_t} = e_{tu}\odot \partial_t e_{tu} + e_{tl}\odot \partial_t e_{tl} + e_{tn}\odot \partial_t e_{tn} = \cL_{\cX_t} h_{\fre_t}
\end{equation*}

\noindent
This, together with equations $\partial_t \frf_t =  \cL_{\cX_t}\frf_t$ and $\partial_t e_{tu} = \mathcal{L}_{\cX_t} e_{tu}$ in the first line of \eqref{eq:auxiliarysystem}, implies that:
\begin{eqnarray*}
 & \partial_t \cX_t = \partial_t (e^{-\frf_t} e^{\sharp_{\fre_t}}_{tu})  =  -  \cL_{\cX_t}\frf_t\, \cX_t + e^{-\frf_t}  \partial_t (e^{\sharp_{\fre_t}}_{tu}) \\
 & =  -  \cL_{\cX_t}\frf_t\, \cX_t + e^{-\frf_t}  ((\partial_t e_{tu})^{\sharp_{\fre_t}} -   (\cL_{\cX_t} e_{tu})^{\sharp_{\fre_t}} + \cL_{\cX_t} e_{tu}^{\sharp_{\fre_t}}) = -  \cL_{\cX_t}\frf_t\, \cX_t + e^{-\frf_t} \cL_{\cX_t} e_{tu}^{\sharp_{\fre_t}} = 0
\end{eqnarray*}

\noindent
Therefore we can set $\cX := \cX_t = e^{-\frf} e^{\sharp}_u$ to be given by its value at $t=0$ for every $t$. Given initial data $(\fre,\frf,\alpha^o,\alpha^{\perp})$, we denote by:
\begin{eqnarray*}
\Phi \colon \cU \to \Sigma
\end{eqnarray*}

\noindent
the maximal flow of the vector field $\cX\in \mathfrak{X}(\Sigma)$. Here $\cU \subset \mathbb{R} \times \Sigma$ is the maximal open tubular neighborhood of $\Sigma$ where the flow is well-defined. Denoting by $t$ the Cartesian coordinate of $\mathbb{R}$, the points of $\cU$ are of the form $(t,p)$, for $t\in\mathbb{R}$ and $p\in \Sigma$. All subsequent flow formulas are understood on their common domains, allowing the time interval to depend on $p$. The first two equations in the first line of \eqref{eq:auxiliarysystem} are Lie transport equations and are immediately solved by:
\begin{eqnarray*}
\frf_t (p) := (\Phi^{\ast}_t\frf) (p)\, , \qquad \  e_{tu}(p) = (\Phi^{\ast}_t e_u)(p) \, .
\end{eqnarray*}

\noindent
The third equation in the first line of \eqref{eq:auxiliarysystem} is again of Lie-transport type but with a zero-order homogeneous term. Its general solution is:
\begin{eqnarray*}
\zeta_t = e^{- i \int_0^t \Phi^{\ast}_{t-s} (\alpha^{\perp}_s (\cX) - e^{-\frf_s} \alpha^o_s)\, \dd s}\, \Phi_t^{\ast} \zeta    
\end{eqnarray*} 

\noindent 
Conversely, given $\alpha^o_t$ and $\alpha^{\perp}_t$, a direct computation shows that every such tuple $(\frf_t , \fre_t)$ solves the evolution equations in \eqref{eq:auxiliarysystem}. The unit complex phase preserves the metric and orientation, so $h_{\fre_t}=\Phi_t^{\ast}h_{\fre}$ and the Hodge operator commutes with these pull-backs. Substituting the solution of the evolution equations into the first constraint equation of \eqref{eq:auxiliarysystem}, we obtain:
\begin{equation}
\label{eq:transversalpullback}
\alpha^{\perp}_t   - \alpha^o_t e_{tu}   = \Phi^{\ast}_t (\alpha^{\perp} - \alpha^o e_{u})   
\end{equation}

\noindent
which immediately implies $\alpha^{\perp}_t (\cX) - e^{-\frf_t} \alpha^o_t = \Phi^{\ast}_t (\alpha^{\perp} (\cX) - e^{-\frf} \alpha^o)$, and thus:
\begin{eqnarray}
\label{eq:zetaux}
\zeta_t = \Phi_t^{\ast} (e^{- i (\alpha^{\perp} (\cX) - e^{-\frf} \alpha^o) t}\, \zeta)    
\end{eqnarray}
The second and third constraint equations in \eqref{eq:auxiliarysystem} can be written as $\dd_{\Sigma}\zeta_t \wedge e_{tu} = 0$. We compute:
\begin{equation*}
0 = \dd_{\Sigma}\zeta_t \wedge e_{tu} = \Phi_t^{\ast} [e^{- i (\alpha^{\perp} (\cX) - e^{-\frf} \alpha^o) t} (- \mathrm{i} t\, \dd_{\Sigma} (\alpha^{\perp} (\cX) - e^{-\frf} \alpha^o) \wedge \zeta \wedge e_u + \dd_{\Sigma}\zeta \wedge e_u) ]
\end{equation*}

\noindent
Since we have $\dd_{\Sigma}\zeta \wedge e_u = 0$ initially, taking real and imaginary parts for sufficiently small nonzero $t$ gives, locally at every point of $\Sigma$:
\begin{equation*}
\dd_{\Sigma} (\alpha^{\perp} (\cX) - e^{-\frf} \alpha^o) \wedge e_u = 0\quad \Rightarrow \quad \dd_{\Sigma} (\alpha^{\perp} (\cX) - e^{-\frf} \alpha^o)  = f e_u
\end{equation*}

\noindent
for a smooth function $f\in C^{\infty}(\Sigma)$. Taking the exterior derivative of this equation and wedging with $e_u$ gives: 
\begin{equation*}
0 = f \,\dd_{\Sigma} e_u \wedge e_u = f\, ( \alpha^{\perp} (e^{\sharp}_u) - \alpha^o)\, e_u\wedge e_l \wedge e_n
\end{equation*}

\noindent
and therefore:
\begin{equation*}
f\, (\alpha^{\perp} (\cX) - e^{-\frf} \alpha^o) = 0
\end{equation*}
Hence:
\begin{equation*}
\frac{1}{2}\dd_{\Sigma} (\alpha^{\perp} (\cX) - e^{-\frf} \alpha^o)^2 = (\alpha^{\perp} (\cX) - e^{-\frf} \alpha^o) \dd_{\Sigma} (\alpha^{\perp} (\cX) - e^{-\frf} \alpha^o) = f\, (\alpha^{\perp} (\cX) - e^{-\frf} \alpha^o) e_u = 0
\end{equation*}
Since $\Sigma$ is connected, continuity implies that $(\alpha^{\perp} (\cX) - e^{-\frf} \alpha^o) \in C^{\infty}(\Sigma)$ is constant on $\Sigma$ and thus $f=0$ identically on $\Sigma$. For the third line in \eqref{eq:auxiliarysystem}, we need to solve the following equation:
\begin{equation*}
 (\cL_{\cX}e_{tn})(e_{tl}^{\sharp}) - (\cL_{\cX}e_{tl})(e_{tn}^{\sharp}) = r - e^{-\frf_t} \alpha^o_t     
\end{equation*}

\noindent
where $ r := \alpha^{\perp} (\cX) - e^{-\frf} \alpha^o \in \mathbb{R}$ is the constant determined above. The rotation by $rt$ in Equation \eqref{eq:zetaux} is constant on each slice, and the left-hand side of the preceding equation is invariant under this rotation. Substituting for $e_{tl}$ and $e_{tn}$ and using the initial scalar constraint, we obtain:
\begin{equation*}
\Phi_t^{\ast} ( e^{-\frf} \alpha^o) =  e^{-\frf \circ \Phi_t} \alpha^o_t         
\end{equation*}

\noindent
from which we conclude that $\alpha_t^o = \Phi_t^{\ast}\alpha^o$. That is, the family of functions $\alpha^o_t$ must necessarily be given as the evolution of $\alpha^o$ via the flow of $\cX = e^{-\frf} e^{\sharp}_u$. This, together with Equation \eqref{eq:transversalpullback}, implies the analog condition for $\alpha^{\perp}_t$, that is:
\begin{equation*}
\alpha^{\perp}_t = \Phi^{\ast}_t \alpha^{\perp}
\end{equation*}

\noindent
Conversely, define:
\begin{eqnarray*}
& \frf_t=\Phi_t^{\ast}\frf\, , \qquad e_{tu}=\Phi_t^{\ast}e_u\, , \qquad
e_{tl}+\mathrm{i}e_{tn}=e^{-\mathrm{i}rt}\Phi_t^{\ast}(e_l+\mathrm{i}e_n)\, , \\
& \alpha_t^o=\Phi_t^{\ast}\alpha^o\, , \qquad \alpha_t^{\perp}=\Phi_t^{\ast}\alpha^{\perp}\, , \qquad \lambda_t=e^{-\frf_t}.
\end{eqnarray*}

\noindent
These variables have the prescribed initial values and verify all equations in \eqref{eq:auxiliarysystem}. By Lemmas \ref{lemma:refor1} and \ref{lemma:refor2}, with $\fra_t=-\frf_t$, these variables reconstruct the required skew-torsion parallelism and spinor. The preceding derivation also proves uniqueness in the gauge $\lambda_t=e^{-\frf_t}$. For clarity in the exposition, we prove the uniqueness of the germ in an arbitrary gauge below in Lemma \ref{lemma:uniqueevol}. The tubular neighborhood $\cU$ contains $\Sigma=\{0\}\times\Sigma$ as a closed, spacelike, acausal and edgeless surface with respect to the reconstructed Lorentzian metric, since $\dd t$ is time-like. The reconstructed metric on $\cU$ need not be globally hyperbolic. Nevertheless, by Theorem 3.45 in \cite{Minguzzi:2019mbe}, $M=\operatorname{Int}D(\Sigma)\subset\cU$ is globally hyperbolic with the induced metric. Since $\Sigma$ is edgeless, $M$ contains $\Sigma$ as a Cauchy hypersurface and gives the desired germ. \medskip

\noindent
For the converse, assume that $(M,g)$ is a globally hyperbolic Lorentzian four-manifold equipped with a skew-torsion parallel spinor, with skew-torsion $\alpha\in\Omega^1(M)$, and let $\Sigma$ be a smooth spacelike Cauchy hypersurface. Applying the null-current flow construction in the proof of Lemma \ref{lemma:uniqueevol} presented below, and composing with the smooth map $((t,p)\mapsto(t,\Phi_t(p))$, where $\Phi_t$ denotes the local flow of $\cX=e^{-\frf}e_u^\sharp$, we can represent the given development in the gauge $\lambda_t=e^{-\frf_t}$, with $\frf_t=\Phi_t^*\frf$. Then, Corollary \ref{cor:Thetaevolutionsystem} and Lemmas \ref{lemma:refor1} and \ref{lemma:refor2}, with $\frf_t=-\fra_t$, now give a solution of \eqref{eq:auxiliarysystem}. Evaluating these equations at $t=0$, gives a solution to all equations in \eqref{eq:constraintIIthm1} and \eqref{eq:constraintIIthm2}, except for the first equation in \eqref{eq:constraintIIthm2}. The preceding calculation, in which we did not assume the first equation of \eqref{eq:constraintIIthm2}, proves that equation as well.
\end{proof} 

\noindent
By the previous theorem, we obtain the following picture for the Cauchy problem posed by a skew-torsion parallel spinor. Let $(M,g)$ be a globally hyperbolic Lorentzian four-manifold equipped with a skew-torsion parallel spinor $\varepsilon$ with torsion $\alpha$. For every smooth spacelike Cauchy hypersurface $\Sigma \subset M$, there exists a tuple $(\fre,\frf,\alpha^{\perp},\alpha^o)$ which is equivalent to the restriction of $(g,\varepsilon,\alpha)$ to $\Sigma$ and satisfies equations \eqref{eq:constraintIIthm1} and \eqref{eq:constraintIIthm2}. Assume for simplicity that $\cX = e^{-\frf} e^{\sharp}_u \in \mathfrak{X}(\Sigma)$ is complete. Then, the germ along $\Sigma$ admits a smooth extension to $\mathbb{R}\times\Sigma$, not necessarily globally hyperbolic or contained in the original $M$, whose metric and skew-torsion read:
\begin{eqnarray*}
g = - e^{- 2 \frf \circ \Phi_t} \dd t \otimes \dd t + \Phi^{\ast}_t h\, , \qquad \alpha = \Phi^{\ast}_t (\alpha^o e^{-\frf})\,\dd t + \Phi^{\ast}_t\alpha^{\perp} 
\end{eqnarray*}
where $\Phi\colon \mathbb{R}\times \Sigma \to \Sigma$ is the flow of $\cX$. Furthermore, by Proposition \ref{prop:Thetaevolution}, the second fundamental form $\Theta_t$ of the embedding of $(\Sigma,h_{\fre_t})$ as the slice $\{t\}\times\Sigma$ in this extension is given by:
\begin{equation*}
\Theta_t = - \frac{1}{2} e^{\frf \circ \Phi_t} \Phi_t^{\ast} \cL_{\cX}h_{\fre}
\end{equation*}

\noindent
where $\cX = e^{-\frf} e^{\sharp}_u$. Interestingly enough, the first equation in \eqref{eq:constraintIIthm2} does not arise from a time-dependent constraint on the reduction of the skew-torsion parallel spinorial system, but rather as a genuine obstruction to solving the evolution and time-dependent constraints starting at given data. Recalling that:
\begin{equation}
\label{eq:skewscreenrotationinvariant}
 \alpha(u^{\sharp_g})  =e^{-\frf_t} (\alpha_t^\perp(e^{\sharp}_{tu})-\alpha_t^o )
\end{equation}

\noindent
it then follows that $\alpha(u^{\sharp_g})$ is constant on $M$. Since the normalization and necessity argument above applies near every spacelike hypersurface patch, this is a local compatibility condition, independent of global hyperbolicity. This gives the following corollary.

\begin{cor}
Let $(g,\varepsilon)$ be a skew-torsion parallel spinor on a connected and oriented four-manifold $M$, and let $u\in\Omega^1(M)$ be the Dirac current of $\varepsilon$. Then $\alpha(u^{\sharp})$ is a real constant on $M$.
\end{cor} 

\noindent
We address now the uniqueness of the germ stated in Theorem \ref{thm:wellposedFinalSkew}.   

\begin{lemma}
\label{lemma:killingreformulation}
A tuple $(\lambda_t,\fre_t,\frf_t)$ satisfies the evolution equations in the first line of \eqref{eq:evolutionconstraintII} if and only if:
\begin{equation}
\label{eq:killingreformulation}
\cL_{u} g = 0\, , \qquad \cL_{u}(u\wedge\zeta) = \mathrm{i}\, r_t\, u\wedge\zeta\, .
\end{equation}

\noindent
where:
\begin{eqnarray*}
& g = -\lambda_t^2\, \dd t\otimes\dd t + h_{\fre_t}\, , \qquad u = e^{-\frf_t}(\lambda_t\, \dd t + e_{tu})\, , \qquad \zeta = e_{tl} + \mathrm{i}\, e_{tn}\, ,\qquad   r_t = e^{-\frf_t} ( \alpha^{\perp}_t(e^{\sharp}_{tu}) - \alpha^o_t ) 
\end{eqnarray*}
In that case, $\cL_{u} u = 0$ and $\iota_u \dd u = 0$.
\end{lemma}

\begin{proof}
For this proof, write:
\begin{equation*}
 [e^{\sharp}_{ta},e^{\sharp}_{tb}] = C^c_{ab}e^{\sharp}_{tc}\, , \qquad   \lambda_a = \dd_{\Sigma}\lambda_t(e^{\sharp}_{ta})\, , \qquad \frf_a = \dd_{\Sigma}\frf_t(e^{\sharp}_{ta}) \, , \qquad \dot{e}_{ab} = (\partial_t e_{ta})(e^{\sharp}_{tb})
\end{equation*}

\noindent
and $\epsilon_{ln} = -\epsilon_{nl} = 1$, $\epsilon_{ab} = 0$ otherwise, for $a,b,c\in\left\{ u,l,n\right\}$. Since:
\begin{equation*}
(\cL_{\lambda_t e^{\sharp}_{tu}}e_{ta})(e^{\sharp}_{tb}) = -\lambda_t C^a_{ub} + \lambda_b\delta_{au}    \, , \qquad (\ast_{\fre_t}e_{ta})(e^{\sharp}_{tu},e^{\sharp}_{tb}) = -\epsilon_{ab}   
\end{equation*}

\noindent
the evolution equations in \eqref{eq:evolutionconstraintII} read:
\begin{equation}
\label{eq:evolutionframecomponents}
\partial_t\frf_t = -\lambda_u\, , \qquad \dot{e}_{ab} = -\lambda_t C^a_{ub} + \lambda_b\delta_{au} - (\lambda_a + \lambda_t\frf_a)\delta_{ub} + \lambda_t ( \alpha^{\perp}_t(e^{\sharp}_{tu}) - \alpha^o_t )\epsilon_{ab}\, .
\end{equation}

\noindent
Set $\mu = e^{-\frf_t}\lambda_t^{-1}$, so that $u^{\sharp_g} = -\mu\, \partial_t + e^{-\frf_t}e^{\sharp}_{tu}$. Differentiating $e_{ta}(e^{\sharp}_{tb}) = \delta_{ab}$ in $t$ gives:
\begin{equation*}
e_{ta}(\partial_t e^{\sharp}_{tb}) = -\dot{e}_{ab}  
\end{equation*}

\noindent
whence:
\begin{equation*}
[u^{\sharp_g},\partial_t] = -\partial_t u^{\sharp_g} = \partial_t\mu\, \partial_t + e^{-\frf_t}\partial_t\frf_t\, e^{\sharp}_{tu} - e^{-\frf_t}\partial_t e^{\sharp}_{tu}
\end{equation*}

\noindent
and:
\begin{equation*}
[u^{\sharp_g} , e^{\sharp}_{tb}] = -\mu\, \partial_te^{\sharp}_{tb} + e^{\sharp}_{tb}(\mu)\partial_t + e^{-\frf_t}C^c_{ub}e^{\sharp}_{tc} + e^{-\frf_t}\frf_b e^{\sharp}_{tu}    
\end{equation*}

\noindent
From this, a calculation gives:
\begin{eqnarray*}
& (\cL_{u}g)(\partial_t,\partial_t) = -2\lambda_t e^{-\frf_t} ( \partial_t\frf_t + \lambda_u )\, , \qquad (\cL_{u}g)(\partial_t,e^{\sharp}_{tb}) = -e^{-\frf_t} ( \dot{e}_{bu} + \lambda_b + \lambda_t\frf_b + \partial_t\frf_t\, \delta_{ub} )\, , \\
& (\cL_{u}g)(e^{\sharp}_{ta},e^{\sharp}_{tb}) = -\mu ( \dot{e}_{ab} + \dot{e}_{ba} ) - e^{-\frf_t} ( C^a_{ub} + C^b_{ua} + \frf_a\delta_{ub} + \frf_b\delta_{ua} )\, .
\end{eqnarray*}

\noindent
Hence $\cL_{u}g = 0$ if and only if:
\begin{eqnarray*}
& \partial_t\frf_t = -\lambda_u \, , \qquad \dot{e}_{bu} = \lambda_u\delta_{ub} - \lambda_b - \lambda_t\frf_b\, , \qquad   \dot{e}_{ab} + \dot{e}_{ba} = -\lambda_t(C^a_{ub} + C^b_{ua}) - \lambda_t(\frf_a\delta_{ub} + \frf_b\delta_{ua}) \\
& \dot{e}_{ab} + \dot{e}_{ba} = -\lambda_t(C^a_{ub} + C^b_{ua}) - \lambda_t(\frf_a\delta_{ub} + \frf_b\delta_{ua})
\end{eqnarray*}

\noindent
which are precisely the first equation in \eqref{eq:evolutionframecomponents} together with the components $\dot{e}_{bu}$, $\dot{e}_{ub}$ and the symmetric part of the second. Assume all these are satisfied. Note that $\cL_{u}(u\wedge\zeta) = u\wedge\cL_{u}\zeta$ since $u$ is isotropic and Killing. Using $\zeta(\partial_t e^{\sharp}_{tb}) = -\dot{e}_{lb} - \mathrm{i}\, \dot{e}_{nb}$ and the brackets above, yields:
\begin{equation*}
(\cL_{u}\zeta)(e^{\sharp}_{tb}) = -\mu ( \dot{e}_{lb} + \mathrm{i}\, \dot{e}_{nb} ) - e^{-\frf_t} ( C^l_{ub} + \mathrm{i}\, C^n_{ub} )\, , \qquad (\cL_{u}\zeta)(\partial_t) = -e^{-\frf_t} ( \dot{e}_{lu} + \mathrm{i}\, \dot{e}_{nu} ) = \lambda_t (\cL_{u}\zeta)(e^{\sharp}_{tu})\, .
\end{equation*}

\noindent
The second identity states that $\cL_{u}\zeta$ evaluated on $\partial_t$ and $e^{\sharp}_{tu}$ is proportional to $u$, which also annihilates $e^{\sharp}_{tl}$ and $e^{\sharp}_{tn}$. Therefore $u\wedge(\cL_{u}\zeta - \mathrm{i}\, r_t\zeta) = 0$ if and only if $(\cL_{u}\zeta)(e^{\sharp}_{tb}) = \mathrm{i}\, r_t\zeta(e^{\sharp}_{tb})$ for $b\in\left\{ l,n\right\}$, that is, if and only if:
\begin{equation*}
\dot{e}_{lb} + \mathrm{i}\, \dot{e}_{nb} = -\lambda_t ( C^l_{ub} + \mathrm{i}\, C^n_{ub} ) - \mathrm{i}\, \lambda_t ( \alpha^{\perp}_t(e^{\sharp}_{tu}) - \alpha^o_t ) ( \delta_{lb} + \mathrm{i}\, \delta_{nb} )\, , \qquad b\in\left\{ l,n\right\}\, ,
\end{equation*}

\noindent
which are exactly the remaining components $\dot{e}_{ab}$, $a,b\in\left\{ l,n\right\}$, of \eqref{eq:evolutionframecomponents}. This proves the equivalence, together with the last assertion.
\end{proof}

\begin{lemma}
\label{lemma:uniqueevol}
The germ of a development of any given initial data $(\fre,\frf,\alpha^o,\alpha^{\perp})$ in Theorem \ref{thm:wellposedFinalSkew} is locally unique, up to an isometry fixing the initial hypersurface.
\end{lemma}
 
\begin{proof}
Let $(M_i,g_i,[u_i,v_i,l_i,n_i],\alpha_i,\iota_i)$, $i = 1,2$, be two developments of the same initial datum $(\fre,\frf,\alpha^o,\alpha^{\perp})$. By the globally hyperbolic reduction of Section \ref{sec:ReductionGloballyHyperbolicparallelisms} and Lemma \ref{lemma:refor2}, after restricting to open neighborhoods of $\iota_i(\Sigma)$ we may assume that $M_i = \cI_i\times\Sigma$, that $\iota_i$ is the inclusion of $\left\{ 0\right\}\times\Sigma$, and that $g_i$, the adapted representative of $[u_i,v_i,l_i,n_i]$ and $\alpha_i$ are reconstructed from a solution $(\lambda^{(i)}_t,\fre^{(i)}_t,\frf^{(i)}_t,\alpha^{o(i)}_t,\alpha^{\perp(i)}_t)$ of \eqref{eq:evolutionconstraintII} as in Lemma \ref{lemma:killingreformulation}, namely:
\begin{equation}
\label{eq:initialcurrentnormalization}
g_i = -(\lambda^{(i)}_t)^2\, \dd t\otimes \dd t + h_{\fre^{(i)}_t}\, , \qquad u_i = e^{-\frf^{(i)}_t} ( \lambda^{(i)}_t\, \dd t + e^{(i)}_{tu} )\, , \qquad l_i + \mathrm{i}\, n_i = e^{(i)}_{tl} + \mathrm{i}\, e^{(i)}_{tn} =: \zeta_i\, ,
\end{equation}

\noindent
with:
\begin{equation*}
\alpha_i = \alpha^{o(i)}_t\lambda^{(i)}_t\, \dd t + \alpha^{\perp(i)}_t    
\end{equation*}

\noindent
and common initial values:
\begin{equation*}
(\fre^{(i)}_0,\frf^{(i)}_0,\alpha^{o(i)}_0,\alpha^{\perp(i)}_0) = (\fre,\frf,\alpha^o,\alpha^{\perp})    
\end{equation*}

\noindent
whereas the lapse functions $\lambda^{(i)}_t$ may be unrelated. Set:
\begin{equation*}
r^{(i)}_t = e^{-\frf^{(i)}_t}(\alpha^{\perp(i)}_t(e^{(i)\sharp}_{tu}) - \alpha^{o(i)}_t)    
\end{equation*}

\noindent
By Lemma \ref{lemma:killingreformulation}, $u_i^{\sharp_{g_i}}$ is a nowhere vanishing null Killing vector field of $g_i$, transverse to $\left\{ 0\right\}\times\Sigma$, whose flow preserves $u_i$ and satisfies:
\begin{equation*}
\cL_{u_i^{\sharp_{g_i}}}(u_i\wedge\zeta_i) = \mathrm{i}\, r^{(i)}_t u_i\wedge\zeta_i    
\end{equation*}

\noindent
Moreover, the first constraint equation in \eqref{eq:evolutionconstraintII} evaluated on $(e^{(i)\sharp}_{tl},e^{(i)\sharp}_{tn})$ gives $r^{(i)}_t = \dd u_i(e^{(i)\sharp}_{tl},e^{(i)\sharp}_{tn})$, which together with $\iota_{u_i^{\sharp_{g_i}}}\dd u_i = 0$ is equivalent to:
\begin{equation*}
u_i\wedge\dd u_i = -r^{(i)}_t\, \ast_{g_i}u_i\, ,
\end{equation*}

\noindent
so that $r^{(i)}_t$ is determined by $g_i$, $u_i$ and the orientation. Its local flow defines a diffeomorphism:
\begin{equation*}
\Psi_i\colon\mathcal V\subset\mathbb R\times\Sigma\longrightarrow M_i,\qquad
\Psi_i(s,p)=\phi_{-s}^{u_i^{\sharp_{g_i}}}(0,p),
\end{equation*}

\noindent
onto an open neighborhood of $\left\{ 0\right\}\times\Sigma$, after restricting to a common flow neighborhood $\mathcal V$ of $\left\{ 0\right\}\times\Sigma$. The time interval within this common neighborhood may depend on $p$, since we do not assume completeness. Since these flows preserve $g_i$ and $u_i$, for $X,Y\in T_p\Sigma$ we obtain:
\begin{eqnarray*}
& (\Psi_i^{\ast}g_i)(X,Y)=h_{\fre}(X,Y),\qquad
(\Psi_i^{\ast}g_i)(\partial_s,X)=-e^{-\frf}e_u(X),\qquad
(\Psi_i^{\ast}g_i)(\partial_s,\partial_s)=0,
\end{eqnarray*}

\noindent
whence:
\begin{equation*}
\Psi_i^{\ast}g_i=h_{\fre}-e^{-\frf}\dd s\odot e_u,\qquad
\Psi_i^{\ast}u_i=e^{-\frf}e_u\, , \qquad \Psi_i^{\ast}u_i^{\sharp_{g_i}} = -\partial_s\, .
\end{equation*}

\noindent 
Hence, we may pull back all the variables to $\mathcal V$ and denote the common metric and null one-form by $g$ and $u$. By the preceding paragraph, $\Psi_1^{\ast}r^{(1)} = \Psi_2^{\ast}r^{(2)} =: r$. The two-forms $\omega_i = \Psi_i^{\ast}(u_i\wedge\zeta_i)$ satisfy $\iota_{\partial_s}\omega_i = 0$, since $u_i$ is isotropic, $\zeta_i(u_i^{\sharp_{g_i}}) = 0$ and, furthermore, the following transport equation is satisfied: 
\begin{equation*}
\cL_{\partial_s}\omega_i = -\mathrm{i}\, r\, \omega_i\, ,
\end{equation*}

\noindent
with the same initial value:
\begin{equation*}
\omega_i|_{\left\{ 0\right\}\times\Sigma} = e^{-\frf}e_u\wedge(e_l + \mathrm{i}\, e_n)
\end{equation*}

\noindent
In a fixed basis of two-forms, this is a linear ordinary differential equation along each line $\mathbb{R}\times\left\{ p\right\}$, whence $\omega_1 = \omega_2$. Taking real parts, we obtain:
\begin{equation*}
  u\wedge l_1 = u\wedge l_2  
\end{equation*}

\noindent
so $l_1 - l_2$ is a multiple of $u$ and therefore the parabolic pairs $(u,[l_i]_u)$ coincide. Since $g$ and the orientation are the same, the isotropic parallelisms coincide, and by Remark \ref{remark:torsiondetermined} so do $\Psi_1^{\ast}\alpha_1$ and $\Psi_2^{\ast}\alpha_2$. Returning to the original developments, $F:=\Psi_2\circ\Psi_1^{-1}$ satisfies:
\begin{equation*}
F\circ\iota_1=\iota_2,\qquad F^{\ast}g_2=g_1,\qquad F^{\ast}\alpha_2=\alpha_1,\qquad
F^{\ast}[u_2,v_2,l_2,n_2]=[u_1,v_1,l_1,n_1].
\end{equation*}

\noindent
The isometry preserves orientation and time orientation by its initial values and the identification of the null currents. This proves uniqueness of the germ independently of the gauge.
\end{proof}  

\begin{remark}
\label{remark:straightening}
For a single development $(M,g,[u,v,l,n],\alpha,\iota)$ of a datum $(\fre,\frf,\alpha^o,\alpha^{\perp})$, the flow of $u^{\sharp_g}$ defines a diffeomorphism $\Psi(s,p) = \phi^{u^{\sharp_g}}_{-s}(\iota(p))$ from a neighborhood $\mathcal{V}$ of $\left\{ 0\right\}\times\Sigma$ in $\mathbb{R}\times\Sigma$ onto a neighborhood of $\iota(\Sigma)$ such that:
\begin{equation}
\label{eq:straightening}
\Psi^{\ast}g = h_{\fre} - e^{-\frf}\dd s\odot e_u\, , \qquad \Psi^{\ast}u = e^{-\frf}e_u\, , \qquad \Psi^{\ast} ( u\wedge(l + \mathrm{i}\, n) ) = e^{-\mathrm{i}\, r s}\, e^{-\frf}e_u\wedge(e_l + \mathrm{i}\, e_n)\, ,
\end{equation}

\noindent
so that $\Psi^{\ast}[u,v,l,n]$ is the isotropic parallelism of the coframe \eqref{eq:timeindependentcoframe} below, and $\Psi^{\ast}\alpha$ is the torsion that Remark \ref{remark:torsiondetermined} associates with it. The argument only uses that the reduction of the development solves \eqref{eq:evolutionconstraintII} and that the data satisfy \eqref{eq:constraintIIthm2}. We do not use either global hyperbolicity or any restriction on the lapse, and the identities \eqref{eq:straightening} hold on the whole domain on which the flow of $u^{\sharp_g}$ starting at $\iota(\Sigma)$ is defined.
\end{remark}

\begin{remark}
\label{remark:evolutionisgauge}
Lemma \ref{lemma:killingreformulation} shows that the evolution equations in \eqref{eq:evolutionconstraintII} carry no dynamical content: in the coordinates $(s,p)$ defined by the flow of $U = u^{\sharp_g}$, every solution becomes the time-independent metric and null one-form \eqref{eq:straightening}, whatever the lapse, and the only residual evolution is the rotation of $u\wedge\zeta$ at the rate $r_t$. In the gauge $\lambda_t = e^{-\frf_t}$ one has $U = -\partial_t + e^{-\frf_t}e^{\sharp}_{tu}$, which is why composing with the flow of $\cX = e^{-\frf}e^{\sharp}_u$ integrates the system explicitly in the proof of Theorem \ref{thm:wellposedFinalSkew}. The entire content of the Cauchy problem therefore resides in the constraint equations and in the obstruction $\dd_{\Sigma}r_t = 0$, in agreement with \eqref{eq:screenrotationinvariant}, which identifies $r_t$ with $\alpha(u^{\sharp_g})$.
\end{remark}

\noindent
We continue with a simple example that shows that the \emph{original} three-dimensional manifold $\Sigma$ may not be a Cauchy hypersurface of the maximal flow tube, even if it is fully contained in it.

\begin{ep} 
Consider $\Sigma=\mathbb{R}^3$ with coordinates $(x,y,z)$ and initial data:
\begin{equation*}
\fre = (\dd x,\dd y,\dd z), \qquad \frf=-\log(1+x^2), \qquad \alpha^o=0, \qquad \alpha^\perp=0.
\end{equation*}

\noindent
These data satisfy the constraint equations. Indeed, $d_\Sigma(e^{-\frf}\dd x)=d_\Sigma((1+x^2)\dd x)=0$, while all the
remaining constraint terms vanish identically. The vector field $\cX \in \mathfrak{X}(\mathbb{R}^3)$ that we use to define the flow $\Phi_t$ is give by:
\begin{equation*}
\cX=e^{-\frf}e_u^\sharp=(1+x^2)\partial_x.
\end{equation*}
Its maximal flow is:
\begin{equation*}
\Phi_t(x,y,z) = (\tan(t+\arctan x),y,z),
\end{equation*}
and its maximal flow domain is
\begin{equation*}
\cU = \left\{ (t,x,y,z)\in\mathbb{R}^4 \ \middle|\ -\frac{\pi}{2}<t+\arctan x<\frac{\pi}{2} \right\}.
\end{equation*}

\noindent
In particular, although $\cU$ contains the entire initial hypersurface $\{0\}\times\Sigma$, it contains no uniform product neighborhood $(-\varepsilon,\varepsilon)\times\Sigma$. Set $\xi=\arctan x$. The explicit evolution gives:
\begin{eqnarray*}
e^{-\frf_t}=\sec^2(t+\xi), \qquad e_{tu} = \frac{\sec^2(t+\xi)}{1+x^2}\,\dd x, \qquad e_{tl}=\dd y, \qquad e_{tn}=\dd z.
\end{eqnarray*}
Consequently, the Lorentzian metric on $\cU$ is:
\begin{equation*}
g = \sec^4(t+\xi)(-\dd t^2+\dd\xi^2 ) + \dd y^2 + \dd z^2.
\end{equation*}
which is a flat domain in Minkowski space. Introduce null coordinates $u=t+\xi$ and $v=t-\xi$. In these coordinates,
the flow domain is:
\begin{equation*}
-\frac{\pi}{2}<u<\frac{\pi}{2}, \qquad -\pi<u-v<\pi .
\end{equation*}
The null curve defined by $v=\pi$, $y=z=0$ and $0<u<\frac{\pi}{2}$ is inextendible in $\cU$ and never intersects $\{0\}\times\Sigma$, since $t=\frac{1}{2}(u+\pi)>0$. Hence $\{0\}\times\Sigma$ is not a Cauchy hypersurface of the entire maximal flow tube. A globally hyperbolic development is obtained by restricting the constructed solution to the domain of dependence $D(\Sigma)$, which can be found to be given by:
\begin{equation*}
D(\Sigma) = \left\{ -\frac{\pi}{2}<u<\frac{\pi}{2},\ -\frac{\pi}{2}<v<\frac{\pi}{2} \right\}\times\mathbb R^2.   
\end{equation*}

\noindent
in the coordinates $(u,v,y,z)$.
\end{ep} 

\noindent
The previous example shows that the maximal flow domain $\cU$ of $\cX = e^{-\frf} e^{\sharp}_u$, equipped with the metric reconstructed from the flow, may fail to be globally hyperbolic. The following corollary realizes the development constructed in the proof of Theorem \ref{thm:wellposedFinalSkew} as an open restriction of a skew-torsion parallel spinor defined on the whole of $\mathbb{R}\times \Sigma$ whose metric and torsion do not depend on $t$, exhibiting in particular the null Killing vector field $u^{\sharp_g}$ of Lemma \ref{lemma:uniqueevol} as $-\partial_t$, and shows that the restriction to the domain of dependence of $\Sigma$ in Theorem \ref{thm:wellposedFinalSkew} is unnecessary when $\Sigma$ is compact.

\begin{cor}
\label{cor:timeindependentdevelopment}
Let $(\fre,\frf,\alpha^o,\alpha^{\perp})$ be a solution of the constraint equations \eqref{eq:constraintIIthm1} and \eqref{eq:constraintIIthm2} on a connected and oriented three-manifold $\Sigma$. Let $\Phi\colon \cU \to \Sigma$ be the maximal flow of $\cX = e^{-\frf} e^{\sharp}_u \in \mathfrak{X}(\Sigma)$ and let $r = \alpha^{\perp}(\cX) - e^{-\frf}\alpha^o \in \mathbb{R}$. Then:
\begin{enumerate}
\item The coframe $(u,v,l,n)$ on $\mathbb{R}\times \Sigma$ given by:
\begin{equation}
\label{eq:timeindependentcoframe}
u = e^{-\frf} e_u \, , \qquad v = \frac{e^{\frf}}{2}\, e_u - \dd t \, , \qquad l + \mathrm{i}\, n = e^{-\mathrm{i}\, r t}\, (e_l + \mathrm{i}\, e_n)
\end{equation}

\noindent
defines a skew-torsion isotropic parallelism $[u,v,l,n]$ on $\mathbb{R}\times \Sigma$ with torsion $H = \ast_g \alpha$, where:
\begin{equation}
\label{eq:timeindependentmetric}
g = u\odot v + l\otimes l + n\otimes n = h_{\fre} - e^{-\frf}\, \dd t \odot e_u \, , \qquad \alpha = \alpha^{\perp} - r\, \dd t
\end{equation}

\noindent
In particular, $g$ and $\alpha$ do not depend on $t$, the vector field $u^{\sharp_g} = -\partial_t$ is a null Killing vector field of $g$, and $\mathbb{R}\times\Sigma$ admits a skew-torsion parallel spinor with torsion $\ast_g\alpha$ regardless of whether $\cX$ is complete or not.

\item The map $\Phi\colon \cU \to \mathbb{R}\times\Sigma$, $(t,p)\mapsto (t,\Phi_t(p))$, is a diffeomorphism onto an open neighborhood of $\Sigma$ in $\mathbb{R}\times\Sigma$ which restricts to the identity on $\Sigma$ and pulls back $(g,\alpha,[u,v,l,n])$ to the solution on $\cU$ constructed in the proof of Theorem \ref{thm:wellposedFinalSkew} in the gauge $\lambda_t = e^{-\frf_t}$.  

\item If the Riemannian metric $e^{2\frf} h_{\fre}$ on $\Sigma$ is complete, in particular if $\Sigma$ is compact, then the map of the previous item is a diffeomorphism of $\mathbb{R}\times\Sigma$ onto itself, and the Lorentzian manifold $(\mathbb{R}\times \Sigma , g)$ is a maximal globally hyperbolic development of the given data.  
\end{enumerate}
\end{cor}

\begin{proof} 
We prove items $(1)$ and $(2)$ together. The forms in \eqref{eq:timeindependentcoframe} constitute a global coframe on $\mathbb R\times\Sigma$, since:
\begin{equation*}
u\wedge v\wedge l\wedge n
=
e^{-\frf}\,\dd t\wedge e_u\wedge e_l\wedge e_n.
\end{equation*}

\noindent
is nowhere vanishing. Their associated Lorentzian metric is:
\begin{equation*}
g = u\odot v+l\otimes l+n\otimes n = h_{\fre} - e^{-\frf}\,\dd t\odot e_u\, .
\end{equation*}

\noindent
In particular, $g$ and $\alpha=\alpha^{\perp}-r\,\dd t$ do not depend on $t$, and $u^{\sharp_g}=-\partial_t$ is a null Killing vector field. To verify the skew-torsion parallelism condition, we first establish the pull-back identities in item $(2)$. The maximal flow property implies that the map:
\begin{equation*}
\Phi\colon\cU\longrightarrow\mathbb R\times\Sigma,
\qquad
(t,p)\longmapsto(t,\Phi_t(p))
\end{equation*}
has image:
\begin{equation*}
\Phi(\cU)
=
\left\{(t,p)\in\mathbb R\times\Sigma
\ \middle|\ (-t,p)\in\cU\right\},
\end{equation*}
and smooth inverse:
\begin{equation*}
\Phi^{-1}(t,p)=(t,\Phi_{-t}(p)).
\end{equation*}

\noindent
Thus $\Phi$ is a diffeomorphism onto an open neighborhood of $\{0\}\times\Sigma$, restricts to the identity on this hypersurface, and preserves orientation. By the construction in the proof of Theorem \ref{thm:wellposedFinalSkew}, the solution on $\cU$ is given by:
\begin{eqnarray*}
& \frf_t=\Phi_t^{\ast}\frf,
\qquad
e_{tu}=\Phi_t^{\ast}e_u,
\qquad
e_{tl}+\mathrm{i}e_{tn}
=
e^{-\mathrm{i}rt}\Phi_t^{\ast}(e_l+\mathrm{i}e_n),\\
& \alpha_t^o=\Phi_t^{\ast}\alpha^o,
\qquad
\alpha_t^{\perp}=\Phi_t^{\ast}\alpha^{\perp},
\qquad
\lambda_t=e^{-\frf_t}.
\end{eqnarray*}

\noindent
For every one-form $\omega$ on $\Sigma$, regarded as a time-independent form on $\mathbb R\times\Sigma$, we have:
\begin{equation*}
\Phi^{\ast}\omega
=
\Phi_t^{\ast}\omega
+
\Phi_t^{\ast} (\omega(\cX) )\,\dd t.
\end{equation*}

\noindent
Since $e_u(\cX)=e^{-\frf}$ and $e_l(\cX)=e_n(\cX)=0$, this gives:
\begin{eqnarray*}
\Phi^{\ast}u = e^{-\frf_t} (e^{-\frf_t}\dd t+e_{tu})\, , \qquad \Phi^{\ast}v = \frac{e^{\frf_t}}{2} (-e^{-\frf_t}\dd t+e_{tu}), \qquad \Phi^{\ast}(l+\mathrm{i}n) =  e^{-\mathrm{i}rt}\Phi_t^{\ast}(e_l+\mathrm{i}e_n) = e_{tl}+\mathrm{i}e_{tn}.
\end{eqnarray*}

\noindent
These are precisely the adapted coframe forms of the reconstructed solution. Consequently:
\begin{equation*}
\Phi^{\ast}g
=
-e^{-2\frf_t}\,\dd t\otimes\dd t+h_{\fre_t},
\end{equation*}
and $\Phi$ pulls back the isotropic parallelism of \eqref{eq:timeindependentcoframe} to that of the reconstructed solution. Moreover, using the constancy of $r$ and its defining identity, we obtain:
\begin{eqnarray*}
\Phi^{\ast}\alpha =
\Phi_t^{\ast}\alpha^{\perp} + (\Phi_t^{\ast} (\alpha^{\perp}(\cX))-r )\dd t = 
\Phi_t^{\ast}\alpha^{\perp} + \Phi_t^{\ast}(e^{-\frf}\alpha^o)\,\dd t =  \alpha_t^{\perp} + e^{-\frf_t}\alpha_t^o\,\dd t.
\end{eqnarray*}

\noindent
This proves the pull-back assertion in item $(2)$. In particular, the solution on $\cU$ and the naturality of the exterior derivative and Hodge operator imply that the coframe \eqref{eq:timeindependentcoframe} satisfies equations \eqref{eq:duvtorsiongeneralII} and \eqref{eq:dlntorsiongeneralII}, with torsion $\ast_g\alpha$, on $\Phi(\cU)$. We now extend this solution to the whole product. For every $s\in\mathbb R$, translation $(t,p)\mapsto(t+s,p)$ preserves $g$, $\alpha$, $u$, $v$ and the orientation, and sends:
\begin{equation*}
l+\mathrm{i}n
\longmapsto
e^{-\mathrm{i}rs}(l+\mathrm{i}n).
\end{equation*}

\noindent
Since $r$ is constant, the exterior system is invariant under the corresponding simultaneous constant rotation of $l+\mathrm{i}n$ and $\kappa+\mathrm{i}\rho$. Its validity therefore extends to every translate of $\Phi(\cU)$. These translates cover $\mathbb R\times\Sigma$, because $\Phi(\cU)$ contains $\{0\}\times\Sigma$. The local auxiliary one-forms $\kappa$ and $\rho$ agree on overlaps by Remark \ref{remark:torsiondetermined}, and hence define global smooth one-forms. Theorem \ref{thm:existencenullcoframeII} now gives the claimed skew-torsion isotropic parallelism and associated parallel spinor on $\mathbb R\times\Sigma$. This proves item $(1)$, including the case in which $\cX$ is not complete.

\noindent
We now consider item $(3)$. We assume that the Riemannian metric:
\begin{equation*}
q:=e^{2\frf}h_{\fre}
\end{equation*}
is complete. Since $|\cX|_q=1$, the vector field $\cX$ is complete. Thus $\cU=\mathbb R\times\Sigma$, and $\Phi$ is a diffeomorphism of $\mathbb R\times\Sigma$ onto itself, identifying the metric constructed in Theorem \ref{thm:wellposedFinalSkew} with $g$. We now prove maximality as a globally hyperbolic development of the prescribed data. For simplicity in the exposition, below we set $M=\mathbb R\times\Sigma$ and:
\begin{equation*}
K:=u^{\sharp_g}=-\partial_t.
\end{equation*}

\noindent
This distinguished null current is complete. Suppose that:
\begin{equation*}
j:(M,g,\alpha,[u,v,l,n]) \longrightarrow (\widetilde M,\widetilde g,\widetilde\alpha, [\widetilde u,\widetilde v,\widetilde l,\widetilde n])
\end{equation*}

\noindent
is an embedding into another globally hyperbolic development of the same initial data, preserving the metric, torsion, isotropic parallelism and initial embedding. In particular, $j(\Sigma)$ is Cauchy in $\widetilde M$, and:
\begin{equation*}
j_*K=\widetilde K|_{j(M)}, \qquad \widetilde K:=\widetilde u^{\sharp_{\widetilde g}}.
\end{equation*}

\noindent
Every maximal integral curve of the nowhere vanishing vector field $\widetilde{K}$ is an inextendible causal curve. It cannot have an endpoint at a finite value of the flow parameter, by the escape lemma, nor at an infinite value, since a limit point of an integral curve as the
parameter tends to $\pm\infty$ is a zero of the vector field. Let $q\in\widetilde M$. The maximal $\widetilde K$-integral curve through $q$ therefore meets $j(\Sigma)$, say at $j(p)$ with $p\in\Sigma$. On the other hand:
\begin{equation*}
s\longmapsto j (\operatorname{Fl}^{K}_s(p) ),
\qquad s\in\mathbb R,
\end{equation*}
is a globally defined $\widetilde K$-integral curve through $j(p)$. Uniqueness for integral curves implies that this is the entire maximal orbit through $j(p)$. Hence $q\in j(M)$. Since $q$ was arbitrary, $j$ is surjective. Thus the constructed development admits no proper globally hyperbolic extension preserving the prescribed solution and initial hypersurface. The same construction also yields the universal property of the maximal development. Let $N$ be any globally hyperbolic development of the initial data, with initial embedding $\iota_N$ and distinguished null vector field $K_N = u_N^{\sharp_{g_N}}$. Its maximal flow from the initial hypersurface defines:
\begin{equation*}
\Psi_N:\mathcal V_N\subset\mathbb R\times\Sigma\longrightarrow N,
\qquad
\Psi_N(t,p)=\operatorname{Fl}^{K_N}_{-t}(\iota_N(p)).
\end{equation*}

\noindent
Every maximal orbit of $K_N$ is an inextendible causal curve and therefore meets $\iota_N(\Sigma)$ exactly once. The absence of closed causal curves then implies injectivity of $\Psi_N$, while the Cauchy property implies surjectivity. Transversality of $K_N$ to the spacelike initial hypersurface shows that $\Psi_N$ is a local diffeomorphism, and hence a diffeomorphism onto $N$. By Remark \ref{remark:straightening}, applied on the whole of $\mathcal V_N$, $\Psi_N$ pulls back $(g_N,\alpha_N,[u_N,v_N,l_N,n_N])$ to the restriction of $(g,\alpha,[u,v,l,n])$ to $\mathcal{V}_N\subset M$. Consequently, $\Psi_N^{-1}$, followed by the inclusion $\mathcal V_N\hookrightarrow M$, embeds developments that fix the initial hypersurface. Thus every globally hyperbolic development of the prescribed data embeds into $M$, which is the unique maximal development up to the prescribed equivalence.
\end{proof}

\begin{remark} 
The time-independence of the metric and torsion does not imply time-independence of the isotropic parallelism under the evolution flow. Indeed, the coframe in \eqref{eq:timeindependentcoframe} satisfies:
\begin{equation*}
l+\mathrm{i}n=e^{-\mathrm{i}rt}(e_l+\mathrm{i}e_n),
\end{equation*}
and therefore:
\begin{equation*}
\mathcal{L}_{\partial_t}l=rn, \qquad \mathcal{L}_{\partial_t}n=-rl.
\end{equation*} 

\noindent
This rotation preserves $l\otimes l+n\otimes n$, and hence the metric, whereas the spinor square satisfies:
\begin{equation*}
\mathcal{L}_{\partial_t} \alpha_{\varepsilon} = \mathcal{L}_{\partial_t}(u+u\wedge l)=r\,u\wedge n.
\end{equation*}

\noindent
Consequently, for $r\neq0$, the isotropic parallelism, equivalently the associated spinor up to sign, is not invariant under the evolution flow, although both $g$ and $\alpha$ are invariant. This rotation cannot be absorbed by changing the representative of the same isotropic parallelism, since such a change leaves $u\wedge l$ unchanged. 
\end{remark}

\noindent
The metric \eqref{eq:timeindependentmetric} admits the nowhere vanishing null Killing vector field $u^{\sharp_g}=-\partial_t$. It therefore belongs to the class admitting a $G_1$ on $N_1$, namely a one-parameter group of isometries with null one-dimensional orbits, in the terminology of \cite[Section 24.4]{StephaniExactSolutions}. Since the Dirac current is null and Killing, the constraints give:
\begin{equation*}
\nabla^g_{u} u = 0\, , \qquad \nabla^g u=\frac{1}{2}\dd u\, , \qquad u\wedge\dd u=r\,u\wedge l\wedge n.
\end{equation*}

\noindent
Thus its null congruence is affinely parametrized, non-expanding and shear-free, with constant optical twist coefficient $r/2$. It is Kundt precisely when $r=0$ \cite{ColeyKundt}, and the cases $r=0$ and $r\neq0$ correspond respectively to the non-twisting and twisting null Killing classes of \cite[Sections 24.4.1 and 24.4.2]{StephaniExactSolutions}. The constancy of $r$ also implies that the metric is \emph{algebraically special}. Indeed, we have $\Psi_0=\Psi_1=0$, where $\Psi_0,\Psi_1$ are the Weyl scalars in a null tetrad adapted to $u^{\sharp_g}$. Hence, wherever the Weyl tensor is nonzero, $u^{\sharp_g}$ is a repeated principal null direction. We conclude that the Petrov type is II, D, III, N, or O at those points where the Weyl tensor vanishes. In particular, the metric is algebraically special at every point. Moreover:
\begin{equation*}
\Ric^g(u^{\sharp_g},u^{\sharp_g})=\frac{r^2}{2}.
\end{equation*}

\noindent
Consequently, the case $r\neq0$ contains no Einstein metrics, even with a cosmological constant. When $r=0$, the metric is locally conformal to a metric for which the same null vector field is parallel. Such a conformally related metric is a \emph{Brinkmann metric} \cite{MehidiZeghib}. The original metric has Levi-Civita parallel Dirac current precisely when:
\begin{equation*}
\nabla^g u=0
\quad\Longleftrightarrow\quad
\dd_{\Sigma}(e^{-\frf}e_u)=0
\quad\Longleftrightarrow\quad
\alpha\wedge u=0,
\end{equation*}

\noindent
by \eqref{eq:duvtorsiongeneralII}. This is stronger than $r=0$ and does not require vanishing skew-torsion. Note that in the previous example the Riemannian metric:
\begin{equation*}
e^{2\frf} h_{\fre} = (1+x^2)^{-2}(\dd x^2 + \dd y^2 + \dd z^2)
\end{equation*}

\noindent
is not complete, in agreement with the third item of the previous corollary. As a corollary to Theorem \ref{thm:wellposedFinalSkew}, we obtain the following result, which provides an initial data characterization of skew-torsion parallel spinors in four dimensions.

\begin{cor}
\label{cor:initialdataskew}
A connected and oriented three-manifold $\Sigma$ occurs as the Cauchy hypersurface of a Lorentzian four-manifold equipped with a skew-torsion parallel spinor if and only if it admits a global coframe $(e_u,e_l,e_n)$ and a function $\frf\in C^{\infty}(\Sigma)$ satisfying the following differential system:
\begin{equation}
\label{eq:constraintsfinal}
\dd_\Sigma e_l\wedge e_u=0,  \qquad \dd_\Sigma e_n\wedge e_u=0, \qquad \dd_{\Sigma}e_u \wedge e_u  = r e^{\frf} e_u\wedge e_l\wedge e_n \, ,
\end{equation}
for a real constant $r\in \mathbb{R}$. If that is the case, the corresponding initial skew-torsion data is uniquely determined by:
\begin{equation}
\label{eq:alphafinal}
\alpha^o= \ast_{\fre} (\dd_{\Sigma} e_u \wedge e_u)- (\cL_{e^{\sharp}_u}e_n)(e^{\sharp}_l) + (\cL_{e^{\sharp}_u}e_l)(e^{\sharp}_n) , \qquad \alpha^\perp =\ast_\fre (\dd_\Sigma e_u-\dd_\Sigma\frf\wedge e_u ) + \alpha^o e_u
\end{equation} 
\end{cor} 

\begin{proof}
Suppose first that $(\fre,\frf,\alpha^o,\alpha^{\perp})$ satisfies equations \eqref{eq:constraintIIthm1} and \eqref{eq:constraintIIthm2}. The first and second equations in \eqref{eq:constraintsfinal} correspond to the second and third equations in \eqref{eq:constraintIIthm1}. Hence, we only need to prove that the third equation in \eqref{eq:constraintsfinal} holds. Applying the Hodge operator $\ast_{\fre}$ to Equation \eqref{eq:constraintIIthm1} and expanding, we obtain:
\begin{equation*}
\alpha^{\perp}-\alpha^o e_u = \ast_{\fre} (\dd_{\Sigma}e_u-\dd_{\Sigma}\frf\wedge e_u)
\end{equation*}

\noindent
which immediately implies: 
\begin{equation*}
\dd_{\Sigma}e_u\wedge e_u = (\alpha^{\perp}(e^{\sharp}_{u})-\alpha^o ) e_u\wedge e_l\wedge e_n = (\alpha^{\perp}(e^{\sharp}_{u})-\alpha^o ) e^{-\frf} e^{\frf} e_u\wedge e_l\wedge e_n.
\end{equation*}

\noindent
Since $\Sigma$ is connected, by the first equation in \eqref{eq:constraintIIthm2}, there is a constant $r\in\mathbb R$ such that $r = (\alpha^{\perp}(e^{\sharp}_{u})-\alpha^o ) e^{-\frf}$, and thus we obtain the third equation in \eqref{eq:constraintsfinal}. The equations in \eqref{eq:alphafinal} are obtained by isolating $\alpha^o$ and $\alpha^{\perp}$ in equations \eqref{eq:constraintIIthm1} and \eqref{eq:constraintIIthm2}. 

\noindent
Conversely, suppose that $(\fre,\frf)$ satisfies the differential system \eqref{eq:constraintsfinal} for a real constant $r$, and define $\alpha^o$ and $\alpha^{\perp}$ by equation \eqref{eq:alphafinal}. Then the second equation in \eqref{eq:alphafinal} is equivalent to the first equation in \eqref{eq:constraintIIthm1}. On the other hand, the second and third equations in \eqref{eq:constraintIIthm1} are equivalent to the first and second equations in \eqref{eq:constraintsfinal}. Combining both equations in \eqref{eq:alphafinal} to eliminate the term $\ast_{\fre}(\dd_{\Sigma} e_u \wedge e_u)$ gives the second equation in \eqref{eq:constraintIIthm2}. Finally, to recover the first equation in \eqref{eq:constraintIIthm2}, we evaluate the second equation in \eqref{eq:constraintIIthm2} in $e^{\sharp}_u$, obtaining:
\begin{equation*}
    \alpha^{\perp}(e^{\sharp}_u) - \alpha^o = \ast_{\fre} (\dd_{\Sigma} e_u \wedge e_u)
\end{equation*}

\noindent
Combining this with the third equation in \eqref{eq:constraintsfinal} gives the first equation in \eqref{eq:constraintIIthm2} and hence we conclude.
\end{proof}

\noindent
By the previous corollary, we will consider the equations in \eqref{eq:constraintsfinal} as the constraint equations of the Cauchy problem posed by a skew-torsion parallel spinor. The solutions to these constraint equations split into two cases, namely, the case $r\neq0$, to which we will refer as the \emph{contact branch}, and the case $r=0$, to which we will refer as the \emph{integrable branch}.

\begin{prop}
\label{prop:geometryconstraintbranches}
Let $\fre = (e_u,e_l,e_n)$ be a positively oriented global coframe on $\Sigma$ and let $\frf\in C^{\infty}(\Sigma)$. Then, $(\fre,\frf)$ satisfies the constraint equations \eqref{eq:constraintsfinal} for a constant $r\in \mathbb{R}$ if and only if:
\begin{equation}
\label{eq:bracketconstraint}
[e^{\sharp}_{l},e^{\sharp}_{n}]=-r e^{\frf}e^{\sharp}_{u}
\end{equation}

\noindent
In this case, $e_u\wedge \dd_{\Sigma} e_u = r e^{\frf} e_u\wedge e_l\wedge e_n$, and the two branches can be characterized as follows.
\begin{enumerate}
\item If $r\neq0$, then $e_u$ is a contact form on $\Sigma$, which is positive, that is, $e_u\wedge\dd_{\Sigma} e_u$ is a positive volume form on $\Sigma$, if and only if $r>0$, and $\Ker(e_u) \subset T\Sigma$ is a cooriented contact distribution trivialized by $(e^{\sharp}_{l},e^{\sharp}_{n})$. Conversely, let $e_u$ be a contact form on $\Sigma$ such that $\Ker(e_u)$ is trivial as a plane bundle, and let $r\in \mathbb{R}$ be a nonzero constant which is positive if and only if $e_u$ is positive. Then, every global frame of $\Ker(e_u)$ is, up to the order of its elements, of the form $(e^{\sharp}_{l},e^{\sharp}_{n})$ for a unique positively oriented global coframe $\fre = (e_u,e_l,e_n)$ and a unique function $\frf\in C^{\infty}(\Sigma)$ such that $(\fre,\frf)$ satisfies \eqref{eq:constraintsfinal} with constant $r$.

\item If $r=0$, then $e^{\sharp}_{l}$ and $e^{\sharp}_{n}$ commute, and $\Ker(e_u) \subset T\Sigma$ is a cooriented integrable distribution whose leaves, equipped with the metric induced by $h_{\fre}$, are oriented flat surfaces, locally isometric to the Euclidean plane through the restrictions of $e_l$ and $e_n$. If $e^{\sharp}_{l}$ and $e^{\sharp}_{n}$ are complete, in particular if $\Sigma$ is compact, their flows define a locally free action of $\mathbb{R}^2$ on $\Sigma$ whose orbits are the leaves of $\Ker(e_u)$, each of which is isometric to the quotient of the Euclidean plane by a discrete group of translations. Conversely, if $(e^{\sharp}_{u},e^{\sharp}_{l},e^{\sharp}_{n})$ is a positively oriented global frame on $\Sigma$ such that $e^{\sharp}_{l}$ and $e^{\sharp}_{n}$ commute, then its dual coframe $\fre$ satisfies \eqref{eq:constraintsfinal} with $r=0$ together with every $\frf \in C^{\infty}(\Sigma)$.
\end{enumerate}

\noindent
In either branch, the compatible skew-torsion data are uniquely reconstructed through Equation \eqref{eq:alphafinal}.
\end{prop}

\begin{proof}
Since $(e^{\sharp}_{u},e^{\sharp}_{l},e^{\sharp}_{n})$ is the frame dual to $\fre$, for every $a = u,l,n$ we have:
\begin{equation*}
\dd_{\Sigma} e_a \wedge e_u = \dd_{\Sigma} e_a (e^{\sharp}_{l},e^{\sharp}_{n})\, e_u\wedge e_l\wedge e_n = - e_a ([e^{\sharp}_{l},e^{\sharp}_{n}])\, e_u\wedge e_l\wedge e_n
\end{equation*}

\noindent
where we have used that $e_a(e^{\sharp}_{l})$ and $e_a(e^{\sharp}_{n})$ are constant. Hence, the first two equations in \eqref{eq:constraintsfinal} are equivalent to $e_l([e^{\sharp}_{l},e^{\sharp}_{n}]) = e_n([e^{\sharp}_{l},e^{\sharp}_{n}]) = 0$, whereas the third equation is equivalent to:
\begin{equation*}
e_u([e^{\sharp}_{l},e^{\sharp}_{n}]) = - r e^{\frf}    
\end{equation*}

\noindent
Together, these conditions are equivalent to \eqref{eq:bracketconstraint}, and
\begin{equation*}
 e_u\wedge \dd_{\Sigma} e_u = r e^{\frf} e_u\wedge e_l\wedge e_n   
\end{equation*}

\noindent
is the third equation in \eqref{eq:constraintsfinal}.  \medskip

\noindent
Suppose that $r\neq 0$. Then $e_u\wedge \dd_{\Sigma} e_u \in \Omega^3(\Sigma)$ is a volume form, which is positive if and only if $r>0$, so $e_u$ is a contact form whose kernel is spanned by $e^{\sharp}_{l}$ and $e^{\sharp}_{n}$ and is cooriented by $e_u$. \medskip

Conversely, let $e_u$, $r$ and a global frame $(X, Y)$ of $\mathrm{Ker}(e_u)$ be as in the statement. Since $e_u(X) = e_u(Y) = 0$, we have:
\begin{equation*}
e_u([X, Y]) = - (\dd _\Sigma e_u) (X, Y)    
\end{equation*}

\noindent
which vanishes nowhere because $\dd _\Sigma e_u\in \Omega^2(\Sigma)$ is nondegenerate on $\mathrm{Ker}(e_u)\subset T\Sigma$ and $e_u$ is a contact form. Hence $[X, Y]\in \mathfrak{X}(\Sigma)$ is everywhere transverse to $\mathrm{Ker}(e_u)$. Dividing the bracket $[X, Y]$ by the nowhere vanishing function $e_u([X, Y])$, we obtain a vector field transverse to $\Ker(e_u)$ on which $e_u$ takes the value one, and which does not change if the order of the two elements of the frame is reversed. Hence, for exactly one ordering of the given frame, this vector field together with the frame form a positively oriented global frame on $\Sigma$, whose dual coframe is of the form $\fre = (e_u,e_l,e_n)$ with $(e^{\sharp}_{l},e^{\sharp}_{n})$ the given frame and $e^{\sharp}_{u}$ the vector field just constructed. By construction: 
\begin{equation*}
e_l([e^{\sharp}_{l},e^{\sharp}_{n}]) = e_n([e^{\sharp}_{l},e^{\sharp}_{n}]) = 0    
\end{equation*}

\noindent
that is, $\fre$ satisfies the first two equations in  \eqref{eq:constraintsfinal}. Moreover, $e_u\wedge \dd_{\Sigma} e_u$ is a nowhere vanishing function with the sign of $r$ times $e_u\wedge e_l\wedge e_n$, whence the third equation in \eqref{eq:constraintsfinal} determines $\frf$ uniquely. Finally, $\fre$ is unique, since by \eqref{eq:bracketconstraint} the vector field $e^{\sharp}_{u}$ of any coframe as in the statement is proportional to $[e^{\sharp}_{l},e^{\sharp}_{n}]$ and satisfies $e_u(e^{\sharp}_{u}) = 1$, which forces it and thus must be equal to the vector field constructed above. \medskip

\noindent
Suppose now that $r = 0$. Then $e^{\sharp}_{l}$ and $e^{\sharp}_{n}$ commute by \eqref{eq:bracketconstraint}, so $\Ker(e_u)$ is integrable by the Frobenius theorem and is cooriented by $e_u$. Since $e_u$ vanishes on the leaves of $\Ker(e_u)$, the metric induced by $h_{\fre}$ on a leaf is the restriction of $e_l\otimes e_l + e_n\otimes e_n$, where the restrictions of $e_l$ and $e_n$ to the leaf are closed because $\dd_{\Sigma} e_l (e^{\sharp}_{l},e^{\sharp}_{n}) = \dd_{\Sigma} e_n (e^{\sharp}_{l},e^{\sharp}_{n}) = 0$. Hence, local primitives of these restrictions define local isometries from the leaf to the Euclidean plane, and the leaf is oriented by $e_l\wedge e_n$. If $e^{\sharp}_{l}$ and $e^{\sharp}_{n}$ are complete, their flows commute and define an action of $\mathbb{R}^2$ on $\Sigma$, which is locally free since $e^{\sharp}_{l}$ and $e^{\sharp}_{n}$ are everywhere linearly independent. The orbits of this action are the leaves of $\Ker(e_u)$, and the orbit map of a leaf pulls back $e_l$ and $e_n$ to the standard coframe of $\mathbb{R}^2$, so it is a local isometry from the Euclidean plane onto the leaf which induces an isometry onto the leaf from the quotient of $\mathbb{R}^2$ by the stabilizer of any of its points, a discrete group of translations. The converse follows from \eqref{eq:bracketconstraint}, since the dual coframe of a positively oriented global frame $(e^{\sharp}_{u},e^{\sharp}_{l},e^{\sharp}_{n})$ with $e^{\sharp}_{l}$ and $e^{\sharp}_{n}$ commuting satisfies \eqref{eq:bracketconstraint} with $r=0$ for every $\frf\in C^{\infty}(\Sigma)$. Finally, the last assertion follows from Equation \eqref{eq:alphafinal}.
\end{proof}

 \begin{remark}
\label{remark:Reebcontactbranch}
In the contact branch, the vector field $e^{\sharp}_u$ is transverse to the contact distribution $\Ker(e_u)$ and satisfies $e_u(e^{\sharp}_u) = 1$, but it is in general not the Reeb vector field of the contact form $e_u$. Indeed, since $e^{\sharp}_u$ is a unit vector field with respect to $h_{\fre}$, we have:
\begin{equation*}
\iota_{e^{\sharp}_u}\dd_{\Sigma} e_u = \cL_{e^{\sharp}_u} e_u = h_{\fre}(\nabla^{h_{\fre}}_{e^{\sharp}_u} e^{\sharp}_u , -)\, ,
\end{equation*}

\noindent
and therefore $e^{\sharp}_u$ is the Reeb vector field of $e_u$ if and only if its integral curves are geodesics of $h_{\fre}$. On the other hand, the first equation in \eqref{eq:constraintIIthm1} gives:
\begin{equation*}
\iota_{e^{\sharp}_u}\dd_{\Sigma} e_u = \ast_{\fre}(\alpha^{\perp}\wedge e_u) - \dd_{\Sigma}\frf + \dd_{\Sigma}\frf(e^{\sharp}_u)\, e_u\, ,
\end{equation*}

\noindent
whence $e^{\sharp}_u$ is the Reeb vector field of $e_u$ if and only if:
\begin{equation*}
\dd_{\Sigma}\frf\wedge e_u = \ast_{\fre}(\alpha^{\perp}\wedge e_u)\wedge e_u\, ,
\end{equation*}

\noindent
that is, if and only if $\dd_{\Sigma}\frf(e^{\sharp}_l) = \alpha^{\perp}(e^{\sharp}_n)$ and $\dd_{\Sigma}\frf(e^{\sharp}_n) = -\alpha^{\perp}(e^{\sharp}_l)$. By Equation \eqref{eq:frfdef} evaluated at $t=0$, this is in turn equivalent to:
\begin{equation*}
\Theta(e^{\sharp}_u , e^{\sharp}_l) = \frac{1}{2}\,\alpha^{\perp}(e^{\sharp}_n)\, , \qquad \Theta(e^{\sharp}_u , e^{\sharp}_n) = -\frac{1}{2}\,\alpha^{\perp}(e^{\sharp}_l)\, ,
\end{equation*}

\noindent
where $\Theta$ denotes the second fundamental form of $\Sigma$ in the development of Theorem \ref{thm:wellposedFinalSkew}. In particular, if $\alpha^{\perp}\wedge e_u = 0$, then $e^{\sharp}_u$ is the Reeb vector field of $e_u$ if and only if it is an eigenvector of the shape operator of $\Sigma$. In general, a direct computation using \eqref{eq:constraintIIthm1} and the third equation in \eqref{eq:constraintsfinal} shows that the Reeb vector field of $e_u$ is given by:
\begin{equation*} 
e^{\sharp}_u + \frac{e^{-\frf}}{r} (  (\alpha^{\perp}(e^{\sharp}_l) + \dd_{\Sigma}\frf(e^{\sharp}_n) ) e^{\sharp}_l +  (\alpha^{\perp}(e^{\sharp}_n) - \dd_{\Sigma}\frf(e^{\sharp}_l) ) e^{\sharp}_n ) \in \mathfrak{X}(\Sigma)
\end{equation*}
This would be the first step towards studying the Reeb dynamics of this class of contact forms.
\end{remark}

\begin{cor}
\label{cor:compactconstraintbranches}
Let $\Sigma$ be a connected, compact and oriented three-manifold without boundary. Then:
\begin{enumerate}
\item $\Sigma$ admits initial data $(\fre,\frf)$ in the contact branch, and the skew-torsion of the corresponding development is necessarily nonzero. More precisely, a cooriented contact distribution on $\Sigma$ is of the form $\Ker(e_u)$ for initial data $(\fre,\frf)$ in the contact branch if and only if it is trivial as a plane bundle, and every such $\Sigma$ admits cooriented contact distributions that are trivial as plane bundles.
\item $\Sigma$ admits initial data in the integrable branch if and only if it is diffeomorphic to an oriented two-torus bundle over the circle.
\end{enumerate}

\noindent
Consequently, every connected, compact and oriented three-manifold without boundary is a Cauchy hypersurface of a globally hyperbolic Lorentzian four-manifold diffeomorphic to $\mathbb{R}\times \Sigma$ and equipped with a skew-torsion parallel spinor, whose skew-torsion is necessarily nonzero unless $\Sigma$ is diffeomorphic to an oriented two-torus bundle over the circle.  
\end{cor}

\begin{remark}
In the integrable case, the connected orbits of the corresponding locally free $\mathbb{R}^2$ action are, when equipped with the induced metric, either flat planes, flat cylinders, or flat tori. This foliation need not coincide with the torus-bundle fibration since nonclosed and dense leaves may occur, as explained in \cite[Proposition 3.15]{Murcia:2020zig}, see also \cite{Murcia:2021dur}.
\end{remark}

\begin{proof}
Every closed oriented three-manifold is parallelizable, so $\Sigma$ admits a global frame, and the plane field spanned by its first two elements is cooriented and trivial as a plane bundle. By the Lutz--Martinet theorem \cite[Theorem 3.1 and Remark 3.9]{GeigesContactGeometry}, every homotopy class of cooriented plane fields on a closed oriented three-manifold contains a positive contact distribution, and homotopic plane fields are isomorphic as plane bundles. Hence, $\Sigma$ admits a cooriented contact distribution which is trivial as a plane bundle. Let now $\Ker(e_u)$ be any such distribution, where $e_u$ is a contact form defining it together with its coorientation, and let $r\in \mathbb{R}$ be a nonzero constant which is positive if and only if $e_u$ is positive. By the first item of Proposition \ref{prop:geometryconstraintbranches}, every global frame of $\Ker(e_u)$ determines initial data $(\fre,\frf)$ in the contact branch with $\fre = (e_u,e_l,e_n)$. The skew-torsion of the corresponding development is nonzero since, by Equation \eqref{eq:skewscreenrotationinvariant} and the definition of $r$, its Dirac current satisfies $\alpha(u^{\sharp_g}) = r \neq 0$ along $\Sigma$. Conversely, if a cooriented contact distribution is of the form $\Ker(e_u)$ for initial data $(\fre,\frf)$ in the contact branch, then it is trivialized by $(e^{\sharp}_{l},e^{\sharp}_{n})$, again by the first item of Proposition \ref{prop:geometryconstraintbranches}. This proves the first item. 

\noindent
Suppose next that $\Sigma$ admits initial data $(\fre,\frf)$ in the integrable branch. Since $\Sigma$ is compact, the commuting vector fields $e^{\sharp}_{l}$ and $e^{\sharp}_{n}$ are complete and, by the second item of Proposition \ref{prop:geometryconstraintbranches}, their flows define a locally free action of $\mathbb{R}^2$ on $\Sigma$. By the theorem of Rosenberg, Roussarie and Weil \cite{RosenbergRoussarieWeil}, a closed oriented three-manifold carrying a locally free action of $\mathbb{R}^2$ is diffeomorphic to an oriented two-torus bundle over the circle. Conversely, every oriented two-torus bundle over the circle is diffeomorphic to the mapping torus of an element $A\in\mathrm{SL}(2,\mathbb{Z})$ acting linearly on $\mathbb{R}^2/\mathbb{Z}^2$, namely to the quotient of $\mathbb{R}^2/\mathbb{Z}^2\times [0,1]$ by the identification $(x,1)\sim(Ax,0)$. Choose a smooth path $B\colon[0,1]\to\mathrm{GL}^+(2,\mathbb{R})$, constant near its endpoints, such that $B(0)=I$ and $B(1)=A^{-1}$, which exists since $\mathrm{GL}^+(2,\mathbb{R})$ is connected. The two vector fields on $\mathbb{R}^2/\mathbb{Z}^2\times [0,1]$ tangent to the torus factor whose coefficient vectors at the point $(x,t)$ are the columns of $B(t)$ commute, since their coefficients depend only on $t$, are everywhere linearly independent, since $B(t)$ is invertible, and descend to the mapping torus, since the identification $(x,1)\sim(Ax,0)$ transforms $B(1) = A^{-1}$ into $AA^{-1} = B(0)$. Together with the vector field tangent to the interval factor, which also descends to the mapping torus, they form a global frame, which is positively oriented after possibly interchanging the two commuting vector fields. By the second item of Proposition \ref{prop:geometryconstraintbranches}, the dual coframe $\fre$ of this frame, together with any $\frf\in C^{\infty}(\Sigma)$, is initial data in the integrable branch. This proves the second item.

\noindent
Finally, by the first item, Equation \eqref{eq:alphafinal} and Corollary \ref{cor:timeindependentdevelopment}, $\Sigma$ is a Cauchy hypersurface of the globally hyperbolic development $\mathbb{R}\times\Sigma$ of initial data in the contact branch, whose skew-torsion is nonzero. If $\Sigma$ is not diffeomorphic to an oriented two-torus bundle over the circle, then by the second item every initial datum on $\Sigma$ belongs to the contact branch, whence the skew-torsion of every such development is nonzero. On the other hand, if the skew-torsion of a development vanishes, then $r = \alpha(u^{\sharp_g}) = 0$ by Equation \eqref{eq:skewscreenrotationinvariant}, so torsion-free initial data belong to the integrable branch and the last assertion follows from the second item.
\end{proof}


\section{The supersymmetric Neveu-Schwarz Cauchy problem}
\label{sec:Globallyhyperbolicsusy}


In this section, we apply the theory of globally hyperbolic skew-torsion isotropic parallelisms developed in the previous section to the globally hyperbolic supersymmetric four-dimensional solutions of Neveu-Schwarz supergravity. This is a celebrated supergravity theory whose supersymmetric conditions on a bosonic configuration involve skew-torsion parallel spinors \emph{coupled} to an abelian bundle gerbe. The key property of the globally hyperbolic reduction of this supergravity theory is that the dilatino equation fixes the normal components of the flux and dilaton curvatures and forces the \emph{screen-rotation} occurring in Theorem~\ref{thm:wellposedFinalSkew} to vanish.  Consequently, in the natural gauge $\lambda_t=e^{-\frf_t}$, the supersymmetric development of the initial data is the pullback development generated by the vector field $e^{-\frf}e^{\sharp}_u$. 


\subsection{The supersymmetric Neveu-Schwarz system}
\label{sec:NSNSsystem}


We begin by introducing the supersymmetric Neveu-Schwarz system in four dimensions, which combines the bosonic sector of four-dimensional Neveu-Schwarz supergravity with its supersymmetry conditions restricted to a bosonic configuration. To do this, we first need to introduce the \emph{geometric object} on which the Neveu-Schwarz system is formulated, which uniquely determines it. In the following, the triple $(\cP, \cA, Y)$ will denote a bundle gerbe \cite{Murray,BehrendXu} defined on an oriented four-manifold $M$. Hence, $Y\to M$ is a smooth submersion, $\cP \to Y\times_M Y$ is a principal $\U(1)$ bundle defined over the fibered product $Y\times_M Y$ of $Y\to M$ with itself, and $\cA$ is connection on $\cP$ satisfying the conditions to be a \emph{connective structure} on $(\cP,Y)$. Given a curving $b\in \Omega^2(Y)$ on $(\cP, \cA, Y)$ we denote its curvature by $H_b\in \Omega^3(M)$. We continue by introducing the following terminology, which relates to the physical jargon used in supergravity and string theory.
 
 \begin{definition}
 \label{def:dilaton}
 A \emph{dilaton} on $M$ is a smooth function $\phi \colon M \to S^1$. The \emph{curvature} $\varphi_{\phi}$ of the dilaton $\phi \colon M \to S^1$ is the pull-back one-form $\varphi_{\phi}  := \phi^{\ast}\mu \in \Omega^1(M)$, where $\mu\in \Omega^1(S^1,\mathbb{R})$ is the Maurer-Cartan form of $S^1$.
 \end{definition}
 
 \noindent
 It follows that the curvature of every dilaton defines an integer class when normalized as follows:
 \begin{eqnarray*}
 	\frac{1}{2\pi} [\varphi_{\phi}] \in H^1_{\mathbb{Z}}(M,\mathbb{R})
 \end{eqnarray*} 
 
 \noindent
 where $H^1_{\mathbb{Z}}(M,\mathbb{R}) \subset H^1(M,\mathbb{R})$ denotes the image of the singular cohomology group $H^1(M,\mathbb{Z})$ in de Rham cohomology of $M$. Since $S^1 = K(\mathbb{Z},1)$ is an Eilenberg-- MacLane space for $\mathbb{Z}$, it follows that we can recover every class in $H^1_{\mathbb{Z}}(M,\mathbb{R})$ as the curvature of a dilaton.  
 
 \begin{remark}
 The dilaton notion introduced above provides a global model for the local dilaton notion used in the physics literature. Let $\phi\colon M\to S^1$ be a dilaton. Since $\varphi_{\phi}\in \Omega^1(M)$ is closed, locally we can write $\varphi_{\phi} = \dd f$ for a local real function $f$, unique modulo a constant, which is what the physics literature refers to as a \emph{dilaton}. We use the dilaton notion introduced in Definition \ref{def:dilaton}, which is locally equivalent to the one used in the physics literature, to construct an appropriate \emph{global} dilaton model.
 \end{remark}
 
 \noindent
 A bundle gerbe with connective structure determines a unique Neveu-Schwarz system on $M$, as prescribed in the following definitions.  
 
 \begin{definition}
 The \emph{configuration space} $\Conf(\cP, \cA, Y)$ of the supersymmetric Neveu-Schwarz system determined by $(\cP, \cA, Y)$ on $M$ is the set of tuples $(g,b,\phi)$ consisting of a Lorentzian metric $g$ on $M$, a curving $b$ on $(\cP, \cA, Y)$ and a dilaton $\phi \colon M \to S^1$.
 \end{definition}

 \noindent
 Given an element $(g,b,\phi)\in \Conf(\cP, \cA, Y)$ we denote by $\nabla^{g,b}$ the unique metric connection on $(M,g)$ with totally skew-symmetric torsion given by $H_b\in \Omega^3(M)$, namely the curvature of $b\in\Omega^2(Y)$. We will denote by $\Ric^{g,b}$ the Ricci curvature of $\nabla^{g,b}$. 
 
 \begin{definition}
 	\label{label:NSNSsystem}
 	The \emph{Neveu-Schwarz system} determined by $(\cP, \cA, Y)$ on $M$ is the following differential system \cite{Ortin,Tomasiello}:
 	\begin{equation}
 		\label{eq:NSNSsystem}
 		\Ric^{g,b} + \nabla^{g,b} \varphi_{\phi}  = 0\, , \qquad  \nabla^{g\ast}\varphi_{\phi} +  \vert \varphi_{\phi}\vert_g^2 = \vert H_b\vert^2_g  
 	\end{equation}
 	
 	\noindent
 	for triples $(g,b,\phi)\in \Conf(\cP, \cA, Y)$.
 \end{definition}
 
 \noindent
 By the previous definition, we can consider the Neveu-Schwarz system as a natural \emph{gauge theoretic} system in Lorentzian geometry with torsion. We will refer to solutions of Equation \eqref{eq:NSNSsystem} as \emph{Neveu-Schwarz solutions}. Given a triple $(g,b,\phi)\in \Conf(\cP, \cA, Y)$, the physics literature refers to $b$ as a \emph{b-field}, a terminology that we will use occasionally. We will denote by $\Sol(\cP, \cA, Y)\subset \Conf(\cP, \cA, Y)$ the solution set of the Neveu-Schwarz system determined by $(\cP, \cA, Y)$ on $M$. The notion of the Neveu-Schwarz system that we have introduced in Definition \ref{label:NSNSsystem} is valid in any dimension. Since in our case $M$ is four-dimensional, we can simplify it accordingly. Given a curving $b\in \Omega^2(Y)$ with curvature $H_b \in \Omega^3(M)$ we set $H_b = \ast_g \alpha_b$ for a uniquely determine one-form. A quick computation gives the following formulae:
 \begin{equation*}
 	H_b\circ_g H_b = \alpha_b\otimes \alpha_b - \vert\alpha_b\vert_g^2 g\, , \quad \vert H_b\vert^2_g  = - \vert \alpha_b\vert^2_g\, , \quad \varphi_{\phi}\lrcorner_g H_b =\ast_g (\alpha_b \wedge \varphi_{\phi})\, , \quad \nabla^{g\ast} H_b = \ast_g \dd\alpha_b 
 \end{equation*}
 
 \noindent
 Hence, the four-dimensional Neveu-Schwarz system is equivalent to the following system of equations:
 \begin{equation}
 	\label{eq:NSNSsystem4d}
 	\Ric^g + \nabla^g \varphi_{\phi} - \frac{1}{2} \alpha_b \otimes \alpha_b + \frac{1}{2} \vert\alpha_b\vert_g^2 g= 0\, , \qquad \dd \alpha_b = \varphi_{\phi}\wedge \alpha_b \, , \qquad    \nabla^{g\ast}\varphi_{\phi} +  \vert \varphi_{\phi}\vert_g^2 + \vert\alpha_b\vert^2_g = 0
 \end{equation} 
 
 \noindent
 for triples $(g,b,\phi)\in \Conf(\cP, \cA, Y)$. Note that, since $H_b = \ast_g\alpha$ is the curvature of a curving on an abelian gerbe, it follows that $\nabla^{g\ast}\alpha = 0$. In the following, when referring to the Neveu-Schwarz system, we will always refer to the differential system \eqref{eq:NSNSsystem4d}.
 
 \begin{definition}
 	A triple $(g,b,\phi)\in \Conf(\cP, \cA, Y)$ is \emph{flux-less} if $b$ is flat, namely if $\alpha_b=0$ identically on $M$, and is \emph{flux} otherwise. A triple $(g,b,\phi)\in \Conf(\cP, \cA, Y)$ is \emph{trivial} if its flux-less and $g$ is flat. 
 \end{definition}
 
\noindent
We have introduced the Neveu-Schwarz system as the equations of motion of the bosonic sector of Neveu-Schwarz supergravity in four dimensions. If that were all, then the Neveu-Schwarz system would simply be a very particular general relativity model, with a very particular \emph{matter content} and energy-momentum tensor. However, and thanks to the underlying supersymmetry of the theory, the Neveu-Schwarz system comes equipped with a first-order \emph{spinorial} differential system that provides \emph{partial integrability} to the full second-order Neveu-Schwarz system and defines the notion of \emph{supersymmetric solution} or \emph{BPS state}. Given an element $(g,b,\phi)\in \Conf(\cP, \cA, Y)$, for ease of notation we denote again by $\nabla^{g,b}$ its lift to any irreducible spinor bundle defined on $(M,g)$.

\begin{definition}
A \emph{supersymmetric Neveu-Schwarz configuration} on $(\cP, \cA, Y)$ is a triple $(g,b,\phi) \in \Conf(\cP, \cA, Y)$ satisfying the following differential system:
\begin{equation}
 \label{eq:KSE4d}
 \nabla^{g,b} \varepsilon = 0\, , \qquad  \varphi_{\phi} \cdot \varepsilon = H_b\cdot \varepsilon
 \end{equation}
 	
\noindent
for a nowhere vanishing section $\varepsilon\in\Gamma(S)$ of an irreducible paired spinor bundle $(S,\Gamma,\cB)$  on $(M,g)$. Such a spinor $\varepsilon\in \Gamma(S)$ is a \emph{supersymmetry parameter} or \emph{supersymmetry generator} for the supersymmetric configuration $(g,b,\phi)$.
\end{definition} 
 
\noindent
Therefore, the supersymmetry parameter of a supersymmetric configuration yields a particular case of a torsion parallel spinor with fully skew-symmetric torsion. We denote by $\Conf_s(\cP, \cA, Y)$ the set of supersymmetric configurations on $(\cP, \cA, Y)$, whose elements consist by definition of tuples $(g,b,\phi,\varepsilon)$ satisfying the differential system \eqref{eq:KSE4d}. Because of their physical origin in the supergravity literature, the first equation in \eqref{eq:KSE4d} is called the \emph{gravitino equation}, whereas the second equation in \eqref{eq:KSE4d} is called the \emph{dilatino equation}. We will occasionally use this terminology below. Since supersymmetric Neveu-Schwarz configurations involve torsion parallel spinors with completely skew-symmetric torsion, we consider the latter separately in the following subsection.


\subsection{The globally hyperbolic reduction of abelian bundle gerbes}
\label{sec:reductionbundlegerbe}


In Section \ref{sec:evolutionskewtorsion}, we obtained the globally hyperbolic reduction of a skew-torsion isotropic parallelism, which plays a fundamental role in reducing globally hyperbolic supersymmetric Neveu-Schwarz solutions. In this section, we address the globally hyperbolic reduction of another fundamental ingredient of Neveu-Schwarz supergravity and its supersymmetric solutions and configurations, namely, an abelian bundle gerbe equipped with a connection and a curving \cite{Murray}. Let $[u,v,l,n]$ be a globally hyperbolic parallelism on $M = \cI \times \Sigma$ with associated metric:
\begin{equation*}
g = u \odot v + l\otimes l + n\otimes n = -\lambda^2_t \dd t\otimes \dd t + h_t\, .
\end{equation*} 

\noindent
Let $(\bar{\cP} , \bar{\cA}, \bar{Y})$ be a bundle gerbe on $M = \cI\times \Sigma$. Since the goal is to apply this reduction to the evolution problem defined by Neveu-Schwarz supergravity, we will assume that the \emph{topological data} contained in $(\cP, \cA, Y)$, namely $(\bar{\cP},\bar{Y})$, is given by the pull-back of a bundle gerbe $(\cP,Y)$ defined over $\Sigma$, while the connective structure $\bar{\cA}$ is allowed to be genuinely time-dependent, not necessarily given by the pull-back of a fixed connective structure on $(\cP,Y)$. We will refer to such bundle gerbes as \emph{reducible} bundle gerbes on $\cI\times \Sigma$. This set-up captures the time-dependence of the evolution problem on a bundle gerbe in its full generality. Since the submersion $\bar{Y} \to M$ of a reducible bundle gerbe on $\cI\times\Sigma$ is, by assumption, a pull-back submersion, it follows that it is isomorphic to:
\begin{equation*}
\bar{Y} = \cI \times Y \to \cI \times \Sigma	
\end{equation*}

\noindent
where the projection on the first factor is the identity map and $Y\to \Sigma$ is the submersion underlying the bundle gerbe $(\cP , Y)$. Similarly:
\begin{equation*}
\bar{Y}\times_M  \bar{Y} = \cI\times Y \times_\Sigma Y \to \cI \times \Sigma
\end{equation*}

\noindent
with the identity map $\cI\to \cI$ acting on the first factor. Since $\bar{\cP}$ is the pull-back of $\cP\to Y\times_{\Sigma} Y$ to $\cI \times Y\times_{\Sigma} Y$, we have:
\begin{eqnarray*}
\bar{\cP}= \cI \times \cP \to \cI \times Y\times_{\Sigma} Y
\end{eqnarray*}

\noindent
again with the identity map $\cI\to \cI$ acting on the first factor. The principal $\U(1)$ action of $\bar{\cP} = \cI \times \cP$ is the one induced by the principal $\U(1)$ action of $\cP$ and the trivial action on the $\cI$ factor. Let $\bar{\cA} \in \Omega^1(\cI \times \cP )$ be a connective structure on $\bar{\cP}$. Then we can write:
\begin{equation*}
\bar{\cA} = \Psi_t \dd t + \cA_t
\end{equation*}  

\noindent
where $\Psi_t$ is a family of functions on $\cP$ and $\cA_t$ is a family of connections on $\cP$. Recall that $t$ denotes the Cartesian coordinate on $\cI$. Since $\bar{\cA}$ is a connection on $\bar{\cP}$, it is in particular invariant under the $\U(1)$ action of the principal bundle $\bar{\cP}$, and since this action is ineffective on $\cI$ we conclude that $\Psi_t$ is in fact a family of invariant functions on $\cP$ whence it descends to a family of invariant functions on $Y\times_{\Sigma} Y$ that we denote by the same symbol for ease of notation. Since $(\bar{\cP} , \bar{\cA}, \bar{Y})$ is a bundle gerbe, it comes equipped with an isomorphism:
\begin{equation*}
 (\bar{\pi}^{\ast}_{12}\bar{\cP}\otimes\bar{\pi}^{\ast}_{23}\bar{\cP} , \bar{\pi}^{\ast}_{12}\bar{\cA}\otimes \bar{\pi}^{\ast}_{23}\bar{\cA})  \xrightarrow{\simeq} (\bar{\pi}^{\ast}_{13} \bar{\cP} , \bar{\pi}^{\ast}_{13}\bar{\cA})
\end{equation*}

\noindent
where $\bar{\pi}\colon \bar{Y} \to M$ and the notation $\pi_{ij} \colon \bar{Y}\times_M \bar{Y}\times_M \bar{Y} \to \bar{Y}\times_M \bar{Y}$ forgets the entry not labelled neither by $i$ nor $j$. This implies $\bar{\delta} \Psi_t = 0$, where $\bar{\delta}$ is the differential of the simplicial manifold determined by $\bar{\pi}\colon \bar{Y} \to M$. Since this differential is exact, it follows that $\Psi_t$ descends through $\bar{\delta}$ to a family of functions on $Y$, unique modulo a choice of a family of elements in $\pi^{\ast}C^{\infty}(\Sigma)$. The preceding discussion thus leads to the following correspondence. 

\begin{prop}
There is a natural one-to-one correspondence between reducible bundle gerbes $(\bar{\cP},\bar{Y},\bar{\cA})$ on $M = \cI\times \Sigma$ and tuples $(\cP,Y,\cA_t , \Psi_t + \pi^{\ast}f_t)$, where  $(\cP,Y)$ is a bundle gerbe on $\Sigma$, $\Psi_t$ is a family of functions on $Y$, $f_t$ is a family of functions on $\Sigma$, and $\cA_t$ is a family of connective structures on $(\cP,Y)$.
\end{prop}

\noindent
For simplicity in the exposition, we will relabel $\Psi_t + \pi^{\ast}f_t$ simply as $\Psi_t$, and we will refer to $(\cP,Y, \cA_t , \Psi_t )$ as the \emph{reduction} of $(\bar{\cP},\bar{Y},\bar{\cA})$. A direct computation gives the following result.

\begin{lemma}
\label{lemma:gerbecurvaturedecomposition}
The curvature $\cF_{\bar{\cA}} \in \Omega^2(\cI \times Y\times_{\Sigma} Y)$ of $\bar{\cA} = \Psi_t \dd t + \cA_t$ satisfies the following equation:
\begin{eqnarray*}
\cF_{\bar{\cA}} = \dd_{Y^{[2]}} \Psi_t \wedge \dd t + \dd t \wedge \partial_t \cA_t + \cF_{\cA_t}
\end{eqnarray*}

\noindent
where $\dd_{Y^{[2]}}\colon \Omega^{\bullet}(Y\times_{\Sigma} Y) \to \Omega^{\bullet}(Y\times_{\Sigma} Y) $ is the exterior derivative on $\Omega^{\bullet}(Y\times_{\Sigma} Y)$ and $\cF_{\cA_t}$ is the family of curvatures of $\cA_t$.
\end{lemma}

\noindent
A curving on $(\bar{\cP},\bar{Y},\bar{\cA})$ is by definition a two-form $\bar{b}\in \Omega^2(\bar{Y})$ satisfying:
\begin{equation*}
\bar{\delta} \bar{b} = F_{\bar{\cA}}
\end{equation*}

\noindent 
Since $\bar{Y}\times_M \bar{Y} = \cI \times Y \times_{\Sigma} Y$, we can write every curving $\bar{b}$ on $(\bar{\cP},\bar{Y},\bar{\cA})$ as follows:
\begin{equation}
\label{eq:decompositionbarb}
\bar{b} = \dd t \wedge a_t + b_t
\end{equation}

\noindent
for uniquely determined families of one forms $a_t$ and two-forms $b_t$ on $Y$.
\begin{lemma}
\label{lemma:reduciblegerbe}
Let $(\bar{\cP},\bar{Y},\bar{\cA})$ be a reducible bundle gerbe with connective structure on $M = \cI \times \Sigma$ and let  $(\cP,Y, \cA_t, \Psi_t)$ be its reduction on $\Sigma$. A two-form $\bar{b}\in \Omega^2(\bar{Y})$ is a curving on $(\bar{\cP},\bar{Y},\bar{\cA})$ if and only if, in the decomposition given in \eqref{eq:decompositionbarb}, $b_t$ is a family of curvings on $(\cP,Y,\cA_t)$ and $a_t$ is a family of one-forms on $Y$ satisfying:
\begin{equation}
\label{eq:conditionat}
\delta a_t = \partial_t \cA_t - \dd_{Y^{[2]}} \Psi_t
\end{equation}

\noindent
for every $t\in \cI$, where $\delta$ is the simplicial differential of the simplicial manifold determined by submersion $\pi \colon Y \to M$.
\end{lemma}

\begin{proof}
Plugging equation \eqref{eq:decompositionbarb} in $\bar{\delta} \bar{b} = F_{\bar{\cA}}$ we have:
\begin{equation*}
\bar{\delta}\bar{b} = \bar{\delta}(\dd t \wedge a_t + b_t) = \dd t \wedge \delta a_t + \delta b_t = \cF_{\bar{\cA}} = \dd_{Y^{[2]}} \Psi_t \wedge \dd t + \dd t \wedge \partial_t \cA_t + \cF_{\cA_t}
\end{equation*}

\noindent
where we have used Lemma \ref{lemma:gerbecurvaturedecomposition}. Isolating by tensor type this gives:
\begin{equation*}
\delta b_t = F_{\cA_t}\, , \qquad \delta a_t = \partial_t \cA_t - \dd_{Y^{[2]}} \Psi_t
\end{equation*}

\noindent
and thus we conclude. 
\end{proof}

\begin{cor}
There is a natural one-to-one correspondence between curvings on a reducible bundle gerbe $(\bar{\cP},\bar{Y},\bar{\cA})$ with reduction $(\cP,Y,\cA_t,\Psi_t)$ and pairs $(a_t,b_t)$ consisting of a family $b_t$ of curvings on $(\cP,Y,\cA_t)$ and a family of one-forms $a_t$ on $Y$ satisfying Equation \eqref{eq:conditionat}.
\end{cor}

\noindent
Let $\bar{b} = \dd t \wedge a_t + b_t$ be a curving on a reducible bundle gerbe over $\cI\times \Sigma$. Then:
\begin{equation*}
\dd_{\bar{Y}}\bar{b} = \dd t \wedge (\partial_t b_t - \dd_Y a_t) + \dd_Y b_t
\end{equation*}

\noindent
where $\dd_{\bar{Y}}$ is the exterior derivative on $\bar{Y}$ and $\dd_Y$ is the exterior derivative on $Y$. 

\begin{lemma}
The family of two-forms $(\partial_t b_t - \dd_Y a_t) \in \Omega^2(Y)$ satisfies $\delta (\partial_t b_t - \dd_Y a_t) = 0$ for every $t\in \cI$, where $\delta$ denotes the simplicial differential of the simplicial manifold determined by $\pi \colon Y\to \Sigma$. 
\end{lemma}

\begin{proof}
We compute:
\begin{equation*}
\delta (\partial_t b_t - \dd_Y a_t) = \partial_t \delta b_t - \dd_{Y^{[2]}} \delta a_t =  \partial_t \cF_{\cA_t}- \dd_{Y^{[2]}} (\partial_t \cA_t - \dd_{Y^{[2]}} \Psi_t) = 0
\end{equation*}

\noindent
where we have used Lemma \ref{lemma:reduciblegerbe}.
\end{proof}

\noindent
Since $\delta$ is exact, by the previous lemma it follows that there exists a unique family of two-forms $c_t \in \Omega^2(\Sigma)$ on $\Sigma$ such that:
\begin{equation*}
\pi^{\ast} c_t = \partial_t b_t - \dd_Y a_t
\end{equation*}

\noindent
Since $\dd_{\bar{Y}} \bar{b} = \bar{\pi}^{\ast} H_{\bar{b}}$ for a uniquely determined three-form $\bar{H}_{\bar{b}}\in \Omega^3(M)$, it follows that:
\begin{equation*}
H_{\bar{b}} = \dd t \wedge c_t + H_{b_t}
\end{equation*}

\noindent
where $H_{b_t}$ is the curvature of the family of curvings $b_t$. Therefore, every curving $\bar{b}$ on a reducible bundle gerbe is equivalent to a pair $(a_t , b_t)$ as introduced in Lemma \ref{lemma:reduciblegerbe}, and this pair defines in turn a family of two-forms $c_t$ on $\Sigma$ as described above. 

\begin{definition}
The family of two-forms $c_t$ is the \emph{derived} family of two-forms of the curving $\bar{b}$ on the reducible bundle gerbe $(\bar{\cP},\bar{\cA},\bar{Y})$.
\end{definition} 


\subsection{Globally hyperbolic Neveu-Schwarz solutions}
\label{sec:NSNSsystemreduction}

Our first task is to characterize the globally hyperbolic Neveu-Schwarz solutions on $(\bar{\cP},\bar{\cA},\bar{Y})$ in terms of families of solutions on $(\cP, Y , \cA_t  , \Psi_t)$.

\begin{definition}
A Neveu-Schwarz configuration $(g ,\bar{b},\bar{\phi})\in \Conf(\bar{\cP},\bar{\cA},\bar{Y})$ is \emph{globally hyperbolic} if $g$ is a globally hyperbolic metric on $\cI\times \Sigma$ of the form $g = - \lambda_t^2 \dd t\otimes \dd t + h_t$ for a family of nowhere-vanishing functions $\lambda_t$ and a family of Riemannian metrics $h_t$ on $\Sigma$.   
\end{definition}

\begin{remark}
By the seminal work of Bernal and S\'anchez \cite{Bernal:2003jb,Bernal:2004gm}, every globally hyperbolic Lorentzian manifold is \emph{globally} isometric to a Lorentzian manifold of the form $(\mathbb{R}\times \Sigma, - \lambda_t^2 \dd t\otimes \dd t + h_t)$. Hence, there is no loss of generality in our previous definition. 
\end{remark}

\noindent
Since $\bar{\phi} \colon M = \cI \times \Sigma \to S^1$ is a circle-valued smooth function on $\cI\times \Sigma$, it follows that it is equivalent to a family of smooth circle-valued functions $\phi_t$ parametrized by $t\in \cI$. Therefore, globally hyperbolic configurations $(g , \bar{b}, \bar{\phi}) \in \Conf(\bar{\cP},\bar{\cA},\bar{Y})$ are in one-to-one correspondence with tuples $( h_t , a_t , b_t , \phi_t , \lambda_t)$ as described above. We will refer to the latter as the \emph{reduction} of the former. 

\noindent
For every globally hyperbolic Neveu-Schwarz configuration, set:
\begin{equation*}
\mathfrak{t}_t:=\lambda_t\dd t\, ,  \qquad D_t:=\lambda_t^{-1}\partial_t\, , \qquad \eta_t:=\lambda_t^{-1}\dd_{\Sigma}\lambda_t\, , \qquad \Theta_t:=-\frac{1}{2}D_t h_t.
\end{equation*}
Here $D_t$ acts coefficientwise in the fixed product splitting, and the sign of $\lambda_t$ is chosen consistently with the normal orientation so that $\mathfrak{t}_t$ is the chosen unit time-like normal one-form. For clarity in the exposition, when understanding $D_t$ as a family of vector fields on $\Sigma$, we will denote by $N_t \in \mathfrak{X}(\Sigma)$. A direct calculation proves the following lemma.  

\begin{lemma}
\label{lemma:hyperbolicderivatives}
Let $g = - \lambda_t^2 \dd t \otimes \dd t + h_t$ be a globally hyperbolic metric on $M = \cI\times \Sigma$. Then, for every family of one-forms $\varphi_{\phi_t} \in \Omega^1(\Sigma)$ we have:
\begin{eqnarray*}
& \nabla^g \mathfrak{t}_t = \Theta_t-\mathfrak{t}_t\otimes\eta_t\, , \\
& \nabla^g\varphi_{\phi_t} = \mathfrak{t}_t\otimes D_t\varphi_{\phi_t} + \mathfrak{t}_t\odot \Theta_t (\beta^{\sharp}_t) + \nabla^{h_t}\varphi_{\phi_t} - \langle\varphi_{\phi_t},\eta_t\rangle_{h_t}\mathfrak{t}_t\otimes\mathfrak{t}_t\\
& \dd_{\nabla^g} \Theta_t = \dd_{\nabla^{h_t}} \Theta_t + \mathfrak{t}_t\wedge D_t\Theta_t + \mathfrak{t}_t\wedge \Theta_t\circ_{h_t} \Theta_t - (\mathfrak{t}_t\wedge\Theta_t(\eta_t^{\sharp}))\otimes\mathfrak{t}_t 
\end{eqnarray*}

\noindent
where $\Theta_t , \Theta_t\circ_{h_t}\Theta_t \in \Gamma(T^{\ast}\Sigma \odot T^{\ast}\Sigma)$ are considered as one-forms taking values on $T^{\ast}\Sigma$ and we have set:
\begin{eqnarray*}
\Theta_t\circ_{h_t}\Theta_t (v_1,v_2) = \langle \Theta_t(v_1) , \Theta_t(v_2)\rangle_{h_t}
\end{eqnarray*} 

\noindent
for every $v_1,v_2 \in T\Sigma$.
\end{lemma}

\noindent
Using the previous lemma, we compute the Ricci and scalar curvatures of the globally hyperbolic metric $g$ in terms of families of objects on $\Sigma$.
\begin{lemma}
\label{lemma:hyperboliccurvatures}
Let $g = - \lambda_t^2 \dd t \otimes \dd t + h_t$ be a globally hyperbolic metric on $M = \cI\times \Sigma$. The Ricci and scalar curvatures of $g$ are given by:
\begin{eqnarray*}
& \Ric^g = (\Tr_{h_t}(D_t\Theta_t) + \vert \Theta_t \vert_{h_t}^2 - \nabla^{h_t\ast}\eta_t + \vert\eta_t\vert_{h_t}^2)\mathfrak{t}_t\otimes\mathfrak{t}_t + \Ric^{h_t} + \Tr_{h_t}(\Theta_t) \Theta_t - 2 \Theta_t\circ_{h_t} \Theta_t\\
& \qquad - D_t\Theta_t - \nabla^{h_t}\eta_t-\eta_t\otimes\eta_t + \mathfrak{t}_t \odot (\dd_{\Sigma} \Tr_{h_t}(\Theta_t) + \nabla^{h_t\ast}\Theta_t)\\
& s^g = s^{h_t} + \Tr_{h_t}(\Theta_t)^2 - 2 \Tr_{h_t}(D_t\Theta_t) - 3 \vert \Theta_t\vert_{h_t}^2 + 2\nabla^{h_t\ast}\eta_t-2\vert\eta_t\vert_{h_t}^2 
\end{eqnarray*}

\noindent
where $\Theta_t\in \Gamma(T^{\ast}\Sigma\odot T^{\ast}\Sigma)$ is understood as a one-form taking values on $T^{\ast}\Sigma$.
\end{lemma}  

\noindent
We decompose now each of the equations of the Neveu-Schwarz system \eqref{eq:NSNSsystem} evaluated on globally hyperbolic configurations in terms of \emph{flow equations} on $\Sigma$. A direct computation proves the following lemma.

\begin{lemma}
Let $(g,\bar{b},\bar{\phi})\in \Conf(\bar{\cP},\bar{\cA},\bar{Y})$ be a globally hyperbolic Neveu-Schwarz configuration with reduction $( h_t , a_t , b_t , \phi_t,\lambda_t)$. Then:
\begin{eqnarray*}
\varphi_{\bar{\phi}} = p_t\mathfrak{t}_t + \varphi_{\phi_t}\, , \qquad p_t:= D_t\phi_t\, ,  
\end{eqnarray*}
	
\noindent
where $\varphi_{\phi_t}$ is the curvature of $\phi_t$. Furthermore, there exist a time-independent closed one-form $\varrho \in \Omega^1(\Sigma)$ such that $\varphi_{\phi_t} = \dd_{\Sigma} \phi_t + \varrho$.
\end{lemma} 

\begin{remark}
For a family of circle-valued dilatons $\phi_t\colon\Sigma\to S^1$, the notation $D_t\phi_t$ is understood through the Maurer--Cartan form $\mu$ of $S^1$, namely:
\begin{equation*}
D_t\phi_t:=\lambda_t^{-1}\mu_{\phi_t}(\partial_t\phi_t).
\end{equation*}

\noindent
This defines a real-valued function without requiring a global real-valued lift of $\phi_t$. In particular, the evolution equation $D_t\phi_t=\varphi_{\phi_t}(e^{\sharp}_{tu})$ is equivalent to the intrinsic transport equation:
\begin{equation*}
\dd\bar{\phi} (\partial_t-\lambda_t e^{\sharp}_{tu} )=0,
\end{equation*}

\noindent
where $\bar{\phi}(t,x)=\phi_t(x)$.
\end{remark} 

\noindent
In terms of the derived family $c_t$ defined in Section \ref{sec:reductionbundlegerbe}, set:
\begin{equation*}
\theta_t:=\lambda_t^{-1}\ast_{h_t}c_t\, , \qquad H_{b_t}:=f_{b_t}\nu_{h_t}.
\end{equation*}
Then $H_{\bar b}=\mathfrak{t}_t\wedge\ast_{h_t}\theta_t+f_{b_t}\nu_{h_t}$. We will refer to $\theta_t$ as the family of \emph{derived one-forms} associated to $(g,\bar{b},\bar{\phi})$ or, equivalently, $( h_t , a_t , b_t , \phi_t,\lambda_t)$.

\begin{lemma}
\label{lemma:Einsteingh}
A globally hyperbolic Neveu-Schwarz configuration $(g,\bar{b},\bar{\phi})\in \Conf(\bar{\cP},\bar{\cA},\bar{Y})$ satisfies the Einstein equation of the Neveu-Schwarz system \eqref{eq:NSNSsystem} if and only if its reduction $( h_t , a_t , b_t , \phi_t,\lambda_t)$ satisfies the following differential system:
\begin{eqnarray*}
& D_tp_t  - \nabla^{h_t\ast}\eta_t + \Tr_{h_t}(D_t\Theta_t) + \vert \Theta_t \vert_{h_t}^2   + \vert\eta_t\vert_{h_t}^2 - \langle\varphi_{\phi_t},\eta_t\rangle_{h_t} - \frac{1}{2} \vert\theta_t\vert_{h_t}^2 = 0\\
& \dd_{\Sigma} \Tr_{h_t}(\Theta_t) + \nabla^{h_t\ast}\Theta_t + D_t\varphi_{\phi_t} +\Theta_t (\varphi_{\phi_t}^{\sharp}) -p_t\eta_t - \frac{1}{2}  f_{b_t} \theta_t  = 0\\
& D_t\Theta_t = \Ric^{h_t} + (\Tr_{h_t}(\Theta_t) + p_t) \Theta_t - 2 \Theta_t\circ_{h_t} \Theta_t + \nabla^{h_t} \varphi_{\phi_t} - \nabla^{h_t}\eta_t-\eta_t\otimes\eta_t  - \frac{1}{2} \theta_t \otimes \theta_t + \frac{1}{2} (\vert \theta_t\vert_{h_t}^2 - f_{b_t}^2) h_t   
\end{eqnarray*}

\noindent
where $\Theta_t$ is the fundamental form of $(\left\{ t\right\} \times \Sigma, h_t)$, $f_{b_t}\in C^{\infty}(\Sigma)$ is the Hodge dual of $H_{b_t}$ and $\theta_t$ is the derived family of one-forms of $( h_t , a_t , b_t , \phi_t,\lambda_t)$. 
\end{lemma}

\begin{remark}
Note that equations:
\begin{equation*}
\varphi_{\phi_t} = \dd_{\Sigma}\psi_t + \rho\, , \qquad D_t f_{b_t}=\Tr_{h_t}(\Theta_t)f_{b_t}-\nabla^{h_t\ast}\theta_t+\langle\eta_t,\theta_t\rangle_{h_t}
\end{equation*}
need not be included in the previous lemma since they are a consequence of $( h_t , a_t , b_t , \phi_t,\lambda_t)$ being the \emph{reduction} of $(g,\bar{b},\bar{\phi})\in \Conf(\bar{\cP},\bar{\cA},\bar{Y})$. 
\end{remark}

\begin{proof}
Let $(g,\bar{b},\bar{\phi})$ be a globally hyperbolic configuration with reduction $( h_t , a_t , b_t , \phi_t,\lambda_t)$. We evaluate the Einstein equation of the Neveu-Schwarz system, namely the first equation in \eqref{eq:NSNSsystem}, on the following configuration:
\begin{eqnarray*}
& g = - \lambda_t^2\dd t\otimes \dd t + h_t \in \Gamma(T^{\ast}M\odot T^{\ast}M) \, , \quad \bar{b} = \dd t\wedge a_t + b_t \in \Omega^2(\cI\times Y)\\
& \bar{\phi}(t,p) = \phi_t(p)\in S^1
\end{eqnarray*}

\noindent
for every $(t,p)\in \cI\times \Sigma$. By Lemma \ref{lemma:hyperboliccurvatures}, the Ricci tensor of $g$ is given by:
\begin{eqnarray*}
& \Ric^{g} (N_t , N_t) =  \Tr_{h_t}(D_t\Theta_t) + \vert \Theta_t\vert^2_{h_t}-\nabla^{h_t\ast}\eta_t+\vert\eta_t\vert_{h_t}^2 \,  ,\quad \Ric^g(N_t)\vert_{T\Sigma} = \dd_{\Sigma} \Tr_{h_t}(\Theta_t) + \nabla^{h_t\ast}\Theta_t\\
& \Ric^{g}\vert_{T\Sigma\times T\Sigma} = \Ric^{h_t} + \Tr_{h_t}(\Theta_t) \Theta_t - 2 \Theta_t\circ_{h_t} \Theta_t - D_t \Theta_t-\nabla^{h_t}\eta_t-\eta_t\otimes\eta_t 
\end{eqnarray*}

\noindent
where $\Theta_t$ is the second fundamental form of the embedding $\left\{ t\right\}\times \Sigma \hookrightarrow \cI \times \Sigma$ as defined in Equation \eqref{eq:secondfundamentalform} with respect to the metric $g$. By the discussion of Section \ref{sec:reductionbundlegerbe} we have:
\begin{equation*}
H_{\bar{b}} = \dd t \wedge c_t + H_{b_t} = \mathfrak{t}_t \wedge \ast_{h_t} \theta_t + H_{b_t}
\end{equation*}

\noindent
from which we compute:
\begin{eqnarray*}
& (H_{\bar{b}}\circ_{g} H_{\bar{b}}) (N_t , N_t) = \vert \theta_t \vert_{h_t}^2\\
& (H_{\bar{b}}\circ_{g} H_{\bar{b}}) (N_t , v) = f_{b_t} \theta_t(v)\, , \quad v\in T\Sigma
\end{eqnarray*}

\noindent
where we have defined the function $f_{b_t} \in C^{\infty}(\Sigma)$ associated to every family of curvings $b_t$ and Riemannian metrics $h_t$ by $H_{b_t} := f_{b_t} \nu_{h_t}$. Furthermore:
\begin{eqnarray*}
& (H_{\bar{b}}\circ_{g} H_{\bar{b}})(v_1,v_2) = \theta_t(v_1)\theta_t(v_2)+(f_{b_t}^2-\vert\theta_t\vert_{h_t}^2)h_t(v_1,v_2)
\end{eqnarray*}

\noindent
for every $v_1 ,v_2 \in T\Sigma$. On the other hand, using Lemma \ref{lemma:hyperbolicderivatives} we obtain:
\begin{eqnarray*}
& \nabla^{g}_{N_t}\varphi_{\bar{\phi}} = (D_tp_t-\langle\varphi_{\phi_t},\eta_t\rangle_{h_t})\mathfrak{t}_t+D_t\varphi_{\phi_t}+\Theta_t(\varphi_{\phi_t}^{\sharp})-p_t\eta_t\\
& \nabla^{g}_v \varphi_{\bar{\phi}} = (\dd_{\Sigma}p_t(v)+\Theta_t(\varphi_{\phi_t}^{\sharp},v))\mathfrak{t}_t+p_t\Theta_t(v)+\nabla^{h_t}_v\varphi_{\phi_t}
\end{eqnarray*}

\noindent
for every $v \in T\Sigma$. With these provisos in mind, we compute:
\begin{eqnarray*}
& \Ric^g(N_t , N_t) + (\nabla^g_{N_t}\varphi_{\bar{\phi}})(N_t) - \frac{1}{2} (H_{\bar{b}}\circ_{g} H_{\bar{b}}) (N_t , N_t) =\\
& \qquad D_tp_t+\Tr_{h_t}(D_t\Theta_t)+\vert\Theta_t\vert_{h_t}^2-\nabla^{h_t\ast}\eta_t+\vert\eta_t\vert_{h_t}^2 - \langle\varphi_{\phi_t},\eta_t\rangle_{h_t}-\frac{1}{2}\vert\theta_t\vert_{h_t}^2 = 0
\end{eqnarray*}

\noindent
for the \emph{time-time} component of the Einstein equation in \eqref{eq:NSNSsystem}. For the mixed \emph{time-space} components we have:
\begin{eqnarray*}
& (\Ric^g(N_t) + (\nabla^g_{N_t}\varphi_{\bar{\phi}})  - \frac{1}{2} (H_{\bar{b}}\circ_{g} H_{\bar{b}}) (N_t ))\vert_{T\Sigma} = \\
& = \dd_{\Sigma} \Tr_{h_t}(\Theta_t) + \nabla^{h_t\ast}\Theta_t+D_t\varphi_{\phi_t}+\Theta_t(\varphi_{\phi_t}^{\sharp})-p_t\eta_t-\frac{1}{2}f_{b_t}\theta_t = 0
\end{eqnarray*}

\noindent
Finally, for the \emph{space-space} components of the Einstein equation we have:
\begin{eqnarray*}
& 0 = (\Ric^g  + (\nabla^g \varphi_{\bar{\phi}})  - \frac{1}{2} H_{\bar{b}}\circ_{g} H_{\bar{b}})\vert_{T\Sigma\times T\Sigma} = \\
& = \Ric^{h_t} + (\Tr_{h_t}(\Theta_t) + p_t) \Theta_t - 2 \Theta_t\circ_{h_t} \Theta_t - D_t \Theta_t + \nabla^{h_t} \varphi_{\phi_t}-\nabla^{h_t}\eta_t-\eta_t\otimes\eta_t\\
& \qquad - \frac{1}{2} \theta_t \otimes \theta_t + \frac{1}{2} (\vert \theta_t\vert_{h_t}^2 - f_{b_t}^2) h_t 
\end{eqnarray*}

\noindent
and hence we conclude. 
\end{proof}

\begin{lemma}
\label{lemma:Maxwellgh}
A globally hyperbolic Neveu-Schwarz configuration $(g,\bar{b},\bar{\phi})\in \Conf(\bar{\cP},\bar{\cA},\bar{Y})$ satisfies the Maxwell equation of the Neveu-Schwarz system \eqref{eq:NSNSsystem} if and only if its reduction $( h_t , a_t , b_t , \phi_t,\lambda_t)$ satisfies the following differential system:
\begin{equation*}
D_t \theta_t = \dd_{\Sigma} f_{b_t} + f_{b_t}\eta_t+p_t\theta_t-f_{b_t}\varphi_{\phi_t}  \, , \qquad \dd_{\Sigma} \theta_t = \varphi_{\phi_t} \wedge \theta_t
\end{equation*}
	
\noindent
where $f_{b_t}\in C^{\infty}(\Sigma)$ is the Hodge dual of $H_{b_t}$ and $\theta_t$ is the derived family of one-forms of $( h_t , a_t , b_t , \phi_t,\lambda_t)$. 
\end{lemma}

\begin{proof}
We consider the Maxwell equation as written in \eqref{eq:NSNSsystem4d}, namely:
\begin{equation*}
\dd \alpha_{\bar{b}} = \varphi_{\bar{\phi}} \wedge \alpha_{\bar{b}}
\end{equation*}

\noindent
We have:
\begin{equation*}
\alpha_{\bar{b}} = \ast_g (\mathfrak{t}_t \wedge \ast_{h_t}\theta_t + f_{b_t} \nu_{h_t}) = - (\theta_t + f_{b_t} \mathfrak{t}_t)
\end{equation*}

\noindent
where we have used equations \eqref{eq:hodgehyperbolic}. Hence, we have:
\begin{eqnarray*}
& \dd \alpha_{\bar{b}} + \alpha_{\bar{b}}\wedge \varphi_{\bar{\phi}} = -\dd (\theta_t + f_{b_t} \mathfrak{t}_t) - (\theta_t + f_{b_t} \mathfrak{t}_t) \wedge (p_t \mathfrak{t}_t + \varphi_{\phi_t})\\
& = \mathfrak{t}_t\wedge (\dd_{\Sigma} f_{b_t} - D_t \theta_t+f_{b_t}\eta_t) - \dd_{\Sigma} \theta_t + \mathfrak{t}_t \wedge (p_t\theta_t - f_{b_t}\varphi_{\phi_t})-\theta_t\wedge\varphi_{\phi_t} = 0
\end{eqnarray*}

\noindent
Isolating by tensor type we conclude.
\end{proof}

\begin{lemma}
\label{lemma:Dilatongh}
A globally hyperbolic Neveu-Schwarz configuration $(g,\bar{b},\bar{\phi})\in \Conf(\bar{\cP},\bar{\cA},\bar{Y})$ satisfies the dilaton equation of the Neveu-Schwarz system \eqref{eq:NSNSsystem} if and only if its reduction $( h_t , a_t , b_t , \phi_t,\lambda_t)$ satisfies the following differential system:
\begin{equation*}
D_tp_t + \nabla^{h_t\ast}\varphi_{\phi_t}  - (\Tr_{h_t}(\Theta_t) + p_t) p_t-\langle\varphi_{\phi_t},\eta_t\rangle_{h_t} - f_{b_t}^2 + \vert \theta_t \vert_{h_t}^2 + \vert \varphi_{\phi_t} \vert_{h_t}^2 = 0
\end{equation*}
 
\noindent
where $\Theta_t$ is the fundamental form of $(\left\{ t\right\} \times \Sigma, h_t)$, $f_{b_t}\in C^{\infty}(\Sigma)$ is the Hodge dual of $H_{b_t}$ and $\theta_t$ is the unit-normal derived family of one-forms of $( h_t , a_t , b_t , \phi_t,\lambda_t)$. 
\end{lemma}

\begin{proof}
We compute:
\begin{eqnarray*}
\vert \varphi_{\bar{\phi}}\vert_g^2 + \vert\alpha_{\bar{b}}\vert^2_g =  \vert \theta_t \vert_{h_t}^2 + \vert \varphi_{\phi_t} \vert_{h_t}^2 - p_t^2 - f_{b_t}^2 
\end{eqnarray*}

\noindent
On the other hand we have:
\begin{equation*}
\nabla^{g\ast}\varphi_{\bar{\phi}}  = \nabla^{g\ast}(p_t\mathfrak{t}_t + \varphi_{\phi_t}) = D_tp_t + \nabla^{h_t\ast}\varphi_{\phi_t}-\Tr_{h_t}(\Theta_t)p_t-\langle\varphi_{\phi_t},\eta_t\rangle_{h_t}
\end{equation*}

\noindent
and hence  we conclude.
\end{proof}

\noindent
The previous lemmata allow us to characterize globally hyperbolic Neveu-Schwarz solutions on a reducible pair $(\bar{\cP},\bar{\cA},\bar{Y})$ in terms of constrained evolution flows for tuples $( h_t , a_t , b_t , \phi_t,\lambda_t)$ on the reduction $(\cP,Y,\cA_t,\Psi_t,\cX)$ of $(\bar{\cP},\bar{\cA},\bar{Y})$.
\begin{prop}
\label{prop:normalhyperbolicflow}
A globally hyperbolic Neveu-Schwarz configuration $(g,\bar{b},\bar{\phi})\in \Conf(\bar{\cP},\bar{\cA},\bar{Y})$ satisfies the Neveu-Schwarz system \eqref{eq:NSNSsystem} if and only if its reduction $( h_t , a_t , b_t , \phi_t,\lambda_t)$ satisfies the following system of evolution equations:
\begin{eqnarray}
& D_t \Theta_t = \Ric^{h_t} + (\Tr_{h_t}(\Theta_t) + p_t) \Theta_t - 2 \Theta_t\circ_{h_t} \Theta_t + \nabla^{h_t} \varphi_{\phi_t} \nonumber \\
& \qquad -\nabla^{h_t}\eta_t-\eta_t\otimes\eta_t - \frac{1}{2} \theta_t \otimes \theta_t + \frac{1}{2} (\vert \theta_t\vert_{h_t}^2 - f_{b_t}^2) h_t \nonumber \\
& D_t \theta_t = \dd_{\Sigma} f_{b_t}+f_{b_t}\eta_t+p_t\theta_t-f_{b_t}\varphi_{\phi_t} \label{eq:normalevolution} \\
& D_tp_t + \nabla^{h_t\ast}\varphi_{\phi_t} = (\Tr_{h_t}(\Theta_t) + p_t)p_t+\langle\varphi_{\phi_t},\eta_t\rangle_{h_t}+f_{b_t}^2-\vert\theta_t\vert_{h_t}^2-\vert\varphi_{\phi_t}\vert_{h_t}^2 \nonumber
\end{eqnarray}
together with the following system of time-dependent constraints:
\begin{eqnarray}
& \cH_t := s^{h_t}  + (\Tr_{h_t}(\Theta_t) + p_t)^2 - \vert\Theta_t\vert^2_{h_t} - \vert\varphi_{\phi_t}\vert_{h_t}^2 - \frac{1}{2}( f_{b_t}^2 + \vert \theta_t \vert^2_{h_t} ) - 2 \nabla^{h_t\ast}\varphi_{\phi_t} = 0 \nonumber \\
&  \mathcal{W}_t := \dd_{\Sigma}(\Tr_{h_t}(\Theta_t)+p_t) + \nabla^{h_t\ast}\Theta_t+\Theta_t(\varphi_{\phi_t}^{\sharp})-\frac{1}{2}f_{b_t}\theta_t = 0\, , \label{eq:normalconstraint}\\
& \mathcal{M}_t := \dd_{\Sigma} \theta_t -\varphi_{\phi_t} \wedge \theta_t = 0 \nonumber
\end{eqnarray}

\noindent
where $\Theta_t$ is the fundamental form of $(\left\{ t\right\} \times \Sigma, h_t)$, $f_{b_t}\in C^{\infty}(\Sigma)$ is the Hodge dual of $H_{b_t}$ and $\theta_t$ is the unit-normal derived family of one-forms of $( h_t , a_t , b_t , \phi_t,\lambda_t)$. 
\end{prop}  

\begin{proof}
Lemmas \ref{lemma:Einsteingh}, \ref{lemma:Maxwellgh} and \ref{lemma:Dilatongh} immediately imply that a globally hyperbolic Neveu-Schwarz configuration $(g,\bar{b},\bar{\phi})\in \Conf(\bar{\cP},\bar{\cA},\bar{Y})$ satisfies the Neveu-Schwarz system \eqref{eq:NSNSsystem} if and only if its reduction $( h_t , a_t , b_t , \phi_t,\lambda_t)$ on $(\cP,Y,\cA_t,\Psi_t,\cX)$ satisfies the following system of evolution equations:
\begin{eqnarray*}
& D_t \Theta_t = \Ric^{h_t} + (\Tr_{h_t}(\Theta_t) + p_t) \Theta_t - 2 \Theta_t\circ_{h_t} \Theta_t + \nabla^{h_t}\varphi_{\phi_t}\\
& \qquad -\nabla^{h_t}\eta_t-\eta_t\otimes\eta_t - \frac{1}{2} \theta_t \otimes \theta_t + \frac{1}{2} (\vert \theta_t\vert_{h_t}^2 - f_{b_t}^2) h_t \\
& D_t \theta_t = \dd_{\Sigma} f_{b_t}+f_{b_t}\eta_t+p_t\theta_t-f_{b_t}\varphi_{\phi_t} \\
& D_tp_t + \nabla^{h_t\ast}\varphi_{\phi_t} = (\Tr_{h_t}(\Theta_t) + p_t)p_t+\langle\varphi_{\phi_t},\eta_t\rangle_{h_t}+f_{b_t}^2-\vert\theta_t\vert_{h_t}^2-\vert\varphi_{\phi_t}\vert_{h_t}^2
\end{eqnarray*}
together with the following system of time-dependent constraints:
\begin{eqnarray*}
& D_tp_t+\Tr_{h_t}(D_t\Theta_t)+\vert\Theta_t\vert_{h_t}^2-\nabla^{h_t\ast}\eta_t+\vert\eta_t\vert_{h_t}^2-\langle\varphi_{\phi_t},\eta_t\rangle_{h_t}-\frac{1}{2}\vert\theta_t\vert_{h_t}^2 = 0\\
& \dd_{\Sigma}\theta_t=\varphi_{\phi_t}\wedge\theta_t \, , \quad \dd_{\Sigma}(\Tr_{h_t}(\Theta_t)+p_t)+\nabla^{h_t\ast}\Theta_t+\Theta_t(\varphi_{\phi_t}^{\sharp})-\frac{1}{2}f_{b_t}\theta_t=0
\end{eqnarray*}

\noindent
Taking the trace of the first evolution equation above and combining it with the third evolution equation, we obtain:
\begin{eqnarray*}
& D_tp_t+\Tr_{h_t}(D_t\Theta_t)=s^{h_t}+(\Tr_{h_t}(\Theta_t)+p_t)^2-2\vert\Theta_t\vert_{h_t}^2-\vert\varphi_{\phi_t}\vert_{h_t}^2-\frac{1}{2}f_{b_t}^2\\
& \qquad -2\nabla^{h_t\ast}\varphi_{\phi_t}+\nabla^{h_t\ast}\eta_t-\vert\eta_t\vert_{h_t}^2+\langle\varphi_{\phi_t},\eta_t\rangle_{h_t}
\end{eqnarray*}

\noindent
Together with the first time-dependent constraint above, this implies:
\begin{eqnarray*}
s^{h_t}+(\Tr_{h_t}(\Theta_t)+p_t)^2-\vert\Theta_t\vert_{h_t}^2-\vert\varphi_{\phi_t}\vert_{h_t}^2-\frac{1}{2}(f_{b_t}^2+\vert\theta_t\vert_{h_t}^2)-2\nabla^{h_t\ast}\varphi_{\phi_t}=0
\end{eqnarray*} 

\noindent
and hence we conclude. 
\end{proof}

\noindent
The previous proposition presents the evolution and constraint equations of the Neveu-Schwarz system in the string frame and serves as the starting point for studying the evolution problem posed by Neveu-Schwarz supergravity. By the seminal results of Oskar Schiller \cite{Schiller:2026vls}, this Cauchy problem for the Neveu-Schwarz system is well-posed. The time-dependent constraints of the previous proposition can be interpreted in terms of the standard general relativity terminology as follows: $\cH_t$ is the \emph{Hamiltonian constraint} of the system, $\mathcal{W}_t$ is the \emph{momentum constraint} of the system, and $\cM_t$ is the constraint equation of the Maxwell equations, or the \emph{Maxwell constraint equations} for short. Hence, by evaluating at the initial time $t=0$, we define:
\begin{eqnarray}
& \cH  := s^{h}  + (\Tr_{h}(\Theta) + \mathfrak{p})^2 - \vert\Theta\vert^2_{h} - \vert\varphi_{\phi}\vert_{h}^2 - \frac{1}{2}( f_{b}^2 + \vert \theta \vert^2_{h} ) - 2 \nabla^{h\ast}\varphi_{\phi}   \nonumber \\
&  \mathcal{W} := \dd_{\Sigma}(\Tr_{h}(\Theta) + \mathfrak{p}) + \nabla^{h\ast}\Theta + \Theta(\varphi_{\phi}^{\sharp})-\frac{1}{2}f_{b}\theta  \, , \label{eq:normalconstraintt0}\\
& \mathcal{M} := \dd_{\Sigma} \theta -\varphi_{\phi} \wedge \theta \nonumber
\end{eqnarray}

\noindent
Hence, by the results of \cite{Schiller:2026vls} adapted to the bundle gerbe case, if $\Sigma$ is an oriented and connected three-manifold equipped with a bundle gerbe $(\cP,\cA,Y)$ and a tuple $(h,\theta,b,\phi,\mathfrak{p})$ satisfying:
\begin{equation*}
 \cH = 0\, , \qquad \mathcal{W} = 0\, , \qquad \mathcal{M} = 0   
\end{equation*}

\noindent
then, there exists a four-dimensional Neveu-Schwarz solution having $(\Sigma,h)$ as its Cauchy hypersurface and restricting to $(h,\theta,b,\phi,\mathfrak{p})$ on $\Sigma$.


\subsection{The supersymmetric Neveu-Schwarz evolution flow}
\label{sec:susyflowCauchy}


In this subsection we consider the evolution flow defined by the globally hyperbolic supersymmetric configurations on a reducible triple $(\bar{\cP},\bar{Y},\bar{\cA})$ with reduction $(\cP,Y,\cA_t,\Psi_t)$. Our first result builds on the theory developed so far for skew-torsion parallelisms and supersymmetric Neveu-Schwarz solutions and characterizes globally hyperbolic supersymmetric solutions in terms of solutions of a constrained evolution system, as explained in the following. Throughout this subsection we set:
\begin{equation*}
\widehat{c}_t:=\lambda_t^{-1}c_t=\ast_{\fre_t}\theta_t\, ,  \qquad i=u,l,n,
\end{equation*}
and write $\widehat{c}_t(e_{ti}):=\iota_{e^{\sharp}_{ti}}\widehat{c}_t$.

\begin{prop}
\label{prop:susyevol}
A globally hyperbolic Neveu-Schwarz configuration $(g,\bar{b},\bar{\phi})\in \Conf(\bar{\cP},\bar{\cA},\bar{Y})$ is supersymmetric if and only if there exists a family of functions $\fra_t$ and a family of coframes $\fre_t$ on $\Sigma$ such that the following differential system is satisfied:
\begin{eqnarray}
& \partial_t \fra_t = \dd_{\Sigma}\lambda_t(e^{\sharp}_{tu})\, , \quad D_t e^t_i + \Theta_t(e^t_i)+\frac{1}{2}\widehat{c}_t(e^t_i)=\delta_{ui}\eta_t-\eta_t(e^{\sharp}_{ti})e_{tu}\, ,  \nonumber \\
& \dd_{\Sigma}e^t_i=\Theta_t(e^t_i)\wedge e_{tu}+\frac{1}{2}f_{b_t}\ast_{h_t}e^t_i-\frac{1}{2}\delta_{ui}\widehat{c}_t\, , \quad i=u,l,n \label{eq:susyevol} \\
& \dd_\Sigma\fra_t=-\Theta_t(e_{tu})+\tfrac12\widehat{c}_t (e_{tu})\, , \quad \varphi_{\phi_t}=p_t e_{tu}+\widehat{c}_t(e_{tu})\, , \quad f_{b_t}=\ast_{h_t}(\widehat{c}_t\wedge e_{tu})\nonumber
\end{eqnarray}

\noindent
where $(h_t,a_t,b_t,\phi_t,\lambda_t)$ is the reduction of $(g,\bar{b},\bar{\phi})$ and $p_t$ and $\varphi_{\phi_t}$ are defined by $\varphi_{\bar\phi}=p_t\mathfrak t_t+\varphi_{\phi_t}$.
\end{prop}

\begin{remark}
In the previous proposition, we are using the Kronecker delta $\delta_{ui}$ for the indices $i = u, l, n$. Hence, $\delta_{uu} = 1$ and $\delta_{ul} = \delta_{un} = 0$. 
\end{remark}

\begin{proof}
By Proposition \ref{prop:Thetaevolutionsystem} together with the discussion in Section \ref{sec:ReductionGloballyHyperbolicparallelisms}, it only remains to consider the equation $H_{\bar{b}}\cdot\varepsilon=\varphi_{\bar\phi}\cdot\varepsilon$. Since the spinorial square of the supersymmetry parameter is represented by $u+u\wedge l$, this equation is equivalent to \cite{Cortes:2019xmk}:
\begin{equation*}
 \varphi_{\bar\phi}\diamond_g(u+u\wedge l)
 =H_{\bar b}\diamond_g(u+u\wedge l)
\end{equation*}
which in turn is equivalent to:
\begin{equation*}
 \varphi_{\bar\phi}\wedge u=\ast_g(\alpha_{\bar b}\wedge u),
 \qquad
 \varphi_{\bar\phi}(u^{\sharp_g})=\alpha_{\bar b}(u^{\sharp_g})=0.
\end{equation*}   

\noindent
On the other hand, we have:
\begin{eqnarray*}
& u=e^{\fra_t}(\mathfrak t_t+e_{tu})\, , \qquad u^{\sharp_g}=e^{\fra_t}(-N_t + e^{\sharp}_{tu})\, ,\\
& \varphi_{\bar\phi}=p_t\mathfrak t_t+\varphi_{\phi_t}\, , \qquad \alpha_{\bar b}=-(f_{b_t}\mathfrak t_t+\theta_t)\, , \qquad \widehat{c}_t =\ast_{h_t}\theta_t.
\end{eqnarray*}

\noindent
as well as:
\begin{eqnarray*}
\varphi_{\bar\phi}(u^{\sharp_g})=e^{\fra_t}(\varphi_{\phi_t}(e^{\sharp}_{tu})-p_t)=0\, , \qquad
\alpha_{\bar b}(u^{\sharp_g})=e^{\fra_t}(f_{b_t}-\theta_t(e^{\sharp}_{tu}))=0
\end{eqnarray*}

\noindent
which are equivalent to:
\begin{equation*}
p_t=\varphi_{\phi_t}(e^{\sharp}_{tu})\, , \qquad f_{b_t}=\theta_t(e^{\sharp}_{tu}).
\end{equation*}

\noindent
Furthermore:
\begin{eqnarray*}
& \varphi_{\bar\phi}\wedge u-\ast_g(\alpha_{\bar b}\wedge u)\\
& =e^{\fra_t}\big(\mathfrak t_t\wedge(p_t e_{tu}-\varphi_{\phi_t}+\ast_{h_t}(\theta_t\wedge e_{tu}))+\varphi_{\phi_t}\wedge e_{tu}-f_{b_t}\ast_{h_t}e_{tu}+\ast_{h_t}\theta_t\big).
\end{eqnarray*}

\noindent
This implies:
\begin{equation*}
p_t e_{tu}-\varphi_{\phi_t}+\ast_{h_t}(\theta_t\wedge e_{tu})=0\, , \qquad
\varphi_{\phi_t}\wedge e_{tu}-f_{b_t}\ast_{h_t}e_{tu}+\ast_{h_t}\theta_t=0.
\end{equation*}

\noindent
or, equivalently:
\begin{equation*}
\varphi_{\phi_t}=p_t e_{tu}+\widehat{c}_t (e_{tu})\, , \qquad
\theta_t=f_{b_t}e_{tu}+\theta_t(e^{\sharp}_{tl})e_{tl}+\theta_t(e^{\sharp}_{tn})e_{tn}.
\end{equation*}

\noindent
Hence, there exists a family of one-forms $w_{b_t}\in\Gamma(e_{tu})^{\perp_{h_t}}$ orthogonal to $e_{tu}$ such that:
\begin{equation}
\label{eq:wbt}
\varphi_{\phi_t}=p_t e_{tu}-\ast_{q_t}w_{b_t}\, , \qquad \theta_t=f_{b_t}e_{tu}+w_{b_t}.
\end{equation}

\noindent
where $\ast_{q_t}$ is the Hodge dual determined on the oriented distribution $(e_{tu})^{\perp_{h_t}}\subset T^{\ast}\Sigma$ by the metric $q_t=e_{tl}\otimes e_{tl}+e_{tn}\otimes e_{tn}$ and the orientation $e_{tl}\wedge e_{tn}$. Similarly, we compute:
\begin{eqnarray*}
& \alpha_{\bar b}(e^{\sharp}_{tl})\ast_g u-\varphi_{\bar\phi}\wedge l\wedge u=e^{\fra_t}(\theta_t(e^{\sharp}_{tl})\nu_{h_t}-\varphi_{\phi_t}\wedge e_{tl}\wedge e_{tu})\\
& \qquad +e^{\fra_t}\mathfrak t_t\wedge(\theta_t(e^{\sharp}_{tl})e_{tl}\wedge e_{tn}-p_t e_{tl}\wedge e_{tu}-\varphi_{\phi_t}\wedge e_{tl})=0.
\end{eqnarray*}

\noindent
which is thus automatically satisfied by \eqref{eq:wbt}. Analogously, the remaining algebraic equation involving $l$ is automatically solved by \eqref{eq:wbt}. Furthermore, we notice that:
\begin{eqnarray*}
& \ast_{h_t}(\alpha_t^{\perp}-\alpha_t^o e_{tu})-\frac{1}{2}e_{tu}\wedge\ast_{h_t}(\alpha_t^{\perp}\wedge e_{tu})\\
& =f_{b_t}\ast_{h_t}e_{tu}-\widehat{c}_t +\frac{1}{2}e_{tu}\wedge\widehat{c}_t (e_{tu})
=\frac{1}{2}(f_{b_t}\ast_{h_t}e_{tu}-\widehat{c}_t ),
\end{eqnarray*}

\noindent
where $\alpha_t^o=-f_{b_t}$ and $\alpha_t^{\perp}=-\theta_t$. Substituting the previous expressions and relations in Proposition \ref{prop:Thetaevolutionsystem}, including the lapse terms, yields the differential system \eqref{eq:susyevol} stated in the proposition. These equations guarantee that the globally hyperbolic reduction of $[u,v,l,n]$ is skew-torsion with torsion given by the curvature of a curving on a reducible bundle gerbe. The converse follows by tracing back the previous steps and hence we conclude.
\end{proof}

\noindent
By the previous proposition, the family of functions $\fra_t\in C^{\infty}(\Sigma)$ that conforms the reduced isotropic parallelism, not to be confused with the family of one-forms $a_t\in\Omega^1(Y)$ associated to the globally hyperbolic reduction of a curving, satisfies $\partial_t\fra_t=\dd_{\Sigma}\lambda_t(e^{\sharp}_{tu})$. This family can be decoupled from the evolution system as the following result shows.

\begin{cor}
\label{cor:susyevolcor}
A globally hyperbolic Neveu-Schwarz configuration $(g,\bar{b},\bar{\phi})\in \Conf(\bar{\cP},\bar{\cA},\bar{Y})$ is supersymmetric if and only if there exists a family of coframes $\fre_t$ on $\Sigma$ such that the following differential system:
\begin{eqnarray}
& D_t e^t_i+\Theta_t(e^t_i)+\frac{1}{2}\widehat{c}_t (e^t_i)=\delta_{ui}\eta_t-\eta_t(e^{\sharp}_{ti})e_{tu}\, , \quad i=u,l,n\nonumber \\
& \dd_{\Sigma}e^t_i=\Theta_t(e^t_i)\wedge e_{tu}+\frac{1}{2}f_{b_t}\ast_{h_t}e^t_i-\frac{1}{2}\delta_{ui}\widehat{c}_t \, , \quad i=u,l,n\label{eq:susyevolcor} \\
& \varphi_{\phi_t}=p_t e_{tu}+\widehat{c}_t (e_{tu})\, , \quad f_{b_t}=\ast_{h_t}(\widehat{c}_t \wedge e_{tu})\nonumber
\end{eqnarray}

\noindent
together with the following conditions:
\begin{eqnarray*}
0=[\Theta_t(e_{tu})-\frac{1}{2}\widehat{c}_t (e_{tu})]\in H^1(\Sigma,\mathbb R)\, , \qquad
\partial_t(\Theta_t(e_{tu})-\frac{1}{2}\widehat{c}_t (e_{tu}))=-\dd_{\Sigma}(\dd_{\Sigma}\lambda_t(e^{\sharp}_{tu})).
\end{eqnarray*}

\noindent
are satisfied, where $(h_t,a_t,b_t,\phi_t,\lambda_t)$ is the reduction of $(g,\bar b,\bar\phi)$.
\end{cor}

\noindent
In practical situations, we are sometimes interested in the existence of globally hyperbolic configurations without necessarily fixing the underlying reducible bundle gerbe and given in terms of the explicit tuple $(\lambda_t,\fre_t,\widehat{c}_t ,f_t,p_t,\varphi_{t})$, where $\varphi_{t}$ is a family of closed one-forms.

\begin{prop}
\label{prop:susyevolII}
A four-manifold of the form $M=\cI\times\Sigma$ admits a supersymmetric globally hyperbolic Neveu-Schwarz configuration with Cauchy hypersurface $\Sigma$ if and only if there exists a tuple $(\lambda_t,\fre_t,\widehat{c}_t ,f_t,p_t,\varphi_{\phi_t})$ on $\Sigma$, with $\lambda_t\neq0$, satisfying the following differential system:
\begin{eqnarray*}
& \partial_t(\Theta_t(e_{tu})-\frac{1}{2}\widehat{c}_t (e_{tu}))=-\dd_{\Sigma}(\dd_{\Sigma}\lambda_t(e^{\sharp}_{tu}))\, ,\\
& D_t e^t_i+\Theta_t(e^t_i)+\frac{1}{2}\widehat{c}_t (e^t_i)=\delta_{ui}\eta_t-\eta_t(e^{\sharp}_{ti})e_{tu}\, , \quad i=u,l,n\nonumber \\
& f_t=\ast_{h_t}(\widehat{c}_t \wedge e_{tu})\, , \quad
\dd_{\Sigma}e^t_i=\Theta_t(e^t_i)\wedge e_{tu}+\frac{1}{2}f_t\ast_{h_t}e^t_i-\frac{1}{2}\delta_{ui}\widehat{c}_t \, , \quad i=u,l,n\\
& D_t(f_t\nu_{h_t})=\dd_{\Sigma}\widehat{c}_t +\eta_t\wedge\widehat{c}_t \, , \qquad \dd_{\Sigma}\varphi_{\phi_t}=0\, ,\\
& D_t\varphi_{\phi_t}=\dd_{\Sigma}p_t+p_t\eta_t\, , \qquad \varphi_{\phi_t}=p_t e_{tu}+\widehat{c}_t (e_{tu}).
\end{eqnarray*}
	
\noindent
together with the following cohomological conditions:
\begin{eqnarray*}
\frac{1}{2\pi}[\varphi_{\phi_t}]\in H^1_{\mathbb Z}(\Sigma,\mathbb R)\, , \quad
\left[\frac{f_t}{2\pi}\nu_{h_t}\right]\in H^3_{\mathbb Z}(\Sigma,\mathbb R)\, , \quad
0=[\Theta_t(e_{tu})-\frac{1}{2}\widehat{c}_t (e_{tu})]\in H^1(\Sigma,\mathbb R).
\end{eqnarray*}

\noindent
for every $t\in \cI$. 
\end{prop}

\noindent
We will refer to the differential system \eqref{eq:susyevol} as occurring in Proposition \ref{prop:susyevol}, or equivalently Corollary \ref{cor:susyevolcor}, as the \emph{Neveu-Schwarz supersymmetry flow}, which is the evolution flow determined by the supersymmetry conditions of Neveu-Schwarz supergravity for families of the form $(\lambda_t,\fre_t,a_t,b_t,\phi_t)$. Given such a family, we obtain a canonical candidate for the Neveu-Schwarz evolution flow given simply by $(h_{\fre_t},a_t,b_t,\phi_t,\lambda_t)$, where $h_{\fre_t}$ is the family of Riemannian metrics on $\Sigma$ given by $h_{\fre_t}=e_{tu}\otimes e_{tu}+e_{tl}\otimes e_{tl}+e_{tn}\otimes e_{tn}$.  
\begin{remark}
\label{remark:covariantevolution}
From equations \eqref{eq:susyevol} or \eqref{eq:susyevolcor} it immediately follows that the covariant derivative of the family of coframes $\fre_t$ belonging to a Neveu-Schwarz supersymmetry flow $(\lambda_t,\fre_t,a_t,b_t,\phi_t)$ is given by:
\begin{eqnarray*}
\nabla^{h_t}e^t_i+\delta_{ui}(\Theta_t+\frac{1}{4}\widehat{c}_t )=\Theta_t(e^t_i)\otimes e_{tu}+\frac{1}{4}f_{b_t}\ast_{h_t}e^t_i-\frac{1}{4}e_{tu}\odot\widehat{c}_t (e^t_i)
\end{eqnarray*}

\noindent
for $i = u, l, n$. The previous covariant derivatives can be equivalently written as follows: 
\begin{eqnarray*}
& \nabla^{h_t}e_{tu}+\Theta_t=\Theta_t(e_{tu})\otimes e_{tu}-\frac{1}{2}e_{tu}\otimes\widehat{c}_t (e_{tu}) \\
& \nabla^{h_t}e_{tl}=\Theta_t(e_{tl})\otimes e_{tu}-\frac{1}{2}f_{b_t}e_{tu}\otimes e_{tn}+\frac{1}{2}\widehat{c}_t (e_{tu},e_{tl})e_{tu}\otimes e_{tu}\\
& \nabla^{h_t}e_{tn}=\Theta_t(e_{tn})\otimes e_{tu}+\frac{1}{2}f_{b_t}e_{tu}\otimes e_{tl}+\frac{1}{2}\widehat{c}_t (e_{tu},e_{tn})e_{tu}\otimes e_{tu}
\end{eqnarray*}
	
\noindent
From this we immediately obtain:
\begin{eqnarray*}
\nabla^{h_t\ast}e^t_i=-\Tr_{h_t}(\nabla^{h_t}e^t_i)=\delta_{ui}\Tr_{h_t}(\Theta_t)-\Theta_t(e_{tu},e^t_i)-\frac{1}{2}\widehat{c}_t (e_{tu},e^t_i)
\end{eqnarray*}

\noindent
or, equivalently:
\begin{eqnarray*}
& \nabla^{h_t\ast}e_{tu}=\Tr_{h_t}(\Theta_t)-\Theta_t(e_{tu},e_{tu})\\
& \nabla^{h_t\ast}e_{tl}=-\Theta_t(e_{tl},e_{tu})-\frac{1}{2}\widehat{c}_t (e_{tu},e_{tl})\, , \qquad
\nabla^{h_t\ast}e_{tn}=-\Theta_t(e_{tn},e_{tu})-\frac{1}{2}\widehat{c}_t (e_{tu},e_{tn})
\end{eqnarray*}
	
\noindent
These equations will be useful in the following.
\end{remark} 

\noindent
It is natural to compare the evolution flow prescribed by the Neveu-Schwarz supersymmetry conditions for $(\lambda_t,\fre_t,a_t,b_t,\phi_t)$ with the Neveu-Schwarz evolution flow for the associated family $(h_{\fre_t},a_t,b_t,\phi_t,\lambda_t)$. For this, we first need to compute the scalar curvature of $h_{\fre_t}$, which we denote by $s_{\fre_t}$.

\begin{lemma}
\label{lemma:curvaturesusyflow}
Let $(\lambda_t,\fre_t, a_t , b_t,\phi_t)$ be a Neveu-Schwarz supersymmetry flow. Then:
\begin{equation*}
s_{\fre_t}=\vert\Theta_t\vert^2_{h_t}-(\Tr_{h_t}\Theta_t)^2+2(\nabla^{h_t\ast}\Theta_t)(e^{\sharp}_{tu})+2\Tr_{h_t}(\nabla^{h_t}_{e^{\sharp}_{tu}}\Theta_t)-(\nabla^{h_t\ast}\widehat{c}_t )(e^{\sharp}_{tu})-\frac{1}{2}\vert\widehat{c}_t (e_{tu})\vert_{h_t}^2
\end{equation*}

\noindent
where $s_{\fre_t}$ denotes the scalar curvature of $h_{\fre_t}$ at a given time $t\in \cI$. 
\end{lemma} 

\begin{proof}
Using equations \eqref{eq:susyevol} together with Remark \ref{remark:covariantevolution}, we compute:
\begin{eqnarray*}
&\nabla^{h_t}_{e^{\sharp}_{ti}}\nabla^{h_t}_{e_{tj}^{\sharp}}e_{tu}=\nabla^{h_t}_{e^{\sharp}_{ti}}(-\Theta_t(e_{tj})+\Theta_t(e_{tu},e_{tj})e_{tu}-\frac{1}{2}\delta_{uj}\widehat{c}_t (e_{tu}))\\
&=-(\nabla^{h_t}_{e^{\sharp}_{ti}}\Theta_t)(e_{tj})-\Theta_t(\nabla^{h_t}_{e^{\sharp}_{ti}}e_{tj})+((\nabla^{h_t}_{e^{\sharp}_{ti}}\Theta_t)(e_{tu},e_{tj})+\Theta_t(\nabla^{h_t}_{e^{\sharp}_{ti}}e_{tu},e_{tj})+\Theta_t(e_{tu},\nabla^{h_t}_{e^{\sharp}_{ti}}e_{tj}))e_{tu}\\
&\quad+\Theta_t(e_{tu},e_{tj})\nabla^{h_t}_{e^{\sharp}_{ti}}e_{tu}-\frac{1}{2}\delta_{uj}((\nabla^{h_t}_{e^{\sharp}_{ti}}\widehat{c}_t )(e_{tu})+\widehat{c}_t (\nabla^{h_t}_{e^{\sharp}_{ti}}e_{tu}))
\end{eqnarray*}

\noindent
On the other hand:
\begin{eqnarray*}
\nabla^{h_t}_{[e^{\sharp}_{ti},e_{tj}^{\sharp}]}e_{tu}=-\Theta_t([e^{\sharp}_{ti},e_{tj}^{\sharp}])+\Theta_t(e^{\sharp}_{tu},[e^{\sharp}_{ti},e_{tj}^{\sharp}])e_{tu}-\frac{1}{2}e_{tu}([e^{\sharp}_{ti},e_{tj}^{\sharp}])\widehat{c}_t (e_{tu})
\end{eqnarray*}

\noindent
From these formulae we obtain, after a tedious calculation, the following expression:
\begin{eqnarray*}
& (\cR^{\fre_t}_{e^{\sharp}_{ti} e^{\sharp}_{tu}}e_{tu})(e^{\sharp}_{ti})=\vert\Theta_t(e_{tu})\vert^2_{h_t}-\Theta_t(e_{tu},e_{tu})\Tr_{h_t}(\Theta_t)+(\nabla^{h_t\ast}\Theta_t)(e^{\sharp}_{tu})+\Tr_{h_t}(\nabla^{h_t}_{e^{\sharp}_{tu}}\Theta_t)\\
&\quad-\frac{1}{2}(\nabla^{h_t\ast}\widehat{c}_t )(e^{\sharp}_{tu})+\frac{1}{2}\widehat{c}_t (\Theta_t(e_{tl}),e_{tl})+\frac{1}{2}\widehat{c}_t (\Theta_t(e_{tn}),e_{tn})-\frac{1}{2}\langle\Theta_t(e_{tu}),\widehat{c}_t (e_{tu})\rangle_{h_t}-\frac{1}{4}\vert\widehat{c}_t (e_{tu})\vert_{h_t}^2
\end{eqnarray*}

\noindent
where summation over repeated indices is assumed. Regarding $e_{tl}$ we have:
\begin{eqnarray*}
& \nabla^{h_t}_{e^{\sharp}_{ti}}\nabla^{h_t}_{e_{tj}^{\sharp}}e_{tl}=\nabla^{h_t}_{e^{\sharp}_{ti}}(\Theta_t(e_{tl},e_{tj})e_{tu}-\frac{1}{2}f_{b_t}\delta_{uj}e_{tn}+\frac{1}{2}\delta_{uj}\widehat{c}_t (e_{tu},e_{tl})e_{tu})\\
&=((\nabla^{h_t}_{e^{\sharp}_{ti}}\Theta_t)(e_{tl},e_{tj})+\Theta_t(\nabla^{h_t}_{e^{\sharp}_{ti}}e_{tl},e_{tj})+\Theta_t(e_{tl},\nabla^{h_t}_{e^{\sharp}_{ti}}e_{tj}))e_{tu}+\Theta_t(e_{tl},e_{tj})\nabla^{h_t}_{e^{\sharp}_{ti}}e_{tu}\\
&\quad+\frac{1}{2}\delta_{uj}((\nabla^{h_t}_{e^{\sharp}_{ti}}\widehat{c}_t )(e_{tu},e_{tl})+\widehat{c}_t (\nabla^{h_t}_{e^{\sharp}_{ti}}e_{tu},e_{tl})+\widehat{c}_t (e_{tu},\nabla^{h_t}_{e^{\sharp}_{ti}}e_{tl}))e_{tu}\\
&\quad+\frac{1}{2}\delta_{uj}(\widehat{c}_t (e_{tu},e_{tl})\nabla^{h_t}_{e^{\sharp}_{ti}}e_{tu}-\dd_{\Sigma}f_{b_t}(e^{\sharp}_{ti})e_{tn}-f_{b_t}\nabla^{h_t}_{e^{\sharp}_{ti}}e_{tn})
\end{eqnarray*}
 
\noindent
as well as:
\begin{eqnarray*}
\nabla^{h_t}_{[e^{\sharp}_{ti},e_{tj}^{\sharp}]}e_{tl}=\Theta_t(e^{\sharp}_{tl},[e^{\sharp}_{ti},e_{tj}^{\sharp}])e_{tu}-\frac{1}{2}f_{b_t}e_{tu}([e^{\sharp}_{ti},e_{tj}^{\sharp}])e_{tn}+\frac{1}{2}\widehat{c}_t (e_{tu},e_{tl})e_{tu}([e^{\sharp}_{ti},e_{tj}^{\sharp}])e_{tu}
\end{eqnarray*}

\noindent
from which we obtain:
\begin{eqnarray*}
& (\cR^{\fre_t}_{e^{\sharp}_{ti} e^{\sharp}_{tl}}e_{tl})(e^{\sharp}_{ti})=(\nabla^{h_t}_{e^{\sharp}_{tu}}\Theta_t)(e_{tl},e_{tl})-(\nabla^{h_t}_{e^{\sharp}_{tl}}\Theta_t)(e_{tl},e_{tu})-\Theta_t(e_{tl},e_{tl})\Tr_{h_t}\Theta_t+\vert\Theta_t(e_{tl})\vert^2_{h_t}\\
&\quad-\frac{1}{2}(\widehat{c}_t (e_{tu},e_{tl})\Theta_t(e_{tu},e_{tl})+(\nabla^{h_t}_{e^{\sharp}_{tl}}\widehat{c}_t )(e_{tu},e_{tl})+\langle\Theta_t(e_{tl}),\widehat{c}_t (e_{tl})\rangle_{h_t})-\frac{1}{4}\widehat{c}_t (e_{tu},e_{tl})^2
\end{eqnarray*}

\noindent
We compute similarly for $e_{tn}$. We have:
\begin{eqnarray*}
& \nabla^{h_t}_{e^{\sharp}_{ti}}\nabla^{h_t}_{e_{tj}^{\sharp}}e_{tn}=\nabla^{h_t}_{e^{\sharp}_{ti}}(\Theta_t(e_{tn},e_{tj})e_{tu}+\frac{1}{2}\delta_{uj}(f_{b_t}e_{tl}+\widehat{c}_t (e_{tu},e_{tn})e_{tu}))=\\
&=((\nabla^{h_t}_{e^{\sharp}_{ti}}\Theta_t)(e_{tn},e_{tj})+\Theta_t(\nabla^{h_t}_{e^{\sharp}_{ti}}e_{tn},e_{tj})+\Theta_t(e_{tn},\nabla^{h_t}_{e^{\sharp}_{ti}}e_{tj}))e_{tu}+\Theta_t(e_{tn},e_{tj})\nabla^{h_t}_{e^{\sharp}_{ti}}e_{tu}\\
&\quad+\frac{1}{2}\delta_{uj}((\nabla^{h_t}_{e^{\sharp}_{ti}}\widehat{c}_t )(e_{tu},e_{tn})+\widehat{c}_t (\nabla^{h_t}_{e^{\sharp}_{ti}}e_{tu},e_{tn})+\widehat{c}_t (e_{tu},\nabla^{h_t}_{e^{\sharp}_{ti}}e_{tn}))e_{tu}\\
&\quad+\frac{1}{2}\delta_{uj}(\widehat{c}_t (e_{tu},e_{tn})\nabla^{h_t}_{e^{\sharp}_{ti}}e_{tu}+\dd_{\Sigma}f_{b_t}(e^{\sharp}_{ti})e_{tl}+f_{b_t}\nabla^{h_t}_{e^{\sharp}_{ti}}e_{tl})
\end{eqnarray*}

\noindent
as well as:
\begin{eqnarray*}
\nabla^{h_t}_{[e^{\sharp}_{ti},e_{tj}^{\sharp}]}e_{tn}=\Theta_t(e^{\sharp}_{tn},[e^{\sharp}_{ti},e_{tj}^{\sharp}])e_{tu}+\frac{1}{2}f_{b_t}e_{tu}([e^{\sharp}_{ti},e_{tj}^{\sharp}])e_{tl}+\frac{1}{2}\widehat{c}_t (e_{tu},e_{tn})e_{tu}([e^{\sharp}_{ti},e_{tj}^{\sharp}])e_{tu}
\end{eqnarray*}

\noindent
from which we obtain:
\begin{eqnarray*}
& (\cR^{\fre_t}_{e^{\sharp}_{ti} e^{\sharp}_{tn}}e_{tn})(e^{\sharp}_{ti})=(\nabla^{h_t}_{e^{\sharp}_{tu}}\Theta_t)(e_{tn},e_{tn})-(\nabla^{h_t}_{e^{\sharp}_{tn}}\Theta_t)(e_{tn},e_{tu})-\Theta_t(e_{tn},e_{tn})\Tr_{h_t}\Theta_t+\vert\Theta_t(e_{tn})\vert^2_{h_t}\\
&\quad-\frac{1}{2}(\widehat{c}_t (e_{tu},e_{tn})\Theta_t(e_{tu},e_{tn})+(\nabla^{h_t}_{e^{\sharp}_{tn}}\widehat{c}_t )(e_{tu},e_{tn})+\langle\Theta_t(e_{tn}),\widehat{c}_t (e_{tn})\rangle_{h_t})-\frac{1}{4}\widehat{c}_t (e_{tu},e_{tn})^2
\end{eqnarray*}

\noindent
The previous formulae yield the scalar curvature $s_{\fre_t}$ of the family of metrics $h_t$ as follows:
\begin{eqnarray*}
s_{\fre_t}=(\cR^{\fre_t}_{e^{\sharp}_{ti} e^{\sharp}_{tu}}e_{tu})(e^{\sharp}_{ti})+(\cR^{\fre_t}_{e^{\sharp}_{ti} e^{\sharp}_{tl}}e_{tl})(e^{\sharp}_{ti})+(\cR^{\fre_t}_{e^{\sharp}_{ti} e^{\sharp}_{tn}}e_{tn})(e^{\sharp}_{ti})
\end{eqnarray*}

\noindent
which, after some manipulations, gives the expression in the statement of the lemma.
\end{proof}

\noindent
We proceed to prove a general \emph{compatibility result} between the Neveu-Schwarz evolution flow and its supersymmetric counterpart.  

\begin{thm}
\label{thm:compatibilityflows}
Let $(\lambda_t,\fre_t, a_t , b_t,\phi_t)$ be a supersymmetric Neveu-Schwarz flow, that is, a solution to the differential system \eqref{eq:susyevol}, and let $\cH_t$, $\cW_t$ and $\cM_t$ be the time-dependent constraints \eqref{eq:normalconstraint} of the Neveu-Schwarz evolution flow evaluated on $(h_{\fre_t}, a_t , b_t,\phi_t,\lambda_t)$. Then, for every $t\in \cI$:
\begin{equation}
\label{eq:constraintidentities0}
\cH_t e_{tu} = 2\, \cW_t  + \ast_{h_t} \cM_t  
\end{equation}

\noindent
In particular, $(h_{\fre_t}, a_t , b_t,\phi_t,\lambda_t)$ satisfies the constraint equations \eqref{eq:normalconstraint} of the Neveu-Schwarz evolution flow if and only if $\cW_t(e^{\sharp}_{tu}) = 0$ and $\cM_t = 0$.
\end{thm}

\begin{proof}
Note that Equation \eqref{eq:constraintidentities0} is equivalent to the following equations:
\begin{equation}
\label{eq:constraintidentities}
\cH_t = 2\, \cW_t(e^{\sharp}_{tu}) + \ast_{h_t} (\cM_t\wedge e_{tu})\, , \qquad \cW_t - \cW_t(e^{\sharp}_{tu})\, e_{tu} =  \frac{1}{2} ( \ast_{h_t} (\cM_t\wedge e_{tu})\, e_{tu} - \ast_{h_t}\cM_t )
\end{equation}

\noindent
which we will obtain in the proof. We fix $t\in \cI$ and write $h_t = h_{\fre_t}$. The identities \eqref{eq:constraintidentities} only involve the constraint equations in \eqref{eq:susyevol}, namely the equations that prescribe $\dd_{\Sigma} e^t_i$ and $\dd_{\Sigma}\fra_t$ together with $\varphi_{\phi_t}=p_t e_{tu}+\widehat{c}_t(e_{tu})$ and $f_{b_t}=\ast_{h_t}(\widehat{c}_t\wedge e_{tu})$. Since $\widehat{c}_t (e_{tu})$ is orthogonal to $e_{tu}$, the last of these equations gives the decomposition:
\begin{equation*}
\widehat{c}_t  = f_{b_t} \ast_{h_t} e_{tu} + e_{tu}\wedge \widehat{c}_t (e_{tu})    
\end{equation*}

\noindent
whence $\theta_t = \ast_{h_t}\widehat{c}_t  = f_{b_t} e_{tu} + \ast_{h_t}(e_{tu}\wedge\widehat{c}_t (e_{tu}))$ and:
\begin{eqnarray}
\label{eq:dilatinoalgebra}
& \theta_t(e^{\sharp}_{tu}) = f_{b_t}\, , \qquad \vert\theta_t\vert^2_{h_t} = f_{b_t}^2 + \vert \widehat{c}_t (e_{tu})\vert^2_{h_t}\, , \qquad \vert\varphi_{\phi_t}\vert^2_{h_t} = p_t^2 + \vert \widehat{c}_t (e_{tu})\vert^2_{h_t} \nonumber \\
& \ast_{h_t}(\varphi_{\phi_t}\wedge\theta_t\wedge e_{tu}) = \vert \widehat{c}_t (e_{tu})\vert^2_{h_t}
\end{eqnarray}

\noindent
Since $\nabla^{h_t\ast}\widehat{c}_t  = \ast_{h_t}\dd_{\Sigma}\theta_t$ on a three-manifold, we obtain:
\begin{equation}
\label{eq:Maxwellnormalcomponent}
\ast_{h_t}(\cM_t\wedge e_{tu}) = (\ast_{h_t}\cM_t)(e^{\sharp}_{tu}) = (\nabla^{h_t\ast}\widehat{c}_t )(e^{\sharp}_{tu}) - \vert \widehat{c}_t (e_{tu})\vert^2_{h_t}
\end{equation}

\noindent
We prove the first identity in \eqref{eq:constraintidentities}. By Lemma \ref{lemma:curvaturesusyflow} and Equation \eqref{eq:dilatinoalgebra}, the Hamiltonian constraint of \eqref{eq:normalconstraint} evaluated on $(h_{\fre_t}, a_t , b_t,\phi_t,\lambda_t)$ reads:
\begin{equation*}
\cH_t = 2(\nabla^{h_t\ast}\Theta_t)(e^{\sharp}_{tu}) + 2\Tr_{h_t}(\nabla^{h_t}_{e^{\sharp}_{tu}}\Theta_t) - (\nabla^{h_t\ast}\widehat{c}_t )(e^{\sharp}_{tu}) - 2\vert \widehat{c}_t (e_{tu})\vert^2_{h_t} + 2 p_t \Tr_{h_t}(\Theta_t) - f_{b_t}^2 - 2\nabla^{h_t\ast}\varphi_{\phi_t}
\end{equation*}

\noindent
Here $\Tr_{h_t}(\nabla^{h_t}_{e^{\sharp}_{tu}}\Theta_t) = \dd_{\Sigma}\Tr_{h_t}(\Theta_t)(e^{\sharp}_{tu})$, whereas by Remark \ref{remark:covariantevolution} we have:
\begin{equation*}
\nabla^{h_t\ast}(p_t e_{tu}) = - \dd_{\Sigma} p_t (e^{\sharp}_{tu}) + p_t (\Tr_{h_t}(\Theta_t) - \Theta_t(e_{tu},e_{tu}))    
\end{equation*}

\noindent
and, since:
\begin{equation*}
\nabla^{h_t\ast}(\widehat{c}_t (e_{tu})) = -(\nabla^{h_t\ast}\widehat{c}_t )(e^{\sharp}_{tu}) - \sum_{i} \widehat{c}_t  (\nabla^{h_t}_{e^{\sharp}_{ti}} e^{\sharp}_{tu}, e^{\sharp}_{ti})    
\end{equation*}

\noindent
for every two-form $\widehat{c}_t $, the formula for $\nabla^{h_t} e_{tu}$ in Remark \ref{remark:covariantevolution} gives:
\begin{equation*}
\nabla^{h_t\ast}(\widehat{c}_t (e_{tu})) = -(\nabla^{h_t\ast}\widehat{c}_t )(e^{\sharp}_{tu}) - \Theta_t (\widehat{c}_t (e_{tu})^{\sharp}, e^{\sharp}_{tu}) - \frac{1}{2}\vert \widehat{c}_t (e_{tu})\vert^2_{h_t}
\end{equation*}

\noindent
Substituting these three identities in the previous expression for $\cH_t$ and using $\theta_t(e^{\sharp}_{tu}) = f_{b_t}$ together with \eqref{eq:Maxwellnormalcomponent}, we obtain:
\begin{eqnarray*}
& \cH_t = 2 ( \dd_{\Sigma}(\Tr_{h_t}(\Theta_t) + p_t)(e^{\sharp}_{tu}) + (\nabla^{h_t\ast}\Theta_t)(e^{\sharp}_{tu}) + \Theta_t(\varphi_{\phi_t}^{\sharp}, e^{\sharp}_{tu}) - \frac{1}{2} f_{b_t}\theta_t(e^{\sharp}_{tu})) \\
& \qquad + (\nabla^{h_t\ast}\widehat{c}_t )(e^{\sharp}_{tu}) - \vert \widehat{c}_t (e_{tu})\vert^2_{h_t} = 2\, \cW_t(e^{\sharp}_{tu}) + \ast_{h_t} (\cM_t\wedge e_{tu})
\end{eqnarray*}

\noindent
For the second identity in \eqref{eq:constraintidentities}, it suffices to prove $\cW_t(e^{\sharp}_{tl}) = -\frac{1}{2}(\ast_{h_t}\cM_t)(e^{\sharp}_{tl})$, since the equations \eqref{eq:susyevol} are invariant under the substitution $(e_{tl},e_{tn}) \mapsto (e_{tn},-e_{tl})$, which transforms this equation into the corresponding equation for $e_{tn}$. We compute $\Ric^{h_t}(e^{\sharp}_{tl},e^{\sharp}_{tu})$ in two ways, using the identity:
\begin{equation*}
\sum_i (\nabla^{h_t}_{e^{\sharp}_{ti}}\nabla^{h_t}\beta)(Y,e^{\sharp}_{ti}) - \dd_{\Sigma}(\Tr_{h_t}(\nabla^{h_t}\beta))(Y) = \Ric^{h_t}(Y,\beta^{\sharp})    
\end{equation*}

\noindent
which is valid for every one-form $\beta$ and vector field $Y$ on $\Sigma$. Applied to $\beta = e_{tu}$ and $Y = e^{\sharp}_{tl}$, the formula for $\nabla^{h_t}e_{tu}$ in Remark \ref{remark:covariantevolution} gives:
\begin{eqnarray*}
& \Ric^{h_t}(e^{\sharp}_{tl},e^{\sharp}_{tu}) = (\nabla^{h_t\ast}\Theta_t)(e^{\sharp}_{tl}) + \dd_{\Sigma}\Tr_{h_t}(\Theta_t)(e^{\sharp}_{tl}) + (\nabla^{h_t}_{e^{\sharp}_{tu}}\Theta_t)(e_{tu},e_{tl}) - (\nabla^{h_t}_{e^{\sharp}_{tl}}\Theta_t)(e_{tu},e_{tu}) \\
& \qquad + \langle \Theta_t(e_{tu}),\Theta_t(e_{tl})\rangle_{h_t} - \Tr_{h_t}(\Theta_t)\Theta_t(e_{tu},e_{tl})
\end{eqnarray*}

\noindent
whereas, applied to $\beta = e_{tl}$ and $Y = e^{\sharp}_{tu}$, the formula for $\nabla^{h_t}e_{tl}$ in Remark \ref{remark:covariantevolution} gives:
\begin{eqnarray*}
& \Ric^{h_t}(e^{\sharp}_{tu},e^{\sharp}_{tl}) = \langle \Theta_t(e_{tu}),\Theta_t(e_{tl})\rangle_{h_t} - \Tr_{h_t}(\Theta_t)\Theta_t(e_{tu},e_{tl}) + \frac{1}{2}\Theta_t(\widehat{c}_t (e_{tu})^{\sharp},e^{\sharp}_{tl}) - \frac{1}{2}\dd_{\Sigma} f_{b_t}(e^{\sharp}_{tn}) \\
& \qquad - \frac{1}{2} f_{b_t}\Theta_t(e_{tn},e_{tu}) - \frac{1}{4} f_{b_t}\widehat{c}_t (e_{tu},e_{tn}) - \frac{1}{2}\widehat{c}_t (e_{tu},e_{tl}) (\Tr_{h_t}(\Theta_t) - \Theta_t(e_{tu},e_{tu}) )
\end{eqnarray*}

\noindent
On the other hand, the equation for $\dd_{\Sigma}\fra_t$ in \eqref{eq:susyevol} gives $\dd_{\Sigma}(\Theta_t(e_{tu})) = \frac{1}{2}\dd_{\Sigma}(\widehat{c}_t (e_{tu}))$, and expanding $\dd_{\Sigma}(\Theta_t(e_{tu}))$ in terms of $\nabla^{h_t}\Theta_t$ and $\nabla^{h_t}e_{tu}$ we obtain:
\begin{equation*}
(\nabla^{h_t}_{e^{\sharp}_{tu}}\Theta_t)(e_{tu},e_{tl}) - (\nabla^{h_t}_{e^{\sharp}_{tl}}\Theta_t)(e_{tu},e_{tu}) = \frac{1}{2}\dd_{\Sigma}(\widehat{c}_t (e_{tu}))(e^{\sharp}_{tu},e^{\sharp}_{tl}) + \frac{1}{2}\Theta_t(\widehat{c}_t (e_{tu})^{\sharp},e^{\sharp}_{tl})
\end{equation*}

\noindent
Furthermore, since $\varphi_{\phi_t}$ is closed, $\dd_{\Sigma}(\widehat{c}_t (e_{tu})) = - \dd_{\Sigma}(p_t e_{tu})$, and evaluating both sides on $(e^{\sharp}_{tu},e^{\sharp}_{tl})$ by means of the equations for $\dd_{\Sigma} e^t_i$ in \eqref{eq:susyevol} we obtain:
\begin{eqnarray*}
& \dd_{\Sigma}(\widehat{c}_t (e_{tu}))(e^{\sharp}_{tu},e^{\sharp}_{tl}) = \dd_{\Sigma}p_t(e^{\sharp}_{tl}) + p_t\Theta_t(e_{tu},e_{tl}) + \frac{1}{2}p_t \widehat{c}_t (e_{tu},e_{tl}) \\
& \qquad = \dd_{\Sigma}(\widehat{c}_t (e_{tu},e_{tl}))(e^{\sharp}_{tu}) - \widehat{c}_t (e_{tu},e_{tl})\Theta_t(e_{tl},e_{tl}) - \widehat{c}_t (e_{tu},e_{tn})\Theta_t(e_{tl},e_{tn}) + \frac{1}{2} f_{b_t}\widehat{c}_t (e_{tu},e_{tn})
\end{eqnarray*}

\noindent
Finally, evaluating $\cM_t = \dd_{\Sigma}\theta_t - \varphi_{\phi_t}\wedge\theta_t$ on $(e^{\sharp}_{tn},e^{\sharp}_{tu})$, again by means of the equations for $\dd_{\Sigma} e^t_i$ in \eqref{eq:susyevol}, we obtain:
\begin{eqnarray*}
& (\ast_{h_t}\cM_t)(e^{\sharp}_{tl}) = \dd_{\Sigma} f_{b_t}(e^{\sharp}_{tn}) + f_{b_t}\Theta_t(e_{tu},e_{tn}) - \dd_{\Sigma}(\widehat{c}_t (e_{tu},e_{tl}))(e^{\sharp}_{tu}) + \widehat{c}_t (e_{tu},e_{tl})\Theta_t(e_{tn},e_{tn}) \\
& \qquad - \widehat{c}_t (e_{tu},e_{tn})\Theta_t(e_{tl},e_{tn}) - f_{b_t}\widehat{c}_t (e_{tu},e_{tn}) + p_t\widehat{c}_t (e_{tu},e_{tl})
\end{eqnarray*}

\noindent
Equating the two expressions for $\Ric^{h_t}(e^{\sharp}_{tl},e^{\sharp}_{tu})$, substituting the previous three identities and using $\Theta_t(\varphi_{\phi_t}^{\sharp},e^{\sharp}_{tl}) = p_t\Theta_t(e_{tu},e_{tl}) + \Theta_t(\widehat{c}_t (e_{tu})^{\sharp},e^{\sharp}_{tl})$ and $\theta_t(e^{\sharp}_{tl}) = -\widehat{c}_t (e_{tu},e_{tn})$, a direct computation gives $\cW_t(e^{\sharp}_{tl}) = -\frac{1}{2}(\ast_{h_t}\cM_t)(e^{\sharp}_{tl})$, which proves the second identity in \eqref{eq:constraintidentities}. 

\noindent
The last assertion follows from \eqref{eq:constraintidentities}: if $\cM_t = 0$ and $\cW_t(e^{\sharp}_{tu}) = 0$, then $\cW_t = 0$ by the second identity and $\cH_t = 0$ by the first.
\end{proof}


\subsection{Supersymmetric initial data and their development}
\label{sec:susyinitialdevelopment}


The aim of this subsection is to solve the Cauchy problem for the first-order supersymmetry equations \eqref{eq:KSE4d} by means of the equivalent differential system given in \eqref{eq:susyevol} for tuples of the form $(\lambda_t,\fre_t,a_t,b_t,\phi_t)$.   

Section \ref{sec:skewtorsionwellposedness} uses the function $\fra_t$ defined by $u=e^{\fra_t}(\mathfrak t_t+e_{tu})$, whereas
Theorem~\ref{thm:wellposedFinalSkew} uses the \emph{integrated} function $\frf_t$. The corresponding spatial identities give
$\dd_\Sigma(\frf_t+\fra_t)=0$.  Since $\Sigma$ is connected, their sum is spatially constant, and the two time-evolution equations show that this constant is independent of $t$.  We fix the remaining additive normalization so that:
\begin{equation}
\label{eq:NSNS-scale-identification}
 \frf_t=-\fra_t.
\end{equation}
Thus $\partial_t\fra_t=\dd_\Sigma\lambda_t(e^{\sharp}_{tu})$ becomes $\partial_t\frf_t + \dd_\Sigma \lambda_t(e^{\sharp}_{tu}) =0$, exactly as in Lemma~\ref{lemma:refor2}. Combining equations \eqref{eq:susyevol} with all equations in \eqref{eq:evolutionconstraintII}, gives the following characterization of globally hyperbolic supersymmetric configurations.

\begin{lemma}
\label{lemma:susyevolutionNS}
Let $(g,\bar b,\bar\phi) \in \Conf(\bar{\cP},\bar{\cA},\bar{Y})$ be a globally hyperbolic configuration with reduction $( h_t , a_t , b_t , \phi_t,\lambda_t)$. Then, $(g,\bar b,\bar\phi)$ is supersymmetric if and only if there exist a family of positively oriented orthonormal coframes $\fre_t=(e_{tu},e_{tl},e_{tn})$ of $(\Sigma , h_t)$ and a family of functions $\frf_t\in C^\infty(\Sigma)$ satisfying the following evolution equations:
\begin{equation}
 \label{eq:NSNSspecializedevolution}
 \partial_t\frf_t+\dd_\Sigma\lambda_t(e^{\sharp}_{tu}) = 0\, , \qquad \partial_t\fre_t
 =\cL_{\lambda_te^{\sharp}_{tu}}\fre_t
 - (\dd_\Sigma\lambda_t(\fre_t^\sharp)
       +\lambda_t\dd_\Sigma\frf_t(\fre_t^\sharp) )e_{tu}\, , \qquad D_t\phi_t = \varphi_{\phi_t}(e^{\sharp}_{tu})
\end{equation}

\noindent
together with the following time-dependent constraints:
\begin{eqnarray}
 \label{eq:NSNSspecializedconstraints}
 & \dd_\Sigma(e^{-\frf_t}e_{tu}) = \varphi_{\phi_t}\wedge e^{-\frf_t}e_{tu}\, , \qquad   \dd_\Sigma e_{tl}\wedge e_{tu} = 0\, ,
 \qquad \dd_\Sigma e_{tn}\wedge e_{tu}=0 \nonumber \\
& \widehat{c}_t = f_{b_t}\ast_{h_t} e_{tu} + e_{tu}\wedge \varphi_{\phi_t}\, , \qquad (\cL_{e^{\sharp}_{tu}}e_{tn})(e^{\sharp}_{tl}) -(\cL_{e^{\sharp}_{tu}}e_{tl})(e^{\sharp}_{tn}) = f_{b_t}
\end{eqnarray}
\end{lemma}

\begin{proof}
By Proposition \ref{prop:susyevol}, $(g,\bar b,\bar\phi)$ is supersymmetric if and only if there exist a family of functions $\fra_t$ and a family of coframes $\fre_t$, orthonormal for $h_t$ and positively oriented since they arise from the globally hyperbolic reduction of the isotropic parallelism of of the sypersymmetry parameter $\varepsilon$, satisfying \eqref{eq:susyevol}. As shown in the proof of Proposition \ref{prop:susyevol}, the system \eqref{eq:susyevol} consists of the differential system of Proposition \ref{prop:Thetaevolutionsystem} for the torsion $\alpha^o_t = - f_{b_t}$, $\alpha^{\perp}_t = -\theta_t$, together with the dilatino conditions $\varphi_{\phi_t}=p_t e_{tu}+\widehat{c}_t(e_{tu})$ and $f_{b_t}=\ast_{h_t}(\widehat{c}_t\wedge e_{tu})$. Since $\widehat{c}_t(e_{tu})$ is orthogonal to $e_{tu}$ and $\ast_{h_t}(\widehat{c}_t\wedge e_{tu})$ is the coefficient of $\ast_{h_t}e_{tu}$ in $\widehat{c}_t$, the dilatino conditions are equivalent to $\widehat{c}_t = f_{b_t}\ast_{h_t} e_{tu} + e_{tu}\wedge \varphi_{\phi_t}$ and $p_t = \varphi_{\phi_t}(e^{\sharp}_{tu})$, where $p_t = D_t \phi_t$. On the other hand, by Lemmas \ref{lemma:refor1} and \ref{lemma:refor2}, the differential system of Proposition \ref{prop:Thetaevolutionsystem} is equivalent to the system \eqref{eq:evolutionconstraintII} for a family of functions $\frf_t$. For the torsion $\alpha^o_t = - f_{b_t}$, $\alpha^{\perp}_t = -\theta_t = -\ast_{h_t}\widehat{c}_t$, the dilatino conditions give $\theta_t(e^{\sharp}_{tu}) = f_{b_t}$ and $\ast_{h_t}\theta_t = \widehat{c}_t = f_{b_t}\ast_{h_t} e_{tu} + e_{tu}\wedge \varphi_{\phi_t}$, whence:
\begin{equation}
\label{eq:NSNSspecializationidentities}
 \alpha_t^\perp(e^{\sharp}_{tu}) = \alpha_t^o\, , \qquad \ast_{h_t}(\alpha_t^\perp-\alpha_t^oe_{tu}) = - e_{tu}\wedge\varphi_{\phi_t}\, ,
 \qquad
 \alpha_t^\perp(e^{\sharp}_{tu})-2\alpha_t^o=f_{b_t}
\end{equation}

\noindent
Substituting \eqref{eq:NSNSspecializationidentities} into \eqref{eq:evolutionconstraintII}, the term proportional to $\alpha^{\perp}_t(e^{\sharp}_{tu}) - \alpha^o_t$ in the evolution equation of $\fre_t$ vanishes, the equation for $\dd_{\Sigma} e_{tu}$ becomes $\dd_{\Sigma} e_{tu} - \dd_{\Sigma}\frf_t\wedge e_{tu} = \varphi_{\phi_t}\wedge e_{tu}$, and the right-hand side of the last equation in \eqref{eq:evolutionconstraintII} becomes $f_{b_t}$. This gives \eqref{eq:NSNSspecializedevolution} and \eqref{eq:NSNSspecializedconstraints}, and hence we conclude.
\end{proof}

\begin{remark}
\label{remark:NSevolutionisgauge}
By Lemma \ref{lemma:killingreformulation} applied in the case $r_t = 0$, the first two equations in \eqref{eq:NSNSspecializedevolution} state that the Dirac current $u^{\sharp_g}$ is a Killing vector field of $g$ whose flow preserves $u$ and the square $u\wedge(l + \mathrm{i}\, n)$, whereas the third one states that $\varphi_{\bar{\phi}}(u^{\sharp_g}) = 0$, since:
\begin{equation*}
\varphi_{\bar{\phi}} =  D_t\phi_t\, \mathfrak{t}_t + \varphi_{\phi_t}    
\end{equation*}

\noindent
that is, that $\bar{\phi}$ is invariant under the same flow. Since $\alpha_{\bar{b}}$ is determined by the isotropic parallelism by Remark \ref{remark:torsiondetermined}, so is $H_{\bar{b}} = \ast_g\alpha_{\bar{b}}$. Hence the supersymmetric Neveu-Schwarz evolution flow transports every field along the Dirac current. This is the mechanism behind the time-independent development of Corollary \ref{cor:timeindependentsusydevelopment}.
\end{remark}

\noindent
The first equation in \eqref{eq:NSNSspecializationidentities} trivializes the obstruction obtained in Theorem \ref{thm:wellposedFinalSkew}, that is:
\begin{equation*} 
 \dd_\Sigma (e^{-\frf_t}
 (\alpha_t^\perp(e^{\sharp}_{tu})-\alpha_t^o) )=0,
\end{equation*}
and removes the additional phase rotation of $e_{tl}+\mathrm{i}e_{tn}$ in the normalized evolution. By Lemma \ref{lemma:susyevolutionNS}, we define the constraint equations of the supersymmetric Neveu-Schwarz evolution flow as the following differential system:
\begin{eqnarray}
 \label{eq:NSNSsusyconstraintsI}
 & e^{\frf} \dd_\Sigma(e^{-\frf}e_{u}) = \varphi_{\phi}\wedge e_{u}\, , \qquad   \dd_\Sigma e_{l}\wedge e_{u} = 0\, ,
 \qquad \dd_\Sigma e_{n}\wedge e_{u}=0 \\ 
& \frc = f_{b}\ast_{\fre} e_{u} + e_{u}\wedge\varphi_{\phi}\, , \qquad (\cL_{e^{\sharp}_{u}}e_{n})(e^{\sharp}_{l}) -(\cL_{e^{\sharp}_{u}}e_{l})(e^{\sharp}_{n}) = f_{b} = \ast_{\fre}H_b.  \label{eq:NSNSsusyconstraintsII}
\end{eqnarray}

\noindent
for tuples $(\fre,\frf,b,\frc,\phi)$ defined on a bundle gerbe $(\cP,\cA,Y)$ over a three-dimensional manifold $\Sigma$, where $\fre$ is a positively oriented global coframe on $\Sigma$, $\frf \in C^{\infty}(\Sigma)$ is a function, $\phi\colon \Sigma \to S^1$ is a dilaton with curvature $\varphi_{\phi}\in \Omega^1(\Sigma)$, $b$ is a curving on $(\cP,\cA,Y)$, and $\frc\in \Omega^2(\Sigma)$ is a two-form. Hence, the admissible initial data for the supersymmetric Neveu-Schwarz evolution flow consists of tuples $(\fre,\frf,b, \frc ,\phi)$ satisfying all equations in \eqref{eq:NSNSsusyconstraintsI} and \eqref{eq:NSNSsusyconstraintsII}. In particular, it immediately follows that every supersymmetric Neveu-Schwarz initial data belongs to the integrable branch $r=0$ of Proposition \ref{prop:geometryconstraintbranches}. Thus $[e^{\sharp}_{l},e^{\sharp}_{n}]=0$ and $\Ker(e_u)\subset T\Sigma$ is tangent to a foliation by flat translation surfaces.  If $\Sigma$ is compact without boundary, Corollary \ref{cor:compactconstraintbranches} implies that $\Sigma$ is an oriented two-torus bundle over the circle.

\begin{thm}
\label{thm:wellposedFinalSusy}
Let $(\cP,\cA,Y)$ be an abelian bundle gerbe on a connected and oriented three-manifold $\Sigma$ equipped with a tuple $(\fre,\frf,b,\frc,\phi)$ satisfying the constraint equations \eqref{eq:NSNSsusyconstraintsI} and \eqref{eq:NSNSsusyconstraintsII}. Then, there exists an oriented and connected four-manifold $M$ equipped with a globally hyperbolic supersymmetric configuration $(g,\bar{b},\bar{\phi},\varepsilon)$ on an abelian bundle gerbe $(\bar{\cP},\bar{\cA},\bar{Y})$ over $M$, together with an isometric embedding $\iota\colon (\Sigma,h_{\fre})\hookrightarrow (M,g)$ such that $(\Sigma,h_{\fre})$ is a Cauchy hypersurface in $(M,g)$, the globally hyperbolic reduction $(\lambda_t,\fra_t,\fre_t)$ of the isotropic parallelism of $\varepsilon$ in the splitting \eqref{eq:globahyp} defined by $\Sigma$ satisfies $(\fra_0,\fre_0) = (-\frf,\fre)$, and:
\begin{equation*}
(\cP,\cA,Y) = (\iota^{\ast}\bar{\cP},\iota^{\ast}\bar{\cA},\iota^{\ast}\bar{Y}) \, , \qquad  b = \iota^{\ast}\bar{b}\, , \qquad \phi = \iota^{\ast}\bar{\phi}\, , \qquad \frc = - \ast_{\fre} \iota^{\ast}\alpha_{\bar{b}}
\end{equation*}

\noindent
where $H_{\bar{b}} = \ast_g \alpha_{\bar{b}}$. The germ along $\Sigma$ of $(g,\bar{\phi},H_{\bar{b}},\varepsilon)$ is unique up to isometries fixing $\Sigma$ and the sign of $\varepsilon$. Conversely, the Cauchy hypersurface $\Sigma$ of every globally hyperbolic supersymmetric configuration on a bundle gerbe reducible over an open tubular neighborhood of $\Sigma$ carries a solution $(\fre,\frf,b,\frc,\phi)$ to the constraint equations \eqref{eq:NSNSsusyconstraintsI} and \eqref{eq:NSNSsusyconstraintsII}.
\end{thm}

\begin{proof}
Let  $(\fre,\frf,b,\frc,\phi)$ be a solution of the constraint equations \eqref{eq:NSNSsusyconstraintsI} and \eqref{eq:NSNSsusyconstraintsII} on a bundle gerbe $(\cP,\cA,Y)$ over $\Sigma$. By Lemma \ref{lemma:susyevolutionNS}, we have $\alpha^o = -f_b$ and $\alpha^{\perp} = -\ast_{\fre}\frc$. By \eqref{eq:NSNSsusyconstraintsII} we have:
\begin{equation*}
 \alpha^{\perp}(e^{\sharp}_u) = -\ast_{\fre}(\frc\wedge e_u) = -f_b = \alpha^o\, , \qquad \ast_{\fre}(\alpha^{\perp} - \alpha^o e_u) = -\frc + f_b \ast_{\fre} e_u = - e_u\wedge\varphi_{\phi}
\end{equation*}

\noindent
whence \eqref{eq:NSNSsusyconstraintsI} and \eqref{eq:NSNSsusyconstraintsII} imply that $(\fre,\frf,\alpha^o,\alpha^{\perp})$ satisfies \eqref{eq:constraintIIthm1} and \eqref{eq:constraintIIthm2}, with vanishing constant $r = \alpha^{\perp}(\cX) - e^{-\frf}\alpha^o$ in the notation of Theorem \ref{thm:wellposedFinalSkew}, where $\cX = e^{-\frf}e^{\sharp}_u$. Let $\Phi\colon \cU\to\Sigma$ be the maximal flow of $\cX$ and define on an open tubular neighborhood $\cU\subset \mathbb{R} \times \Sigma$ of $\Sigma$:
\begin{eqnarray}
& \lambda_t = e^{-\frf_t}\, , \qquad \frf_t = \Phi_t^{\ast}\frf\, , \qquad \fre_t = \Phi_t^{\ast}\fre\, , \qquad \phi_t = \Phi_t^{\ast}\phi \nonumber \\
& a_t = 0\, , \qquad b_t = b + \pi^{\ast}\int_0^t \Phi_s^{\ast}(e^{-\frf}\frc)\, \dd s \label{eq:NSNSexplicitdevelopment}
\end{eqnarray}

\noindent
For each $t$, these families are defined on $\cU_t := \left\{p\in\Sigma\,\,\vert\,\,(t,p)\in\cU\right\}$, with $a_t$ and $b_t$ defined on $\pi^{-1}(\cU_t)$ and all initial fields restricted accordingly. The integral defining $b_t$ is well-defined and smooth because each maximal interval of the flow through $p$ contains the segment from $0$ to $t$ whenever $(t,p)\in\cU$.

\noindent
Let $(\bar{\cP},\bar{\cA},\bar{Y})$ be the pull-back of $(\cP,\cA,Y)$ to $\cU$ through the projection $\cU\to\Sigma$, which we consider as a reducible bundle gerbe with reduction $(\cP,Y,\cA_t,\Psi_t) = (\cP,Y,\cA,0)$, and let $(g,\bar{b},\bar{\phi})$ be the configuration on $(\bar{\cP},\bar{\cA},\bar{Y})$ with reduction $(h_{\fre_t},a_t,b_t,\phi_t,\lambda_t)$ as defined above, that is:
\begin{equation*}
g = -e^{-2\frf_t}\dd t\otimes \dd t + h_{\fre_t}  \, , \qquad \bar{b}(t,p) = b_t(p)  \, , \qquad \bar{\phi}(t,p) = \phi_t(p)
\end{equation*}

\noindent
for every $(t,p)\in \cU$, and in the case of the curving $\bar{b}$, for every:
\begin{equation*}
(t,y)\in \bar{Y}=\{(t,y):(t,\pi(y))\in\mathcal U\},
\end{equation*}

\noindent
Since $\delta\pi^{\ast} = 0$, we have $\delta b_t = \delta b = F_{\cA}$, so $b_t$ is a family of curvings on $(\cP,Y,\cA)$, and $a_t = 0$ satisfies Equation \eqref{eq:conditionat} because $\cA_t = \cA$ and $\Psi_t = 0$. Hence, $\bar{b}$ is a curving on $(\bar{\cP},\bar{\cA},\bar{Y})$ by Lemma \ref{lemma:reduciblegerbe}, with derived family $c_t = \Phi_t^{\ast}(e^{-\frf}\frc) = \lambda_t \Phi_t^{\ast}\frc$, whence $\widehat{c}_t = \Phi_t^{\ast}\frc$. Moreover, since:
\begin{equation*}
e^{-\frf} f_b \ast_{\fre} e_u = \iota_{\cX} H_b    
\end{equation*}

\noindent
it follows that Equations \eqref{eq:NSNSsusyconstraintsI} and \eqref{eq:NSNSsusyconstraintsII} together with the closedness of $\varphi_{\phi}$ give:
\begin{equation*}
\dd_{\Sigma}(e^{-\frf}\frc) = \dd_{\Sigma}(\iota_{\cX}H_b) + \dd_{\Sigma}(e^{-\frf}e_u)\wedge\varphi_{\phi} = \cL_{\cX}H_b + \varphi_{\phi}\wedge e^{-\frf}e_u\wedge\varphi_{\phi} = \cL_{\cX}H_b
\end{equation*}

\noindent
Therefore:
\begin{equation*}
H_{b_t} = H_b + \int_0^t \Phi_s^{\ast}\cL_{\cX}H_b\, \dd s = \Phi_t^{\ast}H_b    
\end{equation*}

\noindent
that is, $f_{b_t} = \Phi_t^{\ast}f_b$. We verify now that $(g,\bar{b},\bar{\phi})$ satisfies all equations in \eqref{eq:NSNSspecializedevolution} and \eqref{eq:NSNSspecializedconstraints} for the families $(\fre_t,\frf_t)$ in \eqref{eq:NSNSexplicitdevelopment}. Since $h_{\fre_t} = \Phi_t^{\ast}h_{\fre}$ and $\cX$ is invariant under its own flow, we have $\lambda_t e^{\sharp}_{tu} = \Phi_t^{\ast}(e^{-\frf}e^{\sharp}_u) = \cX$, whence $\partial_t \frf_t = \cX(\frf_t) = -\dd_{\Sigma}\lambda_t(e^{\sharp}_{tu})$ and $\partial_t\fre_t = \cL_{\cX}\fre_t = \cL_{\lambda_t e^{\sharp}_{tu}}\fre_t$, which are the first two equations in \eqref{eq:NSNSspecializedevolution} due to $\dd_{\Sigma}\lambda_t = -\lambda_t\dd_{\Sigma}\frf_t$. Similarly, $\varphi_{\phi_t} = \Phi_t^{\ast}\varphi_{\phi}$ and:
\begin{equation*}
D_t\phi_t = \lambda_t^{-1}\partial_t\phi_t = \lambda_t^{-1}\varphi_{\phi_t}(\cX) = \varphi_{\phi_t}(e^{\sharp}_{tu})    
\end{equation*}

\noindent
which is the third equation in \eqref{eq:NSNSspecializedevolution}. Finally, since $\Phi_t$ is orientation preserving and $h_{\fre_t} = \Phi_t^{\ast}h_{\fre}$, we have $\ast_{h_{\fre_t}}\Phi_t^{\ast} = \Phi_t^{\ast}\ast_{\fre}$, so the constraints \eqref{eq:NSNSspecializedconstraints} are the pull-backs through $\Phi_t$ of \eqref{eq:NSNSsusyconstraintsI} and \eqref{eq:NSNSsusyconstraintsII}, using $\widehat{c}_t = \Phi_t^{\ast}\frc$ and $f_{b_t} = \Phi_t^{\ast}f_b$. By Lemma \ref{lemma:susyevolutionNS}, whose proof can be applied locally, $(g,\bar{b},\bar{\phi})$ is supersymmetric on $\cU$, with supersymmetry parameter $\varepsilon$ whose isotropic parallelism has globally hyperbolic reduction $(\lambda_t,\fra_t,\fre_t) = (e^{-\frf_t},-\frf_t,\fre_t)$ by \eqref{eq:NSNS-scale-identification}. In fact, $(g,\varepsilon)$ is the development of $(\fre,\frf,\alpha^o,\alpha^{\perp})$ constructed in the proof of Theorem \ref{thm:wellposedFinalSkew} in the gauge $\lambda_t = e^{-\frf_t}$, since $r = 0$ in Equation \eqref{eq:zetaux}. As in the proof of Theorem \ref{thm:wellposedFinalSkew}, the open interior $M\subset \cU$ of the domain of dependence of $\Sigma$ is globally hyperbolic with Cauchy hypersurface $\Sigma$, and we restrict $(\bar{\cP},\bar{\cA},\bar{Y})$, $(g,\bar{b},\bar{\phi})$ and $\varepsilon$ to $M$. Since $\Phi_0$ is the identity, we have:
\begin{equation*}
 (\fra_0,\fre_0) = (-\frf,\fre)\, , \qquad   \iota^{\ast}\bar{b} = b_0 = b\, , \qquad  \iota^{\ast}\bar{\phi} = \phi_0 = \phi
\end{equation*}

\noindent
the restriction of $(\bar{\cP},\bar{\cA},\bar{Y})$ to $\Sigma$ is $(\cP,\cA,Y)$ and, since $\alpha_{\bar{b}} = -(f_{b_t}\frt_t + \theta_t)$ as shown in the proof of Proposition \ref{prop:susyevol}, we obtain:
\begin{equation*}
\iota^{\ast}\alpha_{\bar{b}} = -\theta_0 = -\ast_{\fre}\widehat{c}_0 = -\ast_{\fre}\frc    
\end{equation*}

\noindent
This proves existence. If $e^{2\frf}h_{\fre}$ is complete, in particular if $\Sigma$ is compact, then $\cU = \mathbb{R}\times\Sigma$ and $(\cU,g)$ is already globally hyperbolic with Cauchy hypersurface $\Sigma$ by Corollary \ref{cor:timeindependentdevelopment}, so that no restriction is needed. 

\noindent
Regarding uniqueness, let $(M_i,g_i,\bar{b}_i,\bar{\phi}_i,\varepsilon_i,\iota_i)$, $i = 1,2$, be two tuples as in the statement. By the gravitino equation, $(g_i,\varepsilon_i)$ is a globally hyperbolic skew-torsion parallel spinor with torsion:
\begin{equation*}
\alpha_{\bar{b}_i} = -(f_{b_t}\frt_t + \theta_t)    
\end{equation*}

\noindent
see the proof of Proposition \ref{prop:susyevol}, whose restriction to $\Sigma$ in the sense of Theorem \ref{thm:wellposedFinalSkew} is, by the conditions in the statement, given by:
\begin{equation*}
(\fre_0,\frf_0,-f_{b_0},-\theta_0) = (\fre,\frf,-f_b,-\ast_{\fre}\frc) = (\fre,\frf,\alpha^o,\alpha^{\perp})    
\end{equation*}

\noindent
for $i = 1,2$. Hence, by Lemma \ref{lemma:uniqueevol} and Remark \ref{remark:straightening}, after restricting $M_1$ and $M_2$ to open neighborhoods of $\iota_i(\Sigma)$, the flows of the Dirac currents $u_i^{\sharp_{g_i}}$ define diffeomorphisms $\Psi_i\colon \mathcal{V}\to M_i$ from a common neighborhood $\mathcal{V}$ of $\left\{ 0\right\}\times\Sigma$ in $\mathbb{R}\times\Sigma$ which pull back $(g_i,\alpha_{\bar{b}_i},[u_i,v_i,l_i,n_i])$ to the same metric, torsion and isotropic parallelism, and which satisfy $\Psi_i^{\ast}u_i^{\sharp_{g_i}} = -\partial_s$. In particular, $F := \Psi_2\circ\Psi_1^{-1}$ is an isometry with $F\circ\iota_1 = \iota_2$ and $F^{\ast}\alpha_{\bar{b}_2} = \alpha_{\bar{b}_1}$, whence $F^{\ast}H_{\bar{b}_2} = H_{\bar{b}_1}$, and $F$ identifies $\varepsilon_1$ and $\varepsilon_2$ up to a sign. On the other hand, the dilatino equation implies:
\begin{equation*}
\varphi_{\bar{\phi}_i}(u_i^{\sharp_{g_i}}) = 0\, ,
\end{equation*}

\noindent
as shown in the proof of Proposition \ref{prop:susyevol}, so that $\Psi_i^{\ast}\bar{\phi}_i$ is constant along every segment $\mathcal{V}\cap(\mathbb{R}\times\left\{ p\right\})$ and therefore equals $\phi(p)$ there, because $\iota_i^{\ast}\bar{\phi}_i = \phi$. Hence $F^{\ast}\bar{\phi}_2 = \bar{\phi}_1$. This proves uniqueness.

\noindent
For the converse, let $(g,\bar{b},\bar{\phi})$ be a globally hyperbolic supersymmetric configuration on a bundle gerbe $(\bar{\cP},\bar{\cA},\bar{Y})$ and let $\Sigma \subset (M,g)$ be a smooth Cauchy hypersurface. We choose a small enough tubular neighborhood $\cU$ of $\Sigma$ such that $\cU\simeq \cI\times \Sigma$ and the restriction of $(g,\bar{b},\bar{\phi})$ to $\cU$ becomes a reducible bundle gerbe. Denote by $(h_t,a_t,b_t,\phi_t,\lambda_t)$ the reduction of $(g,\bar{b},\bar{\phi})$, and let $(\fre_t,\frf_t)$ be as in Lemma \ref{lemma:susyevolutionNS}. Then $(\fre,\frf,b,\frc,\phi) := (\fre_0,\frf_0,b_0,\widehat{c}_0,\phi_0)$ satisfies \eqref{eq:NSNSsusyconstraintsI} and \eqref{eq:NSNSsusyconstraintsII}, which are the constraints \eqref{eq:NSNSspecializedconstraints} at $t = 0$. 
\end{proof}

\noindent 
The following corollary is the supersymmetric analog of Corollary \ref{cor:timeindependentdevelopment}, and its proof goes along the same lines: it identifies the metric, the curvature of the curving, the dilaton and the isotropic parallelism of the development of Theorem \ref{thm:wellposedFinalSusy} with those of an open restriction of a supersymmetric configuration defined on the whole $\mathbb{R}\times\Sigma$ whose metric, curving and dilaton do not depend on $t$.

\begin{cor}
\label{cor:timeindependentsusydevelopment}
Let $(\fre,\frf,b,\frc,\phi)$ be a solution of the constraint equations \eqref{eq:NSNSsusyconstraintsI} and \eqref{eq:NSNSsusyconstraintsII} of the supersymmetric Neveu-Schwarz evolution flow on a bundle gerbe $(\cP,\cA,Y)$ over a connected and oriented three-manifold $\Sigma$, and let $\Phi\colon \cU\to\Sigma$ be the maximal flow of $\cX = e^{-\frf}e^{\sharp}_u\in\mathfrak{X}(\Sigma)$. Then:
\begin{enumerate}
\item The pull-back of $(\cP,\cA,Y)$ to $\mathbb{R}\times\Sigma$ through the projection is a reducible bundle gerbe with reduction $(\cP,Y,\cA,0)$, and the triple $(g,\bar{b},\bar{\phi})$ on it given by:
\begin{equation}
\label{eq:timeindependentsusyconfiguration}
g = h_{\fre} - e^{-\frf}\, \dd t\odot e_u\, , \qquad \bar{b} = \dd t\wedge \pi^{\ast}(e^{-\frf}e_u) + b\, , \qquad \bar{\phi} = \phi
\end{equation}

\noindent
where $\fre$, $\frf$, $b$ and $\phi$ are considered as objects on $\mathbb{R}\times\Sigma$ and $\mathbb{R}\times Y$ that do not depend on $t$, is a supersymmetric configuration with supersymmetry parameter the skew-torsion parallel spinor associated to the isotropic parallelism $[u,v,l,n]$ of Equation \eqref{eq:timeindependentcoframe} with $r = 0$, and:
\begin{equation}
\label{eq:timeindependentsusyflux}
H_{\bar{b}} = H_b + \dd t\wedge e^{-\frf}e_u\wedge\varphi_{\phi}\, , \qquad \alpha_{\bar{b}} = -\ast_{\fre}\frc\, , \qquad \varphi_{\bar{\phi}} = \varphi_{\phi}
\end{equation}

\noindent
In particular, $g$, $\bar{b}$ and $\bar{\phi}$ do not depend on $t$, the vector field $u^{\sharp_g} = -\partial_t$ is a null Killing vector field of $g$ preserving $H_{\bar{b}}$ and $\varphi_{\bar{\phi}}$, and $\mathbb{R}\times\Sigma$ carries a supersymmetric configuration on the pull-back of $(\cP,\cA,Y)$.

\item The map $\cU\to\mathbb{R}\times\Sigma$, $(t,p)\mapsto (t,\Phi_t(p))$, is a diffeomorphism onto an open neighborhood of $\left\{ 0 \right\} \times \Sigma$ in $\mathbb{R}\times\Sigma$ which restricts to the identity on $\Sigma$ and pulls back $g$, $H_{\bar{b}}$, $\bar{\phi}$ and $[u,v,l,n]$ to the metric, the curvature of the curving, the dilaton and the isotropic parallelism of the supersymmetry parameter of the supersymmetric configuration \eqref{eq:NSNSexplicitdevelopment} constructed in the proof of Theorem \ref{thm:wellposedFinalSusy}.

\item If the Riemannian metric $e^{2\frf}h_{\fre}$ on $\Sigma$ is complete, then the map of the previous item is a diffeomorphism of $\mathbb{R}\times\Sigma$ onto itself, and the Lorentzian four-manifold $(\mathbb{R}\times\Sigma , g)$ is a maximal globally hyperbolic development of the given initial data.
\end{enumerate}
\end{cor} 

\noindent
The choice $\lambda_t=e^{-\frf_t}$ is used to integrate the supersymmetric Neveu-Schwarz evolution flow explicitly.  By Lemma \ref{lemma:uniqueevol}, changing to an arbitrary positive lapse changes only the representative of the same local germ.  Theorem \ref{thm:wellposedFinalSusy} proves well-posedness of the gravitino and dilatino equations without imposing the Einstein, Maxwell, or dilaton equations of motion, nor their Hamiltonian and momentum constraints. We address their compatibility with the supersymmetric development below.


\subsection{Compatibility of flows and initial data}
\label{sec:compatibleinitialdata}


Theorem \ref{thm:wellposedFinalSusy} is concerned with the Neveu-Schwarz supersymmetry flow, not with the Neveu-Schwarz flow determined by the bosonic equations of the theory. However, as we will prove below in Theorem \ref{thm:susyflowcompatibility}, these flows are compatible when starting with initial data admissible to both flows, that is, satisfying both sets of constraint equations.  

\begin{thm}
\label{thm:susyflowcompatibility}
Let $(\lambda_t,\fre_t, a_t , b_t,\phi_t)$ be a Neveu-Schwarz supersymmetry flow and let $(g,\bar{b},\bar{\phi})$ be the globally hyperbolic supersymmetric configuration with reduction $(h_{\fre_t}, a_t , b_t,\phi_t,\lambda_t)$. Then, $(g,\bar{b},\bar{\phi})$ is a solution of the Neveu-Schwarz system if and only if:
\begin{equation}
\label{eq:initialdatacondition}
\cW_t(e^{\sharp}_{tu})\vert_{t=0} = 0 \, , \qquad  \cM_t\vert_{t=0} = 0
\end{equation}
where $\mathcal{W}_t$ and $\mathcal{M}_t$ respectively denote the time-dependent momentum and Maxwell constraints. 
\end{thm}

\begin{proof}
By Proposition \ref{prop:normalhyperbolicflow}, such $(g,\bar{b},\bar{\phi})$ is a solution of the Neveu-Schwarz system \eqref{eq:NSNSsystem4d} if and only if its reduced tuple $(h_{\fre_t}, a_t , b_t,\phi_t,\lambda_t)$ satisfies the evolution equations \eqref{eq:normalevolution} and the time-dependent constraint equations \eqref{eq:normalconstraint} for every $t\in\cI$. Evaluating the constraint equations at $t=0$ yields the equations in the statement.

\noindent
For the converse, we need to prove that if $(\lambda_t,\fre_t, a_t , b_t,\phi_t)$ is a supersymmetric Neveu-Schwarz flow satisfying \eqref{eq:initialdatacondition}, then its associated tuple $(h_{\fre_t}, a_t , b_t,\phi_t,\lambda_t)$ is a Neveu-Schwarz flow, that is, satisfies equations  \eqref{eq:normalevolution} and \eqref{eq:normalconstraint}, which is equivalent to $(g,\bar{b},\bar{\phi})$ being a solution of the Neveu-Schwarz equations of motion. In order to do this, we first rewrite the evolution equations \eqref{eq:normalevolution} in terms of the time-dependent constraints $\cH_t$, $\cW_t$, and $\cM_t$. For simplicity in the exposition, we set $h_t=h_{\fre_t}$ and write: 
\begin{eqnarray*}
& \cE_t
:=\Ric^{h_t}
+(\Tr_{h_t}(\Theta_t)+p_t)\Theta_t
-2\Theta_t\circ_{h_t}\Theta_t
+\nabla^{h_t}\varphi_{\phi_t}
-\nabla^{h_t}\eta_t-\eta_t\otimes\eta_t\\
& \qquad
-\frac{1}{2}\theta_t\otimes\theta_t
+\frac{1}{2}
 (\vert\theta_t\vert_{h_t}^2-f_{b_t}^2 )h_t
-D_t\Theta_t\, ,\\
& \cD_t
:=D_tp_t+\nabla^{h_t\ast}\varphi_{\phi_t}
-(\Tr_{h_t}(\Theta_t)+p_t)p_t
-\langle\varphi_{\phi_t},\eta_t\rangle_{h_t}
-f_{b_t}^2+\vert\theta_t\vert_{h_t}^2
+\vert\varphi_{\phi_t}\vert_{h_t}^2\, ,\\
& \cB_t
:=D_t\theta_t-\dd_{\Sigma}f_{b_t}
-f_{b_t}\eta_t-p_t\theta_t+f_{b_t}\varphi_{\phi_t}.
\end{eqnarray*}

\noindent
We will respectively refer to $\cE_t$, $\cD_t$ and $\cB_t$ as the Einstein, dilaton, and Maxwell \emph{defects}. We will prove that, when evaluated on a Neveu-Schwarz supersymmetry flow, these equations can be written as follows:
\begin{equation}
\label{eq:reducedevolutiondefects}
\cE_t
=e_{tu}\odot\cW_t
-\frac{1}{2}\cH_t\,e_{tu}\otimes e_{tu}\, , \qquad
\cD_t=-\ast_{h_t}(\cM_t\wedge e_{tu})\, , \qquad
\cB_t=\iota_{e^{\sharp}_{tu}}\cM_t.
\end{equation}

\noindent
We establish these identities directly from the reduced equations \eqref{eq:susyevol}, without imposing the initial conditions \eqref{eq:initialdatacondition}. We first recall the following identities:
\begin{eqnarray*}
& \dd_{\Sigma}\varphi_{\phi_t}=0\, , \qquad
D_t\varphi_{\phi_t}=\dd_{\Sigma}p_t+p_t\eta_t\, ,\\
& D_tf_{b_t}
=\Tr_{h_t}(\Theta_t)f_{b_t}
-\nabla^{h_t\ast}\theta_t
+\langle\eta_t,\theta_t\rangle_{h_t}.
\end{eqnarray*}

\noindent
They follow from the closedness of the curvatures of the configuration. Moreover, $\dd_{\Sigma}\eta_t=0$, since
$\eta_t=\lambda_t^{-1}\dd_{\Sigma}\lambda_t$.
The algebraic supersymmetry equations give:
\begin{eqnarray*}
& \varphi_{\phi_t}=p_t e_{tu}+\widehat{c}_t(e_{tu})\, , \qquad
p_t=\varphi_{\phi_t}(e^{\sharp}_{tu})\, , \qquad
f_{b_t}=\theta_t(e^{\sharp}_{tu})\, ,\\
& \widehat{c}_t
=f_{b_t}\ast_{h_t}e_{tu}
+e_{tu}\wedge\varphi_{\phi_t}\, , \qquad
\theta_t
=f_{b_t}e_{tu}
+\ast_{h_t}(e_{tu}\wedge\varphi_{\phi_t}).
\end{eqnarray*}

\noindent
We begin with the dilaton defect. Differentiating
$p_t=\varphi_{\phi_t}(e^{\sharp}_{tu})$,
using $D_th_t=-2\Theta_t$, the temporal equation for $e_{tu}$ in \eqref{eq:susyevol} and the curvature identity for $D_t\varphi_{\phi_t}$, yields:
\begin{equation*}
D_tp_t
=\dd_{\Sigma}p_t(e^{\sharp}_{tu})
+\Theta_t(\varphi_{\phi_t}^{\sharp},e^{\sharp}_{tu})
-\frac{1}{2}\vert\widehat{c}_t(e_{tu})\vert_{h_t}^2
+\langle\varphi_{\phi_t},\eta_t\rangle_{h_t}.
\end{equation*}

\noindent
On the other hand, the divergence identities obtained in the proof of Theorem \ref{thm:compatibilityflows}, together with
$\varphi_{\phi_t}=p_t e_{tu}+\widehat{c}_t(e_{tu})$, give:
\begin{eqnarray*}
\nabla^{h_t\ast}\varphi_{\phi_t}
=-\dd_{\Sigma}p_t(e^{\sharp}_{tu})
+p_t\Tr_{h_t}(\Theta_t)
-\Theta_t(\varphi_{\phi_t}^{\sharp},e^{\sharp}_{tu})-(\nabla^{h_t\ast}\widehat{c}_t)(e^{\sharp}_{tu})
-\frac{1}{2}\vert\widehat{c}_t(e_{tu})\vert_{h_t}^2.
\end{eqnarray*}

\noindent
Substituting these expressions into $\cD_t$ and using \eqref{eq:dilatinoalgebra}, we obtain:
\begin{equation*}
\cD_t =\vert\widehat{c}_t(e_{tu})\vert_{h_t}^2 -(\nabla^{h_t\ast}\widehat{c}_t)(e^{\sharp}_{tu}) =-\ast_{h_t}(\cM_t\wedge e_{tu}),
\end{equation*}

\noindent
where the last equality is \eqref{eq:Maxwellnormalcomponent}. We next consider the Maxwell defect. The variation of the Hodge operator under $D_th_t=-2\Theta_t$ gives:
\begin{equation*}
D_t\theta_t = \ast_{h_t}(D_t\widehat{c}_t) + \Tr_{h_t}(\Theta_t)\theta_t - 2\Theta_t(\theta_t^{\sharp},-).
\end{equation*}

\noindent
Differentiate $\widehat{c}_t=f_{b_t}\ast_{h_t}e_{tu} +e_{tu}\wedge\varphi_{\phi_t}$ and substitute the reduced curvature identities and:
\begin{equation*}
D_te_{tu} =-\Theta_t(e^{\sharp}_{tu},-) -\frac{1}{2}\widehat{c}_t(e_{tu}) +\eta_t-\eta_t(e^{\sharp}_{tu})e_{tu}.
\end{equation*}

\noindent
Computing the spatial Lie derivative by means of Remark \ref{remark:covariantevolution}, the resulting expression simplifies to:
\begin{equation*}
D_t\theta_t-\cL_{e^{\sharp}_{tu}}\theta_t-f_{b_t}\eta_t =-\ast_{h_t}\dd_{\Sigma}\varphi_{\phi_t}=0.
\end{equation*}

\noindent
In components, the difference on the left evaluates on $e^{\sharp}_{tu}$, $e^{\sharp}_{tl}$ and $e^{\sharp}_{tn}$ respectively as:
\begin{equation*}
-\dd_{\Sigma}\varphi_{\phi_t} (e^{\sharp}_{tl},e^{\sharp}_{tn})\, , \qquad \dd_{\Sigma}\varphi_{\phi_t}
(e^{\sharp}_{tu},e^{\sharp}_{tn})\, , \qquad -\dd_{\Sigma}\varphi_{\phi_t} (e^{\sharp}_{tu},e^{\sharp}_{tl}).
\end{equation*}

\noindent
Thus Cartan's formula, together with $\theta_t(e^{\sharp}_{tu})=f_{b_t}$, gives:
\begin{eqnarray*}
D_t\theta_t =  \dd_{\Sigma}f_{b_t} + \iota_{e^{\sharp}_{tu}}\dd_{\Sigma}\theta_t
+f_{b_t}\eta_t = \dd_{\Sigma}f_{b_t}
+f_{b_t}\eta_t+p_t\theta_t-f_{b_t}\varphi_{\phi_t}
+\iota_{e^{\sharp}_{tu}}\cM_t,
\end{eqnarray*}

\noindent
which proves $\cB_t=\iota_{e^{\sharp}_{tu}}\cM_t$. It remains to compute the Einstein evolution defect. We differentiate the spatial coframe equations and the equation for $\dd_{\Sigma}\fra_t$ in \eqref{eq:susyevol}, using:
\begin{eqnarray*}
D_t(\dd_{\Sigma}e^t_i) =\dd_{\Sigma}(D_te^t_i)+\eta_t\wedge D_te^t_i\, ,\qquad  D_t(\dd_{\Sigma}\fra_t)
=\dd_{\Sigma} (\eta_t(e^{\sharp}_{tu}) ) +\eta_t(e^{\sharp}_{tu})\eta_t.
\end{eqnarray*}

\noindent
Substituting the temporal coframe equations determines the components of $D_t\Theta_t$. The component in the $e^{\sharp}_{tu}$ direction is:
\begin{eqnarray*}
(D_t\Theta_t+\nabla^{h_t}\eta_t + \eta_t\otimes\eta_t) (e^{\sharp}_{tu},e^{\sharp}_{tu}) =-\vert\Theta_t(e^{\sharp}_{tu},-)\vert_{h_t}^2 +\frac{1}{4}\vert\widehat{c}_t(e_{tu})\vert_{h_t}^2.
\end{eqnarray*}

\noindent
For $i=l,n$, the mixed components are:
\begin{eqnarray*}
&  (D_t\Theta_t+\nabla^{h_t}\eta_t
+\eta_t\otimes\eta_t )
(e^{\sharp}_{tu},e^{\sharp}_{ti})\\
& \qquad
=-(\Theta_t\circ_{h_t}\Theta_t)
(e^{\sharp}_{tu},e^{\sharp}_{ti})
+\frac{1}{2}\Theta_t
(\widehat{c}_t(e_{tu})^{\sharp},e^{\sharp}_{ti})
+\frac{1}{2}\dd_{\Sigma}p_t(e^{\sharp}_{ti})
-\frac{p_t}{2}
(\nabla^{h_t}_{e^{\sharp}_{tu}}e_{tu})(e^{\sharp}_{ti}).
\end{eqnarray*}

\noindent
Finally, for $i,j=l,n$, we obtain:
\begin{eqnarray*}
&  (D_t\Theta_t+\nabla^{h_t}\eta_t
+\eta_t\otimes\eta_t )
(e^{\sharp}_{ti},e^{\sharp}_{tj})\\
& \qquad
=(\nabla^{h_t}_{e^{\sharp}_{tu}}\Theta_t)
(e^{\sharp}_{ti},e^{\sharp}_{tj})
-(\nabla^{h_t}_{e^{\sharp}_{ti}}\Theta_t)
(e^{\sharp}_{tu},e^{\sharp}_{tj})
-(\Theta_t\circ_{h_t}\Theta_t)
(e^{\sharp}_{ti},e^{\sharp}_{tj})\\
& \qquad
+\frac{1}{2}
 (\nabla^{h_t}_{e^{\sharp}_{ti}}
(\widehat{c}_t(e_{tu})) )(e^{\sharp}_{tj})
+\frac{1}{4}
\widehat{c}_t(e_{tu},e_{ti})
\widehat{c}_t(e_{tu},e_{tj}).
\end{eqnarray*}

\noindent 
The intrinsic Ricci tensor is computed from the covariant derivatives in Remark \ref{remark:covariantevolution}. Its components involving $e^{\sharp}_{tu}$ are those obtained in the proofs of Lemma \ref{lemma:curvaturesusyflow} and Theorem \ref{thm:compatibilityflows}. For $i,j=l,n$, the same computation gives:
\begin{eqnarray*}
& \Ric^{h_t}(e^{\sharp}_{ti},e^{\sharp}_{tj})
=(\nabla^{h_t}_{e^{\sharp}_{tu}}\Theta_t)
(e^{\sharp}_{ti},e^{\sharp}_{tj})
-(\nabla^{h_t}_{e^{\sharp}_{ti}}\Theta_t)
(e^{\sharp}_{tu},e^{\sharp}_{tj})\\
& \qquad
-\Tr_{h_t}(\Theta_t)\Theta_t
(e^{\sharp}_{ti},e^{\sharp}_{tj})
+(\Theta_t\circ_{h_t}\Theta_t)
(e^{\sharp}_{ti},e^{\sharp}_{tj})\\
& \qquad
-\frac{1}{2}
 (\nabla^{h_t}_{e^{\sharp}_{ti}}
(\widehat{c}_t(e_{tu})) )(e^{\sharp}_{tj})
-\frac{1}{4}
\widehat{c}_t(e_{tu},e_{ti})
\widehat{c}_t(e_{tu},e_{tj}).
\end{eqnarray*}

\noindent
Substituting these Ricci and evolution formulas into the definition of $\cE_t$, and using:
\begin{eqnarray*}
& (\nabla^{h_t}\varphi_{\phi_t})
(e^{\sharp}_{ti},e^{\sharp}_{tj})
=-p_t\Theta_t(e^{\sharp}_{ti},e^{\sharp}_{tj})
+ (\nabla^{h_t}_{e^{\sharp}_{ti}}
(\widehat{c}_t(e_{tu})) )(e^{\sharp}_{tj})\, ,\\
& \theta_t(e^{\sharp}_{ti})\theta_t(e^{\sharp}_{tj})
=\vert\widehat{c}_t(e_{tu})\vert_{h_t}^2\delta_{ij}
-\widehat{c}_t(e_{tu},e_{ti})
\widehat{c}_t(e_{tu},e_{tj})
\qquad\text{for }i,j=l,n,
\end{eqnarray*}

\noindent
we obtain:
\begin{eqnarray*}
& \cE_t(e^{\sharp}_{tu},e^{\sharp}_{tu})
=\cW_t(e^{\sharp}_{tu})
-\frac{1}{2}(\ast_{h_t}\cM_t)(e^{\sharp}_{tu})\, ,\\
& \cE_t(e^{\sharp}_{tu},e^{\sharp}_{ti})
=-\frac{1}{2}(\ast_{h_t}\cM_t)(e^{\sharp}_{ti})
\qquad\text{for }i=l,n\, ,\\
& \cE_t(e^{\sharp}_{ti},e^{\sharp}_{tj})=0
\qquad\text{for }i,j=l,n.
\end{eqnarray*}

\noindent
Here $\cE_t$ is symmetric because
$\dd_{\Sigma}\varphi_{\phi_t}=\dd_{\Sigma}\eta_t=0$.
By the compatibility identities \eqref{eq:constraintidentities}, its first component can be rewritten as:
\begin{equation*}
\cE_t(e^{\sharp}_{tu},e^{\sharp}_{tu})
=2\cW_t(e^{\sharp}_{tu})-\frac{1}{2}\cH_t,
\end{equation*}

\noindent
whereas its mixed components satisfy:
\begin{equation*}
\cE_t(e^{\sharp}_{tu},e^{\sharp}_{ti})
=\cW_t(e^{\sharp}_{ti})
\qquad\text{for }i=l,n.
\end{equation*}

\noindent
These component identities are precisely:
\begin{equation*}
\cE_t=e_{tu}\odot\cW_t
-\frac{1}{2}\cH_t\,e_{tu}\otimes e_{tu},
\end{equation*}

\noindent
which completes the proof of \eqref{eq:reducedevolutiondefects}. In particular, at every time at which the constraints
$\cH_t$, $\cW_t$ and $\cM_t$ vanish, all three reduced evolution equations hold. Hence, to prove that $(h_{\fre_t}, a_t , b_t,\phi_t,\lambda_t)$ is a Neveu-Schwarz flow, and therefore determines a solution of the Neveu-Schwarz equations of motion, it is enough to prove that $\cH_t=\cW_t=\cM_t = 0$ for every $t\in \cI$. For this, we study the \emph{propagation} of these constraints. By the identity $\cB_t=\iota_{e^{\sharp}_{tu}}\cM_t$ in \eqref{eq:reducedevolutiondefects}, the definition of $\cM_t$ and the algebraic relations $\theta_t(e^{\sharp}_{tu})=f_{b_t}$ and $\varphi_{\phi_t}(e^{\sharp}_{tu})=p_t$, we have:
\begin{eqnarray*}
D_t\theta_t = \dd_{\Sigma}f_{b_t} +f_{b_t}\eta_t+p_t\theta_t-f_{b_t}\varphi_{\phi_t}
+\iota_{e^{\sharp}_{tu}}\cM_t = \dd_{\Sigma}f_{b_t} +f_{b_t}\eta_t
+\iota_{e^{\sharp}_{tu}}\dd_{\Sigma}\theta_t.
\end{eqnarray*}

\noindent Multiplying by $\lambda_t$ and using
$\lambda_t\eta_t=\dd_{\Sigma}\lambda_t$, Cartan's formula gives:
\begin{equation*}
\partial_t\theta_t
=\dd_{\Sigma}(\lambda_t f_{b_t})
+\iota_{\lambda_t e^{\sharp}_{tu}}\dd_{\Sigma}\theta_t
=\cL_{\lambda_t e^{\sharp}_{tu}}\theta_t.
\end{equation*}

\noindent
Similarly, the reduced curvature identities for $\varphi_{\phi_t}$ imply:
\begin{equation*}
\partial_t\varphi_{\phi_t}
=\dd_{\Sigma}(\lambda_t p_t)
=\dd_{\Sigma} (
\varphi_{\phi_t}(\lambda_t e^{\sharp}_{tu})
 )
=\cL_{\lambda_t e^{\sharp}_{tu}}\varphi_{\phi_t}.
\end{equation*}

\noindent
Since the spatial exterior derivative commutes with the Lie derivative, differentiating the Maxwell constraint yields:
\begin{eqnarray*}
\partial_t\cM_t
&=& \dd_{\Sigma}(\partial_t\theta_t)
-(\partial_t\varphi_{\phi_t})\wedge\theta_t
-\varphi_{\phi_t}\wedge\partial_t\theta_t\\
&=& \cL_{\lambda_t e^{\sharp}_{tu}}
 (\dd_{\Sigma}\theta_t-\varphi_{\phi_t}\wedge\theta_t ).
\end{eqnarray*}

\noindent
Thus we obtain the following equation, which is of \emph{transport type}:
\begin{equation}
\label{eq:Maxwellconstrainttransport}
\partial_t\cM_t
=\cL_{\lambda_t e^{\sharp}_{tu}}\cM_t.
\end{equation}

\noindent
Its characteristics are the integral curves of $\partial_t-\lambda_t e^{\sharp}_{tu}$. Since the time coordinate is itself the parameter along each characteristic, a future or past endpoint $p$ of a maximal characteristic would be reached at the finite parameter value $t(p)$, which contradicts maximality by the escape lemma. Hence every maximal characteristic is an inextendible null curve and, $\Sigma$ being a Cauchy hypersurface, it meets $\{t = 0\}$ exactly once. Along each characteristic, \eqref{eq:Maxwellconstrainttransport} is a homogeneous linear ordinary differential system for the components of $\cM_t$. Therefore the initial condition $\cM_0=0$ implies:
\begin{equation*}
\cM_t=0
\qquad\text{for every }t\in\cI.
\end{equation*}

\noindent
By \eqref{eq:reducedevolutiondefects}, it follows that $\cB_t=\cD_t=0$, so the Maxwell and dilaton evolution equations hold for every $t$. Furthermore, Theorem \ref{thm:compatibilityflows} and the identity for $\cE_t$ give:
\begin{equation*}
\cH_t=2\cW_t(e^{\sharp}_{tu})\, , \qquad
\cW_t=\frac{1}{2}\cH_t\,e_{tu}\, , \qquad
\cE_t=\frac{1}{2}\cH_t\,e_{tu}\otimes e_{tu}.
\end{equation*}

\noindent
Since $\cM_t = 0$, all remaining defects are determined by the single function or scalar constraint  $\cH_t = 2\cW_t(e_{tu}^{\sharp})$ as follows:
\begin{equation*}
\mathcal{W}_t = \tfrac12 \cH_t\, e_{tu}\, , \qquad \cE_t = \tfrac12 \cH_t\, e_{tu} \otimes e_{tu}    
\end{equation*}

\noindent
It therefore suffices to show that $\cH_t$ vanishes for every $t \in I$, which we do by deriving a transport equation for it. We derive its evolution by differentiating the expression for $\cH_t$ in \eqref{eq:normalconstraint}. The Riemannian variation formulas associated to
$\partial_t h_t=-2\lambda_t\Theta_t$ are:
\begin{eqnarray*}
& \partial_t s^{h_t}
=2\lambda_t\langle\Theta_t,\Ric^{h_t}\rangle_{h_t}
-2\nabla^{h_t\ast}\nabla^{h_t\ast}(\lambda_t\Theta_t)
-2\nabla^{h_t\ast}
\dd_{\Sigma} (\lambda_t\Tr_{h_t}(\Theta_t) )\, ,\\
& \partial_t(\nabla^{h_t\ast}\varphi_{\phi_t})
=\nabla^{h_t\ast}(\partial_t\varphi_{\phi_t})
-2\lambda_t
\langle\Theta_t,\nabla^{h_t}\varphi_{\phi_t}\rangle_{h_t}\\
& \qquad
+\left\langle
2\nabla^{h_t\ast}(\lambda_t\Theta_t)
+\dd_{\Sigma} (\lambda_t\Tr_{h_t}(\Theta_t) ),
\varphi_{\phi_t}
\right\rangle_{h_t}.
\end{eqnarray*}

\noindent
Moreover, taking the $h_t$-trace in the definition of $\cE_t$, using
\begin{equation*}
D_t\Tr_{h_t}(\Theta_t)
=\Tr_{h_t}(D_t\Theta_t)+2\vert\Theta_t\vert_{h_t}^2
\end{equation*}

\noindent
and adding the dilaton evolution equation gives:
\begin{eqnarray*}
& D_t (\Tr_{h_t}(\Theta_t)+p_t )
=\cH_t+\vert\Theta_t\vert_{h_t}^2
+\frac{1}{2}\vert\theta_t\vert_{h_t}^2
+\nabla^{h_t\ast}\eta_t-\vert\eta_t\vert_{h_t}^2\\
& \qquad
+\langle\eta_t,\varphi_{\phi_t}\rangle_{h_t}
-\Tr_{h_t}(\cE_t).
\end{eqnarray*}

\noindent
Substituting these formulas into the derivative of $\cH_t$, together with the Maxwell evolution equation, the reduced curvature identities and the definition of $\cE_t$, gives:
\begin{eqnarray*}
& D_t\cH_t
=2 (\Tr_{h_t}(\Theta_t)+p_t )\cH_t
-2\nabla^{h_t\ast}\cW_t
+2\langle2\eta_t-\varphi_{\phi_t},\cW_t\rangle_{h_t}\\
& \qquad
+2\langle\Theta_t,\cE_t\rangle_{h_t}
-2 (\Tr_{h_t}(\Theta_t)+p_t )
\Tr_{h_t}(\cE_t).
\end{eqnarray*}

\noindent
In this computation, the terms containing spatial derivatives of
$\Theta_t$ and $p_t$ are grouped using:
\begin{equation*}
\nabla^{h_t\ast}\Theta_t
+\dd_{\Sigma} (\Tr_{h_t}(\Theta_t)+p_t )
=\cW_t-\Theta_t(\varphi_{\phi_t}^{\sharp},-)
+\frac{1}{2}f_{b_t}\theta_t.
\end{equation*}

\noindent
The terms involving derivatives of $f_{b_t}$ and $\theta_t$ cancel by the Maxwell evolution equation and the reduced curvature identity for $D_tf_{b_t}$. The givenevolution identity for $\cH_t$ retains the Einstein defect $\cE_t$ and therefore does not assume the Einstein evolution equation. Since $\cW_t=\frac{1}{2}\cH_t e_{tu}$, Remark \ref{remark:covariantevolution} gives:
\begin{equation*}
\nabla^{h_t\ast}\cW_t
=-\frac{1}{2}\dd_{\Sigma}\cH_t(e^{\sharp}_{tu})
+\frac{1}{2}
 (\Tr_{h_t}(\Theta_t)
-\Theta_t(e^{\sharp}_{tu},e^{\sharp}_{tu}) )\cH_t.
\end{equation*}

\noindent
Furthermore:
\begin{equation*}
\Tr_{h_t}(\cE_t)=\frac{1}{2}\cH_t\, , \qquad
\langle\Theta_t,\cE_t\rangle_{h_t}
=\frac{1}{2}\Theta_t(e^{\sharp}_{tu},e^{\sharp}_{tu})\cH_t.
\end{equation*}

\noindent
Substituting these expressions into the evolution identity for $\cH_t$ and using
$p_t=\varphi_{\phi_t}(e^{\sharp}_{tu})$, we obtain:
\begin{equation*}
D_t\cH_t-\dd_{\Sigma}\cH_t(e^{\sharp}_{tu})
=2 (
\eta_t(e^{\sharp}_{tu})
+\Theta_t(e^{\sharp}_{tu},e^{\sharp}_{tu})
 )\cH_t.
\end{equation*}

\noindent
The equations for $\fra_t$ in \eqref{eq:susyevol} imply:
\begin{equation*}
\partial_t\fra_t
=\lambda_t\eta_t(e^{\sharp}_{tu})\, , \qquad
\dd_{\Sigma}\fra_t(e^{\sharp}_{tu})
=-\Theta_t(e^{\sharp}_{tu},e^{\sharp}_{tu}),
\end{equation*}

\noindent
and hence:
\begin{equation*}
 (\partial_t-\lambda_t e^{\sharp}_{tu} )\fra_t
=\lambda_t (
\eta_t(e^{\sharp}_{tu})
+\Theta_t(e^{\sharp}_{tu},e^{\sharp}_{tu})
 ).
\end{equation*}

\noindent
Combining these identities yields:
\begin{equation*}
 (\partial_t-\lambda_t e^{\sharp}_{tu} )
 (e^{-2\fra_t}\cH_t )=0.
\end{equation*}

\noindent
Since $\cH_t=2\cW_t(e^{\sharp}_{tu})$, this is equivalently the second transport equation:
\begin{equation}
\label{eq:momentumconstrainttransport}
 (\partial_t-\lambda_t e^{\sharp}_{tu} )
 (e^{-2\fra_t}\cW_t(e^{\sharp}_{tu}) )=0.
\end{equation}

\noindent
Thus $e^{-2\fra_t}\cW_t(e^{\sharp}_{tu})$ is constant along the same characteristics used above. Every such characteristic meets $t=0$, where this function vanishes by \eqref{eq:initialdatacondition}. Consequently:
\begin{equation*}
\cW_t(e^{\sharp}_{tu})=0
\qquad\text{for every }t\in\cI.
\end{equation*}

\noindent
The compatibility identities \eqref{eq:constraintidentities} now imply $\cH_t=\cW_t=0$, and \eqref{eq:reducedevolutiondefects} gives $\cE_t=\cD_t=\cB_t=0$. Therefore the reduced tuple satisfies all evolution equations
\eqref{eq:normalevolution} and constraint equations \eqref{eq:normalconstraint}. By Proposition \ref{prop:normalhyperbolicflow}, $(g,\bar{b},\bar{\phi})$ is a solution of the Neveu-Schwarz system, and thus we conclude.
\end{proof}  

\noindent
Since the Cauchy problem for the supersymmetric Neveu-Schwarz flow is well-posed by Theorem \ref{thm:wellposedFinalSusy}, the previous proposition implies that the supersymmetric Neveu-Schwarz flow defines a Neveu-Schwarz flow when starting from common admissible data. Hence, when starting on common admissible data, the Neveu-Schwarz supersymmetry flow provides \emph{first-order integrability} of the Neveu-Schwarz flow. Similarly, the constraint equations of the supersymmetric Neveu-Schwarz flow provide partial first-order integrability of the constraint equations of the second-order Neveu-Schwarz supergravity flow, as codified in equation \eqref{eq:initialdatacondition}. This leads in turn to the following initial data characterization of supersymmetry solutions in four-dimensional Neveu-Schwarz supergravity.

\begin{cor}
\label{cor:NSNSsusyinitialdata}
Let $(\fre,\frf,b,\frc,\phi)$ be a solution of the differential system \eqref{eq:NSNSsusyconstraintsI} and \eqref{eq:NSNSsusyconstraintsII} on a bundle gerbe $(\cP,\cA,Y)$ over a connected and oriented three-manifold $\Sigma$, and let $(g,\bar{b},\bar{\phi})$ be the supersymmetric configuration \eqref{eq:NSNSexplicitdevelopment} constructed in the proof of Theorem \ref{thm:wellposedFinalSusy}. Then, $(g,\bar{b},\bar{\phi})$ is a solution of the Neveu-Schwarz system \eqref{eq:NSNSsystem4d} if and only if:
\begin{eqnarray}
\label{eq:NSNSsusysolutionconstraints}
& \dd_{\Sigma}\ast_{\fre}\frc = \varphi_{\phi}\wedge\ast_{\fre}\frc \nonumber \\
& \dd_{\Sigma} (\Tr_{h_{\fre}}(\Theta) + \varphi_{\phi}(e^{\sharp}_u))(e^{\sharp}_u) + (\nabla^{h_{\fre}\ast}\Theta)(e^{\sharp}_u) + \Theta(\varphi_{\phi}^{\sharp},e^{\sharp}_u) = \frac{1}{2} f_b^2
\end{eqnarray}

\noindent
where:
\begin{equation*}
\Theta = -\frac{1}{2} (\cL_{e^{\sharp}_u}h_{\fre} - \dd_{\Sigma}\frf\odot e_u)    
\end{equation*}

\noindent
is the second fundamental form of $\Sigma$ with respect to $g$.  
\end{cor}
  
\noindent
Since the Cauchy problem for supersymmetric Neveu-Schwarz configurations is well-posed by Theorem \ref{thm:wellposedFinalSusy}, the previous corollary can be more explicitly reformulated as follows.

\begin{cor}
\label{cor:NSNSsusyinitialdataII}
A connected and oriented three-dimensional manifold $\Sigma$ occurs as the Cauchy hypersurface of a supersymmetric Neveu-Schwarz solution if and only if it admits a bundle gerbe $(\cP,\cA,Y)$ and a tuple $(\fre,\frf,b,\frc,\phi)$  satisfying the following differential system:
\begin{eqnarray} 
 & e^{\frf} \dd_\Sigma(e^{-\frf}e_{u}) = \varphi_{\phi}\wedge e_{u}\, , \qquad   \dd_\Sigma e_{l}\wedge e_{u} = 0\, ,
 \qquad \dd_\Sigma e_{n}\wedge e_{u}=0\nonumber \\ 
& \frc = f_{b}\ast_{\fre} e_{u} + e_{u}\wedge\varphi_{\phi}\, , \qquad (\cL_{e^{\sharp}_{u}}e_{n})(e^{\sharp}_{l}) -(\cL_{e^{\sharp}_{u}}e_{l})(e^{\sharp}_{n}) = f_{b} = \ast_{\fre}H_b \label{eq:combinedconstraints}\\
& \dd_{\Sigma}\ast_{\fre}\frc = \varphi_{\phi}\wedge\ast_{\fre}\frc\, , \qquad  \dd_{\Sigma} (\Tr_{h_{\fre}}(\Theta) + \varphi_{\phi}(e^{\sharp}_u) )(e^{\sharp}_u) + (\nabla^{h_{\fre}\ast}\Theta)(e^{\sharp}_u) + \Theta(\varphi_{\phi}^{\sharp},e^{\sharp}_u) = \frac{1}{2} f_b^2\nonumber
\end{eqnarray}

\noindent
where:
\begin{equation*}
\Theta = -\frac{1}{2}(\cL_{e^{\sharp}_u}h_{\fre} - \dd_{\Sigma}\frf\odot e_u)    
\end{equation*}

\noindent
is the second fundamental form of $\Sigma$ with respect to $g$.  
\end{cor}

\noindent
By the previous corollary, we will refer to all equations in \eqref{eq:combinedconstraints} as the \emph{combined constraint equations} of the Neveu-Schwarz system. This system subsumes the necessary and sufficient conditions for a given datum to define, via evolution, a globally hyperbolic supersymmetric Neveu-Schwarz solution, and therefore provides an initial-data characterization of supersymmetry. \medskip

\noindent
Interestingly, the globally hyperbolic development of supersymmetric initial data satisfying all equations in \eqref{eq:combinedconstraints} on an oriented compact three-manifold without boundary satisfies a strong rigidity condition: the corresponding four-dimensional supersymmetric Neveu-Schwarz solution is flat and flux-less. The main difficulty in proving this result is that the flat leaves determined by the initial data may be noncompact and even dense, so normalized integration over an individual leaf is not always available. We address this situation using a leafwise Bochner--Hodge argument. The flat Weitzenböck identity identifies weakly closed and coclosed leafwise one-forms with forms whose coefficients have vanishing tangential derivatives. Orthogonal projection in $L^2$ then decomposes each leafwise closed form into a parallel part and an $L^2$ limit of leafwise exact forms \cite{AlvarezLopezTondeur,BruningLesch}. Integration by parts shows that the latter do not contribute to the pairing in the following lemma.  \medskip

\noindent
Here and below, $L^2(\Sigma,\omega)$ denotes the Hilbert space of real-valued square-integrable functions with respect to the volume form $\omega$, with respect to the standard $L^2$ inner product: 
\begin{equation*}
\langle F,G\rangle=\int_{\Sigma}FG\,\omega    
\end{equation*}

\noindent
For $F\in L^2(\Sigma,\omega)$ and a complete vector field $X$ whose flow $\Phi^X$ preserves $\omega$, the weak derivative $X(F)$ is the distribution:
\begin{equation*}
 \phi\mapsto -\int_{\Sigma}F\,X(\phi)\,\omega   
\end{equation*}

\noindent
on $C^{\infty}(\Sigma)$. It vanishes if and only if $F$ is invariant under the flow of $X$. This is because:
\begin{equation*}
\int_{\Sigma}F\,(\phi\circ\Phi^X_t)\,\omega=\int_{\Sigma}(F\circ\Phi^X_{-t})\,\phi\,\omega    
\end{equation*}

\noindent
and the derivative of the left-hand side with respect to $t$ is precisely $\int_{\Sigma}F\,X(\phi\circ\Phi^X_t)\,\omega$.

\begin{lemma}
\label{lemma:leafwiseaveraging}
Let $\Sigma$ be a connected, compact, and oriented three-manifold without boundary, equipped with a closed, nowhere-vanishing one-form $\beta$ and two commuting vector fields $e_l^{\sharp},e_n^{\sharp}$ trivializing $\Ker(\beta) \subset T\Sigma$. Let $e_l,e_n\in\Omega^1(\Sigma)$ satisfy $e_a(e^{\sharp}_b)=\delta_{ab}$ for $a,b\in\{l,n\}$, such that the volume form $\omega=\beta\wedge e_l\wedge e_n$ is positive. Let:
\begin{equation*}
\Pi\colon L^2(\Sigma,\omega)\longrightarrow L^2_{l,n}(\Sigma,\omega)
\end{equation*}
	
\noindent
be the orthogonal projection onto the closed subspace of functions whose weak derivatives with respect to $e_l^{\sharp}$ and $e_n^{\sharp}$ vanish. For $\sigma\in\Omega^1(\Sigma)$, write $\sigma_a=\sigma(e_a^{\sharp})$, $a=l,n$, and extend this notation by:
\begin{equation*}
\Pi(\sigma)=\Pi(\sigma_l)e_l+\Pi(\sigma_n)e_n.
\end{equation*}
	
\noindent
If $\sigma,\sigma^{\prime}\in\Omega^1(\Sigma)$ satisfy:
\begin{equation*}
\beta\wedge\dd_{\Sigma}\sigma=0, \qquad \beta\wedge\dd_{\Sigma} \sigma^{\prime} =0,
\end{equation*}
	
\noindent
then:
\begin{equation}
\label{eq:leafwiseaveraging}
\int_{\Sigma}\psi\,\beta\wedge\sigma\wedge\sigma^{\prime} =
\int_{\Sigma}\psi\,\beta\wedge\Pi(\sigma)\wedge\Pi(\sigma^{\prime})
\end{equation}
	
\noindent
for every $\psi\in C^{\infty}(\Sigma)$ satisfying $\dd_{\Sigma}\psi \wedge \beta=0$.  
\end{lemma}

\begin{remark}
Note that the functions whose weak derivatives along both $e^{\sharp}_l$ and $e^{\sharp}_n$ vanish are precisely those orthogonal to every function of the form $e^{\sharp}_l(\psi)+e^{\sharp}_n(\psi^{\prime})$ with $\psi,\psi^{\prime}\in C^{\infty}(\Sigma)$. They therefore form a closed subspace $L^2_{l,n}(\Sigma,\omega)\subset L^2(\Sigma,\omega)$, whose orthogonal complement is the $L^2$-closure of the space of such functions.
\end{remark}

\begin{proof}
Two choices of $e_l,e_n$ as in the statement differ by multiples of $\beta$, so neither $\omega$ nor the integrands in \eqref{eq:leafwiseaveraging} depend on the choice. Since $[e^{\sharp}_l,e^{\sharp}_n]=0$, the two-forms $\dd_{\Sigma}e_l , \dd_{\Sigma}e_n \in \Omega^2(\Sigma)$ vanish on the pair $(e^{\sharp}_l,e^{\sharp}_n)$, whence $\beta \wedge \dd_{\Sigma}e_l = \beta \wedge \dd_{\Sigma}e_n=0$. Consequently, the two-forms:
\begin{equation*}
\iota_{e^{\sharp}_l}\omega = -\beta\wedge e_n \, , \qquad \iota_{e^{\sharp}_n}\omega = \beta\wedge e_l 
\end{equation*}

\noindent
are closed, so the flows of $e^{\sharp}_l$ and $e^{\sharp}_n$ preserve $\omega$. Moreover:
	\begin{equation*}
		\beta\wedge\dd_{\Sigma}\sigma= (e^{\sharp}_l(\sigma_n)-e^{\sharp}_n(\sigma_l) )\omega\, ,
		\qquad
		\beta\wedge\sigma\wedge\sigma^{\prime}= (\sigma_l\sigma^{\prime}_n-\sigma_n\sigma^{\prime}_l )\omega\, .
	\end{equation*}
	
	\noindent
	The restrictions of $\sigma,\sigma^{\prime}$ to the leaves are closed, and the pairing depends only on these restrictions.

    Moreover, since $\beta$ vanishes nowhere, the condition $\dd _\Sigma \psi \wedge \beta = 0$ implies that $\dd _\Sigma \psi$ vanishes on $\mathrm{Ker}(\beta)$, so $\psi$ is constant along the leaves of $\mathrm{Ker}(\beta)$, that is, $\psi$ is basic. In components, Equation \eqref{eq:leafwiseaveraging} reads:
	\begin{equation*}
		\int_{\Sigma}\psi\,(\sigma_l\sigma^{\prime}_n-\sigma_n\sigma^{\prime}_l)\,\omega
		=\int_{\Sigma}\psi\, (\Pi(\sigma_l)\Pi(\sigma^{\prime}_n)-\Pi(\sigma_n)\Pi(\sigma^{\prime}_l) )\omega\, .
	\end{equation*}

\noindent
The hypotheses on $\sigma,\sigma^{\prime}$ and the fact that the flows of $e^{\sharp}_l$ and $e^{\sharp}_n$ commute, imply:
\begin{equation*}
e_l^{\sharp}(\sigma_n)-e_n^{\sharp}(\sigma_l)=0\, ,
\qquad
e_l^{\sharp}(\sigma_n^{\prime})-e_n^{\sharp}(\sigma_l^{\prime})=0.
\end{equation*}

\noindent
Equip each leaf with the flat metric for which $(e_l^{\sharp},e_n^{\sharp})$ is a parallel orthonormal frame. Since the two fields preserve $\omega$, integration by parts and commutativity of the flows give, for every $F,G\in C^{\infty}(\Sigma)$:
\begin{eqnarray}
\label{eq:leafwisebochner}
\Vert e_l^{\sharp}(F)+e_n^{\sharp}(G)\Vert_{L^2}^2
+\Vert e_l^{\sharp}(G)-e_n^{\sharp}(F)\Vert_{L^2}^2 = \Vert e_l^{\sharp}(F)\Vert_{L^2}^2
+\Vert e_n^{\sharp}(F)\Vert_{L^2}^2
+\Vert e_l^{\sharp}(G)\Vert_{L^2}^2
+\Vert e_n^{\sharp}(G)\Vert_{L^2}^2
\end{eqnarray}

\noindent
All norms are taken with respect to $\omega$. Here the mixed terms cancel because:
\begin{equation*}
\int_{\Sigma}e_l^{\sharp}(F)e_n^{\sharp}(G)\,\omega
=-\int_{\Sigma}F e_l^{\sharp}e_n^{\sharp}(G)\,\omega
=\int_{\Sigma}e_n^{\sharp}(F)e_l^{\sharp}(G)\,\omega.
\end{equation*}

\noindent
Equation \eqref{eq:leafwisebochner} is the integrated Weitzenb\"ock identity for leafwise one-forms in the parallel frame. Its left-hand side is the sum of the squared norms of the leafwise codifferential and differential, and its right-hand side is the squared norm of the leafwise covariant derivative. In $L^2(\Sigma,\omega)\oplus L^2(\Sigma,\omega)$, we subtract from $(\sigma_l,\sigma_n)$ its orthogonal projection onto the closed subspace:
\begin{equation*}
\overline{\left\{
 (e_l^{\sharp}(\alpha),e_n^{\sharp}(\alpha) )
: \alpha\in C^{\infty}(\Sigma)
\right\}}^{\,L^2},
\end{equation*}

\noindent
and denote the remainder by $(\tau_l,\tau_n)$. Orthogonality to every pair of tangential derivatives gives:
\begin{equation*}
e_l^{\sharp}(\tau_l)+e_n^{\sharp}(\tau_n)=0
\end{equation*}

\noindent
in the distributional sense. Moreover, every pair in this closed subspace satisfies the leafwise closedness equation distributionally, by commutation of the fields and continuity of differentiation from $L^2$ to distributions. Hence:
\begin{equation*}
e_l^{\sharp}(\tau_n)-e_n^{\sharp}(\tau_l)=0.
\end{equation*}

\noindent 
To apply \eqref{eq:leafwisebochner} to this weak remainder, we regularize it along the leaves. Since $\Sigma$ is compact, $e_l^{\sharp}$ and $e_n^{\sharp}$ are complete. Denote by $\Phi^l_x$ and $\Phi^n_y$ their flows at \emph{times} $x,y\in\mathbb{R}$. Since $e_l^{\sharp}$ and $e_n^{\sharp}$ commute, so do their flows, and:
\begin{equation*}
\Psi\colon\mathbb{R}^2\times\Sigma\to\Sigma\, , \qquad (x,y,p)\mapsto\Psi_{x,y}(p) = \Phi^l_x(\Phi^n_y(p))\, ,
\end{equation*}

\noindent
is an action of the additive group $\mathbb{R}^2$ on $\Sigma$, which preserves $\omega$, $e_l^{\sharp}$ and $e_n^{\sharp}$. For $F\in L^2(\Sigma,\omega)$ and $(x,y)\in\mathbb{R}^2$, set:
\begin{equation*}
U_{x,y}F = F\circ\Psi_{x,y}\, .
\end{equation*}

\noindent
Since $\Psi_{x,y}$ preserves $\omega$, every $U_{x,y}$ is a unitary operator on $L^2(\Sigma,\omega)$, and furthermore:
\begin{equation*}
U_{x,y}U_{x^{\prime},y^{\prime}} = U_{x+x^{\prime},y+y^{\prime}}    
\end{equation*}

\noindent
Moreover, $U_{x,y}F$ converges to $F$ in $L^2(\Sigma,\omega)$ as $(x,y)\to(0,0)$. For continuous $F$ the convergence is uniform, since $\Sigma$ is compact, and the general case follows because continuous functions are dense and the operators $U_{x,y}$ are unitary. For $\rho\in C^{\infty}_c(\mathbb{R}^2)$ and $F\in L^2(\Sigma,\omega)$, define:
\begin{equation*}
F_{\rho} = \int_{\mathbb{R}^2}\rho(x,y)\, U_{x,y}F\,\dd x\wedge\dd y\in L^2(\Sigma,\omega)\, ,
\end{equation*}

\noindent
as the Riemann integral of a continuous, compactly supported function with values in $L^2(\Sigma,\omega)$. Since:
\begin{equation*}
U_{t,0}F_{\rho} = \int_{\mathbb{R}^2}\rho(x-t,y)\, U_{x,y}F\,\dd x\wedge\dd y    
\end{equation*}

\noindent
and furthermore:
\begin{equation*}
 U_{0,t}F_{\rho} = \int_{\mathbb{R}^2}\rho(x,y-t)\, U_{x,y}F\,\dd x\wedge \dd y   
\end{equation*}

\noindent
the derivatives of $F_{\rho}$ along the two flows exist in $L^2(\Sigma,\omega)$ and are given by:
\begin{equation*}
D_lF_{\rho} := \lim_{t\to 0}\frac{U_{t,0}F_{\rho}-F_{\rho}}{t} = -F_{\partial_x\rho}\, , \qquad D_nF_{\rho} := \lim_{t\to 0}\frac{U_{0,t}F_{\rho}-F_{\rho}}{t} = -F_{\partial_y\rho}\, .
\end{equation*}

\noindent
Iterating, $F_{\rho}$ admits derivatives of all orders along both flows, and:
\begin{equation*}
D_lD_nF_{\rho} = D_nD_lF_{\rho} = F_{\partial_x\partial_y\rho}    
\end{equation*}

\noindent
For smooth functions $\phi$, the same limits exist uniformly and $D_a\phi = e_a^{\sharp}(\phi)$. Differentiating $\langle U_{t,0}V,U_{t,0}W\rangle = \langle V,W\rangle$ and $\langle U_{0,t}V,U_{0,t}W\rangle = \langle V,W\rangle$ at $t=0$ gives:
\begin{equation*}
\langle D_aV,W\rangle = -\langle V,D_aW\rangle    
\end{equation*}

\noindent
for all $V,W\in L^2(\Sigma,\omega)$ whose derivatives $D_aV$ and $D_aW$ exist. The proof of \eqref{eq:leafwisebochner} only uses this skew-symmetry and the commutation of the derivatives, so it applies verbatim to the pair $(F,G) = ((\tau_l)_{\rho},(\tau_n)_{\rho})$. Moreover, for every $\phi\in C^{\infty}(\Sigma)$, unitarity of $U_{x,y}$ and the invariance of $e_l^{\sharp}$ and $e_n^{\sharp}$ under $\Psi$ together with the the weak co-closedness of $(\tau_l,\tau_n)$ give:
\begin{eqnarray*}
& \langle D_l(\tau_l)_{\rho} + D_n(\tau_n)_{\rho} , \phi\rangle = -\langle\tau_l , e_l^{\sharp}(\check{\phi})\rangle - \langle\tau_n , e_n^{\sharp}(\check{\phi})\rangle = 0\, , \\
& \check{\phi} = \int_{\mathbb{R}^2}\rho(x,y)\, U_{-x,-y}\phi\,\dd x\wedge \dd y\in C^{\infty}(\Sigma)\, ,
\end{eqnarray*}

\noindent
and similarly $D_l(\tau_n)_{\rho} - D_n(\tau_l)_{\rho} = 0$ by its weak closedness. The left-hand side of \eqref{eq:leafwisebochner} therefore vanishes, so $D_a(\tau_b)_{\rho} = 0$ and furthermore:
\begin{equation*}
\langle(\tau_b)_{\rho} , e_a^{\sharp}(\phi)\rangle = -\langle D_a(\tau_b)_{\rho},\phi\rangle = 0    
\end{equation*}

\noindent
for all $a,b\in\left\{ l,n\right\}$ and $\phi\in C^{\infty}(\Sigma)$. Finally, let $\rho_k\geq 0$ have unit integral and support in the disc of radius $1/k$. Then:
\begin{equation*}
\Vert(\tau_b)_{\rho_k} - \tau_b\Vert_{L^2}\leq\int_{\mathbb{R}^2}\rho_k(x,y)\,\Vert U_{x,y}\tau_b - \tau_b\Vert_{L^2}\,\dd x\,\dd y\, ,
\end{equation*}

\noindent
which tends to zero as $k\to\infty$. Passing to the limit, we obtain: 
\begin{equation*}
e_a^{\sharp}(\tau_b)=0,
\qquad a,b\in\{l,n\},
\end{equation*}

\noindent
distributionally. Every tangential derivative of a smooth function is orthogonal to the subspace of functions with vanishing weak tangential derivatives, and the same holds for its $L^2$ limits. Therefore:
\begin{equation*}
\tau_l=\Pi(\sigma_l),
\qquad
\tau_n=\Pi(\sigma_n).
\end{equation*}

\noindent
By the definition of the closed subspace above, there exist $\alpha_j\in C^{\infty}(\Sigma)$ such that:
\begin{equation*}
 (e_l^{\sharp}(\alpha_j),e_n^{\sharp}(\alpha_j) )
\longrightarrow
 (\sigma_l-\Pi(\sigma_l),\sigma_n-\Pi(\sigma_n) )
\quad\text{in }L^2.
\end{equation*}

\noindent
Since $\psi$ is constant along the leaves, integration by parts and leafwise closedness of $\sigma^{\prime}$ give:
\begin{eqnarray*}
\int_{\Sigma}\psi ( e_l^{\sharp}(\alpha_j)\sigma_n^{\prime}
-e_n^{\sharp}(\alpha_j)\sigma_l^{\prime})\,\omega = -\int_{\Sigma}\alpha_j\psi ( e_l^{\sharp}(\sigma_n^{\prime}) - e_n^{\sharp}(\sigma_l^{\prime}))\,\omega=0.
\end{eqnarray*}

\noindent
Passing to the limit by the Cauchy--Schwarz inequality yields:
\begin{equation*}
\int_{\Sigma}\psi (\sigma_l\sigma_n^{\prime}-\sigma_n\sigma_l^{\prime} )\omega
=\int_{\Sigma}\psi (\Pi(\sigma_l)\sigma_n^{\prime}-\Pi(\sigma_n)\sigma_l^{\prime} )\omega.
\end{equation*}

\noindent
Finally, $\psi\,\Pi(\sigma_l)$ and $\psi\,\Pi(\sigma_n)$ belong to $L^2(\Sigma,\omega)$ and have vanishing weak tangential derivatives. Orthogonality of $\Pi$ therefore allows us to replace $\sigma_n^{\prime}$ and $\sigma_l^{\prime}$ in the last integral by their projection and hence we conclude.
\end{proof}

\begin{remark}
If the foliation determined by the kernel of $\beta$ admits a compact leaf $L$, then it must be a flat torus with respect to the restriction of the Riemannian metric induced by $(e_l,e_n)$. The identity \eqref{eq:leafwiseaveraging} then has the usual Hodge-theoretic interpretation, namely the integral $\int_L \sigma \wedge \sigma^{\prime}$ depends only on the cohomology classes of the closed restrictions, and their harmonic representatives have constant coefficients given by normalized integration over $L$. The proof above extends this argument to noncompact and dense leaves by taking the $L^2$-closure of the space of leafwise exact forms on the compact ambient manifold. 
\end{remark}

\begin{thm}
\label{thm:compactrigidity}
Every globally hyperbolic supersymmetric solution of the Neveu-Schwarz system that admits a compact Cauchy hypersurface is trivial, that is, its metric is flat, and it is flux-less. Moreover, the Dirac current of the supersymmetry parameter is parallel and proportional by a real constant to the dilaton curvature.  
\end{thm}

\begin{proof}
Let $(g,\bar{b},\bar{\phi},\varepsilon)$ be a supersymmetric Neveu-Schwarz solution on a bundle gerbe $(\bar{\cP},\bar{Y},\bar{\cA})$, and let $\Sigma\subset(M,g)$ be a compact spacelike Cauchy hypersurface. We restrict the bundle gerbe $(\bar{\cP},\bar{Y},\bar{\cA})$ and the supersymmetric solution $(g,\bar{b},\bar{\phi},\varepsilon)$ to $\Sigma$, and we denote by $(\fre,\frf,b,\frc,\phi)$ be the corresponding initial data. The tuple $(\fre,\frf,b,\frc,\phi)$ satisfies \eqref{eq:combinedconstraints}, namely equations \eqref{eq:NSNSsusyconstraintsI} and \eqref{eq:NSNSsusyconstraintsII} by the converse part of Theorem \ref{thm:wellposedFinalSusy}, and equations \eqref{eq:NSNSsusysolutionconstraints} by Theorem \ref{thm:susyflowcompatibility}, see Corollary \ref{cor:NSNSsusyinitialdata}. We first prove that the time-independent development of Corollary \ref{cor:timeindependentsusydevelopment} is flat and fluxless, and then transfer the conclusion to the original solution.

\medskip
\noindent
\emph{Step 1: the Maxwell constraint.} Set:
\begin{equation*}
\beta : =e^{-\frf}e_u.
\end{equation*}
	
\noindent
The supersymmetry conditions together with the Maxwell constraint give:
\begin{equation*}
\dd_{\Sigma}\beta=\varphi_{\phi}\wedge\beta, \qquad \frc=f_b\ast_{\fre}e_u+e_u\wedge\varphi_{\phi}, \qquad \dd_{\Sigma}\ast_{\fre}\frc = \varphi_{\phi}\wedge\ast_{\fre}\frc.
\end{equation*}
	
\noindent
Since $\ast_{\fre}\frc=f_be_u+\ast_{\fre}(e_u\wedge\varphi_{\phi})$ and $e_u\wedge\beta=0$, a direct computation yields:
\begin{equation*}
\dd_{\Sigma}(\ast_{\fre}\frc\wedge\beta) = - 2 \ast_{\fre}\frc\wedge\varphi_{\phi}\wedge\beta = 2e^{-\frf} \vert e_u\wedge\varphi_{\phi}\vert_{h_{\fre}}^2 \,\nu_{h_{\fre}}.
\end{equation*}
	
\noindent
Since $\Sigma$ is compact, integrating this equation we obtain $e_u\wedge\varphi_{\phi}=0$. Consequently, we can write:
\begin{equation*}
\dd_{\Sigma}\beta=0, \qquad \varphi_{\phi}=P\beta, \qquad \ast_{\fre}\frc=Q\beta,
\end{equation*}
	
\noindent
where we have defined:
\begin{equation*}
P := e^{\frf}\varphi_{\phi}(e_u^{\sharp}), \qquad Q := e^{\frf}f_b.
\end{equation*}
	
\noindent
The fact that $\varphi_{\phi} \in \Omega^1(\Sigma)$ is closed together with the previous calculation implies then:
\begin{equation*}
\dd_{\Sigma}P\wedge\beta = \dd_{\Sigma}Q\wedge\beta=0.
\end{equation*}
	
\noindent
Thus, $P,Q$ are basic functions on $\Sigma$ with respect to $\beta$, that is, constant along the leaves of $\Ker(\beta)$.
	
\medskip
\noindent
\emph{Step 2: the leafwise geometry.}
For clarity in the exposition, we introduce the following notation for the remainder of this proof:
\begin{equation*}
Z=e^{\frf}e_u^{\sharp}, \qquad \eta=e_l\wedge e_n, \qquad \omega=\beta\wedge\eta=e^{-\frf}\nu_{h_{\fre}}, \qquad
		A_{ab}=e^{\frf}\Theta(e_a^{\sharp},e_b^{\sharp}),
		\qquad a,b\in\{l,n\},
\end{equation*}
	
\noindent
where:
\begin{equation*}
\Theta = -\frac{1}{2} (\cL_{e_u^{\sharp}}h_{\fre} -\dd_{\Sigma}\frf\odot e_u).
\end{equation*}
	
\noindent
Thus $A$ is a symmetric matrix-valued function. In the following, vector fields acting on matrix-valued functions are understood as acting entrywise in the prescribed frame, and matrix norms are assumed to be Frobenius norms. We also write:
\begin{equation*}
\Delta=(e_l^{\sharp})^2+(e_n^{\sharp})^2, \qquad \nabla_L F = e_l^{\sharp}(F)e_l^{\sharp} + e_n^{\sharp}(F)e_n^{\sharp}
\end{equation*}
	
\noindent
for the leafwise Laplacian and gradient, and set:
\begin{equation*}
\sigma_1=A_{ll}e_l+A_{ln}e_n, \qquad \sigma_2=A_{ln}e_l+A_{nn}e_n.
\end{equation*}
	
\noindent
The remaining supersymmetry constraints together with the definition of $\Theta$ given above yield:
\begin{equation}
\label{eq:compactrigiditystructure}
\dd_{\Sigma}e_l = -\beta\wedge (\sigma_1+\frac{1}{2}Qe_n ), \qquad
\dd_{\Sigma}e_n = -\beta\wedge (\sigma_2-\frac{1}{2}Qe_l ).
\end{equation}
	
\noindent
Together with $\dd_{\Sigma}\beta=0$, these identities imply:
\begin{equation*}
[e_l^{\sharp},e_n^{\sharp}]=0, \qquad\cL_Z\beta=0.
\end{equation*}
	
\noindent
The \emph{screen frame} $(e^{\sharp}_l, e^{\sharp}_n)$ is therefore parallel and orthonormal for the flat metric induced on each leaf, and the flow of $Z$ maps leaves to leaves. As in Lemma \ref{lemma:leafwiseaveraging}, the corresponding screen flows preserve $\omega$. Moreover, a computation gives:
\begin{equation}
\label{eq:compactrigiditydivergence}
\dd_{\Sigma}\eta=-\Tr(A)\omega.
\end{equation}
	
\noindent
Since $\iota_Z\omega=\eta$, this is equivalent to $\mathrm{div}_{\omega}(Z)=-\Tr(A)$. Applying the exterior derivative to \eqref{eq:compactrigiditystructure}, and using $\dd_{\Sigma}Q\wedge\beta=0$, gives:
\begin{equation}
\label{eq:compactrigiditycodazzi}
\beta\wedge\dd_{\Sigma}\sigma_1 = \beta\wedge\dd_{\Sigma}\sigma_2=0.
\end{equation}
	
\noindent
These are the leafwise Codazzi equations:
\begin{equation*}
e_l^{\sharp}(A_{ln})=e_n^{\sharp}(A_{ll}), \qquad e_l^{\sharp}(A_{nn})=e_n^{\sharp}(A_{ln}).
\end{equation*}
	
\noindent
In particular, the rows of $A$ define leafwise closed one-forms, and therefore we obtain:
\begin{equation*}
\beta\wedge\sigma_1\wedge\sigma_2=\det(A)\omega.
\end{equation*}
	
\medskip
\noindent
\emph{Step 3: the scalar constraint.} The spatial Levi-Civita connection satisfies:
\begin{equation*}
\nabla^{h_{\fre}}_{e_u^{\sharp}}e_u^{\sharp} = -\sum_{a=l,n}e_a^{\sharp}(\frf)e_a^{\sharp}, \qquad \nabla^{h_{\fre}}_{e_a^{\sharp}}e_u^{\sharp} = -e^{-\frf}\sum_{b=l,n}A_{ab}e_b^{\sharp}, \qquad \nabla^{h_{\fre}}_{e_a^{\sharp}}e_b^{\sharp} = e^{-\frf}A_{ab}e_u^{\sharp},
		\qquad a,b\in\{l,n\}.
\end{equation*}
	
\noindent
The first identity follows from $\dd_{\Sigma}e_u=\dd_{\Sigma}\frf\wedge e_u$, the second from the definition of $\Theta$ and integrability of the distribution $\Ker(e_u)\subset T\Sigma$, and the third from parallelity of the screen frame along the leaves. Since:
\begin{equation*}
\Theta(e_u^{\sharp},e_u^{\sharp})=e_u^{\sharp}(\frf), \qquad \Theta(e_u^{\sharp},e_a^{\sharp})=e_a^{\sharp}(\frf),
\end{equation*}
	
\noindent
a computation gives:
\begin{eqnarray*}
(\nabla^{h_{\fre}\ast}\Theta)(e_u^{\sharp}) = -(e_u^{\sharp})^2(\frf)-\Delta\frf -2\sum_{a=l,n} (e_a^{\sharp}(\frf) )^2 + e^{-\frf}\Tr(A)e_u^{\sharp}(\frf) -e^{-2\frf}\vert A\vert^2.
\end{eqnarray*}
	
\noindent
Substituting this expression, together with:
\begin{equation*}
\Tr_{h_{\fre}}(\Theta) = e_u^{\sharp}(\frf)+e^{-\frf}\Tr(A), \qquad \Theta(\varphi_{\phi}^{\sharp},e_u^{\sharp}) = e^{-\frf}P\,e_u^{\sharp}(\frf),
\end{equation*}
	
\noindent
into the scalar equation in \eqref{eq:NSNSsusysolutionconstraints}, and multiplying by $e^{2\frf}$, gives:
\begin{equation}
\label{eq:compactrigidityscalar} 
Z (P+\Tr(A) ) = \vert A\vert^2 + \frac{1}{2}\Delta(e^{2\frf}) + \frac{1}{2} Q^2.
\end{equation}
	
\noindent
Here we have used:
\begin{equation*}
\frac{1}{2}\Delta(e^{2\frf}) = e^{2\frf} ( \Delta\frf + 2\sum_{a=l,n} e_a^{\sharp}(\frf)^2 ).
\end{equation*}
	
\noindent
For the following, consider:
\begin{equation*}
W := \Tr(A)Z-\frac{1}{2}\nabla_L(e^{2\frf}).
\end{equation*}
	
\noindent
Since $\mathrm{div}_{\omega}(Z)=-\Tr(A)$ and $\mathrm{div}_{\omega}(e_a^{\sharp})=0$, we have:
\begin{equation*}
\mathrm{div}_{\omega}(W) = Z(\Tr(A))-(\Tr(A))^2-\frac{1}{2}\Delta(e^{2\frf}).
\end{equation*}
	
\noindent
Using $\vert A\vert^2=(\Tr(A))^2-2\det(A)$, the scalar constraint reduces to:
\begin{equation}
\label{eq:compactrigidityraychaudhuri}
\dd_{\Sigma}(\iota_W\omega) = \mathrm{div}_{\omega}(W)\omega = (\frac{1}{2}Q^2-Z(P)-2\det(A) )\omega.
\end{equation}

\noindent
which we will use below.

\medskip
\noindent
\emph{Step 4: vanishing flux.}
We first identify the projection $\Pi$ of Lemma \ref{lemma:leafwiseaveraging}. Let $k$ be a primitive of the lift of $\beta$ to the universal cover of $\Sigma$. The vector field $Z$ is complete by compactness, and its lift is again complete. Since $Z(k)=1$, its lifted flow identifies the universal cover with $\mathbb{R}\times k^{-1}(0)$. Connectedness of the universal cover therefore implies connectedness of every level of $k$. Let $\Gamma\subset\mathbb{R}$ be the period group of $\beta$. It is nonzero, since a nowhere vanishing one-form on compact $\Sigma$ cannot be exact. Deck transformations translate $k$ by the elements of $\Gamma$. If $\Gamma$ is discrete, $k$ modulo $\Gamma$ defines a fibration:
\begin{equation*}
\Sigma\longrightarrow\mathbb{R}/\Gamma
\end{equation*}
	
\noindent
whose fibers are the connected compact leaves, hence flat tori. Each fiber is the image of a connected level of the function $k$, and local trivializations are provided by the flow of $Z$. If $\Gamma$ is dense, the inverse image of each leaf is the union of the levels whose values lie in one coset of $\Gamma$, so every leaf is dense. A function in the range of $\Pi$ lifts to a locally square-integrable function with vanishing weak derivatives in directions tangent to the levels of $k$. In local foliation coordinates, it therefore depends only on $k$, almost everywhere. Connectedness of the levels makes these local representatives agree, giving a locally integrable function. In the compact-leaf case, the range of $\Pi$ consequently consists of functions pulled back from the base circle. Since $\omega=\beta\wedge\eta$, orthogonal projection is normalized integration over the fibers:
\begin{equation*}
\Pi(F)\vert_L = \frac{\int_L F\,\eta}{\int_L\eta}.
\end{equation*}
	
\noindent
In the dense-leaf case, deck invariance implies that the lifted function is invariant under translations by $\Gamma$, as a locally integrable function. The continuity of translations and the density of $\Gamma$ imply invariance under every translation, so the function is constant almost everywhere. Hence:
\begin{equation*}
\Pi(F) = \frac{\int_{\Sigma}F\,\omega}{\int_{\Sigma}\omega}
\end{equation*}
	
\noindent
In the compact-leaf case, let:
\begin{equation*}
a(k)=\int_{L_k}\eta,
\end{equation*}
	
\noindent
viewed also as a basic function on $\Sigma$. It is smooth and strictly positive, and the norm of a basic function is:
\begin{equation*}
\Vert F\Vert_{L^2(\Sigma,\omega)}^2 = \int_{\mathbb{R}/\Gamma}|F(r)|^2a(r)\,\dd r.
\end{equation*}
	
\noindent
Since $a$ is bounded above and away from zero, smooth functions are dense in this weighted $L^2$ space. In the dense-leaf case, the range of $\Pi$ consists of constants. Thus smooth basic functions are dense in the range of $\Pi$ in both cases. Moreover, $\Pi$ sends smooth functions to smooth basic functions, by smooth fiber integration in the first case and global averaging in the second. For every smooth basic function $\psi$, we have:
\begin{equation*}
\dd_{\Sigma}\psi=Z(\psi)\beta.
\end{equation*}
	
\noindent
Since $\cL_Z\beta=0$, the function $Z(\psi)$ is again basic. Applying Lemma \ref{lemma:leafwiseaveraging} to $\sigma_1,\sigma_2$ gives:
\begin{equation}
\label{eq:compactrigidityaveraging}
\int_{\Sigma}\psi\det(A)\,\omega = \int_{\Sigma}\psi\det(\Pi(A))\,\omega,
\end{equation}
	
\noindent
where $\Pi$ acts entrywise on matrices. Since $\Pi(A)$ is smooth and basic, its determinant is also smooth and basic. Density therefore reduces this equation to:
\begin{equation*}
\Pi(\det(A))=\det(\Pi(A)).
\end{equation*}
	
\noindent
Multiplying \eqref{eq:compactrigidityraychaudhuri} by $\psi$ and applying Stokes' theorem, with $W(\psi)=\Tr(A)Z(\psi)$, yields:
\begin{equation}
\label{eq:compactrigidityintegrated}
\int_{\Sigma}\psi\, ( \frac{1}{2}Q^2-Z(P)-2\det(\Pi(A)) )\omega
= -\int_{\Sigma}Z(\psi)\Tr(A)\,\omega.
\end{equation}
	
\noindent
Similarly, integration by parts using $\mathrm{div}_{\omega}(Z)=-\Tr(A)$ gives:
\begin{equation*}
\int_{\Sigma}\psi\, Z(P)\,\omega = \int_{\Sigma}P (\psi\,\Tr(A)-Z(\psi) )\omega.
\end{equation*}
	
\noindent
We now choose a smooth positive basic function $\psi$ satisfying:
\begin{equation*}
Z(\psi)=\psi\,\Pi(\Tr(A)).
\end{equation*}
	
\noindent
If the leaves are compact, $\iota_Z\eta=0$ and \eqref{eq:compactrigiditydivergence} give:
\begin{equation*}
\cL_Z\eta = \iota_Z\dd_{\Sigma}\eta = -\Tr(A)\eta.
\end{equation*}
	
\noindent
The flow of $Z$ transports fibers to fibers, so differentiation of their areas gives:
\begin{equation*}
Z(a) = -\int_{L_r}\Tr(A)\eta = -a\,\Pi(\Tr(A)).
\end{equation*}
	
\noindent
Thus $\psi=1/a$ has the required property. If the leaves are dense, integrating we obtain:
\begin{equation*}
\Pi(\Tr(A)) = -\frac{\int_{\Sigma}\dd_{\Sigma}\eta}{\int_{\Sigma}\omega} =0,
\end{equation*}
	
\noindent
and we take $\psi=1$. For this choice of weight, orthogonality of $\Pi$ gives:
\begin{eqnarray*}
& \int_{\Sigma}\psi Z(P)\,\omega = \int_{\Sigma}\psi P \, (\Tr(A)-\Pi(\Tr(A)) )\, \omega =0, \\ 
& \int_{\Sigma}Z(\psi)\Tr(A)\,\omega =  \int_{\Sigma}\psi (\Pi(\Tr(A)) )^2\omega,
\end{eqnarray*}
	
\noindent
because $\psi\, P$ and $\psi\,\Pi(\Tr(A))$ are basic. Since $\Pi(A)$ is symmetric and $\Tr(\Pi(A))=\Pi(\Tr(A))$, equation \eqref{eq:compactrigidityintegrated} becomes:
\begin{equation*}
\int_{\Sigma}\psi\, (\vert\Pi(A)\vert^2+\frac{1}{2}Q^2 )\,\omega=0.
\end{equation*}
	
\noindent
Positivity of $\psi,\omega$ and smoothness imply $Q=0$ and $\Pi(A)=0$ everywhere. Consider now an arbitrary smooth basic test function $\psi$ in \eqref{eq:compactrigidityintegrated}. The vanishing just obtained gives:
\begin{equation*}
\int_{\Sigma}\psi Z(P)\,\omega = \int_{\Sigma}Z(\psi)\Tr(A)\,\omega
= \int_{\Sigma}Z(\psi)\Pi(\Tr(A))\,\omega =0 .
\end{equation*}
	
\noindent
Taking $\psi=Z(P)$, which is smooth and basic, yields $Z(P)=0$. Since $\dd_{\Sigma}P=Z(P)\beta$, we have proved:
\begin{equation}
\label{eq:compactrigiditymeanvanishing}
Q=0, \qquad \dd_{\Sigma}P=0, \qquad \Pi(A)=0.
\end{equation}
	
\noindent
By \eqref{eq:timeindependentsusyflux}, the torsion one-form of the time-independent development is:
\begin{equation*}
\alpha_{\bar b}=-\ast_{\fre}\frc=-Q\beta=0.
\end{equation*}
	
\noindent
Therefore, its flux vanishes.
	
\medskip
\noindent
\emph{Step 5: vanishing curvature.}
We consider differential forms $\Sigma$ in terms of their time-independent pullbacks to $\mathbb{R} \times \Sigma$. The isotropic coframe of Corollary \ref{cor:timeindependentsusydevelopment} reads:
\begin{equation*}
u=\beta, \qquad v=\frac{1}{2}e^{2\frf}\beta-\dd t, \qquad l=e_l, \qquad n=e_n.
\end{equation*}
	
\noindent
Since $Q=0$, the one-forms:
\begin{equation*}
\kappa = \sigma_1+\frac{1}{2}e_l^{\sharp}(e^{2\frf})\beta, \qquad \rho = \sigma_2 + \frac{1}{2}e_n^{\sharp}(e^{2\frf})\beta
\end{equation*}
	
\noindent
satisfy:
\begin{equation*}
\dd u=0, \qquad \dd l=\kappa\wedge u, \qquad \dd n=\rho\wedge u, \qquad \dd v=-\kappa\wedge l-\rho\wedge n.
\end{equation*}
	
\noindent
These are the torsion-free Cartan equations for the metric-compatible connection described in \eqref{eq:existencenullcoframegeneral}, which is therefore the Levi-Civita connection. By Proposition \ref{prop:curvaturetorsion}, its curvature is determined by the two-forms $\dd\kappa,\dd\rho \in \Omega^2(\Sigma)$. Pulling back the one-forms $\kappa,\rho\in \Omega^1(\Sigma)$ from $\Sigma$ to $\mathbb{R}\times \Sigma$ and applying Equation \eqref{eq:compactrigiditycodazzi}, we obtain:
\begin{equation*}
\beta \wedge \dd\kappa = 0\, , \qquad \beta \wedge \dd\rho=0    
\end{equation*}

\noindent
Hence there is a uniquely determined smooth matrix-valued function $C$ such that:
\begin{equation}
\label{eq:compactrigiditycurvaturematrix}
\dd\kappa = \beta\wedge(C_{ll}e_l+C_{ln}e_n), \qquad \dd\rho = \beta \wedge (C_{nl}e_l+C_{nn}e_n).
\end{equation}
	
\noindent
Consider the matrix of two-forms on $\Sigma$:
\begin{equation*}
\mathcal{M} =
\begin{pmatrix}
\kappa\wedge e_n & -\kappa\wedge e_l\\
\rho\wedge e_n & -\rho\wedge e_l
\end{pmatrix}.
\end{equation*}
	
\noindent
Using $\dd_{\Sigma}e_l=\kappa\wedge\beta$, $\dd_{\Sigma}e_n=\rho\wedge\beta$ and:
\begin{equation*}
\beta\wedge\kappa\wedge\rho = \beta\wedge\sigma_1\wedge\sigma_2 =
\det(A)\omega,
\end{equation*}
	
\noindent
we obtain:
\begin{equation}
\label{eq:compactrigiditymatrix}
\beta\wedge\mathcal{M}=A\omega, \qquad \dd_{\Sigma}\mathcal{M} =
(C-\det(A)\operatorname{Id} )\omega.
\end{equation}
	
\noindent
The four entries of the second identity are:
\begin{eqnarray*}
& \dd_{\Sigma} (\kappa\wedge e_n) = (C_{ll}-\det(A))\omega, \qquad \dd_{\Sigma}(-\kappa\wedge e_l)=C_{ln}\omega, \\
& \dd_{\Sigma}(\rho\wedge e_n)=C_{nl}\omega, \qquad \dd_{\Sigma}(-\rho\wedge e_l) = (C_{nn}-\det(A))\omega.
\end{eqnarray*}
	
\noindent
These identities yield the three properties of $C$ needed below. First:
\begin{equation*}
\Tr(\mathcal{M}) = \Tr(A)\eta + \frac{1}{2}e_l^{\sharp}(e^{2\frf})\beta\wedge e_n -\frac{1}{2}e_n^{\sharp}(e^{2\frf})\beta\wedge e_l =  \iota_W\omega.
\end{equation*}
	
\noindent
Taking the trace in \eqref{eq:compactrigiditymatrix} and comparing with \eqref{eq:compactrigidityraychaudhuri} gives:
\begin{equation*}
\Tr(C)=\frac{1}{2}Q^2-Z(P)=0.
\end{equation*}
	
\noindent
Second:
\begin{equation*}
\mathcal{M}_{ln}-\mathcal{M}_{nl} = -\kappa\wedge e_l-\rho\wedge e_n =
-\frac{1}{2}\beta\wedge\dd_{\Sigma}(e^{2\frf}).
\end{equation*}
	
\noindent
This two-form is closed because $\dd_{\Sigma}\beta=0$. Equation \eqref{eq:compactrigiditymatrix} therefore gives $C_{ln}=C_{nl}$, so $C$ is symmetric. Third, for every smooth basic $\psi$, Stokes' theorem and $\dd_{\Sigma}\psi=Z(\psi)\beta$ imply:
\begin{equation*}
\int_{\Sigma}\psi\, (C-\det(A)\operatorname{Id} )\omega = -\int_{\Sigma}Z(\psi)A\,\omega.
\end{equation*}
	
\noindent
The right-hand side vanishes because $Z(\psi)$ is basic and $\Pi(A)=0$. Equation \eqref{eq:compactrigidityaveraging} also gives:
\begin{equation*}
\int_{\Sigma}\psi\,\det(A)\,\omega=0    
\end{equation*}

\noindent
Consequently:
\begin{equation*}
\int_{\Sigma}\psi\, C\,\omega=0
\end{equation*}
	
\noindent
for every smooth basic $\psi$. Their density in the range of $\Pi$ implies, entrywise, that:
\begin{equation*}
\Pi(C)=0.
\end{equation*}
	
\noindent
Finally, applying $\dd$ to \eqref{eq:compactrigiditycurvaturematrix} gives:
\begin{equation*}
\beta\wedge\dd_{\Sigma}(C_{ll}e_l+C_{ln}e_n)=0, \qquad \beta\wedge\dd_{\Sigma}(C_{nl}e_l+C_{nn}e_n)=0.
\end{equation*}
	
\noindent
Therefore, the rows of $C$ define leafwise closed one-forms. Applying Lemma \ref{lemma:leafwiseaveraging} to these forms with $\Psi=1$, and using $\Pi(C)=0$, yields:
\begin{equation*}
0 = \int_{\Sigma} \beta\wedge(C_{ll}e_l+C_{ln}e_n) \wedge(C_{nl}e_l+C_{nn}e_n) = \int_{\Sigma}\det(C)\,\omega = -\frac{1}{2}\int_{\Sigma}\vert C\vert^2\omega.
\end{equation*}
	
\noindent
The last equality follows because $C$ is symmetric and trace-free. Positivity of $\omega$ and smoothness give $C=0$, and hence $\dd\kappa=\dd\rho=0$. Proposition \ref{prop:curvaturetorsion} therefore implies that the time-independent development is flat.
	
\medskip
\noindent
\emph{Step 6: the original solution and the dilaton.}
Since $\Sigma$ is compact, the Riemannian metric $e^{2\frf}h_{\fre}$ is complete. The universal property established in the proof of Corollary \ref{cor:timeindependentdevelopment} therefore embeds the original development into the time-independent development of the same initial data, preserving the metric, the torsion and the isotropic parallelism. Pulling back its vanishing curvature and torsion gives:
\begin{equation*}
\cR^g=0, \qquad H_{\bar b}=0
\end{equation*}
	
\noindent
on $M$. Thus the original solution is flat and flux-less. Since $H_{\bar b}=0$, the gravitino equation in \eqref{eq:KSE4d} makes $\varepsilon$ Levi-Civita parallel, and hence its Dirac current $u$ is parallel. By the proof of Proposition \ref{prop:susyevol}, the dilatino equation reduces to:
\begin{equation*}
\varphi_{\bar{\phi}}\wedge u=0.
\end{equation*}
	
\noindent
Since $u$ vanishes nowhere, $\varphi_{\bar{\phi}}=c\,u$ for a smooth function $c$. The Einstein equation in \eqref{eq:NSNSsystem4d}, with $\cR^g=0$ and $\alpha_{\bar b}=0$, gives $\nabla^g\varphi_{\bar{\phi}}=0$. Therefore:
\begin{equation*}
\dd c\otimes u=0,
\end{equation*}
	
\noindent
and hence we conclude that $c$ is constant. If $\bar{\phi}=e^{\mathrm{i}F}$ for a real-valued function $F$ on that component, then $\dd F=c\,u$. Since $u\vert_{\Sigma}=\beta$, restriction to $\Sigma$ gives:
\begin{equation*}
\dd_{\Sigma}(F\vert_{\Sigma})=c\,\beta.
\end{equation*}
	
\noindent
The function $F\vert_{\Sigma}$ has a critical point by compactness. Since $\beta$ vanishes nowhere, this forces $c=0$, and hence $\dd F=0$.  
\end{proof}

\noindent
The following example shows that the second equation in \eqref{eq:NSNSsusysolutionconstraints} cannot be omitted, even for compact initial data satisfying the first one.

\begin{ep}
\label{ep:compactremainingscalar}
Let $\Sigma = \mathbb{R}^3/(\mathbb{Z}\oplus\mathbb{Z}\oplus 2\pi\mathbb{Z})$ with coordinates $(x,y,z)$, let $n\in\mathbb{Z}\setminus\left\{ 0\right\}$ and consider the initial data:
\begin{equation*}
\fre =  ( \dd z\, , \,\, \cos(nz)\dd x + \sin(nz)\dd y\, , \,\, -\sin(nz)\dd x + \cos(nz)\dd y )\, , \qquad \frf = 0\, , \qquad \frc = -2n\, \dd x\wedge\dd y
\end{equation*}

\noindent
together with a constant dilaton $\phi$ and a curving $b$ with curvature $H_b = -2n\, \dd x\wedge\dd y\wedge\dd z$ on a bundle gerbe $(\cP,\cA,Y)$ over $\Sigma$, which exists since:
\begin{equation*}
\frac{1}{2\pi}[H_b]\in H^3_{\mathbb{Z}}(\Sigma,\mathbb{R})    
\end{equation*}

\noindent
due to the fact that:
\begin{equation*}
 \frac{1}{2\pi}\int_{\Sigma}H_b = -2n   
\end{equation*}

\noindent
We have $h_{\fre} = \dd x^2 + \dd y^2 + \dd z^2$, $\nu_{h_{\fre}} = \dd x\wedge\dd y\wedge\dd z$, $f_b = -2n$, $\varphi_{\phi} = 0$ and $\ast_{\fre}e_u = e_l\wedge e_n = \dd x\wedge\dd y$, so $\frc = f_b\ast_{\fre}e_u + e_u\wedge\varphi_{\phi}$. Moreover, $\dd_{\Sigma}e_l = n\, \dd z\wedge e_n$ and $\dd_{\Sigma}e_n = -n\, \dd z\wedge e_l$, whence $\dd_{\Sigma}e_l\wedge e_u = \dd_{\Sigma}e_n\wedge e_u = 0$ and, since $e^{\sharp}_u = \partial_z$, $\cL_{e^{\sharp}_u}e_l = n\, e_n$ and $\cL_{e^{\sharp}_u}e_n = -n\, e_l$, so that:
\begin{equation*}
 (\cL_{e^{\sharp}_u}e_n)(e^{\sharp}_l) - (\cL_{e^{\sharp}_u}e_l)(e^{\sharp}_n) = -2n = f_b   
\end{equation*}

\noindent
Since $\dd_{\Sigma}e_u = 0$, these data satisfy \eqref{eq:NSNSsusyconstraintsI} and \eqref{eq:NSNSsusyconstraintsII}. Furthermore, $\cX = e^{-\frf}e^{\sharp}_u = \partial_z$ is complete, $2 \Theta = - \cL_{\partial_z}h_{\fre} = 0$ and $\ast_{\fre}\frc = -2n\, \dd z$ is closed, so the first equation in \eqref{eq:NSNSsusysolutionconstraints} holds, whereas the second one reads $0 = 2n^2$ and fails. By Corollaries \ref{cor:timeindependentsusydevelopment} and \ref{cor:NSNSsusyinitialdata}, the supersymmetric configuration \eqref{eq:NSNSexplicitdevelopment} determined by these data is defined on the whole $\mathbb{R}\times\Sigma$, is globally hyperbolic with compact Cauchy hypersurface $\Sigma$, and is not a solution of the Neveu-Schwarz system. Explicitly, $\Phi_t(x,y,z) = (x,y,z+t)$, $\lambda_t = 1$, $h_{\fre_t} = h_{\fre}$, $\phi_t = \phi$ and $b_t = b - 2nt\, \pi^{\ast}(\dd x\wedge\dd y)$, whence:
\begin{equation*}
g = -\dd t^2 + \dd x^2 + \dd y^2 + \dd z^2\, , \qquad H_{\bar{b}} = -2n\, (\dd t + \dd z)\wedge\dd x\wedge\dd y\, , \qquad \varphi_{\bar{\phi}} = 0
\end{equation*}

\noindent
and the Dirac current of the supersymmetry parameter is $u = \dd t + \dd z$, with $\alpha_{\bar{b}} = 2n\, u$, in agreement with the first step of the proof of Theorem \ref{thm:compactrigidity}, which only uses the first equation in \eqref{eq:NSNSsusysolutionconstraints}.

Since $\alpha_{\bar b}$ is closed and null, $(g,\bar b,\bar\phi)$ satisfies the Maxwell and dilaton equations in \eqref{eq:NSNSsystem4d}, whereas the left-hand side of the Einstein equation equals:
\begin{equation*}
-\tfrac12\,\alpha_{\bar b}\otimes\alpha_{\bar b} = -2n^2\,u\otimes u    
\end{equation*}

\noindent
This agrees with the proof of Theorem \ref{thm:susyflowcompatibility}, since $2\cW(e^{\sharp}_u) = - f_b^2 = -4n^2$ and $\cH = -4n^2$. Hence, $(g,\bar{b},\bar{\phi})$ is a globally hyperbolic supersymmetric configuration with a compact Cauchy hypersurface and non-zero flux that satisfies every equation of the Neveu-Schwarz system except the component of the Einstein equation along $u\otimes u$. Here, as in Theorem \ref{thm:existencenullcoframeII}, the supersymmetry parameter is the spinor associated to the isotropic parallelism, for the spin structure selected by the latter.
\end{ep}

\noindent
On the other hand, compactness does not force the dilaton of a supersymmetric solution to be constant, as the following example shows. 

\begin{ep}
\label{ep:compactnontrivialdilaton}
Let $\Sigma = \mathbb{R}^3/(\mathbb{Z}\oplus\mathbb{Z}\oplus 2\pi\mathbb{Z})$ with coordinates $(x,y,z)$, let $m\in\mathbb{Z}$ and consider the initial data:
\begin{equation*}
\fre = (\dd z\, , \,\, \dd x\, , \,\, \dd y)\, , \qquad \frf = 0\, , \qquad \frc = 0\, , \qquad \phi = e^{\mathrm{i}\, m z}
\end{equation*}

\noindent
together with a flat curving $b$ on a bundle gerbe $(\cP,\cA,Y)$ over $\Sigma$, so that $H_b = 0$ and $f_b = 0$. We have:
\begin{equation*}
 h_{\fre} = \dd x^2 + \dd y^2 + \dd z^2\, , \qquad  \varphi_{\phi} = \phi^{\ast}\mu = m\, \dd z = \varphi_{\phi}(e^{\sharp}_u)\, e_u  
\end{equation*}

\noindent
whence $\varphi_{\phi}\wedge e_u = 0$ and $\iota_{e^{\sharp}_u}\frc = 0$. Since the coframe $\fre$ is closed and $e^{\sharp}_u = \partial_z$, these data satisfy \eqref{eq:NSNSsusyconstraintsI} and \eqref{eq:NSNSsusyconstraintsII}, and, since $\Theta = 0$ and $\ast_{\fre}\frc = 0$, they also satisfy \eqref{eq:NSNSsusysolutionconstraints}. By Corollaries \ref{cor:timeindependentsusydevelopment} and \ref{cor:NSNSsusyinitialdata}, the supersymmetric configuration \eqref{eq:NSNSexplicitdevelopment} determined by these data is a globally hyperbolic supersymmetric solution of the Neveu-Schwarz system on $\mathbb{R}\times\Sigma$ with compact Cauchy hypersurface $\Sigma$. We have:
\begin{equation*}
g = -\dd t^2 + \dd x^2 + \dd y^2 + \dd z^2\, , \qquad H_{\bar{b}} = 0\, , \qquad \bar{\phi} = e^{\mathrm{i}\, m (t + z)}\, , \qquad \varphi_{\bar{\phi}} = m\, (\dd t + \dd z) = m\, u
\end{equation*}

\noindent
where $u = \dd t + \dd z$ is the Dirac current of the supersymmetry parameter, in agreement with Theorem \ref{thm:compactrigidity}. Hence, $(g,\bar{b},\bar{\phi})$ is flux-less with flat metric and parallel null dilaton curvature, but its dilaton is not constant when $m\neq 0$, since the class:
\begin{equation*}
\frac{1}{2\pi}[\varphi_{\phi}] = \frac{m}{2\pi}[\dd z]\in H^1_{\mathbb{Z}}(\Sigma,\mathbb{R})   
\end{equation*}

\noindent
does not vanish, so $\varphi_{\phi}$ is not exact and $\phi$ does not lift to a real-valued function on $\Sigma$. In accordance with Theorem \ref{thm:compactrigidity}, the following expression:
\begin{equation*}
e^{\frf}\varphi_{\phi}(e^{\sharp}_u) = m    
\end{equation*}

\noindent
is constant, and the example shows that the vanishing of $\varphi_{\phi}$ asserted there for exact dilaton curvatures does not extend to arbitrary circle-valued dilatons. These data realize the nonzero parallel dilaton curvature allowed by the compact rigidity result of Theorem \ref{thm:compactrigidity}.
\end{ep}
  
The following example shows that Theorem 4.41 constrains only the spacetime, that is, the Cauchy data of a trivial supersymmetric solution may still be non-flat, with a non-trivial dilaton as a function into $S^1$ and a non-trivial $b$-field holonomy.

\begin{example}
Let $\Sigma = (\mathbb{R}/2\pi\mathbb{Z})^3$ with angular coordinates $(x,y,z)$, and let $M = \mathbb{R}\times\Sigma$ carry the flat, globally hyperbolic metric $g = -\dd t^2 + \dd x^2 + \dd y^2 + \dd z^2$ and the parallel isotropic coframe:
\begin{equation*}
u = \dd t + \dd z, \qquad v = \tfrac12 (\dd z - \dd t), \qquad
l = \dd x, \qquad n = \dd y.
\end{equation*}

\noindent
Fix $m \in \mathbb{Z}\setminus\{0\}$ and $\theta \in \mathbb{R}\setminus 2\pi\mathbb{Z}$, and define the dilaton $\bar\phi := e^{\mathrm{i}m(T+z)}$ together with the flat curving:
\begin{equation*}
 \bar b = \theta(2\pi)^{-2}\,\dd x\wedge\dd y   
\end{equation*}

\noindent
on the topologically trivial bundle gerbe over $M$. Then $\alpha_{\bar b} = 0$ and $\varphi_{\bar\phi} = m\,u$ is null and parallel, so $(g,\bar b,\bar\phi)$ is a trivial Neveu-Schwarz solution. Since the coframe is parallel and $\varphi_{\bar\phi}\wedge u = 0 = \varphi_{\bar\phi}(u^{\sharp_g})$, it is supersymmetric for the spin structure determined by $[u,v,l,n]$. Let:
\begin{equation*}
F(x,y) = a(\cos x + \cos y)    
\end{equation*}

\noindent
with $0 < a < 1/\sqrt2$. We introduce the hypersurface $\Sigma_F := \{t = F\}\subset \mathbb{R} \times \Sigma$, which is diffeomorphic to $\Sigma$. Since $|\dd F|^2 \le 2a^2 < 1$, $\Sigma_F$ is spacelike with induced metric $h = \dd x^2 + \dd y^2 + \dd z^2 - \dd F \otimes \dd F$, and it is a Cauchy hypersurface. By the Gauss equation, the sectional curvature of $h$ on $\mathrm{span}\{\partial_x,\partial_y\}$ at $x = y = 0$ equals $-a^2 \neq 0$. Moreover, the restricted dilaton $e^{\mathrm{i}m(F+z)}$ has curvature class $m[\dd z] \neq 0$, so it is non-constant and admits no real-valued lift. Furthermore, the restricted curving has holonomy $e^{\mathrm{i}\theta} \neq 1$ on the tori $\{z = \mathrm{const}\}$.
\end{example}

\end{document}